\documentclass[english,fleqn,10pt]{article}
\usepackage[T1]{fontenc}
\usepackage[utf8]{inputenc}
\usepackage{morefloats}
\usepackage[letterpaper]{geometry}
\usepackage{float}
\usepackage{amsmath}
\usepackage[title]{appendix}

\usepackage{graphicx}
\usepackage{color}
\usepackage{adjustbox}
\usepackage{diagbox}
\usepackage{multirow}
\usepackage{array}
\usepackage{mathrsfs}
\usepackage{amsmath, amssymb, amsthm}
\usepackage{subfigure}
\usepackage{subfiles}
\usepackage{tikz}
\usepackage{bm}
\usepackage{comment}

\usepackage[hidelinks]{hyperref}
\usepackage{url}
\usepackage[lined,boxed,linesnumbered,ruled]{algorithm2e}
\usepackage{natbib}
\usepackage{enumitem}
\setlist{nolistsep}

\usepackage{epstopdf}
\graphicspath{{./figures/}}

\newtheorem{theorem}{Theorem} 
\newtheorem{definition}[theorem]{Definition}

\newtheorem{remark}[theorem]{Remark}

\usepackage{amsfonts}

\usepackage{color}
\usepackage{adjustbox}
\usepackage{diagbox}
\definecolor{black}{rgb}{0,0,0}

\definecolor{red}{rgb}{1,0,0}

\definecolor{blue}{rgb}{0,0,1}

\newcommand{\mbs}[1]{\mathbf{#1}}
\newcommand{\mbb}[1]{\mathbb{#1}}

\newcommand{\cA}{\mathcal{A}}
\newcommand{\cB}{\mathcal{B}}
\newcommand{\cE}{\mathcal{E}}

\newcommand{\cL}{\mathcal{L}}
\newcommand{\cT}{\mathcal{T}}

\newcommand{\cN}{\mathcal{N}}

\newcommand{\bc}{\mathbf{c}}

\newcommand{\bp}{\mathbf{p}}

\newcommand{\bbR}{\mathbb{R}}

\newcommand{\bbM}{\mathbb{M}}

\newcommand{\ins}{\text{in}}
\newcommand{\bcd}{\text{bc}}

\usepackage{multirow}

\makeatother

\usepackage{babel}
\usepackage{authblk}

\usepackage{graphicx}
\def\al#1\eal{\begin{align}#1\end{align}}
\def\als#1\eals{\begin{align*}#1\end{align*}}
\def\eq#1\eeq{\begin{eqnarray}#1\end{eqnarray}}
\def\eqs#1\eeqs{\begin{eqnarray*}#1\end{eqnarray*}}
\def\eqn#1\eeqn{\begin{equation}#1\end{equation}}
\def\eqns#1\eeqns{\begin{equation*}#1\end{equation*}}
\def\ad#1\ead{\begin{aligned}#1\end{aligned}}
\def\itm#1\eitm{\begin{itemize}#1\end{itemize}}
\def\df#1\edf{\begin{definition}#1\end{definition}}
\def\pm#1\epm{\begin{pmatrix}#1\end{pmatrix}}
\def\sl#1\esl{\begin{solution}#1\end{solution}}
\def\ex#1\eex{\begin{example}#1\end{example}}
\def\tab#1\etab{\begin{table}[!htbp] \centering#1\end{table}}
\graphicspath{{./figures/}}

\title{
An Extreme Learning Machine-Based Method for Computational PDEs in Higher Dimensions
}
\author{Yiran Wang, Suchuan Dong\thanks{Author of correspondence. Emails: wang2335@purdue.edu (Y. Wang), sdong@purdue.edu (S. Dong)} \\
Center for Computational \& Applied Mathematics \\
Department of Mathematics \\
Purdue University, USA
}

\date{(September 13, 2023)}

\begin{document}
	
\maketitle
\begin{abstract}

We present two effective methods  for solving  high-dimensional partial differential equations (PDE) based on randomized neural networks. Motivated by the universal approximation property of this type of networks, both methods extend the extreme learning machine (ELM) approach from low to high dimensions. 
With the first method the unknown solution field  in $d$ dimensions is represented by a randomized feed-forward neural network, in which the hidden-layer parameters are randomly assigned and fixed while the output-layer parameters are trained. The PDE and the boundary/initial conditions, as well as the continuity conditions (for the local variant of the method), are enforced on a set of random interior/boundary collocation points. The resultant linear or nonlinear algebraic system, through its least squares solution, provides the trained values for the network parameters. 
With the second method the high-dimensional PDE problem is reformulated through a constrained expression based on an Approximate variant of the Theory of Functional Connections (A-TFC), which avoids the exponential growth in the number of terms of TFC as the dimension increases. The free field function in the A-TFC constrained expression is represented by a randomized  neural network and is trained by a procedure analogous to the first method. 
We present ample numerical simulations for a number of high-dimensional  linear/nonlinear stationary/dynamic PDEs to demonstrate their performance. These methods can produce accurate solutions to high-dimensional PDEs,  in particular with their errors reaching levels not far from the machine accuracy for relatively lower dimensions.
Compared with the physics-informed neural network (PINN) method, the current method is both cost-effective and more accurate for high-dimensional PDEs.
\end{abstract}

\vspace{0.05in}
 
 \section{Introduction}


This work concerns the numerical approximation of partial differential equations (PDEs) in higher dimensions (typically beyond three). Mathematical models describing natural and physical processes or phenomena are usually expressed in  PDEs. In a number of fields and domains, including physics, biology and finance, the models are naturally formulated in terms of high-dimensional PDEs. Well-known examples include the Schrodinger equation for many-body problems in quantum mechanics, the Black-Scholes equation for the price evolution of financial derivatives, and the Hamilton-Jacobi-Bellman (HJB) equation in dynamic programming and game theory~\cite{han2018solving,Ruthottoetal2020}.
Development of computational techniques for  PDEs  is a primary thrust  in scientific computing. In low dimensions, traditional numerical methods such as the finite difference, finite element (FEM), finite volume, and spectral type methods (and their variants), which are typically grid- or mesh-based, have achieved a tremendous success and are routinely used in computational science and engineering applications. 
For high-dimensional PDEs, on the other hand, these mesh-based approaches encounter severe challenges owing to the curse of dimensionality, because the computational effort/complexity involved therein grows exponentially with increasing problem dimension~\cite{bellman2010dynamic, darbon2016algorithms,hutzenthaler2019multilevel,han2020derivative}.

In the past few years deep neural networks (DNN or NN) have emerged as a promising approach to alleviate or  overcome the curse of dimensionality for solving high-dimensional PDEs~\cite{Beck2019Machine,Berner2020Analysis,HutzenthalerJKN2020,JentzenSW2021}.
DNN-based methods usually compute the PDE solution in a mesh-free manner by transforming the PDE problem into an optimization problem. The PDE and the boundary/initial conditions are encoded into the loss function by penalizing their residual norms on a set of sampling points. The differential operators involved therein are typically computed by automatic differentiation. The loss function is minimized by an optimizer, usually based on some flavor of gradient descent type algorithms~\cite{Werbos1974,Haykin1999}. Early works on NN-based methods for differential equations can be traced to the 1990s (see~\cite{LeeK1990,MeadeF1994,MeadeF1994b,DissanayakeP1994,YentisZ1996,LagarisLF1998}).
More recent prominent methods in this area include the physics-informed neural network (PINN) method~\cite{RaissiPK2019}, deep Galerkin method (DGM)~\cite{sirignano2018dgm}, deep Ritz method~\cite{EY2018}, deep Nitsche method~\cite{LiaoW2021}, deep mixed residual method~\cite{LyuZCC2022}, weak adversarial network method~\cite{zang2020weak}, as well as other related approaches, variants and extensions (see e.g.~\cite{zhu2019physics,nabian2019deep,KharazmiZK2019,LiTWL2020,JagtapKK2020,JagtapK2020,CyrGPPT2020,WangYP2020,Karniadakisetal2021,lu2021deepxde,KrishnapriyanGZKM2021,nakamura2021adaptive,lu2021priori,abueidda2021meshless,Penwardenetal2023}, among others).
%
Another approach for solving high-dimensional PDEs is to reformulate the problem using stochastic differential equations, thus casting the PDE  problem into a learning problem. Representative techniques of this type include the deep backward stochastic differential equation (Deep BSDE)~\cite{EHJ2017,han2018solving} and the forward-backward stochastic neural network method~\cite{Raissi2018}. Temporal difference learning has been employed in~\cite{ZengCZ2022,LuGYZ2023} for solving high-dimensional parabolic PDEs and partial integro-differential equations, which discretizes the problem in time and represents the solution by a neural network at each time step.
A data-driven method is developed in \cite{nakamura2021adaptive} to approximate the semi-global solutions to the HJB equations for high-dimensional nonlinear systems and to compute the optimal feedback controls. In \cite{lu2021priori}  the generalization error bounds are derived for two-layer neural networks in the framework of deep Ritz method for solving two  elliptic PDEs, and it is shown that the  errors are independent of the problem dimension. We would also like to refer the reader to~\cite{weinan2021algorithms} for a recent review of NN-based techniques for  high-dimensional PDEs.

%

For the neural network-based techniques reviewed above for high-dimensional PDEs, all the weight/bias parameters in the neural network are trained and determined by an optimizer, which in most cases is Adam, L-BFGS or some related variant. Unlike these methods, in the current work we consider another type of neural networks for the computation of high-dimensional PDEs, referred to as randomized  neural networks (or random-weight neural networks), in which a subset of the network parameters is assigned to random values and fixed (not trainable) while  the rest of the network parameters are trained. 


Randomness has long been exploited in neural networks~\cite{ScardapaneW2017}. Randomized neural networks can be traced  to the un-organized machine by Turing~\cite{Webster2012} and the perceptron by Rosenblatt~\cite{Rosenblatt1958} in the 1950s.
Since the early 1990s, methods based on randomized NNs have witnessed a strong resurgence and expansion~\cite{SuganthanK2021}, with  prominent techniques  widely applied  and exerting a profound influence over a variety of areas~\cite{ScardapaneW2017,FreireRB2020}. 

A simple strategy underlies randomized neural networks. Since it is extremely hard and expensive to optimize the full set of weight/bias parameters in the neural network, it seems sensible if a subset of the network parameters is randomly assigned and fixed, so that the resultant optimization problem of network training can become simpler, and in certain cases linear, hopefully without severely sacrificing the network's achievable approximation capacity~\cite{DongY2022rm,NiD2023}.
When applied to different types of neural networks or under different configurations, randomization gives rise to several techniques, including the random vector functional link (RVFL) network~\cite{PaoT1992,PaoPS1994,IgelnikP1995}, the extreme learning machine (ELM)~\cite{HuangZS2006,HuangCS2006},
and the echo-state network~\cite{JaegerLPS2007,LukoseviciusJ2009}, among others.


We consider the extreme learning machine (ELM) approach for high-dimensional PDE problems. The original  work on ELM was~\cite{HuangZS2004,HuangZS2006}, developed for linear classification and regression problems with single hidden-layer feed-forward neural networks. This method has since found widespread applications in many fields~\cite{HuangHSY2015,Alabaetal2019}. ELM is characterized by two ideas, randomly-assigned  non-trainable (fixed) hidden-layer parameters, and trainable linear output-layer parameters determined by linear least squares method or by the pseudo-inverse of coefficient matrix~\cite{VermaM1994,PaoPS1994,BraakeS1995,GuoCS1995}. Randomized neural networks of the ELM type and its close cousin RVFL type, with a single hidden layer, are universal function approximators. Their universal approximation ability has been established by the theoretical studies of~\cite{IgelnikP1995,LiCIP1997,HuangCS2006,NeedellNSS2020}. In particular, the expected rate of convergence for approximating Lipschitz continuous functions has been provided by~\cite{IgelnikP1995,RahimiR2008,NeedellNSS2020} (see also Section~\ref{sec_conv} below).

The adoption of ELM for scientific computing, in particular for the numerical solution of differential equations,  occurs only fairly recently. The existing works in this area have all been confined to PDEs in low dimensions (primarily one or two spacial dimensions) or ordinary differential equations (ODEs) so far.
Early works in this regard~\cite{YangHL2018,Sunetal2019,LiuXWL2020} have used  polynomials (e.g.~Chebyshev, Legendre, Bernstein) as activation functions  for solving linear ODEs/PDEs. Subsequent contributions have explored other types of  functions and made advances on a variety of fronts. While many studies are  confined to linear ODE/PDE problems (see e.g.~\cite{PanghalK2020,DwivediS2020,LiuHWC2021,CalabroFS2021,DwivediS2022,LiLX2023,QuanH2023}), ELM-based methods for nonlinear PDEs/ODEs have been developed in e.g.~\cite{DongL2021,DongL2021bip,Schiassietal2021,GalarisFCSS2021,FabianiCRS2021,DongY2022rm,Schiassietal2022,NiD2023,DongW2023,FabianiGRS2023,FlorioSCF2023} (among others). 
%
As has become clear from these  studies, the ELM technique can produce highly accurate solutions to linear and nonlinear PDEs in low dimensions (and ODEs) with a competitive computational cost. For smooth  solutions the ELM errors decrease exponentially as the number of degrees of freedom (number of training  points, or number of trainable parameters) increases, and the errors can reach the level of machine accuracy as the degrees of freedom become large~\cite{DongL2021,DongY2022rm}. In the presence of local complex features (e.g.~sharp gradient) in the solution field, a combination of domain decomposition and ELM, referred to as local ELM (or locELM) in~\cite{DongL2021}, will be critical to achieving a high accuracy~\cite{NiD2023}. ELM-based methods have been compared extensively with the traditional numerical methods (e.g.~classical FEM, high-order finite elements) and with the dominant DNN-based  solvers (e.g.~PINN/DGM) for low-dimensional PDE problems; see e.g.~\cite{DongY2022rm,DongL2021}. ELM far outperforms the classical FEM, and also outperforms the high-order FEM markedly when the problem size is not very small~\cite{DongY2022rm}. With a small problem size, the performance of ELM and  high-order FEM is comparable, with the latter being slightly better~\cite{DongY2022rm}.
Here ``outperform'' refers to the ability of a method to achieve a better accuracy under the same computational cost or to incur a lower computational cost for the same accuracy. 
ELM also considerably outperforms DGM and PINN for low-dimensional problems~\cite{DongL2021}. Very recently it has been shown by~\cite{DongW2023}  that the ELM-based method exhibits a spectral accuracy for solving inverse PDE problems (in low dimensions) if the measurement data is noise-free, when the network is trained by nonlinear least squares or the variable projection algorithm~\cite{GolubP1973,DongY2022}.


In the current paper we focus on the computation of high-dimensional PDEs with the ELM-based approach. To the best of the authors' knowledge, there is very little (or none) investigation in this aspect and no method of this type seems to be available in the literature for solving high-dimensional PDEs  so far. We are especially interested in the following question:
\begin{itemize}
\item Is the ELM-type randomized neural network approach effective for computational PDEs in high dimensions?

\end{itemize}
The objective of this paper is to present two ELM-based methods for solving high-dimensional PDEs, and to demonstrate with numerical simulations that these methods  provide a positive answer to the above question, at least for the range of problem dimensions studied in this paper.


The first method (termed simply ELM herein) extends the ELM technique and its local variant locELM developed in~\cite{DongL2021} (for low-dimensional problems) to linear and nonlinear PDEs in high dimensions. The solution field to the high-dimensional PDE problem is represented by a randomized feed-forward neural network, with its hidden-layer coefficients randomly assigned and fixed and its output-layer coefficients trained. Enforcing the PDE, the boundary and initial conditions on a random set of collocation points from the domain interior and domain boundaries gives rise to a linear or nonlinear algebraic system of equations about the trainable NN parameters. We seek a least squares solution to this algebraic system, attained by either linear or nonlinear least squares method, which provides the trained values for the network parameters. In the local variant of this method, the high-dimensional domain is decomposed along a maximum of $\mathcal M$ ($\mathcal M=2$ herein) directions, and the solution field on each sub-domain is represented by an ELM-type randomized neural network. We enforce the PDE, the boundary/initial conditions and appropriate continuity conditions across  sub-domains on a set of random collocation points from each sub-domain, from the domain boundaries and from the shared sub-domain boundaries. The resultant linear or nonlinear algebraic system  yields, by its least squares solution, the trained values for the network parameters of the local NNs. 

The second method (termed ELM/A-TFC herein) combines the ELM approach and an approximate variant of the theory of functional connections (TFC) for solving high-dimensional PDEs. TFC~\cite{Mortari2017,MortariL2019} provides a systematic approach for enforcing the boundary/initial conditions through a constrained expression (see e.g.~\cite{Schiassietal2021,LeakeJM2022}). However, the number of terms in  TFC constrained expressions grows exponentially with respect to the problem dimension, rendering TFC infeasible for high-dimensional problems. By noting a hierarchical decomposition of the constrained expression, we introduce an approximate variant of TFC (referred to as A-TFC herein) that retains only the dominant terms therein. A-TFC avoids the exponential growth in the number of terms of TFC and is suitable for high-dimensional problems. On the other hand, since A-TFC is an approximation of TFC, its constrained expression does not  satisfy the boundary conditions unconditionally for an arbitrary free function contained therein. However, the conditions for the free function of the A-TFC constrained expression in general involve functions of a simpler form, which is effectively a linearized form of those of the original boundary/initial conditions. A-TFC represents a trade-off. It carries a level of benefit of TFC for enforcing the boundary/initial conditions and is simultaneously suitable for high-dimensional PDEs. 
The ELM/A-TFC method uses the A-TFC constrained expression to reformulate the given high-dimensional PDE problem into a transformed problem about the free function contained in the expression. This free function is then represented by an ELM-type randomized neural network, and the reformulated PDE problem is enforced on a set of random collocation points. The least squares solution to the resultant algebraic system provides the trained values for the network parameters, thus leading to the solution for the free function. The solution to the original high-dimensional PDE problem is then computed based on the A-TFC constrained expression. 


Ample numerical simulations are presented to test these methods for a number of high-dimensional PDEs that are linear or nonlinear,  stationary or time-dependent. The current method has also been compared with the PINN method  for a range of problem dimensions. The numerical results show that the current methods exhibit a clear sense of convergence with respect to the number of training parameters and the number of  boundary  collocation points for high-dimensional PDEs. The rate of convergence is close to exponential for an initial range of parameter values (before saturation). These methods can capture the solutions to high-dimensional PDEs quite accurately, in particular with their errors reaching levels not far from the machine accuracy for comparatively lower dimensions. Compared with PINN, the current ELM method can achieve a significantly better accuracy under a markedly lower computational cost (network training time) for solving high-dimensional PDEs.


The contributions of this paper lie in the ELM method and the ELM/A-TFC method presented herein for computing high-dimensional  PDE problems. To the best of our knowledge, this seems to be the first time that a technique based on ELM-type randomized neural networks is developed for solving high-dimensional PDEs.


The methods presented in this paper are implemented in Python based on the Tensorflow and Keras libraries. The linear and nonlinear least squares methods are based on routines from the Scipy library.

The rest of this paper is organized as follows. In Section~\ref{sec_2} we first briefly recall  the theoretical result on ELM-type randomized NNs for function approximations in high dimensions, and then describe the ELM method and the ELM/A-TFC method for solving high-dimensional PDEs.
In Section~\ref{sec:tests} we present extensive numerical simulations to test these two methods with several linear and nonlinear, stationary and dynamic PDEs for a range of problem dimensions. The current method is also compared with  PINN.
Section~\ref{sec:summary}  concludes the presentation with a summary of the results and some further remarks.

\section{Extreme Learning Machine for High-Dimensional PDEs}
\label{sec_2}

 Suppose $\Omega = \Omega_1\times \Omega_2\times\cdots \times\Omega_d$ is a domain in $\bbR^d$ ($d$ being a positive integer) with boundary $\partial \Omega$,
 where $\Omega_i=[a_i,b_i]$ for given constants $a_i$ and $b_i$ ($1\leqslant i\leqslant d$). We consider the boundary value problem below,
 \begin{subequations}\label{problem}
	\begin{align}
			\mathcal Lu(x)+\mu \mathcal N(u(x))&=Q(x),\quad x\in \Omega, 
   \label{eq_1a}
   \\
			\mathcal B u(x) &= H(x),\quad x\in \partial \Omega.
   \label{eq_1b}
	\end{align}
 \end{subequations}
	Here $\mathcal L$ and $\mathcal N$ are linear and nonlinear differential operators, respectively.
 $u(x)\in \mbb R$ is the unknown field  to be computed.
 $\mathcal B$ is a linear differential or algebraic operator, and equation~\eqref{eq_1b} represents the boundary conditions.
 $Q$ and $H$  are given functions, and $\mu$ is a constant. If $\mu=0$, the problem is linear.
We assume that $\mathcal L$ may contain time derivatives (e.g.~$\frac{\partial}{\partial t}$ or $\frac{\partial^2}{\partial t^2}$, with $t$ being the time variable). In this case, the problem is time dependent, and
we will treat $t$ in the same fashion as $x$.
More specifically, we will treat this as a $(d+1)$-dimensional problem with $t$ as the last dimension, $x=(x_1,\dots,x_d,x_{d+1}=t)$, where $x_i$ ($1\leqslant i\leqslant d+1$) denotes the components of $x$. Accordingly, in this case we will assume that~\eqref{eq_1b} contains appropriate initial conditions with respect to $t$.
The point here is that the problem~\eqref{problem} may represent an initial/boundary value problem, and we will not distinguish this case in the following discussions unless necessary.

 In what follows we present two methods for solving the system~\eqref{problem}. The first is an extension to high dimensions of the ELM technique originally developed in~\cite{DongL2021} for low-dimensional problems.
  The second method is a combination of ELM with an approximate variant of the Theory of Functional Connections, termed A-TFC, which avoids the exponential growth in the computational effort of  TFC in high dimensions.

 \subsection{Randomized Feed-Forward Neural Networks}
 \label{sec:rfnn}


We consider the approximation of the solution field $u(x)$ to system~\eqref{problem} by a randomized feed-forward neural network.
A feed-forward neural network (FNN) having $(L+1)$ ($L\geqslant 2$)  layers represents a parameterized function $\mathcal G(x;\theta)$ given by (for the input $x$ and parameter $\theta$)~\cite{GoodfellowBC2016},
	\begin{equation}\label{eq_2}
	\mathcal{G}(x ; \theta)=W_L\cdot\sigma\left(W_{L-1} \cdot \sigma\left(\cdots \sigma\left(W_2\cdot \sigma\left(W_1 \cdot x+b_1\right)+b_2\right) \cdots\right)+b_{L-1}\right)+b_L,
	\end{equation}
where $W_i$ and $b_i$ ($1\leqslant i\leqslant L$) are the weight  and bias in the $i$-th layer,  $\theta=\left(W_1, \ldots, W_L, b_1, \ldots, b_L\right)$, and $\sigma: \mathbb R\rightarrow\mathbb R$ is the activation function.
Layer $0$ (input layer) contains $d$ nodes, representing the components of $x$, and layer $L$ (output layer) contains a single node, representing $u(x)$. The layers in between are the hidden layers.
Note that the output layer in~\eqref{eq_2} is linear, with no activation function applied. In the current paper we  further assume that the output layer has zero bias, i.e.~$b_L=0$.

A randomized feed-forward neural network is an FNN in which a subset of the network parameters $\theta$ are  assigned to random values and fixed (non-trainable), while only the rest of the network parameters are trained. Extreme learning machine (ELM)~\cite{HuangZS2006} is one type of randomized neural networks, in which all the hidden-layer coefficients are randomly assigned and fixed and only the output-layer coefficients are trained~\cite{DongL2021}.
In the current work we approximate the solution field $u(x)$ to the system~\eqref{problem} by ELM, and  assign the network coefficients in all the hidden layers to uniform random values from the interval $[-R_m,R_m]$, where $R_m$ is a constant. 

 \subsection{Randomized NNs for High-Dimensional Function Approximation}
\label{sec_conv}

Extreme learning machines (or RVFL networks) are universal function approximators; see e.g.~\cite{IgelnikP1995,HuangCS2006}. The universal approximation theorems~\cite{IgelnikP1995,HuangCS2006}
basically state that any given continuous function can be approximated by a randomized NN having a single hidden layer, in which the hidden-layer coefficients are randomly assigned and fixed and the output-layer coefficients can be adjusted/trained, to any desired degree of accuracy, if the number of hidden units is sufficiently large.

We next recall the result from~\cite{IgelnikP1995} concerning the convergence rate of randomized NNs for function approximations in high dimensions, which motivates the development of the current methods for high-dimensional computational PDEs.

	Define $I^d:=[0,1]^d$ and consider a continuous function $f\in C(I^d)$ satisfying the Lipschitz condition, i.e., there exists a constant $\gamma>0$ such that for any $x,y\in I^d$,
	\begin{align}
		|f(x)-f(y)|\leq \gamma \|x-y\|,\label{lip}
	\end{align}
	where $\|x-y\|:=\sum_{i=1}^d|x_i-y_i|.$
	To approximate $f$, we construct a sequence  functions $\{f_{\omega_n}\}$  as follows,
	\begin{align}\label{eq_4}
		f_{\omega_n}(x) = \sum\limits_{j=1}^n a_jg(w_j\cdot x+b_j),
	\end{align}
	where $\omega_n$ is defined as $\omega_n:=(n,a_1,\ldots, a_n,b_1,\ldots, b_n,w_1,\ldots,w_n)$. In particular, $\lambda_n:=\{b_1,\ldots, b_n,w_1,\ldots,w_n\}$ denotes a set of random parameters from some probabilistic space $S_n(\Omega,\alpha)$, where $\alpha$ is a parameter. The corresponding  probability measure $\mu_{n,\Omega,\alpha}$ is specified as follows. Suppose  $\hat{w}_0=(\hat{w}_{01},\ldots, \hat{w}_{0d})$, $y_0=(y_{01},\ldots, y_{0d})$ and $u_0$ are independent and uniformly distributed in $V^d=[0,\Omega]\times \cdots \times [-\Omega, \Omega]$, $I^d$ and $[-2\Omega,2\Omega]$, respectively. Then $w_0=\alpha \hat{w}_0$ and $b_0 = -w_0\cdot y_0-u_0$.  $(w_1,\ldots, w_n)$ and $(b_1,\ldots,b_n)$ are two sets of samples of the random variables $w_0$ and $b_0$. $g$ in~\eqref{eq_4} is the activation function,  chosen to be absolutely integrable, i.e., 
	\begin{align}
		\int_{\mathbb{R}} g^2(x)dx< +\infty.\label{abs_int}
	\end{align}
	We further restrict $g$ on some compact support $\prod_{i=1}^{d}[-\beta w_i,\beta w_i]$ to get $g_{\beta}$ ($\beta$ denotes a parameter). Then the following result holds.
	\begin{theorem}
		\cite{IgelnikP1995} For any $f\in C(I^d)$ satisfying \eqref{lip}, any  compact $K\subset I^d$ that is a proper subset of $I^d$, and any activation function $g_{\beta}$ satisfying \eqref{abs_int}, there exists a sequence of $\{f_{\omega_n}\}$ and probability measure $\mu_{n,\Omega,\alpha}$ such that 
		\begin{align}
			E\int_{K}|f(x)-f_{\omega_n}(x)|^2 dx\leq \frac{C_{f,g,\Omega,\alpha,\beta,d}}{n},
			\label{approx_err}
		\end{align}
  for some constant $C_{f,g,\Omega,\alpha,\beta,d}$ independent of $n$.
	\end{theorem}
	\begin{remark}
		The theorem can be generalized when $I^d$ is replaced by $[a,b]^d$  by a change of variables. 
  We omit this detail and consider the generic situation.
	\end{remark}
	\begin{remark}
		It is notable that the approximation error is of the order $1/\sqrt{n}$ as the number of  basis functions $n$ increases, irrespective of the dimension $d$. 
  This indicates that the approximation~\eqref{eq_4} with random basis functions can be effective for high dimensions.
  On the other hand, one notes from \cite{barron1993universal} that  the approximation by linear combinations of deterministic and fixed bases leads to an approximation error on the order of $\mathcal O(1/n^{\frac{1}{d}})$. 
  In other words, if  deterministic and fixed basis functions are used, it is impossible to avoid the exponential growth in $d$ for the number of basis functions. 
  The randomized bases (such as  in ELM and RVFL), however, can be effective  for high-dimensional function approximations  in the sense of the expectation.
  
	\end{remark}
	
 \subsection{Solving High-Dimensional PDEs with ELM}
\label{sec_elm}

Adopting ELM for computational PDEs is characterized by two ideas:
(i) The hidden-layer coefficients are assigned to random values and fixed, and only the output-layer coefficients are trainable, as already mentioned previously.
(ii) The trainable network parameters are determined by the linear or nonlinear least squares method~\cite{DongL2021}, not by the gradient descent type algorithms.
This means that in equation~\eqref{eq_2} the coefficients $(W_i,b_i)$ for $1\leqslant i\leqslant L-1$ will be assigned to uniform random values from $[-R_m,R_m]$ and fixed, while only $W_L$ is trained (noting that we set $b_L=0$). 

Let $N_u$ denote the number of nodes in the last hidden layer of the neural network, and $V_j(x)$ ($1\leqslant j\leqslant N_u$) denote the output fields of the last hidden layer. Then equation~\eqref{eq_2} can be written into,
\begin{equation}
u(x) = \sum_{i=1}^{N_u}\phi_j V_j(x)
=\mbs V(x)\bm\Phi,
\quad x\in\Omega,
\end{equation}
where 
$\mbs V(x)=(V_1(x),\dots,V_{N_u}(x))
=\sigma\left(W_{L-1} \cdot \sigma\left(\cdots \sigma\left(W_2\cdot \sigma\left(W_1 \cdot x+b_1\right)+b_2\right) \cdots\right)+b_{L-1}\right)
$, and
$\bm\Phi = (\phi_1,\dots,\phi_{N_u})^T=W_L$.
Note that $\mbs V(x)$ is fixed once the hidden-layer coefficients are randomly assigned and $\phi_j$ ($1\leqslant j\leqslant N_u$) are the output-layer coefficients (trainable parameters) of ELM.

The residual function of the system~\eqref{problem} is,
\begin{equation}\label{eq_8}
R(x,y,\bm\Phi) 
=\begin{bmatrix}
R^{\text{pde}}(x,\bm\Phi)\\
R^{\text{bc}}(y,\bm\Phi)
\end{bmatrix}
=\begin{bmatrix}
\mathcal L\mbs V(x)\bm\Phi
+\mu\mathcal{N}(\mbs V(x)\bm\Phi)
-Q(x), \ x\in\Omega \\
\mathcal B\mbs V(y)\bm\Phi - H(y), \ y\in\partial\Omega
\end{bmatrix}
\end{equation}
where $R^{\text{pde}}$ and $R^{\text{bc}}$ are the residuals corresponding to the PDE and the boundary/initial conditions, respectively.

We next choose a set of collocation points from the domain interior and the domain boundaries, and enforce the residual function~\eqref{eq_8} to be zero on these collocation points.
For solving high-dimensional PDEs we will employ a set of random collocation points from the interior/boundary of the domain in this work. We note that the regular grid points, which have been used with ELM in e.g.~\cite{DongL2021,DongL2021bip,DongY2022rm,NiD2023,DongW2023} as the sampling points for low-dimensional problems, are not feasible for high-dimensional PDEs, because of the exponential growth in the number of  points with increasing dimension.

More specifically, the collocation points are set as follows.
Let $N_{\ins}$ denote the number of random collocation points in the interior of $\Omega$, and $N_{\bcd}$ denote the number of random collocation points on each hyperface of $\partial\Omega$.
 To generate the interior collocation points we set 
 \begin{equation}
 P_i =(x_{1,i},\ldots, x_{d,i}),
 \quad \text{for}\ i\in[N_{\ins}],
 \end{equation} 
 where $x_{j,i}\in (a_j,b_j)$ is a uniform random value for $j\in[d]$.  Here we have used the notation in combinatorics that $[d]=\{1,\ldots,d\}$. 
 For the boundary collocation points, we choose $N_{\bcd}$ random points on each hyperface  of $\partial\Omega$, which is a  $(d-1)$-dimensional  hyperplane, and there are $2d$ hyperfaces  in total. The total number of boundary collocation points  is $N_{\bcd}^{\text{tot}}=2dN_{\bcd}$. For $i\in [N_{\bcd}]$, $j\in [d]$ and $l\in\{0,1\}$, we set 
 \begin{equation}
 P_{i,j,l} = (x_{1,i},\ldots, x_{j-1,i}, x_{j,i}=a_j\delta_{l,0}+b_j\delta_{l,1},x_{j+1,i},\ldots, x_{d,i})
 \end{equation}
 as the boundary collocation points, in which $x_{k,i}\in (a_k,b_k)$ is a uniform random value if $k\neq j$. 
 Here $\delta_{i,j}$ denotes the Kronecker delta, $\delta_{i,j}=1$ if $i=j$ and $0$ otherwise. Overall, the total number of collocation points is $N_c = N_{\ins}+N_{\bcd}^{\text{tot}}$.
%
Let $N_a = N_c + N_{\bcd}^{\text{tot}}$, and
\begin{equation}\label{eq_11}
\mbs x_{\ins} = \begin{bmatrix}
\vdots\\
P_i\\ \vdots
\end{bmatrix} \in \mbb M^{N_{\ins}\times d},
\quad
\mbs y = \begin{bmatrix}
\vdots\\
P_{i,j,l}\\ \vdots
\end{bmatrix} \in \mbb M^{N_{\bcd}^{\text{tot}}\times d},
\quad
\mbs x=\begin{bmatrix}
\mbs x_{\ins}\\
\mbs y
\end{bmatrix}
\in \mbb M^{N_c\times d},
\end{equation}
where $\bbM^{a\times b}$ denotes the set of  matrices with shape $a\times b$.

Enforcing the residual function~\eqref{eq_8} to be zero on all the collocation points gives rise to the following  system,
\begin{equation}\label{eq_12}
\mbs 0 = \mbs R(\bm\Phi)=
R(\mbs x,\mbs y,\bm\Phi) 
=\begin{bmatrix}
R^{\text{pde}}(\mbs x,\bm\Phi)\\
R^{\text{bc}}(\mbs y,\bm\Phi)
\end{bmatrix}
=\begin{bmatrix}
\mathcal L\mbs V(\mbs x)\bm\Phi
+\mu\mathcal{N}(\mbs V(\mbs x)\bm\Phi)
-Q(\mbs x) \\
\mathcal B\mbs V(\mbs y)\bm\Phi - H(\mbs y)
\end{bmatrix}
=\begin{bmatrix}
\mbs R^{\text{pde}}(\bm\Phi)\\
\mbs R^{\text{bc}}(\bm\Phi)
\end{bmatrix}.
\end{equation}
Here $\mbs R(\bm\Phi)\in\mbb M^{N_a\times 1}$, $\mbs R^{\text{pde}}(\bm\Phi)\in\mbb M^{N_{c}\times 1}$,
$\mbs R^{\text{bc}}(\bm\Phi)\in\mbb M^{N_{\bcd}^{\text{tot}}\times 1}$,
$\mbs V(\mbs x)\in\mbb M^{N_{c}\times N_u}$, and $\mbs V(\mbs y)\in\mbb M^{N_{\bcd}^{\text{tot}}\times N_u}$.

The  system~\eqref{eq_12} is an algebraic system about $\bm\Phi$,
containing $N_a$ equations with $N_u$ unknowns. We seek a least squares solution to this system.
When $\mu=0$, this system is linear about $\bm\Phi$ and can be written as
\begin{equation}\label{eq_13}
\begin{bmatrix}
\mathcal L\mbs V(\mbs x) \\
\mathcal B\mbs V(\mbs y)
\end{bmatrix} \bm\Phi
=
\begin{bmatrix}
Q(\mbs x) \\
 H(\mbs y)
\end{bmatrix}.
\end{equation}
In this case we compute $\bm\Phi$ by solving the system~\eqref{eq_13} using the linear least squares method (with a minimum norm if the coefficient matrix is rank deficient)~\cite{Bjorck1996}.

When $\mu\neq 0$, the algebraic system~\eqref{eq_12} is nonlinear with respect to $\bm\Phi$. In this case we compute $\bm\Phi$ by solving this system using the nonlinear least squares method with perturbations (NLLSQ-perturb) from~\cite{DongL2021,DongY2022rm,DongW2023}. The nonlinear least squares method~\cite{Bjorck1996} as in~\cite{DongL2021,DongY2022rm,DongW2023} represents a  Gauss-Newton method 
combined with a trust-region strategy.
The NLLSQ-perturb algorithm requires, for an arbitrary given $\bm\Phi\in\mbb R^{N_u}$,
the computation of $\mbs R(\bm\Phi)$ and its Jacobian matrix $\mbs J(\bm\Phi)=\frac{\partial\mbs R}{\partial\bm\Phi}\in\mbb M^{N_a\times N_u}$ for the Gauss-Newton iterations (see~\cite{DongL2021,DongY2022rm} for details). The Jacobian matrix for~\eqref{eq_12} is given by
\begin{equation}
\mbs J(\bm\Phi) = \begin{bmatrix}
\mathcal L\mbs V(\mbs x) + \mu\mathcal N'(\mbs V(\mbs x)\bm\Phi)\mbs V(\mbs x) \\
\mathcal B\mbs V(\mbs y)
\end{bmatrix},
\end{equation}
where $\mathcal N'(u)=\frac{\partial\mathcal N}{\partial u}$.


\begin{remark}
In our implementation, the input data to the ELM neural network consists of $\mbs x$ (coordinates of  all  collocation points), and the output data consists of $u(\mbs x)\in\mbb M^{N_c\times 1}$ (the  solution field evaluated on the collocation points). After the hidden-layer coefficients are randomly assigned, $V(\mbs x)$ is computed by a forward evaluation of a sub-network, implemented in Keras as a sub-model of the original NN, whose input is $\mbs x$  and  output is the last hidden layer of the original NN. 
The differential operators involved in $\mathcal LV(\mbs x)$ and $\mathcal BV(\mbs x)$ are then computed by a forward-mode automatic differentiation with this sub-model. The linear least squares method is based on the routine scipy.linalg.lstsq from the SciPy library. 
The NLLSQ-perturb algorithm for the nonlinear least squares method is based on the routine scipy.optimize.least\_squares from the SciPy library (see~\cite{DongL2021} or the Appendix A of~\cite{DongW2023} for more details).

\end{remark}

\begin{remark}\label{rem_a5}
When the solution field contains local features (e.g.~sharp gradient), it would be  preferable to combine ELM with domain decomposition for approximating the  solution, thus leading to the local ELM  approach (called locELM in~\cite{DongL2021}). In this case, we represent the solution on each sub-domain by a local randomized FNN, and impose $C^k$ (with $k$ related to the PDE order) continuity conditions across the shared sub-domain boundaries.

For high-dimensional PDE problems, if domain decomposition is performed in every direction, the number of sub-domains would increase exponentially with increasing problem dimension. To avoid the exponential growth in the number of sub-domains in local ELM, we require that the domain should be decomposed only along a maximum of $\mathcal M$ directions, where $\mathcal M$ is a fixed small integer   ($0\leqslant\mathcal M\leqslant d$). The specific directions in which the domain is decomposed can be any subset of the $d$ directions with a size not exceeding $\mathcal M$. We use  $\mathcal M=2$ in the current work, i.e.~domain decomposition in a maximum of two directions with local ELM.


In the presence of domain decomposition, the residual function~\eqref{eq_8} for the system~\eqref{problem} needs to be modified accordingly to account for the $C^k$ continuity conditions across sub-domain boundaries. Let $N_{\Omega}$ denote the number of sub-domains, and $\Omega_i$ ($1\leqslant i\leqslant N_{\Omega}$) denote the $i$-th sub-domain.
Symbolically, the modified residual function can be written as,
\begin{equation}\label{eq_15}
R(x,y,z,\bm\Phi) 
=\begin{bmatrix}
R^{\text{pde}}(x,\bm\Phi)\\
R^{\text{bc}}(y,\bm\Phi)\\
R^{\text{ck}}(z,\bm\Phi)
\end{bmatrix}
=\begin{bmatrix}
\mathcal L\mbs V_i(x)\bm\Phi_i
+\mu\mathcal{N}(\mbs V_i(x)\bm\Phi_i)
-Q(x), \ x\in\Omega_i, \ 1\leqslant i\leqslant N_{\Omega} \\
\mathcal B\mbs V_i(y)\bm\Phi_i - H(y), \ y\in\partial\Omega_i\cap \partial\Omega, \ 1\leqslant i\leqslant N_{\Omega} \\
\mathcal C\mbs V_i(z)\bm\Phi_i - \mathcal C\mbs V_j(z)\bm\Phi_j,
\ z\in\partial\Omega_i\cap\partial\Omega_j,
\ \forall\ \text{adjacent}\ (\Omega_i, \Omega_j)
\end{bmatrix},
\end{equation}
where $\mbs V_i(x)=(V_{i1}(x),\dots,V_{iN_u}(x))$ denotes the output fields of the last hidden layer of the local NN on $\Omega_i$,  $\bm\Phi_i=(\phi_{i1},\dots\phi_{iN_u})^T$ denotes the vector of output-layer coefficients of the local network for $\Omega_i$, and $\bm\Phi=(\bm\Phi_1^T,\dots,\bm\Phi_{N_{\Omega}}^T)^T$ denotes the set of training parameters of the overall problem. The solution field $u(x)$, when restricted to $\Omega_i$ ($1\leqslant i\leqslant N_{\Omega}$), is given by
\begin{equation}
u_i(x) = \sum_{j=1}^{N_u}\phi_{ij}V_{ij}(x)
=\mbs V_i(x)\bm\Phi_i,
\quad x\in\Omega_i.
\end{equation}
In~\eqref{eq_15}, $R^{\text{ck}}(z,\bm\Phi)=\mathcal C u_i(z)-\mathcal C u_j(z)$ denotes the residual corresponding to the continuity conditions between $u_i(z)$ and $u_j(z)$ on their shared sub-domain boundary, and $\mathcal C$ denotes a differential (or algebraic) operator corresponding to the continuity conditions. If the PDE order is $m$ ($m\geqslant 1$) with respect to $x_i$,  we will in general impose $C^{m-1}$ conditions across the sub-domain boundaries along the $x_i$ direction.

Accordingly, we choose a set of random collocation points on the interior and on the boundaries of each sub-domain, and enforce the residual function~\eqref{eq_15} to be zero on these collocation points. This leads to a linear or nonlinear algebraic system about $\bm\Phi$, which is solved by the linear or nonlinear least squares method to attain a least squares solution for the training parameters $\bm\Phi$.

	\begin{figure}[tb]
			\centering
	\includegraphics[height=1.5in]{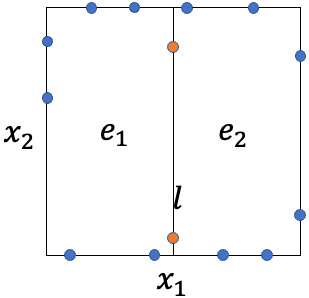}
			\caption{Illustration of random boundary collocation points on two adjacent sub-domains in two dimensions.  
   }
			\label{fg_1}
		\end{figure}

Let us next further comment on enforcing the continuity conditions across sub-domains and illustrate it with an example.
In order to impose the continuity conditions on the common boundary between two adjacent sub-domains $(\Omega_i,\Omega_j)$, the random boundary collocation points for $\Omega_i$ and the random boundary collocation points for $\Omega_j$, when restricted to their shared boundary, must be identical.
%
We illustrate this point using
   Figure \ref{fg_1}, which shows two adjacent sub-domains in 2D.
   Two random collocation points ($N_{\bcd}=2$) are generated on each boundary of each sub-domain.
Note that the  collocation points for different sub-domains are generated independently. But we need to make sure that the points for the common face $l$ (highlighted in orange in Figure~\ref{fg_1}) in both sub-domains sharing $l$ should be identical.

This requirement is implemented by the following procedure in our implementation. For each sub-domain we arrange the boundary collocation points in the following order: those random points on the left face, followed by those on the right face, the bottom face, and the top face. We first generate the random boundary  points for all sub-domains independently. Then we go through all sub-domains, and the four boundaries (left, right, bottom, top) on each sub-domain, successively. 
If the boundary being examined is a shared boundary and the ID of the current sub-domain is higher than that of the neighboring sub-domain, then we replace the collocation points for the this  boundary of the current sub-domain by those  collocation points for the same  boundary in  the neighboring sub-domain. 
As shown in Figure~\ref{fg_1}, 
one can see that the random collocation points for the left boundary of the second sub-domain will be replaced by (thus identical to) those for the right boundary of the first sub-domain.
   
\end{remark}

\begin{remark}\label{rem_a6}

As mentioned previously, the hidden-layer coefficients of the ELM neural network are assigned to uniform random values generated on the interval $[-R_m,R_m]$ in this paper. We observe from numerical simulations that the value of the constant $R_m$ influences the ELM accuracy. In  the numerical simulations of Section~\ref{sec:tests} below, we have used $R_m=R_{m0}$, where $R_{m0}$ is a value determined by the following procedure for a given problem. 
For a  PDE problem of a given dimension $d$, we first fix the number of training parameters in the NN and the number of  collocation points to some chosen values. Then we perform  preliminary simulations of the given problem using this fixed network setting, and a set of different $R_m$ values for generating the random hidden-layer coefficients. We record the errors of the computed solution (when the exact solution is available), or the norm of the residual vector $\mbs R(\bm\Phi)$ corresponding to the computed solution (when the exact solution is unavailable), for this set of  $R_m$ values. We choose the value with the lowest error or the lowest residual norm, and denote it by $R_{m0}$. Then we fix $R_m=R_{m0}$ for generating the random hidden-layer coefficients for the subsequent simulations of this PDE problem in a given dimension $d$, while some other  simulation parameters (e.g.~number of training parameters or collocation points) are varied.
We observe that the $R_{m0}$ as determined above in general decreases with respect to the  problem dimension $d$.

\end{remark}

\subsection{Solving High-Dimensional PDEs by Combined ELM and Approximate Theory of Functional Connections (A-TFC)}
\label{sec:atfc}

 In this subsection we present an alternative strategy for approximating the solution field $u(x)$, by combining ELM and an approximate variant of the Theory of Functional Connections (TFC). TFC provides a systematic technique for handling linear constraints, in particular the boundary or initial conditions~\cite{Mortari2017,MortariL2019,Schiassietal2021,LeakeJM2022}.  The number of terms in TFC constrained expression, however, increases exponentially with increasing problem dimension, rendering TFC  infeasible for high-dimensional problems. The approximate variant of TFC, termed A-TFC here, 
 avoids the exponential growth  in the number of terms of TFC and  is  suitable for computational PDEs in high dimensions.

 \subsubsection{TFC and Approximate TFC} 
	\label{sec_TFC}

Consider the domain $\Omega=[a,b]^d$, where $a$ and $b$ are constants, and a function $u(x)$ defined on $\Omega$ satisfying the  condition,
\begin{equation}\label{eq_17}
u|_{\partial\Omega} = C(x), \quad
x\in \partial\Omega,
\end{equation}
where $C(x)$ is a prescribed function  on $\partial\Omega$. Then the general form of $u(x)$ is given by~\cite{LeakeJM2022},
\begin{equation}\label{eq_18}
u(x) = g(x) - \mathcal T_{\Omega}g(x) + \mathcal T_{\Omega}C(x), \quad x\in\Omega,
\end{equation}
where $g$ is an arbitrary (free) function on $\Omega$. $\mathcal T_{\Omega}$ is a linear operator
satisfying the property that, for any function $f$ defined on $\partial\Omega$, $\mathcal T_{\Omega}f(x)$ is a function defined on $\Omega$ and,
\begin{equation}\label{eq_19}
\left.\mathcal T_{\Omega}f\right|_{\partial\Omega} = f(x),
\quad x\in\partial\Omega.
\end{equation}
One can verify that $u(x)$ as given by~\eqref{eq_18} satisfies~\eqref{eq_17} for an arbitrary $g$.
The expression~\eqref{eq_18} is called the constrained expression of TFC~\cite{LeakeJM2022}.

For a function $f$ defined on $\partial\Omega$, $\mathcal T_{\Omega}f(x)$ can be constructed as follows.
We first use a 2D ($d=2$) example to illustrate the  procedure and then discuss the general $d$-dimensional case.
	Suppose $\Omega=[a,b]^2$ and $f(x)$ is a function defined on $\partial\Omega$. 
 Let $c_i^1\in\{a,b\}$ for $i\in\{1,2\}$, and define $f_i^{c_i^1}(x) := f(x)|_{x_i=c_i^1}$. Let  $c_{1,2}^i\in \{a,b\}$ for $i\in\{1,2\}$, and define $f_{1,2}^{c_{1,2}^1,c_{1,2}^2}(x) := f(x_1=c_{1,2}^1,x_2=c_{1,2}^2)$.  
 We further define
	\begin{align}
		\phi_{a,b}^{c_i^1} (x)= \delta_{a,c_i^1}\frac{b-x_i}{b-a}+\delta_{b,c_i^1}\frac{x_i-a}{b-a},
    \quad
    \phi_{a,b}^{c_{1,2}^i} (x)= \delta_{a,c_{1,2}^i}\frac{b-x_i}{b-a}+\delta_{b,c_{1,2}^i}\frac{x_i-a}{b-a},
    \quad i\in\{1,2 \}.
	\end{align}
	Then 
 \begin{equation}
 \cT_{\Omega} f(x) = T_{\Omega}^1 f(x)-T_{\Omega}^2 f(x),  
 \end{equation}
 where 
	\begin{eqnarray}\left\{
		\begin{aligned}
			\cT_{\Omega}^1 f(x) =& \sum\limits_{ i=1}^2 \sum\limits_{ c_i^1\in\{a,b\}} \phi_{a,b}^{c_i^1}(x)f_i^{c_i^1}(x)= \sum\limits_{ i=1}^2 \left[\frac{b-x_i}{b-a}f(x)|_{x_i=a}+\frac{x_i-a}{b-a}f(x)|_{x_i=b}\right],
   \\
			\cT_{\Omega}^2 f(x)=&\sum\limits_{ c_{1,2}^1,c_{1,2}^2\in\{a,b\}}\phi_{a,b}^{c_{1,2}^1}(x)\phi_{a,b}^{c_{1,2}^2}(x)f_{1,2}^{c_{1,2}^1,c_{1,2}^2}(x) \\
   =& \frac{(b-x_1)(b-x_2)}{(b-a)^2}f(a,a)+\frac{(b-x_1)(x_2-a)}{(b-a)^2}f(a,b)+\frac{(x_1-a)(b-x_2)}{(b-a)^2}f(b,a)\\
   &+\frac{(x_1-a)(x_2-a)}{(b-a)^2}f(b,b).
		\end{aligned}
  \right.
	\end{eqnarray}
	
 For the general case, one can construct $\cT_{\Omega}f$ in a similar fashion. The domain $\Omega = [a,b]^d$ is a $d$-dimensional hypercube (referred to as $d$-cube hereafter). We need to consider all $j$-cubes contained on the boundary $\partial \Omega$, for $j=0,\cdots, d-1$. To this end, we define the following notations.
	For $j=1,\ldots, d$, we define $\cE_j$ as the collection of all $j$ tuples $(k_1,\ldots, k_j)$ with $1\leq k_1< k_2<\cdots< k_j\leq d$. The cardinality of $\cE_j$ is $\begin{bmatrix}d\\j\end{bmatrix}$.
	Hence, there are $2^j\begin{bmatrix}d\\j\end{bmatrix}$ $(d-j)$-cubes on $\partial \Omega$.  We define 
 \begin{eqnarray}\label{tuple}
 \left\{
	\begin{aligned}
&f_{\bp}^{c_{\bp}^1,\ldots,c_{\bp}^{j}}(x) = \left.f(x)\right|_{(x_{p_1} =c_{\bp}^1, x_{p_2} =c_{\bp}^2,\ldots, x_{p_j} =c_{\bp}^j)},
\quad \bp=(p_1,\dots,p_j)\in \cE_j; \\
&\phi_{a,b}^{c_\bp^i}(x) = \delta_{a,c_\bp^i}\frac{b-x_{p_i}}{b-a}
+\delta_{b,c_\bp^i}\frac{x_{p_i}-a}{b-a}, \quad 
\bp\in\cE_j, \quad i\in[j],
	\end{aligned}
 \right.
 \end{eqnarray}
 where $c_{\bp}^i\in\{a,b \}$ for $1\leqslant i\leqslant j$.
	Then $\cT_{\Omega} f(x)$ is given by 
	\begin{align}
		\cT_{\Omega} f(x) = \sum\limits_{i=1}^d 	(-1)^{i-1}\cT_{\Omega}^i f(x),
  \label{TFC_expan}
	\end{align}
	where
	\begin{align}
		\cT_{\Omega}^i f(x) = \sum\limits_{\bp \in \cE_{i}}\text{ }\sum\limits_{c_{\bp}^1,\ldots,c_{\bp}^i\in \{a,b\}}\phi_{a,b}^{c_{\bp}^1}(x)\cdots\phi_{a,b}^{c_{\bp}^i}(x)f_{\bp}^{c_{\bp}^1,\ldots,c_{\bp}^{i}}(x), \quad\forall  i\in [d].\label{eq_25}
	\end{align}

While the constrained expression~\eqref{eq_18} satisfies the condition~\eqref{eq_17} exactly, the number of terms contained therein grows exponentially with respect to the dimension $d$.
The expression~\eqref{eq_25} contains $2^i\begin{bmatrix}d\\i\end{bmatrix}$ terms, giving rise to a total number of $(3^d-1)$ terms in~\eqref{TFC_expan}. Therefore, the constrained expression~\eqref{eq_18} contains $\mathcal O(3^d)$ terms in $d$ dimensions, rendering TFC infeasible for high-dimensional PDE problems.


To devise a TFC-like approximation of the solution field suitable for high-dimensional problems, we notice that the expression~\eqref{TFC_expan} represents a hierarchical decomposition of $\cT_{\Omega}f(x)$ in some sense, in which $\cT_{\Omega}^if(x)$ represents the contributions of $f(x)$ from the $(d-i)$-dimensional hyperplanes on the boundary $\partial\Omega$. Taking $d=3$ as an example, one can note that $\cT_{\Omega}^1f$, $\cT_{\Omega}^2f$ and $\cT_{\Omega}^3f$ represent the contributions of $f(x)$ from the faces, the edges, and the vertices of the cube $\Omega=[a,b]^3$.

This observation inspires the following strategy for approximating $\cT_{\Omega}f(x)$. We can truncate the expression~\eqref{TFC_expan} by keeping only the leading terms, in a spirit analogous to the truncation in Taylor expansion. 
Specifically, by choosing a number $k\in[d-1]$ as the cut-off value, we retain all the terms $\cT_{\Omega}^jf$ with $1\leqslant j\leqslant k$ in~\eqref{TFC_expan} and discard the rest of the terms. 
In the current paper we choose $k=1$ for simplicity. In other words, for a function $f(x)$ defined on $\partial\Omega$, we define   $\cA_{\Omega}f(x)$ by
\begin{equation}\label{eq_26}
\cA_{\Omega}f(x) := \cT_{\Omega}^1f(x),
\quad x\in\Omega.
\end{equation}
Then we approximate the function $u(x)$ satisfying the condition~\eqref{eq_17} by,
\begin{equation}\label{eq_27}
u(x) = g(x) - \cA_{\Omega}g(x) + \cA_{\Omega}C(x),
\end{equation}
where $g(x)$ is a function to be determined. We refer to the approximation~\eqref{eq_27} as the approximate TFC (or A-TFC) constrained expression.


The terms $\cA_{\Omega}g$ and $\cA_{\Omega}C$ in the A-TFC expression~\eqref{eq_27} both  contain $2d$ terms. The computation of~\eqref{eq_27} is thus feasible for large $d$, and A-TFC  avoids the exponential growth in the number of terms of the TFC constrained expression~\eqref{eq_18}. However, there is a trade-off with A-TFC.  Specifically, the A-TFC expression~\eqref{eq_27} for $u(x)$ does not satisfy the boundary condition~\eqref{eq_17} unconditionally for an arbitrary function $g(x)$, because $\cA_{\Omega}f$ is  an approximation of $\cT_{\Omega}f$. In general, $g(x)$ needs to satisfy a certain condition on $\partial\Omega$ in order for $u(x)$ to satisfy~\eqref{eq_17}.

We first illustrate this point using the 2D case. Suppose $\Omega = [a,b]^2$, and we substitute the expression~\eqref{eq_27} into~\eqref{eq_17}. On the boundary $x_2=a$, we have
	\begin{eqnarray}
		\begin{aligned}
			C(x_1,a) &= \cA C(x_1,a)+g(x_1,a)-\cA g(x_1,a)\\
			&=C(x_1,a)+[C(a,a)-g(a,a)]\frac{b-x_1}{b-a}+[C(b,a)-g(b,a)]\frac{x_1-a}{b-a}.
		\end{aligned}
	\end{eqnarray}
%
%
 This leads to the following condition for $g$, 
\begin{equation}\label{eq_29}
 [C(a,a)-g(a,a)]\frac{b-x_1}{b-a}+[C(b,a)-g(b,a)]\frac{x_1-a}{b-a}=0, \quad \forall x_1\in[a,b].
 \end{equation}
Similarly, by considering the other boundaries we attain the following conditions,  
	\begin{eqnarray}
 \left\{
		\begin{aligned}
		&[C(a,b)-g(a,b)]\frac{b-x_1}{b-a}+[C(b,b)-g(b,b)]\frac{x_1-a}{b-a}=0, \quad \forall x_1\in[a,b];\\
			&[C(a,a)-g(a,a)]\frac{b-x_2}{b-a}+[C(a,b)-g(a,b)]\frac{x_2-a}{b-a}=0, \quad \forall x_2\in[a,b];\\
			&[C(b,a)-g(b,a)]\frac{b-x_2}{b-a}+[C(b,b)-g(b,b)]\frac{x_2-a}{b-a}=0, \quad \forall x_2\in[a,b].
\label{eq_30}		
  \end{aligned}
  \right.
	\end{eqnarray}
	In this 2D case, it is straightforward to choose $g(a,a) = C(a,a)$, $g(a,b) = C(a,b)$, $g(b,a) = C(b,a)$ and $g(b,b) = C(b,b)$ as the boundary conditions for $g$. However, in high-dimensional cases the conditions given by~\eqref{eq_29} and~\eqref{eq_30} are easier to implement with randomly chosen collocation points. 

For the general $d$-dimensional case,  we again use the notation in \eqref{tuple}. The conditions for $g$ are then given by,
\begin{subequations}\label{eq_31}
 \begin{align}
			&\sum\limits_{j\neq i}\left[\left(C_{i,j}^{a,a}-g_{i,j}^{a,a}\right)\frac{b-x_j}{b-a}+\left(C_{i,j}^{a,b}-g_{i,j}^{a,b}\right)\frac{x_j-a}{b-a}\right] = 0,
   \quad \text{on}\ x_i=a,
   \quad \forall i\in[d];
   \\
			&\sum\limits_{j\neq i}\left[\left(C_{i,j}^{b,a}-g_{i,j}^{b,a}\right)\frac{b-x_j}{b-a}+\left(C_{i,j}^{b,b}-g_{i,j}^{b,b}\right)\frac{x_j-a}{b-a}\right]=0,
   \quad \text{on}\ x_i=b,
   \quad \forall i\in[d].
	\end{align}
 \end{subequations}
	We define  $\tilde{\cB}_{\Omega}$ as 
	\begin{align}\label{eq_32}
		\tilde{\cB}_{\Omega}f(x) =
		\begin{bmatrix}
  \vdots\\
			\sum\limits_{j\neq i}\left(f_{i,j}^{a,a}\frac{b-x_j}{b-a}+f_{i,j}^{a,b}\frac{x_j-a}{b-a}\right)\\
			\sum\limits_{j\neq i}\left(f_{i,j}^{b,a}\frac{b-x_j}{b-a}+f_{i,j}^{b,b}\frac{x_j-a}{b-a}\right) \\
   \vdots
		\end{bmatrix},
  \quad x\in\partial\Omega.
	\end{align}
	Then one can rewrite~\eqref{eq_31} as 
	\begin{align}
		\tilde{\cB}_{\Omega}g(x) = s_{\Omega}(x),
  \quad x\in\partial\Omega,
  \label{bc_g}
	\end{align}
	where $s_{\Omega}(x):=\tilde{\cB}_{\Omega}C(x)$ denotes the boundary data. 
 This is the boundary condition for $g(x)$ with A-TFC.

We can make the following observation from the above discussions. With A-TFC, while the constrained expression does not satisfy the condition~\eqref{eq_17} exactly for an arbitrary $g$, the condition that $g$ needs to satisfy generally involves functions of a simpler  form than the original condition. Specifically, there are $2(d-1)$ terms  in the condition for $g$ on each boundary, 
and each term involves the product of a linear function and the $g$ evaluated on a $(d-2)$-dimensional hyperplane of the boundary. Therefore, A-TFC can simplify the functional forms involved in the condition in some sense. Let us again use 2D as an example to illustrate this point.  Assume a boundary distribution $C(x_1,a)=\sin(x_1)$ on the boundary $x_2=a$. Then, without A-TFC, the boundary condition for $u(x)$ on $x_2=a$ is given by, 
\begin{equation*}
u(x_1,a) = \sin(x_1), \quad \forall x_1\in[a,b].
\end{equation*}
In contrast, with A-TFC, the condition for $g$ on $x_2=a$ is reduced to (see~\eqref{eq_29}),
\begin{equation*}
 [\sin(a)-g(a,a)]\frac{b-x_1}{b-a}+[\sin(b)-g(b,a)]\frac{x_1-a}{b-a}=0, \quad \forall x_1\in[a,b].
 \end{equation*}
The linear function involved in the condition for $g$  
with A-TFC is obviously simpler than that of the original condition without A-TFC.

\subsubsection{A-TFC Embedded ELM}
\label{sec_tfc_elm}

Let us now consider how to combine ELM and A-TFC for solving the system~\eqref{problem}. In this sub-section we will assume that $\cB = \mbs I$ (identity operator) in~\eqref{eq_1b}, i.e.~the problem has Dirichlet boundary conditions.

Based on A-TFC, we perform the following transformation,
\begin{equation}\label{eq_34}
u(x) = g(x) - \cA_{\Omega}g(x) + \cA_{\Omega}H(x),
\quad x\in\Omega,
\end{equation}
where $g(x)$ is an unknown  function, $\cA_{\Omega}$ is defined in~\eqref{eq_26}, and $H$ is the boundary distribution in~\eqref{eq_1b}. The system~\eqref{problem} is accordingly transformed into,
\begin{subequations}\label{ELM_TFC_1}
 \begin{align}
		&\cL g(x)-\cL\cA_{\Omega}g(x)
   +\mu\cN(\cA_{\Omega}H(x)+g(x)-\cA_{\Omega}g(x))=Q(x)-\cL\cA_{\Omega} H(x)
   =Q_1(x),
   \quad x\in\Omega,
   \\
			&\tilde{\cB}_{\Omega}g(x) = 
   \tilde{\cB}_{\Omega}H(x) = S_{\Omega}(x),
   \quad x\in\partial\Omega,
\end{align}
\end{subequations}
where $\tilde{\cB}_{\Omega}$ is defined in~\eqref{eq_32}. $g(x)$ is the field function to be solved for in~\eqref{ELM_TFC_1}.


We represent $g(x)$ by an ELM-type randomized neural network, following the settings as outlined in Section~\ref{sec_elm}.
In particular, the input layer of the network contains $d$ nodes (representing $x$), and the output layer contains a single node (representing $g(x)$) with zero bias and no activation function. The hidden-layer coefficients are set to uniform random values from the interval $[-R_m,R_m]$. Let $N_g$ denote the number of nodes in the last hidden layer of the neural network, and $V_j(x)$ ($1\leqslant j\leqslant N_g$) denote the output fields of the last hidden layer. Then the network logic gives rise to,
\begin{align}\label{eq_36}
g(x) = \sum_{j=1}^{N_g}\phi_j V_j(x) 
= \mbs V(x)\bm\Phi,
\end{align}
where $\mbs V(x)=(V_1(x),\dots,V_{N_g}(x))$, $\phi_j$ ($1\leqslant j\leqslant N_g$) are the output-layer coefficients (training parameters), and $\bm\Phi=(\phi_1,\dots,\phi_{N_g})^T$.

We determine the training parameters $\bm\Phi$ in~\eqref{eq_36} in a  fashion similar to in Section~\ref{sec_elm}. We again employ random collocation points inside the domain and on the domain boundaries, and let $N_{\ins}$ and $N_{\bcd}$ denote the number of interior collocation points and the number of boundary collocation points on each hyperface,  respectively. Enforcing the residual function of the system~\eqref{ELM_TFC_1} to be zero on these collocation points leads to the following algebraic system,
 \begin{equation}\label{eq_37}
 \mbs R(\bm\Phi) =
 \begin{bmatrix}
		\cL (\mbs V(\mbs x)-\cA_{\Omega}\mbs V(\mbs x))\bm\Phi
   +\mu\cN(\cA_{\Omega}H(\mbs x)+\mbs V(\mbs x)\bm\Phi-\cA_{\Omega}\mbs V(\mbs x)\bm\Phi)-Q_1(\mbs x)
   \\
	\tilde{\cB}_{\Omega}\mbs V(\mbs y)\bm\Phi - S_{\Omega}(\mbs y)
   \end{bmatrix}
   =\mbs 0,
\end{equation}
where $\mbs x$ and $\mbs y$ are defined in~\eqref{eq_11}, and $\mbs R(\bm\Phi)\in\mbb R^{N_a}$ denotes the residual vector of the system~\eqref{ELM_TFC_1} on the collocation points, where $N_a=N_{\ins}+2N_{\bcd}^{\text{tot}} = N_{\ins}+4dN_{\bcd}$.

We seek a least squares solution to the system~\eqref{eq_37}, and solve this system for $\bm\Phi$ by the linear least squares method, if $\mu=0$, and by the nonlinear least squares method with perturbations (NLLSQ-perturb)~\cite{DongL2021}, if $\mu\neq 0$. 
When $\mu=0$, the system~\eqref{eq_37} is reduced to,
\begin{equation}
\begin{bmatrix}
\cL (\mbs V(\mbs x)-\cA_{\Omega}\mbs V(\mbs x)) \\
\tilde{\cB}_{\Omega}\mbs V(\mbs y)
\end{bmatrix}\bm\Phi
=\begin{bmatrix}
Q_1(\mbs x)\\ S_{\Omega}(\mbs y)
\end{bmatrix},
\end{equation}
which is linear, and $\bm\Phi$ can be computed by the linear least squares method. When $\mu\neq 0$, the system~\eqref{eq_37} is nonlinear. The Jacobian matrix of this system is given by,
\begin{equation}
\mbs J(\bm\Phi) = \frac{\partial\mbs R}{\partial\bm\Phi}=\begin{bmatrix}
\cL (\mbs V(\mbs x)-\cA_{\Omega}\mbs V(\mbs x))
   +\mu\cN'(\cA_{\Omega}H(\mbs x)+\mbs V(\mbs x)\bm\Phi-\cA_{\Omega}\mbs V(\mbs x)\bm\Phi)(\mbs V(\mbs x) - \cA_{\Omega}\mbs V(\mbs x))\\
   \tilde{\cB}_{\Omega}\mbs V(\mbs y)
\end{bmatrix},
\end{equation}
which is needed by the NLLSQ-perturb method~\cite{DongL2021} for solving~\eqref{eq_37}. 

Upon solving the system~\eqref{eq_37} by the linear or nonlinear least squares method, we set the output-layer coefficients of the neural network by the least squares solution $\bm\Phi$. Afterwards, a forward evaluation of the neural network provides $g(x)$, and the solution field to  system~\eqref{problem} is computed by~\eqref{eq_34}.

\section{Numerical Examples}
\label{sec:tests}

In this section we test the performance of the proposed methods using several high-dimensional PDEs. These include the linear and nonlinear Poisson equations, which are time-independent, and the heat, Korteweg-de Veries (KdV), and the advection diffusion equations, which are time-dependent.  For each problem we investigate the error convergence with respect to the number of training parameters and number of  collocation points for a range of problem dimensions. 

A single hidden layer has been employed in the neural network with both the ELM and the ELM/A-TFC methods for all  test problems.
In  subsequent discussions we denote the neural network architecture by the following vector (or list) of positive integers, referred to as the architectural vector henceforth,
\begin{equation}
\mbs M_{\text{arch}} = [m_{\ins}, M, m_{\text{out}}]
\end{equation}
where $m_{\ins}$, $M$ and $m_{\text{out}}$ are the number of nodes in the input layer, the hidden layer, and the output layer, respectively. $m_{\text{out}}=1$ in all tests of this section, and $M$ equals the number of training parameters of ELM.
The hidden-layer coefficients  are assigned to uniform random values generated on the interval $[-R_m,R_m]$ with $R_m=R_{m0}$ in the simulations below, where $R_{m0}$ is determined by the procedure discussed in Remark~\ref{rem_a6}. The $\tanh$ activation function has been employed in the hidden layer for all the numerical tests of this section.

In all numerical experiments,  after the NN with a particular setting (given number of training parameters, and given number of training collocation points) is trained, we compute the errors of the  network solution as follows. We generate a set of test points on the $d$-dimensional domain $\Omega$: $N_{\ins}^{(v)}$ random points on the interior of $\Omega$, and $N_{\bcd}^{(v)}$ random points on each boundary of $\Omega$ (($d-1$)-dimensional hypercube). There is a total of $N_c^{(v)}=N_{\ins}^{(v)}+2dN_{\bcd}^{(v)}$ random test points. They  are different from the random collocation points used in the network training, and the  number is much larger than that of the latter. We evaluate the trained NN on these test points to obtain the NN solution data, and evaluate the exact solution to the problem on the same test points. We compare the data of the NN solution and the exact solution on these test points to compute the maximum error ($e_{l^{\infty}}$) and the root-mean-squares (rms) error ($e_{l^2}$) as follows,
\begin{equation}
e_{l^{\infty}} = \max\left\{\ |u(\mbs x_i) - u_{ex}(\mbs x_i)| \ \right\}_{i=1}^{N_c^{(v)}},
\quad
e_{l^2} = \sqrt{\frac{1}{N_c^{(v)}}
\sum_{i=1}^{N_c^{(v)}} \left[u(\mbs x_i) - u_{ex}(\mbs x_i) \right]^2
},
\end{equation}
where $\mbs x_i$ ($1\leqslant i\leqslant N_c^{(v)}$) denote the test points, and $u(\mbs x_i)$ and $u_{ex}(\mbs x_i)$ denote the NN solution and the exact solution, respectively.
We refer to the computed $e_{l^{\infty}}$ and $e_{l^2}$ as the  errors associated with the given NN setting  that is used for training the network.
In all the numerical simulations of this section we employ a fixed $N_{\bcd}^{(v)}=100$ and $N_{\ins}^{(v)}=7000$ for computing the  $e_{l^{\infty}}$ and $e_{l^2}$ errors.



In the numerical tests below, $d$ refers to the dimension of the spatial domain. For time-dependent problems, the time variable $t$ is not included in $d$ when we talk about the problem in $d$ dimensions. In other words, for a $d$-dimensional time-dependent problem, the input layer of the neural network would contain $m_{\ins}=d+1$ nodes (space and time). To make the reported results in this section exactly re-producible, we have set the seed to a fixed value $1$ in the random number generators of the Numpy and Tensorflow libraries for  all the numerical examples. The test results with the ELM method are presented in Section~\ref{elm_test} first, and those obtained with the combined ELM/A-TFC method are  discussed in Section~\ref{atfc_test}.
We then compare the current ELM method with the PINN method~\cite{RaissiPK2019} for selected test problems in Section~\ref{sec_pinn}.


\subsection{Numerical Tests with the ELM Method }
\label{elm_test}

\subsubsection{Poisson Equation}
\label{sec_poisson}

\begin{table}[tb]
\centering\small
\begin{tabular}{ccccccc}
	\hline
	$R_m$ &0.01&0.05&0.1&0.15&0.25&0.5 \\ \hline 
	$e_{l^\infty}$ &1.76E-6&1.30E-8&1.04E-7&1.23E-6& 4.44E-5&8.33E-4 \\ 
	$e_{l^{2}}$ &1.96E-7&6.68E-10&6.91E-9&8.48E-8&2.28E-6&9.03E-5 \\ 
	\hline
\end{tabular}
\caption{\small Poisson equation: Determining the $R_{m0}$ based on the procedure of Remark~\ref{rem_a6} for $d=5$. 
NN architecture: $[d,2000,1]$; $(N_{\bcd},N_{\ins})=(100,1000)$.
}
\label{tab:rm_5}
\end{table}

\begin{table}
\centering\small
\begin{tabular}{l| ccc}
\hline
$d$ & 5 & 7 & 15 \\ \hline
$R_{m0}$ & 0.05 & 0.05 & 0.001 \\
\hline
\end{tabular}
\caption{\small Poisson equation: $R_{m0}$ determined by the procedure from Remark~\ref{rem_a6} for several problem dimensions.
}
\label{tab_2}
\end{table}

We  consider the Poisson equation on the domain $\Omega=[-1,1]^d$, 
\begin{subequations}\label{eq_42}
\begin{align}
	-\Delta u&=f(x), \quad  x\in \Omega,\\
	u&=h(x), \quad x\in \partial \Omega,
\end{align}
\end{subequations}
where $\Delta = \sum_{i=1}^d\frac{\partial^2}{\partial x_i^2}$,   $h(x)=\left(\frac{1}{d}\sum_{i=1}^{d} x_i\right)^2+\sin\left(\frac{1}{d}\sum_{i=1}^{d} x_i\right)$, and $f(x)=\frac{1}{d}\left(\sin\left(\frac{1}{d}\sum_{i=1}^{d} x_i\right)-2\right)$. The exact solution to this system is
$	u(x)=\left(\frac{1}{d}\sum_{i=1}^{d} x_i\right)^2+\sin\left(  \frac{1}{d}\sum_{i=1}^{d} x_i\right).
$

We solve the system~\eqref{eq_42} by the ELM method with an NN architecture $\mbs M_{\text{arch}}=[d,M,1]$, where the $d$ input nodes represent $x$, the single output node represents $u(x)$, and the hidden-layer width $M$ (i.e.~number of training parameters) is varied systematically. 
The neural network is trained by the algorithm from Section~\ref{sec_elm} on a set of random collocation points, consisting of $N_{\ins}$ points from the interior of $\Omega$ and $N_{\bcd}$ points on each of the $2d$ boundaries of $\partial\Omega$, where $(N_{\bcd},N_{\ins})$ are varied in the tests. For the network training, the input data consists of the coordinates of  all the collocation points.
After the NN is trained for each case, as discussed previously, the maximum and rms errors of the NN solution are computed on another set of random test points, characterized by $(N_{\bcd}^{(v)},N_{\ins}^{(v)})=(100,7000)$ for the test points on each boundary and in the interior of $\Omega$. In addition, on selected 2D cross sections of $\Omega$, such as the $x_i$-$x_j$ plane ($1\leqslant i<j\leqslant d$), we have evaluated the network solution on a set of regular $Q^{(v)}\times Q^{(v)}$ grid points ($Q^{(v)}=800$), and compared with the exact solution on the same set of points to study the point-wise  errors of the NN solution.

We first determine the $R_{m0}$ based on the procedure from Remark~\ref{rem_a6} for generating the random hidden-layer coefficients.
Table \ref{tab:rm_5} shows the maximum and rms errors of the ELM solution for dimension $d=5$, corresponding to several $R_m$ values for generating the random hidden-layer coefficients. These are obtained using a network architecture $[d,2000,1]$ and the collocation points $(N_{\bcd},N_{\ins})=(100,1000)$. The ELM method produces accurate results in a range of $R_m$ values around $R_{m}\approx 0.05$, thus leading to  $R_{m0}=0.05$ for $d=5$.
Table~\ref{tab_2} lists the $R_{m0}$ for several problem dimensions with the Poisson equation, which are obtained with the same network architecture and the same number of collocation points as in Table~\ref{tab:rm_5}. We observe that $R_{m0}$ tends to decrease as the  dimension $d$ increases. This seems to be a  characteristic common to all the test problems  in this work.
In the subsequent  simulations, we set $R_m=R_{m0}$ in ELM for generating the hidden-layer coefficients, while the other simulation parameters are varied.


\begin{figure}[tb]
	\centering
	\subfigure[d=5]
	{ \includegraphics[width=0.3\textwidth]{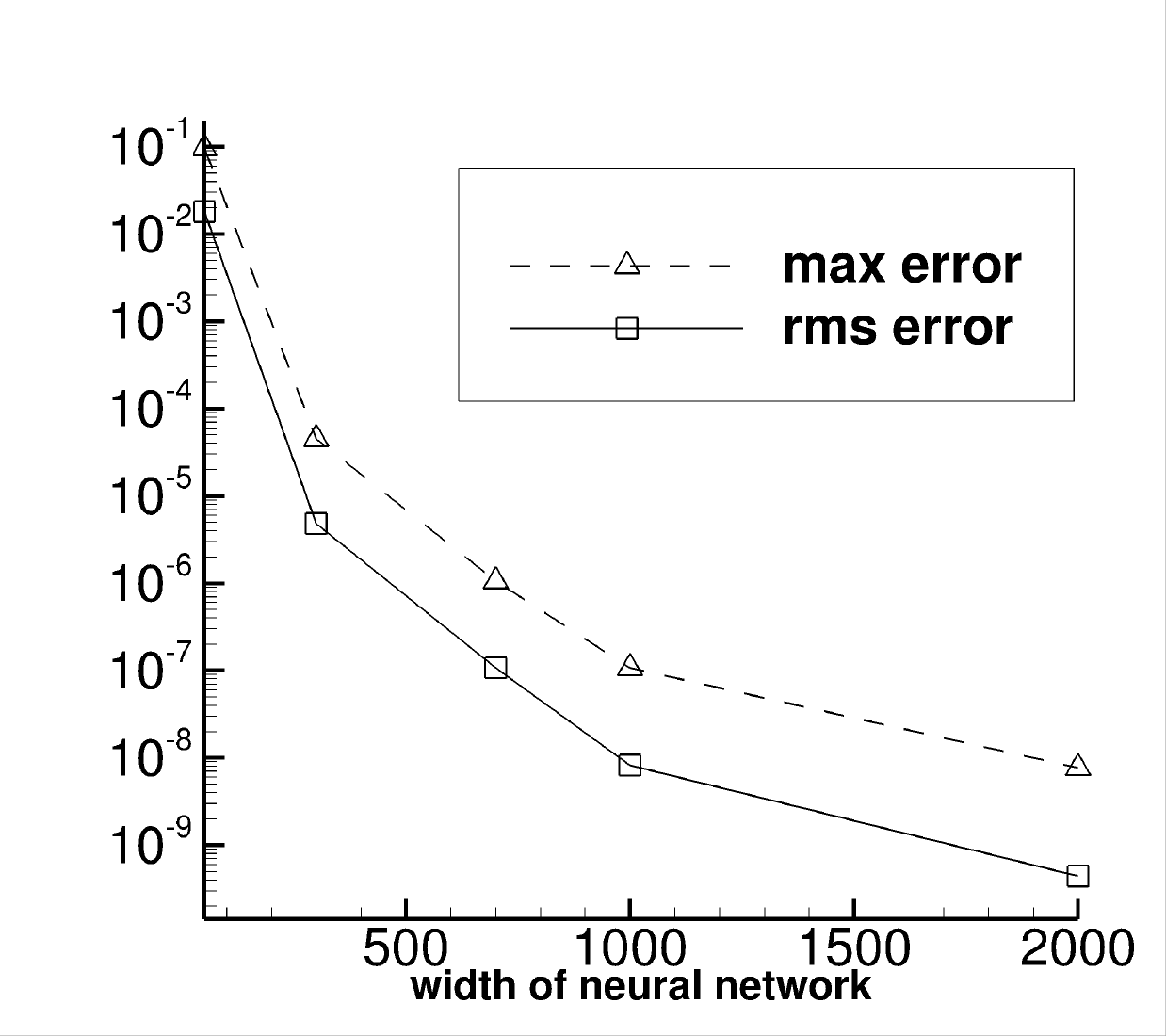}}
	\subfigure[d=7]
	{ \includegraphics[width=0.3\textwidth]{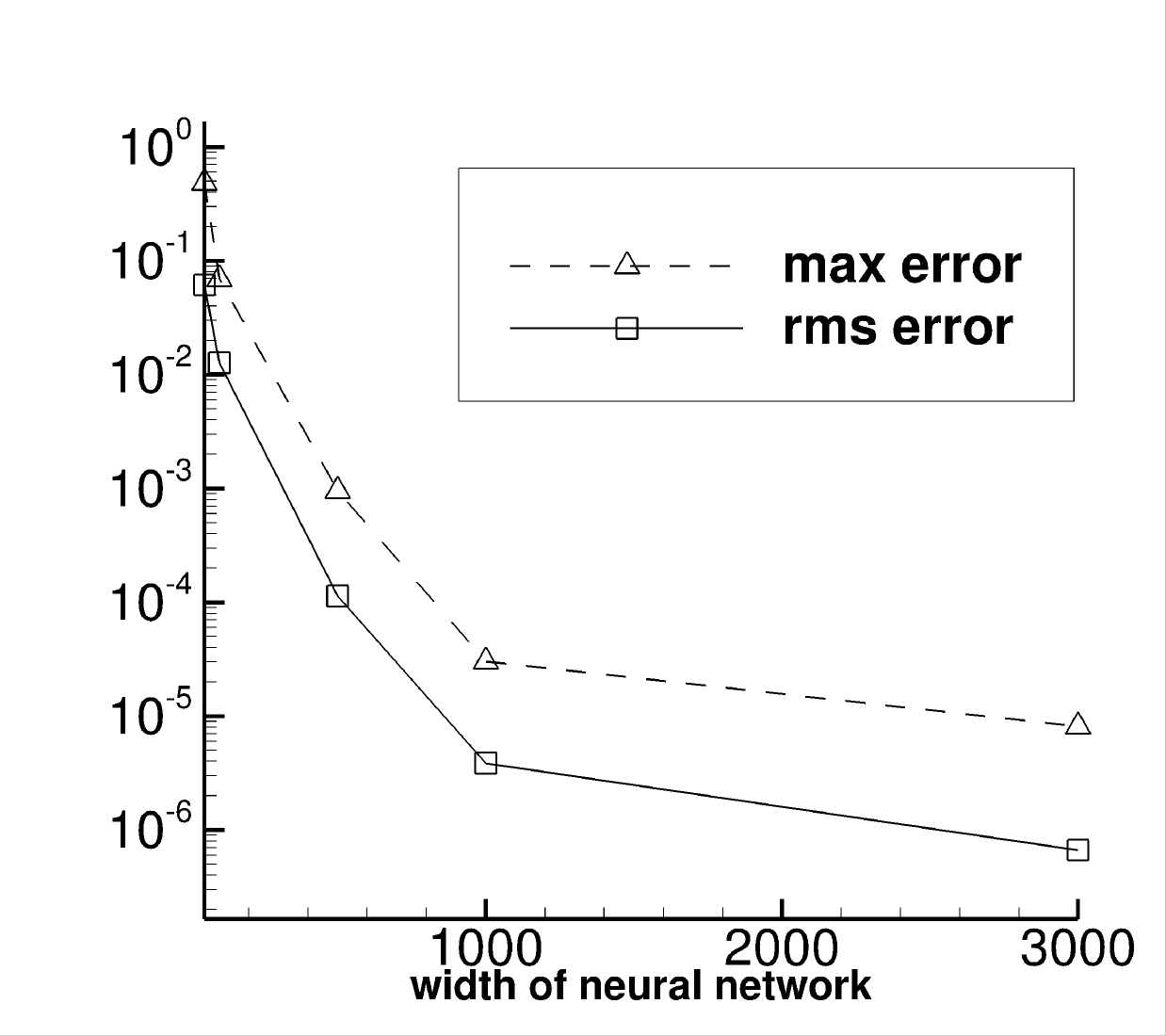}}
	\subfigure[d=15]
	{ \includegraphics[width=0.3\textwidth]{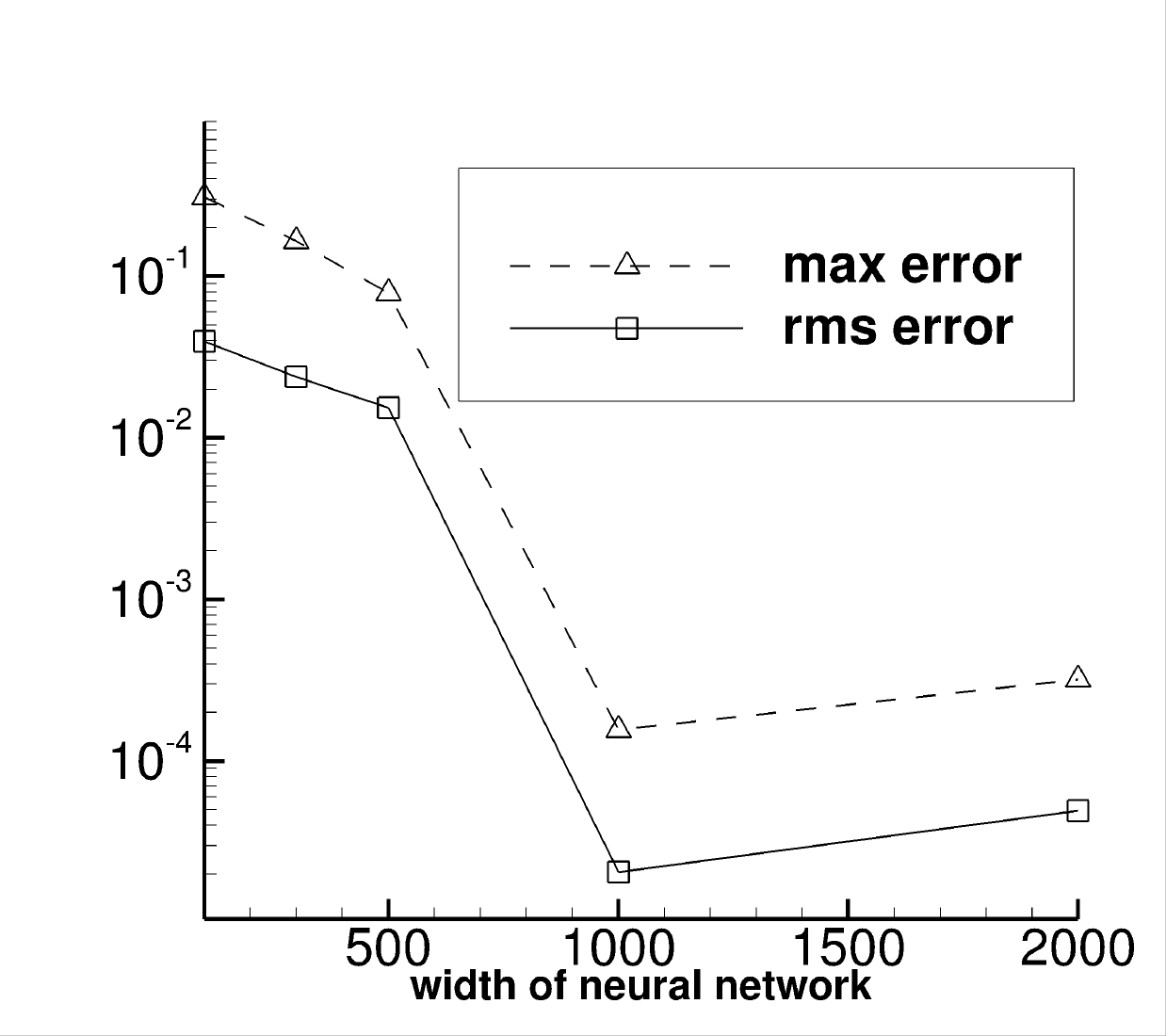}}
	\caption{\small Poisson equation: $e_{l^{\infty}}$ and $e_{l^2}$ errors versus the number of training parameters ($M$) for dimensions (a) $d=5$, (b) $d=7$, (c) $d=15$.  
 $(N_{\bcd},N_{\ins})=(100,2000)$ in (a), $(100,1000)$ in (b), and $(120,2000)$ in (c).
 }
 \label{poisson_err_width}
\end{figure}

Figure \ref{poisson_err_width} illustrates the effect of the number of training parameters  on the ELM accuracy. It shows the maximum and rms errors of ELM versus the number of training parameters ($M$) for the Poisson equation in $d=5$, $7$ and $15$ dimensions. The other crucial simulation parameters are listed in the figure caption. It is observed that the errors decrease quasi-exponentially with increasing number of training parameters (when $M\lesssim 1000$), but appears to stagnate at a certain level as $M$ increases further. The stagnant error level tends to be larger with a higher  dimension. For example, the $e_{l^2}$ error level is on the order of $10^{-10}$, $10^{-7}$ and $10^{-5}$ for the dimensions $d=5$, $7$ and $15$, respectively.


\begin{figure}[tb]
	\centering
	\subfigure[d=5]
	{ \includegraphics[width=0.3\textwidth]{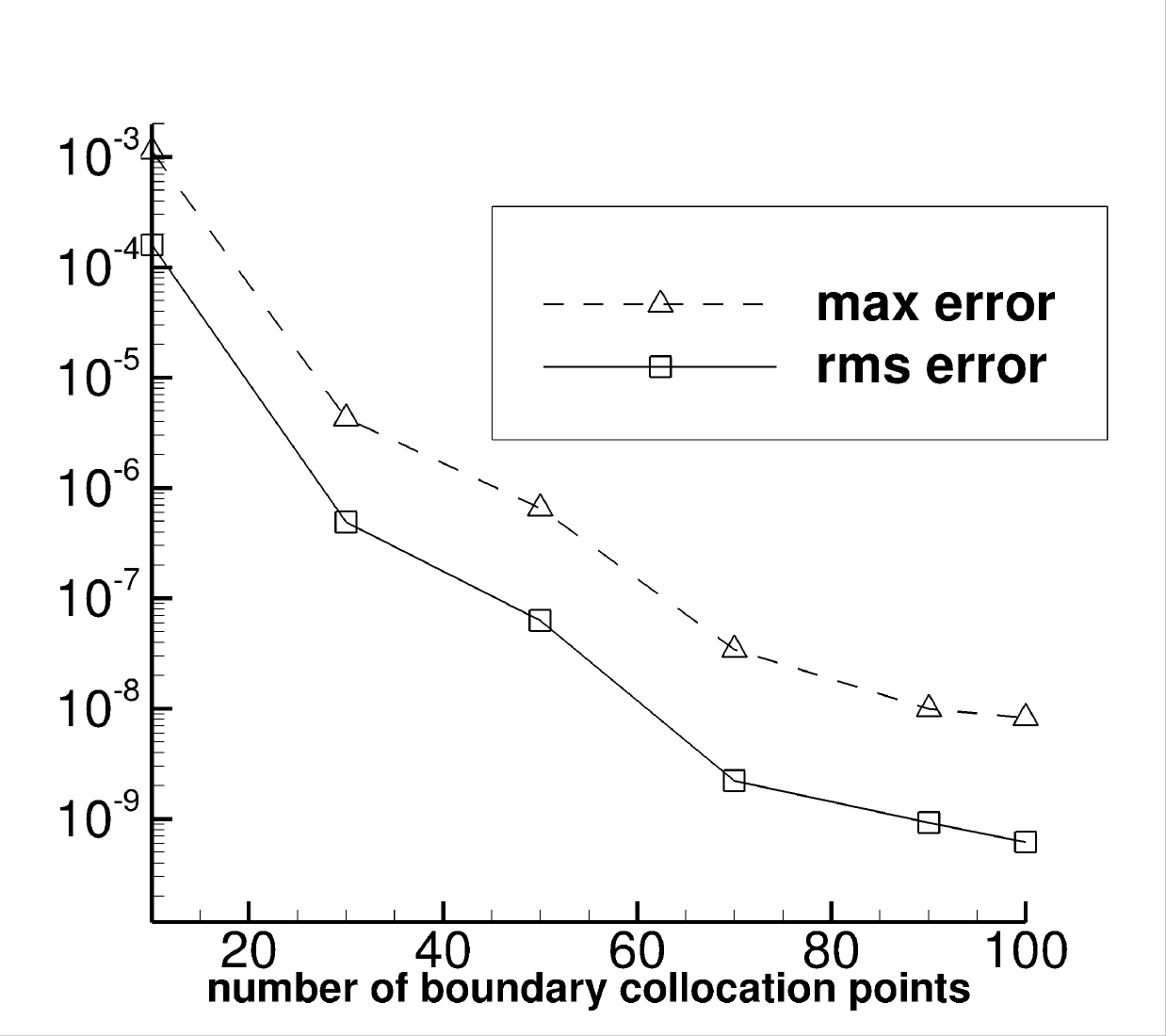}}
	\subfigure[d=7]
	{ \includegraphics[width=0.3\textwidth]{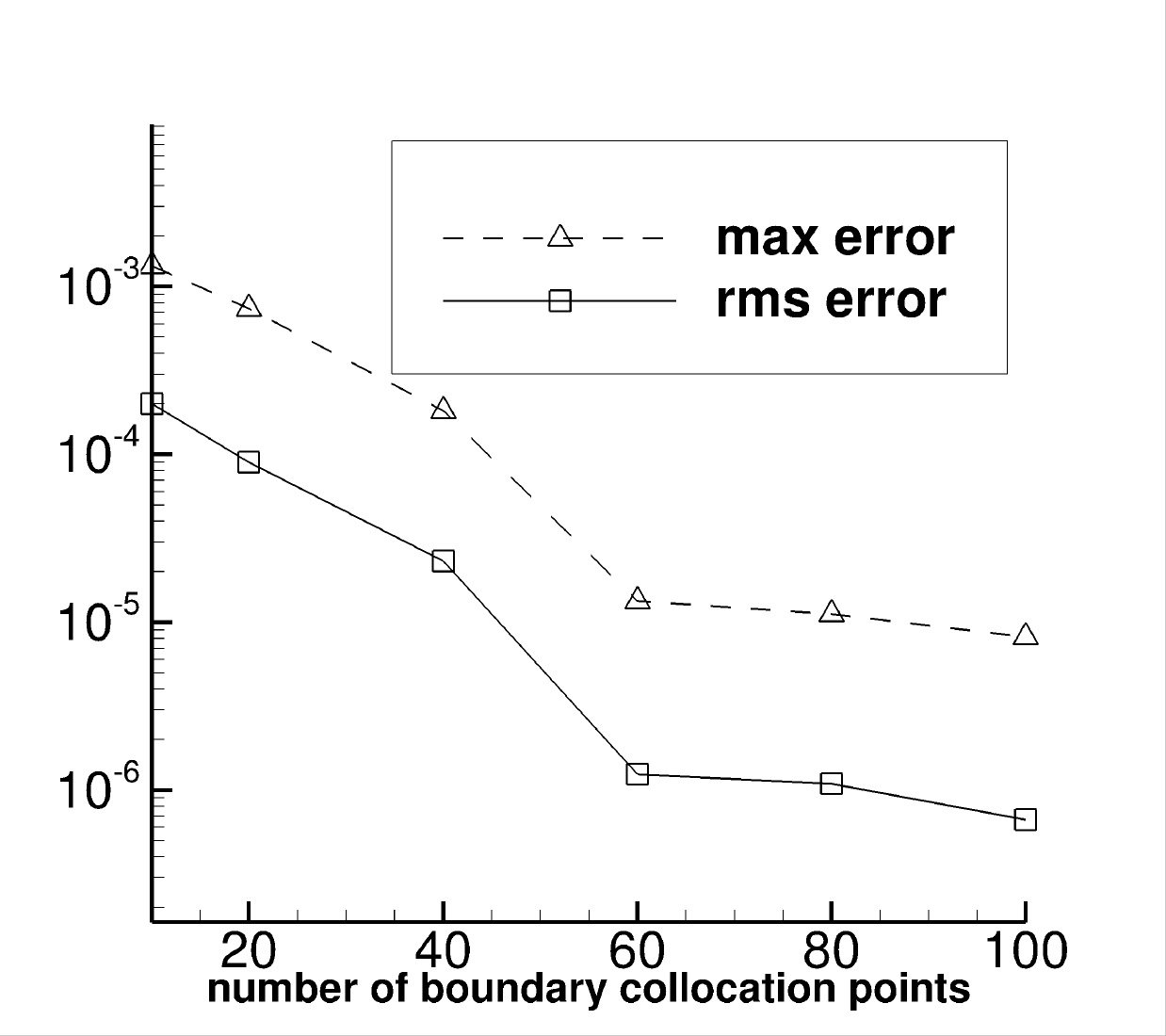}}
		\subfigure[d=15]
	{\includegraphics[width=0.3\textwidth]{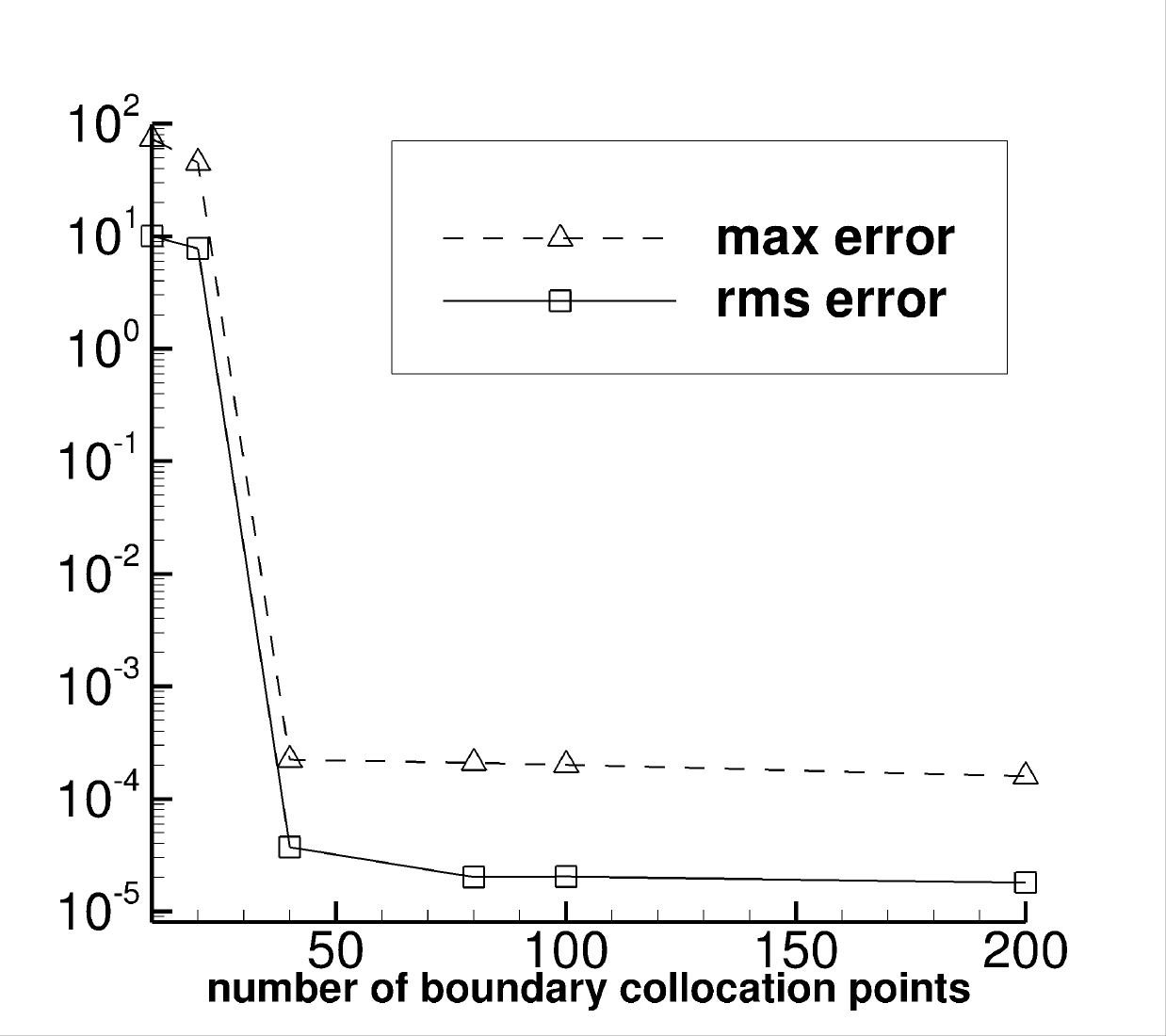}}
 \caption{\small Poisson equation: $e_{l^{\infty}}$ and $e_{l^2}$ errors versus the number of boundary collocation points ($N_{\bcd}$) for dimensions (a) $d=5$, (b) $d=7$, and (c) $d=15$. NN architecture: $[d,2000,1]$ in (a,b,c); $N_{\ins} = 500$ in (a), $1000$ in (b), and $50$ in (c).
 }
 \label{fg_3}
 \end{figure}
 
 \begin{figure}[tb]
			\subfigure[d=5]
	{ \includegraphics[width=0.3\textwidth]{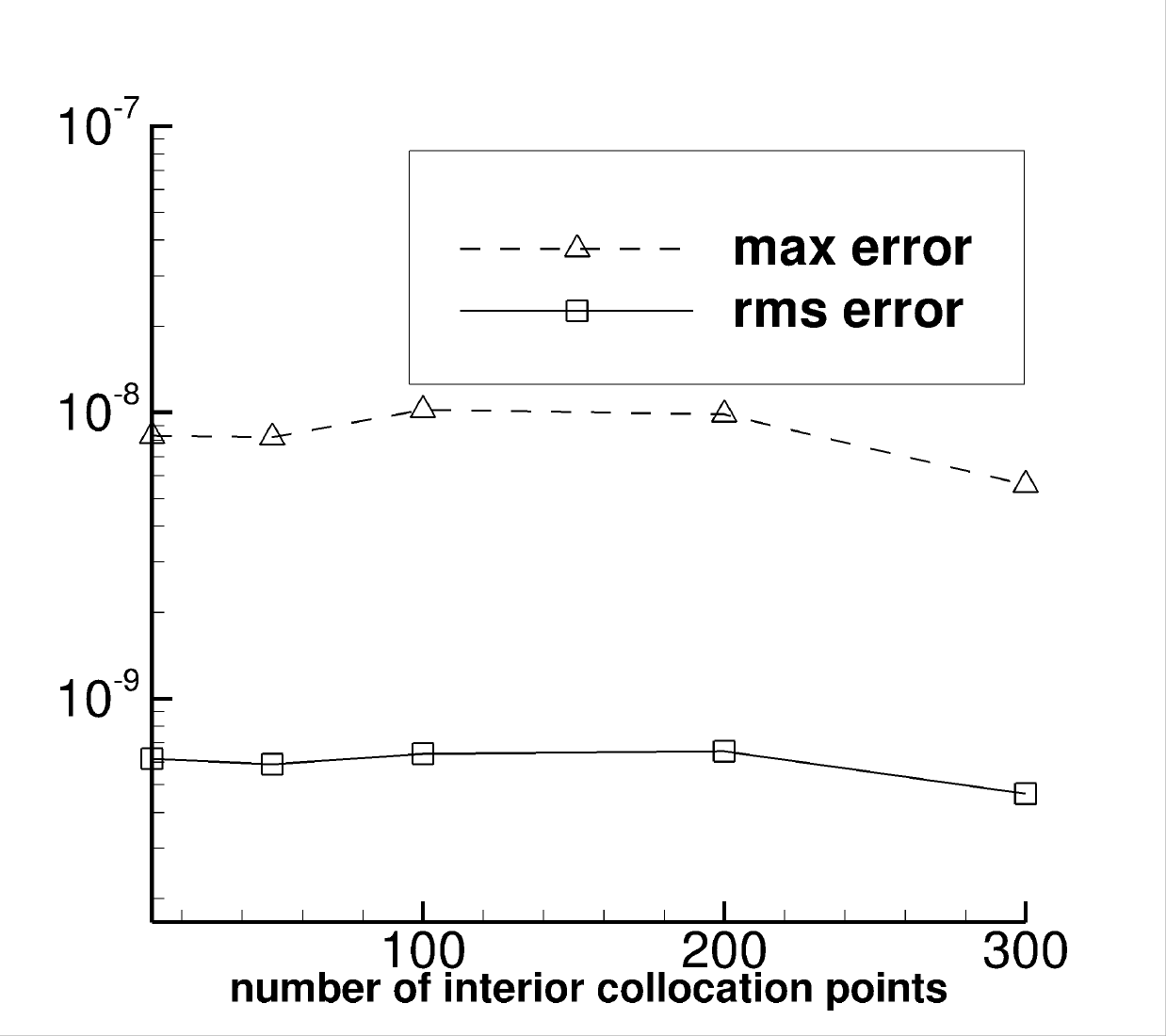}}
	\subfigure[d=7]
	{ \includegraphics[width=0.3\textwidth]{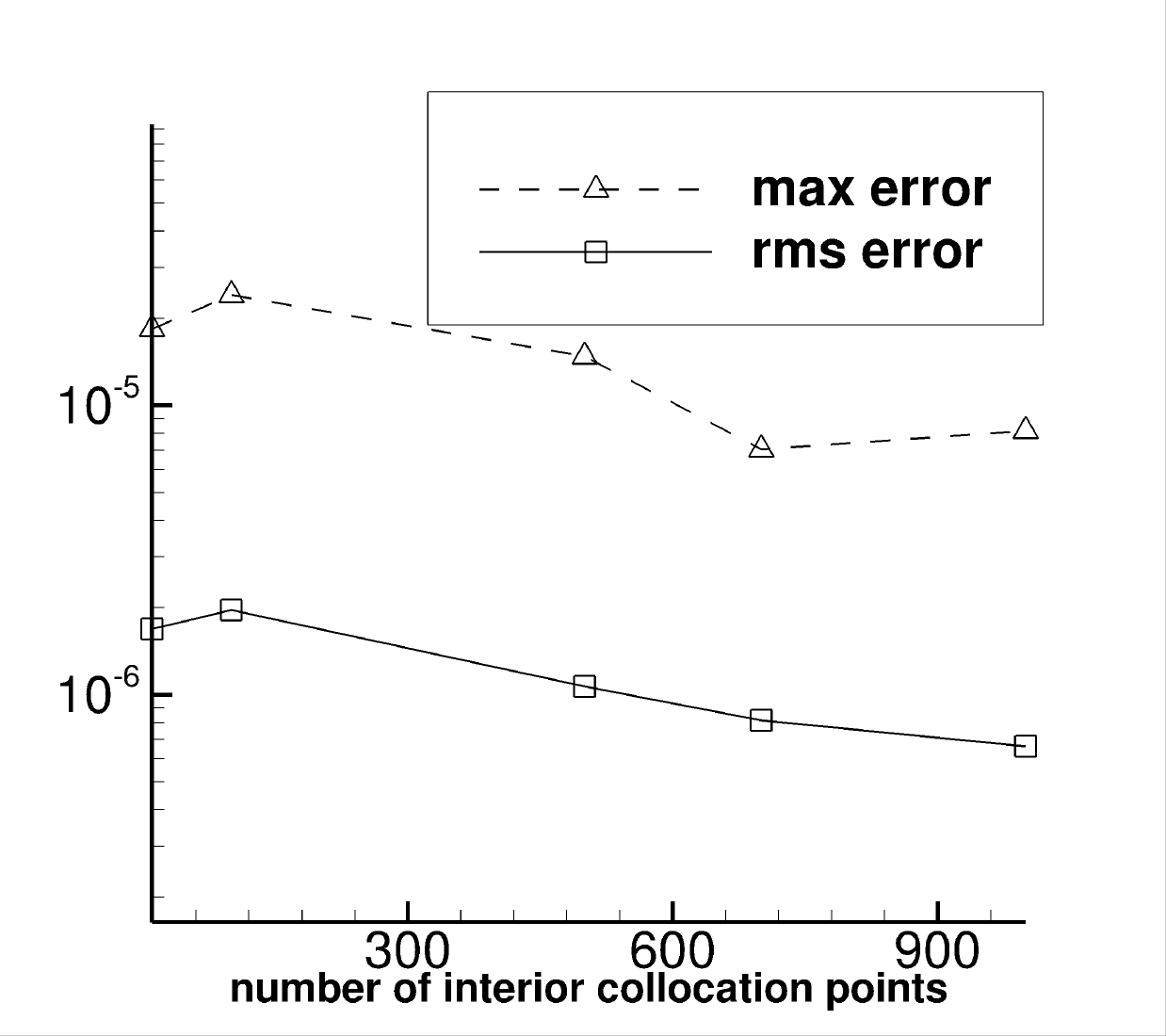}}
	\subfigure[d=15]
	{ \includegraphics[width=0.3\textwidth]{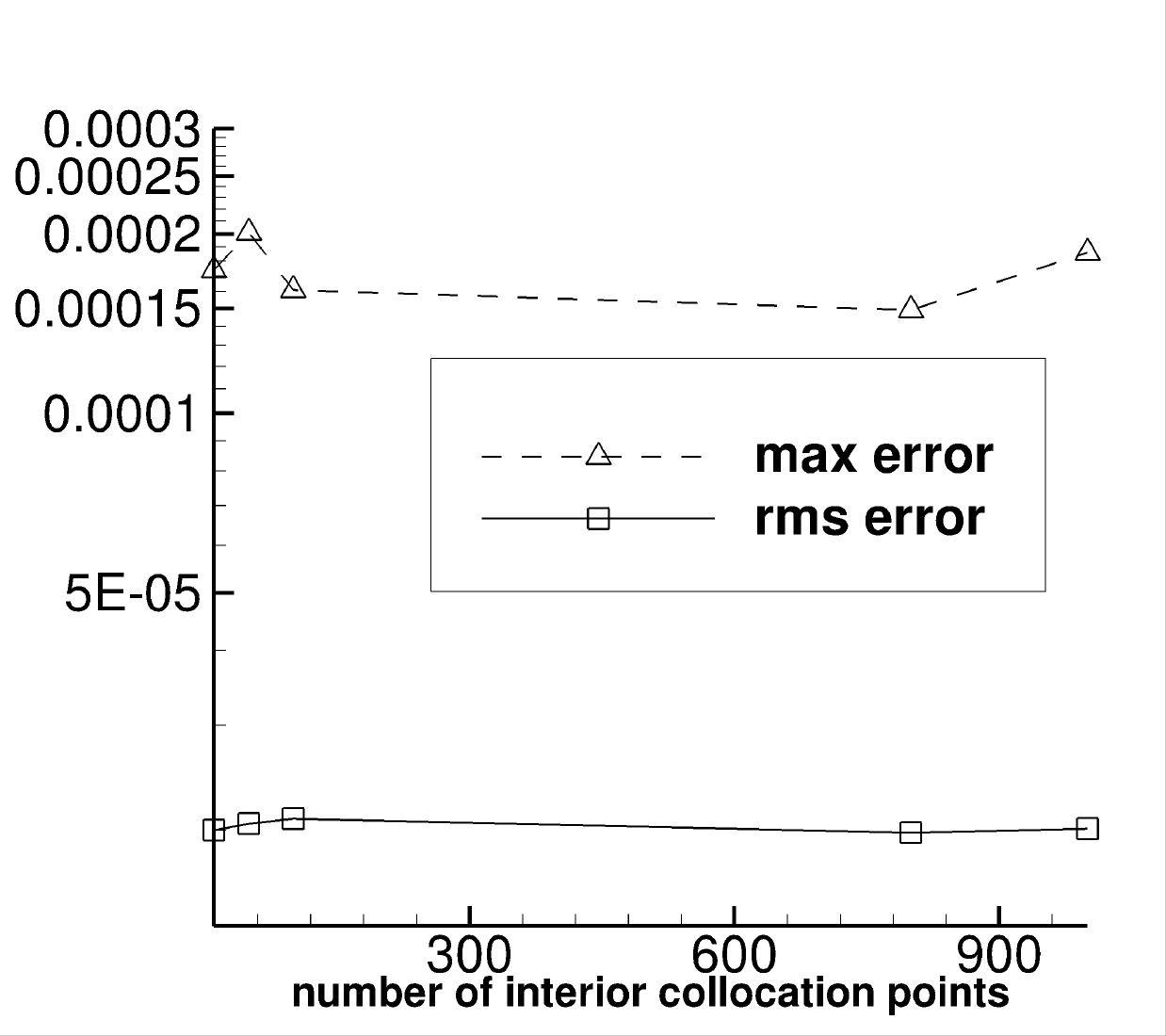}}
	\caption{\small Poisson equation: $e_{l^{\infty}}$ and $e_{l^2}$ errors versus the number of interior collocation points ($N_{\ins}$) for dimensions (a) $d=5$, (b) $d=7$, and (c) $d=15$. NN architecture: $[d,2000,1]$ in (a,b,c); $N_{\bcd} = 100$ in (a,b,c).
 }
 \label{poisson_err_nc}
\end{figure}

Figures~\ref{fg_3} and~\ref{poisson_err_nc} illustrate the effect of the training collocation points on the ELM accuracy. Figure~\ref{fg_3} shows the ELM errors ($e_{l^{\infty}}$ and $e_{l^2}$) as a function of the number of  collocation points on each boundary ($N_{\bcd}$) for dimensions $d=5$, $7$ and $15$. Figure~\ref{poisson_err_nc} shows the ELM errors as a function of the number of interior collocation points ($N_{\ins}$). The other simulation parameters are fixed in the tests and their values are provided in the captions of these figures. Increasing the number of boundary collocation points ($N_{\bcd}$) improves the ELM accuracy significantly. The ELM maximum/rms  errors decrease approximately exponentially with increasing $N_{\bcd}$ for $d=5$, and also for $d=7$ when $N_{\bcd}\lesssim 60$ or $d=15$ when $N_{\bcd}\lesssim 40$. The errors stagnate when $N_{\bcd}$ increases beyond around $60$ for $d=7$ and beyond $40$ for $d=15$.
On the other hand, varying the number of interior collocation points $N_{\ins}$ appears to have little effect on the ELM accuracy for all three  dimensions, which is evident from Figure~\ref{poisson_err_nc}. With increasing problem dimension, the surface of the hypercube (and hence the boundary collocation points) becomes more dominant, while the interior (hence the interior collocation points) becomes less important. Therefore, a small number of interior collocation points will typically suffice in higher dimensions.

 \begin{figure}[!htb]
	\centerline{
	 \includegraphics[width=0.26\textwidth]{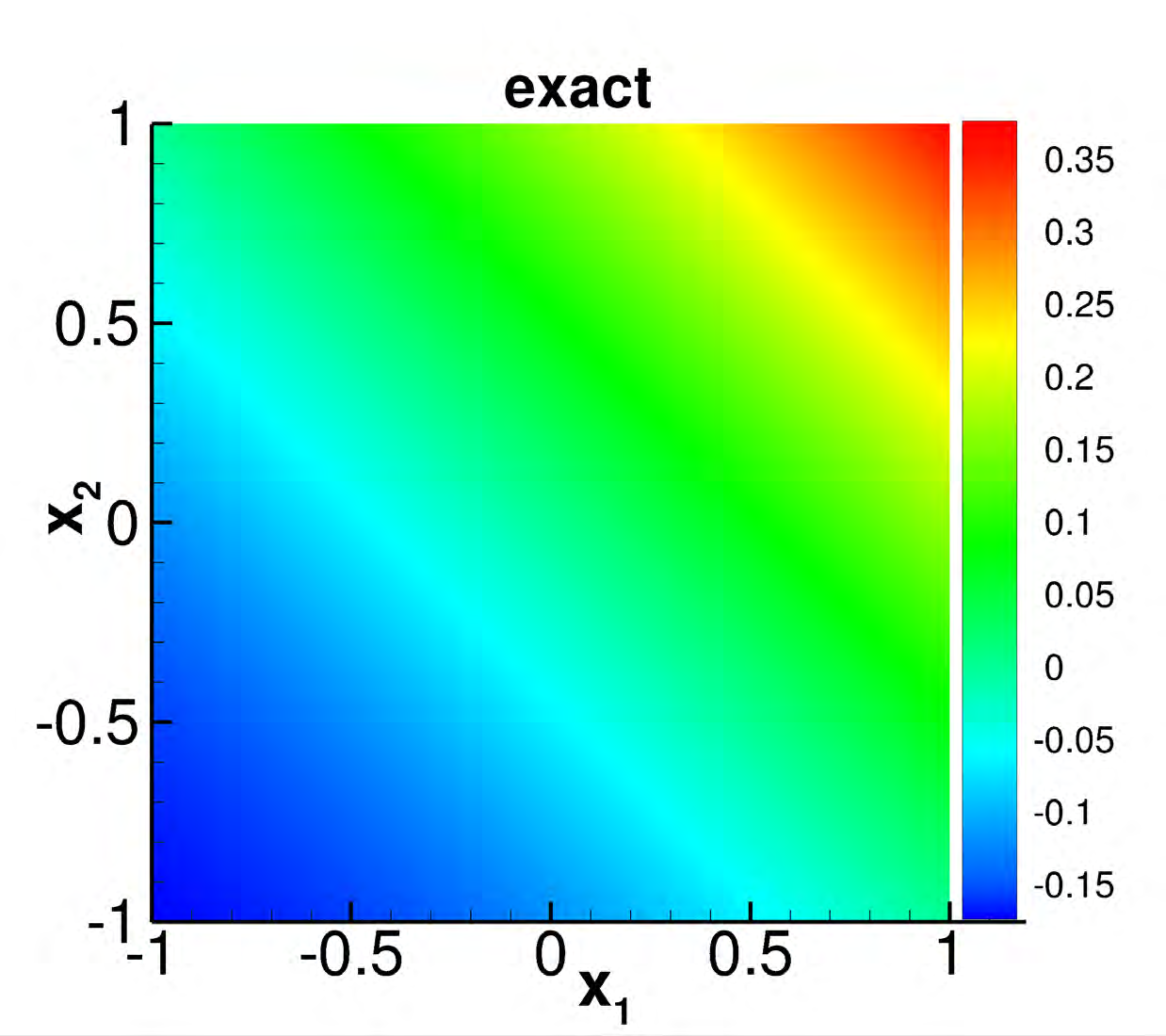}(a)
	 \includegraphics[width=0.26\textwidth]{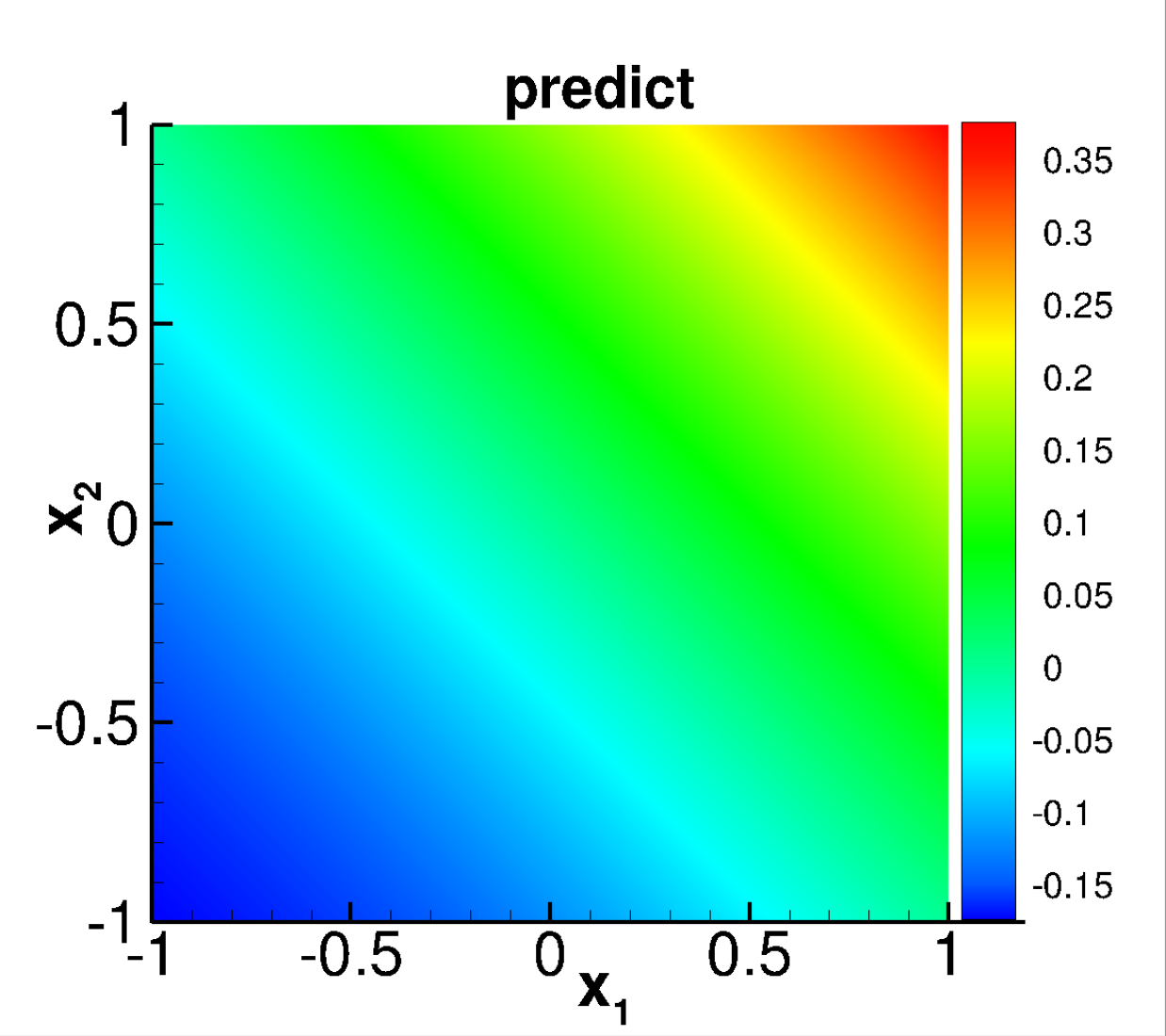}(b)
	 \includegraphics[width=0.26\textwidth]{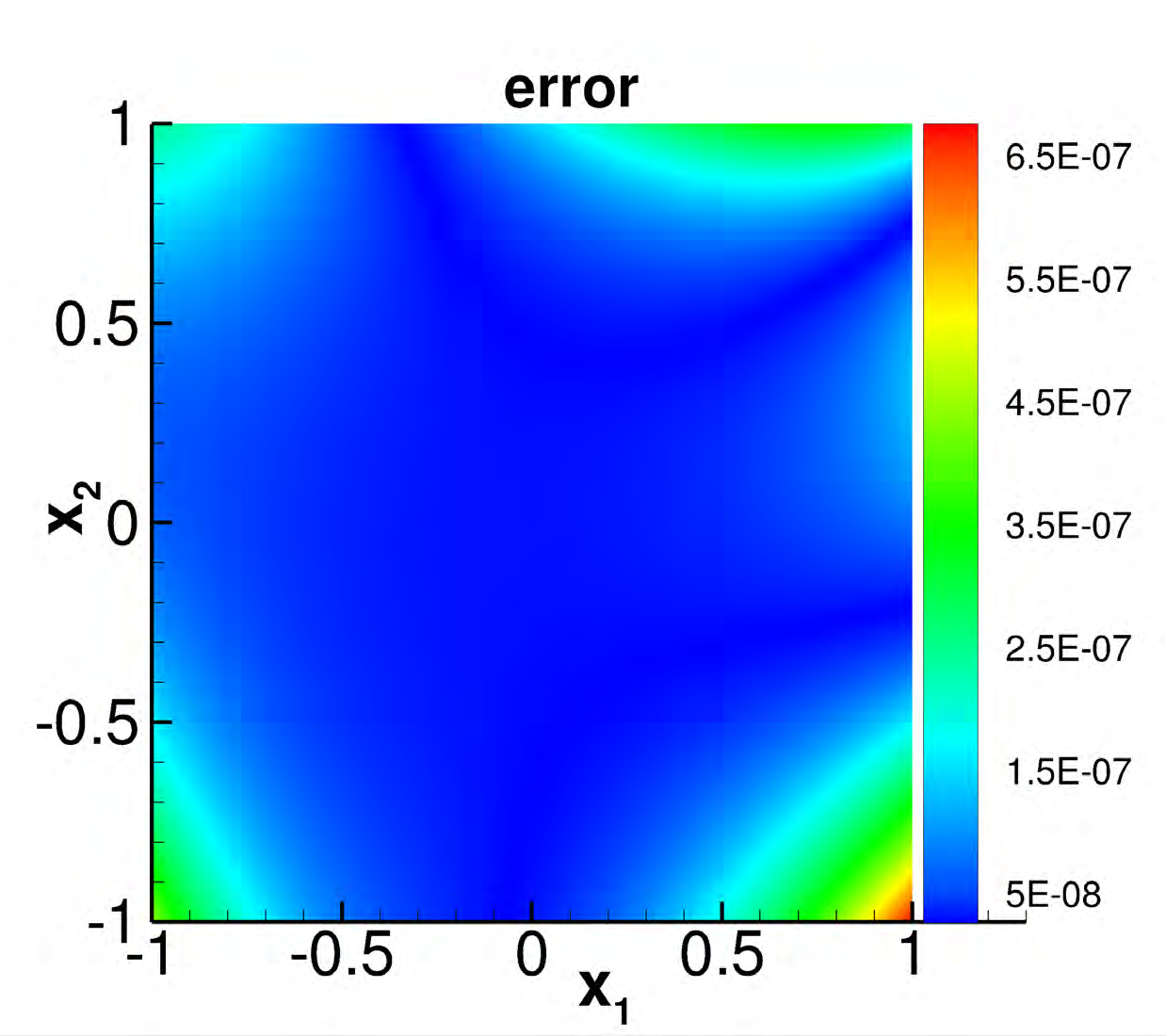}(c)
 }
 \centerline{
	 \includegraphics[width=0.26\textwidth]{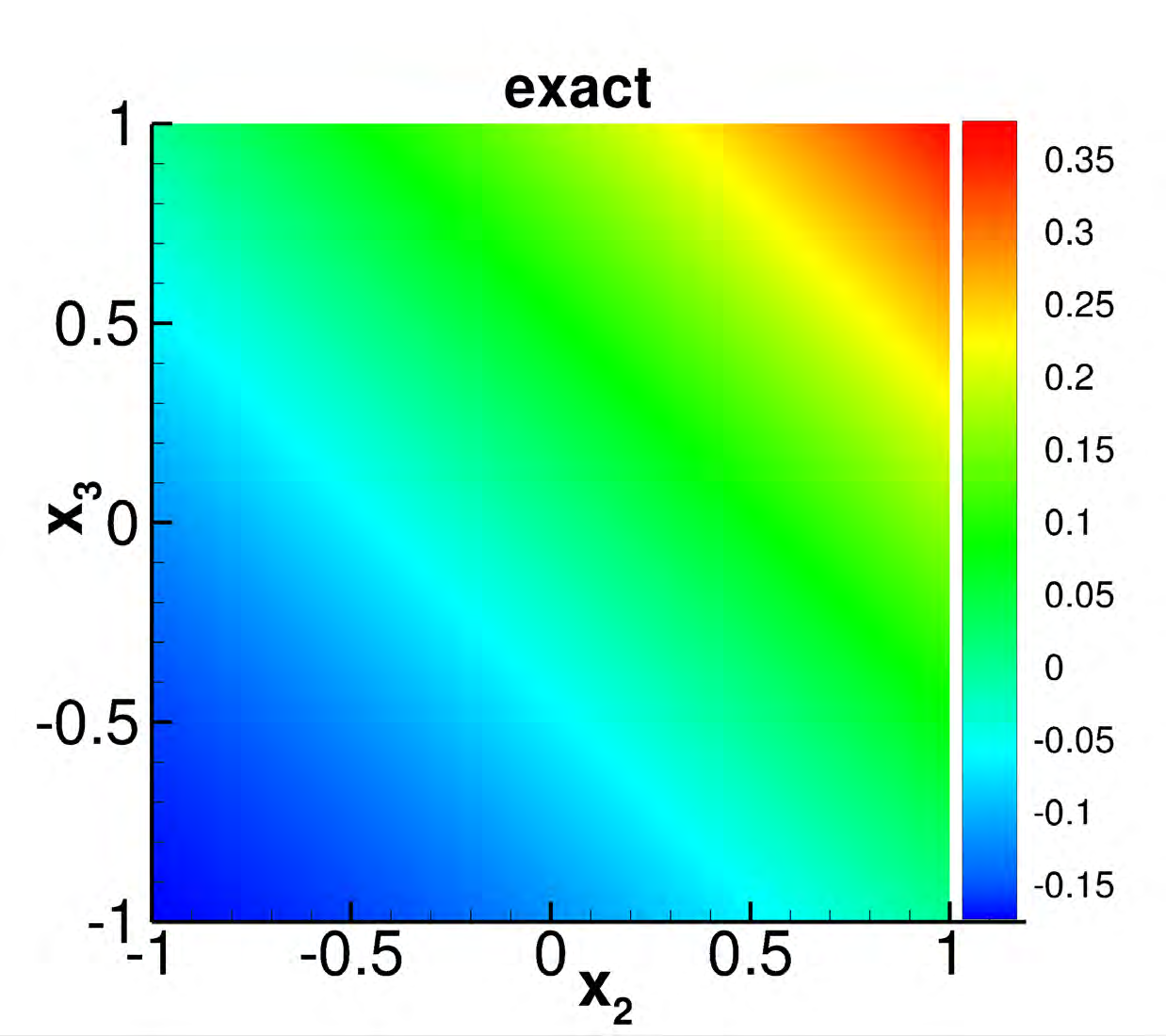}(d)
	 \includegraphics[width=0.26\textwidth]{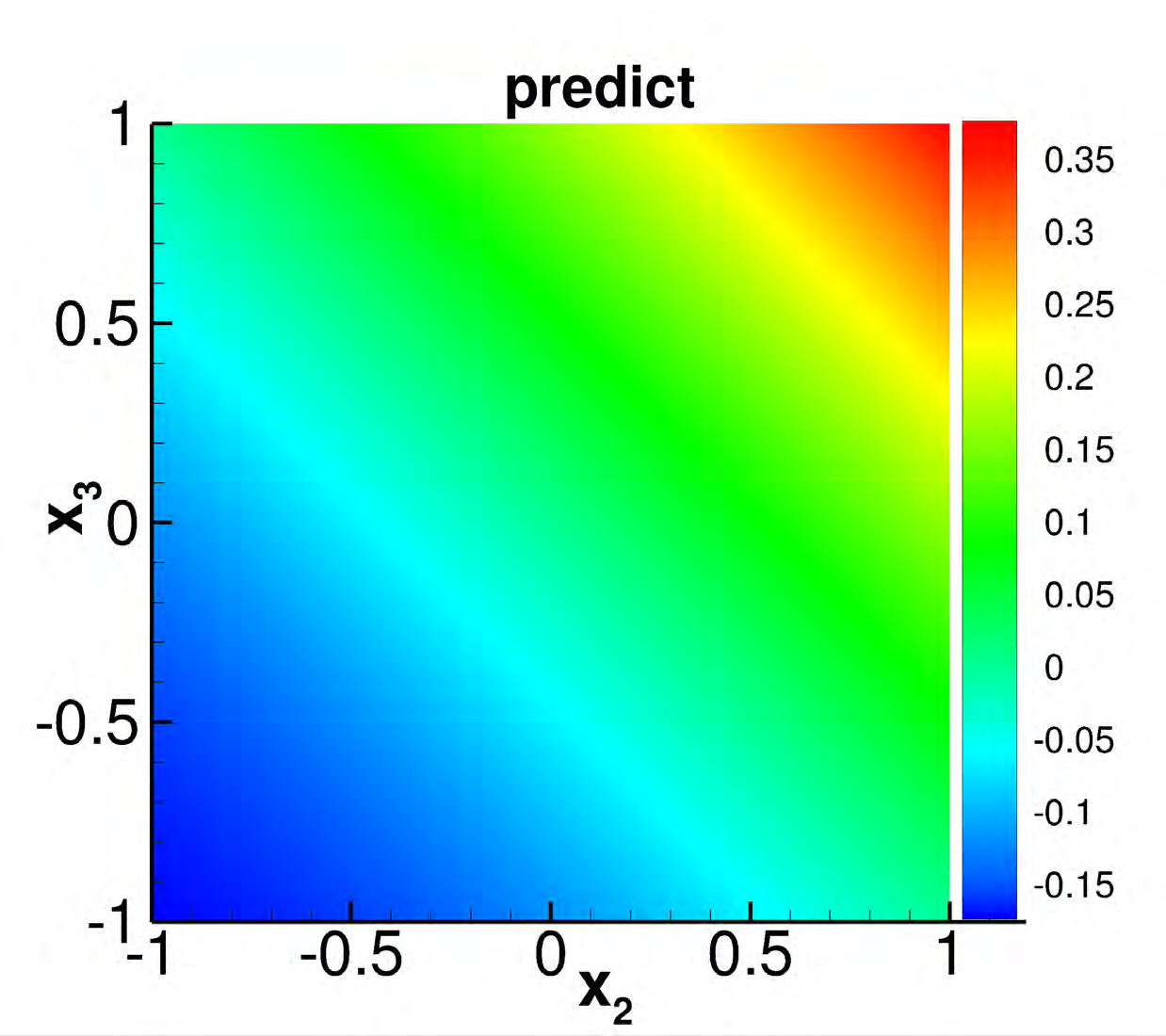}(e)
	\includegraphics[width=0.26\textwidth]{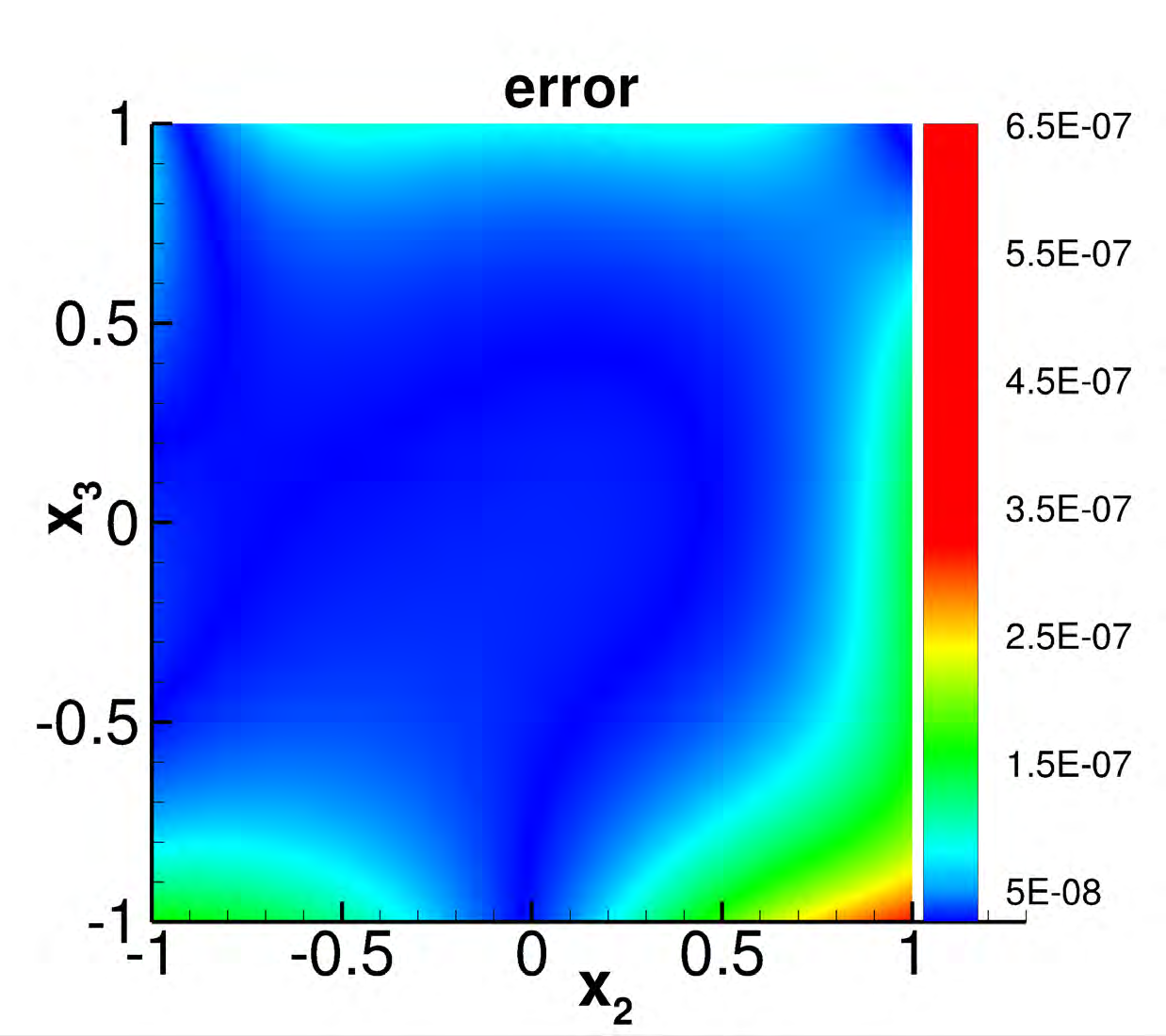}(f)
 }
 \centerline{
	\includegraphics[width=0.26\textwidth]{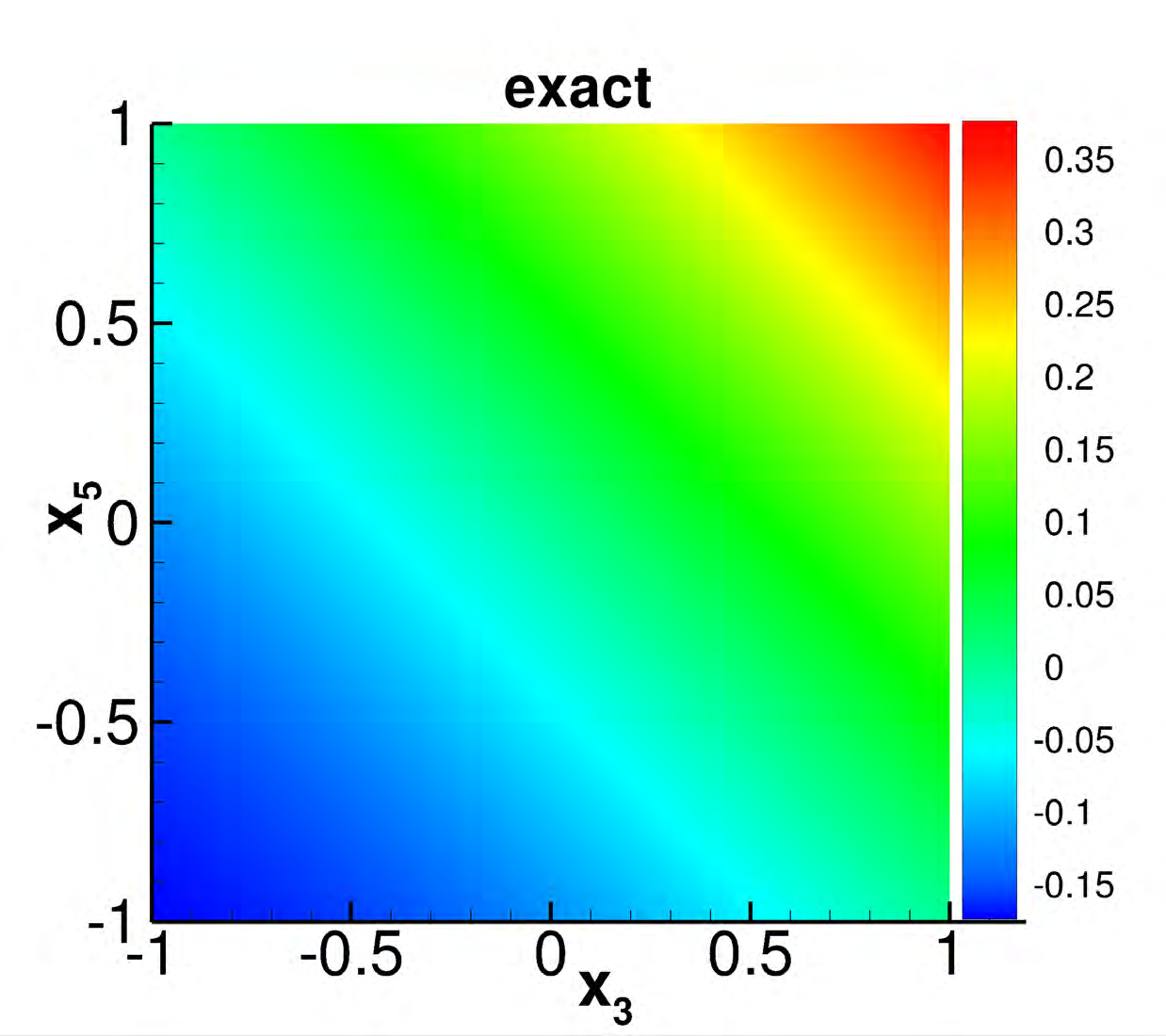}(g)
	\includegraphics[width=0.26\textwidth]{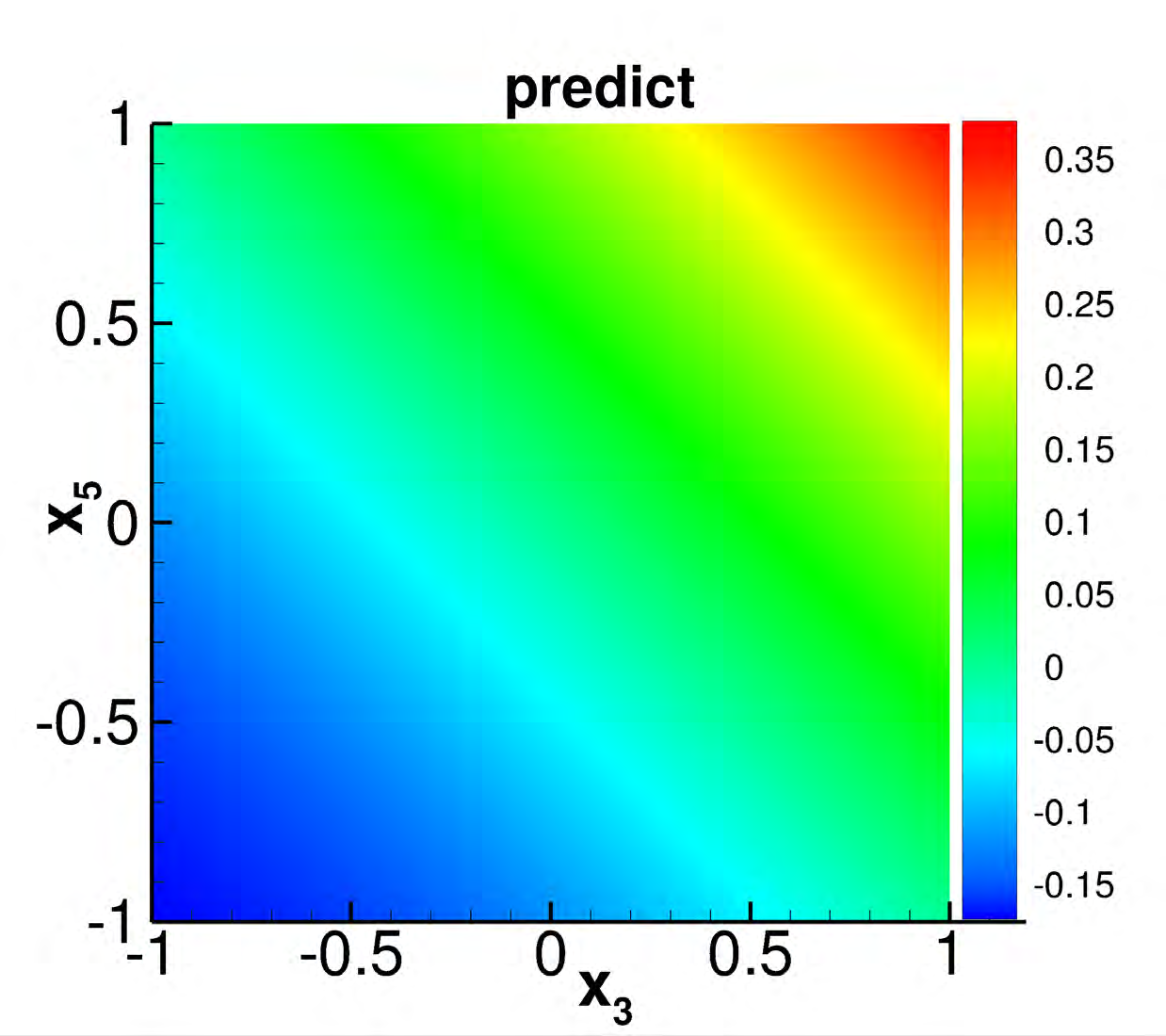}(h)
	 \includegraphics[width=0.26\textwidth]{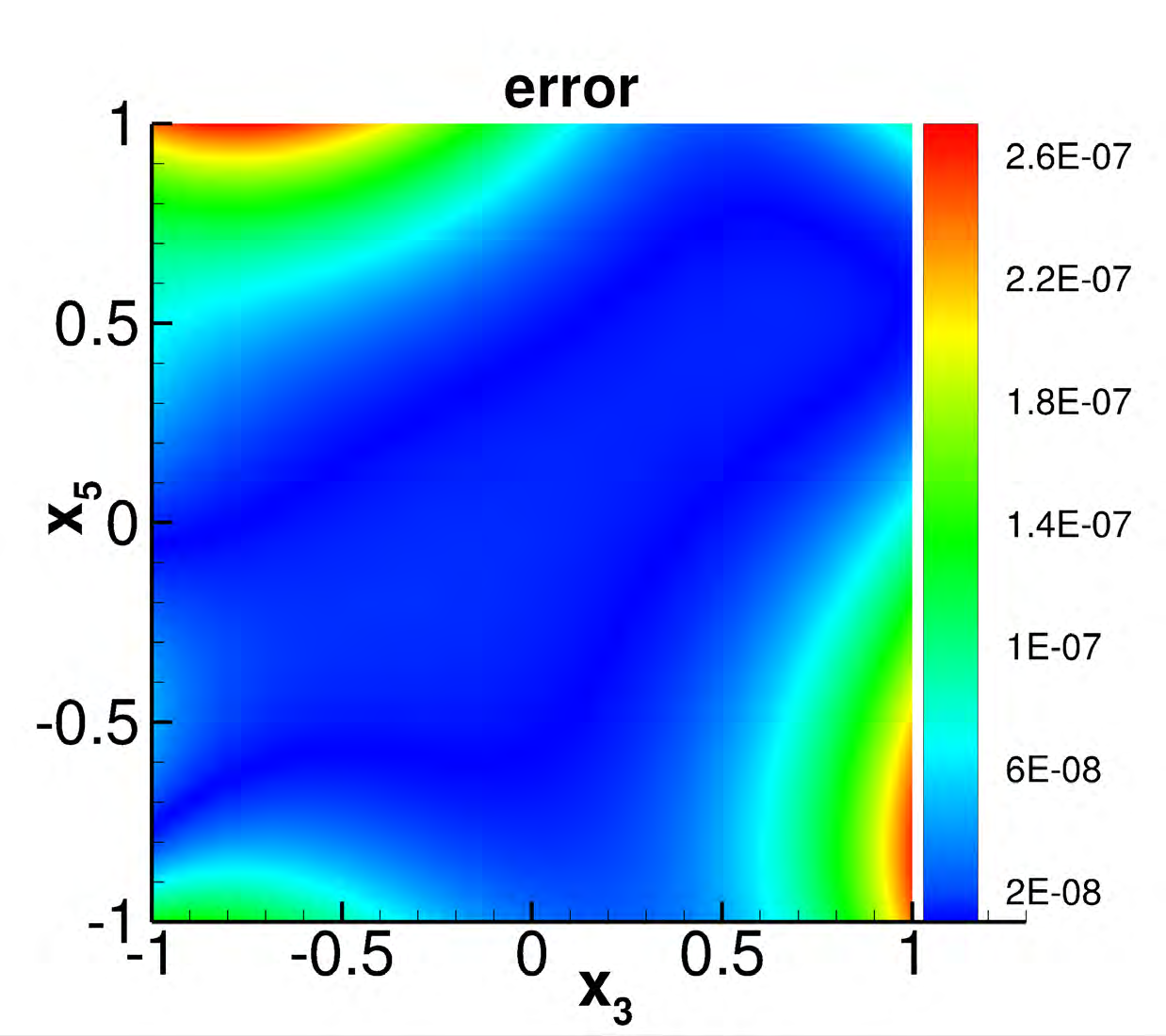}(i)
 }
\caption{\small Poisson equation ($d=7$): Distributions of the exact solution (left column), the ELM solution (middle column), and the point-wise absolute error of ELM (right column), in selected cross sections of the domain. Top row, $x_1$-$x_2$ plane; Middle row, $x_2$-$x_3$ plane; Bottom row, $x_3$-$x_5$ plane.
NN architecture: $[d,2000,1]$; $(N_{\bcd},N_{\ins})=(100,1000)$.
 }\label{poisson_dim7}
\end{figure}

In all the studies so far, the $e_{l^{\infty}}$ and $e_{l^2}$ errors are evaluated on a finite set of random test points from the domain, characterized by $(N_{\bcd}^{(v)},N_{\ins}^{(v)})=(100,7000)$. We have observed the error levels on the order of $10^{-10}$ to $10^{-5}$ for the problem dimensions from $d=5$ to $15$. We would like to consider the following question. Are these error levels representative of the ELM solution error on the entire domain $\Omega$?

To answer this question, ideally one would generate a regular set of  grid points on $\Omega$, with a sufficiently large number of grid points in each direction, and then evaluate and visualize the point-wise errors of the ELM solution on these grid points. This is feasible for low dimensions, but immediately becomes impractical  when the  dimension $d$ increases  to even a moderate value.
On the other hand, we note that it is possible to extract/compute the ELM solution error on certain low-dimensional hyper-planes (e.g.~2D cross sections) in high  dimensions. By looking into the point-wise error distributions in selected cross sections, one can gain a general sense of the representative error levels in the domain.

Figure~\ref{poisson_dim7} is an illustration of the ELM error for the Poisson equation in dimension $d=7$, using  cross sections. It shows distributions of the exact solution, the ELM solution, and the point-wise absolute error of ELM, on three cross sections of the domain ($x_1$-$x_2$ plane, $x_2$-$x_3$ plane, and $x_3$-$x_5$ plane). For a selected $x_i$-$x_j$ plane, the other coordinates of this plane has been set to zero, $x_k=0$ for $k\neq i,j$ (i.e.~the middle of the domain). For each cross section, the ELM solution, the exact solution, and the ELM error are evaluated on a uniform set of $Q^{(v)}\times Q^{(v)}$ ($Q^{(v)}=800$) grid points. The other simulation parameters are listed in the caption of this figure. One can observe that the point-wise error levels  as shown in these cross sections are comparable to (or consistent with) those observed in the convergence studies for $d=7$. This suggests that the $e_{l^{\infty}}$ and $e_{l^2}$ errors computed on the random set of test points with $(N_{\bcd}^{(v)},N_{\ins}^{(v)})=(100,7000)$ indeed seems to reflect well the ELM error on the domain $\Omega$.

\subsubsection{Nonlinear Poisson Equation}
\label{sec_312}

We consider the domain $\Omega=[-1,1]^d$ and the following problem on $\Omega$,
\begin{subequations}\label{eq_44}
\begin{align}
-\nabla \cdot( a(u)\nabla u)&=f(x), \quad x\in \Omega,\\
u&=g(x), \quad  x\in \partial \Omega,
\end{align}
\end{subequations}
where $a(u)=u^2-u$, $g(x) = \exp\left(-\frac{1}{d}\sum_{i=1}^{d}x_i\right)$, and
$
f(x) = \frac{1}{d}\left[ 
-3\exp\left(-\frac{3}{d}\sum_{i=1}^d x_i \right)
+2\exp\left(-\frac{2}{d}\sum_{i=1}^d x_i \right)
\right].
$
This system has an exact solution $u(x)=\exp\left(-\frac{1}{d}\sum_{i=1}^{d}x_i\right)$.

\begin{table}[tb]
\centering\small
\begin{tabular}{ccccccc}
	\hline
	$R_m$ &0.01&0.05&0.1&0.5&1&2 \\  \hline
	$e_{l^\infty}$ &5.44E-6&4.10E-9&2.59E-10&1.65E-8&5.92E-6&4.59E-4 \\
	$e_{l^{2}}$ &9.07E-7&7.57E-10&3.16E-11&3.79E-9&1.25E-6&9.01E-5\\
	\hline
\end{tabular}
\caption{\small Nonlinear Poisson equation: Determining $R_{m0}$ by the procedure from Remark~\ref{rem_a6} for $d=3$.
NN architecture: $[d,1000,1]$; $(N_{\bcd},N_{\ins})=(100,100)$.
}
\label{tab_3}
\end{table}

\begin{table}[tb]
\centering\small
\begin{tabular}{l| ccc}
\hline
$d$ & 3 & 5 & 9 \\ \hline
$R_{m0}$ & 0.1 & 0.05 & 0.001 \\
\hline
\end{tabular}
\caption{\small Nonlinear Poisson equation: $R_{m0}$ determined by the procedure from Remark~\ref{rem_a6} for several  dimensions.
}
\label{tab_4}
\end{table}

We employ the same notations as in Section~\ref{sec_poisson}. The simulation parameters include the network architecture $\mbs M_{\text{arch}}=[d,M,1]$, the random collocation points characterized by $(N_{\bcd},N_{\ins})$, the random test points characterized by $(N_{\bcd}^{(v)},N_{\ins}^{(v)})=(100,7000)$, the set of $Q^{(v)}\times Q^{(v)}$ uniform grid points with $Q^{(v)}=800$ on selected cross sections of the domain for evaluating the ELM and exact solutions, and  $R_m=R_{m0}$ for generating the random hidden-layer coefficients in the ELM neural network.

We first employ a fixed network architecture $[d,1000,1]$ and a set of collocation points characterized by  $(N_{\bcd},N_{\ins})=(100,100)$ to determine $R_{m0}$ by the procedure from Remark~\ref{rem_a6}. Table~\ref{tab_3} lists the $e_{l^{\infty}}$ and $e_{l^2}$ errors corresponding to several $R_m$ values for generating the random hidden-layer coefficients in ELM for  $d=3$, which leads to $R_{m0}\approx 0.1$. Table~\ref{tab_4} lists the $R_{m0}$ values determined by this procedure  for several  dimensions ranging from $d=3$ to $d=9$, using the same simulation parameters (network architecture, collocation points) as in Table~\ref{tab_3}. One can again observe that $R_{m0}$ decreases with increasing $d$.
We  employ $R_{m}=R_{m0}$ when generating the random hidden-layer coefficients in ELM in the subsequent tests.

\begin{figure}[tb]
	\centerline{
	\subfigure[d=3]
	{ \includegraphics[width=0.3\textwidth]{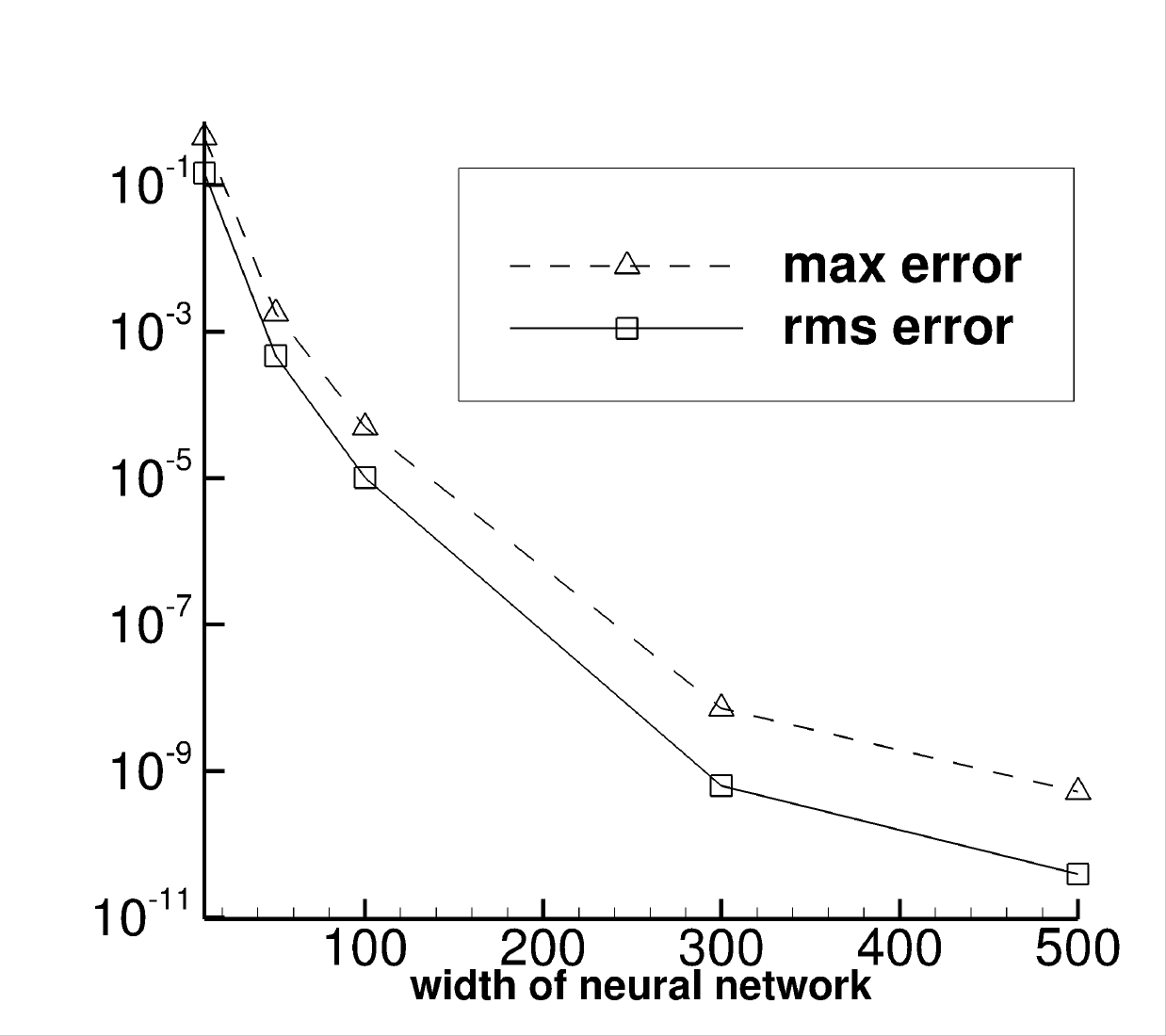}}
	\subfigure[d=5]
	{ \includegraphics[width=0.3\textwidth]{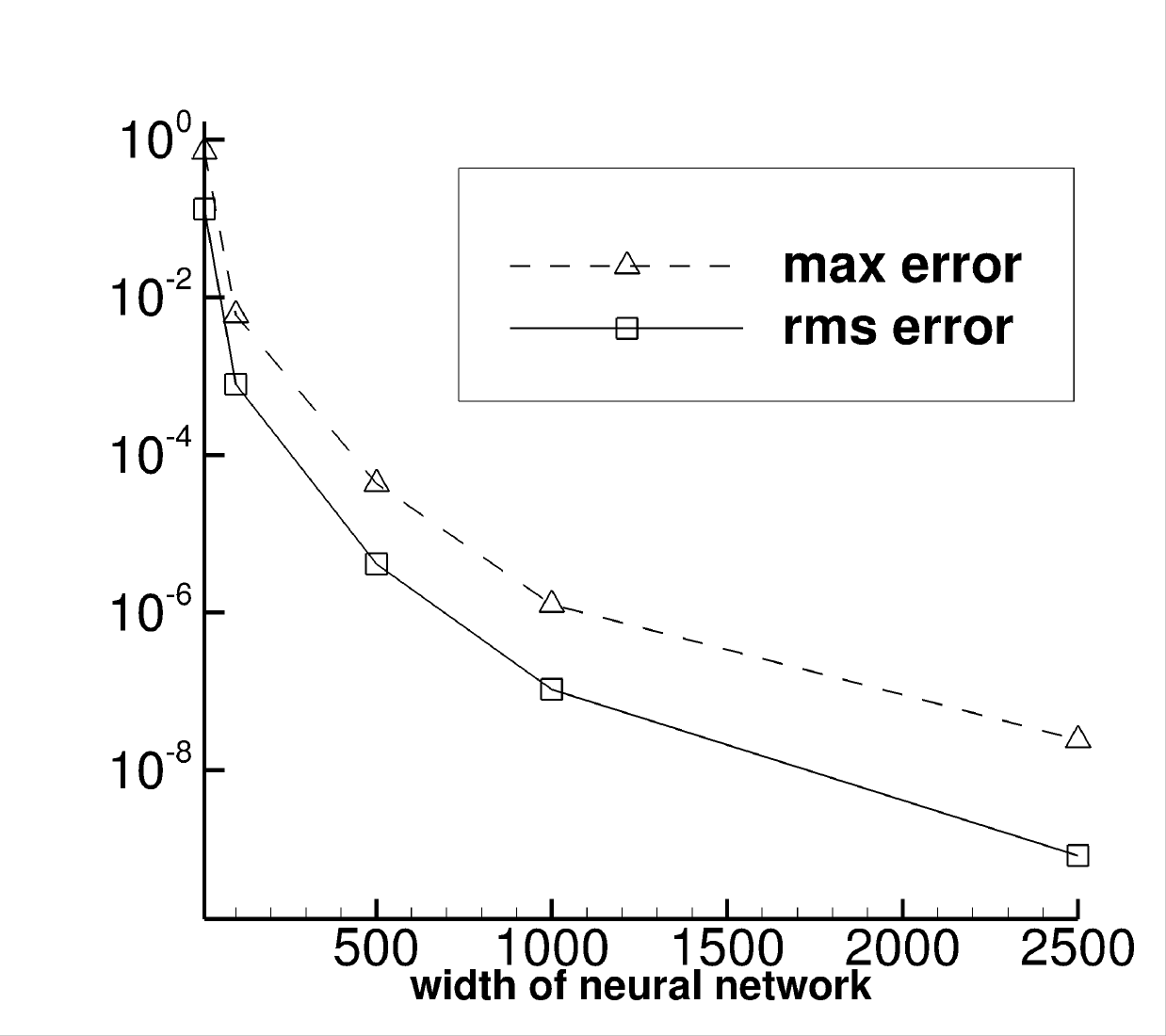}}
	\subfigure[d=9]
	{ \includegraphics[width=0.3\textwidth]{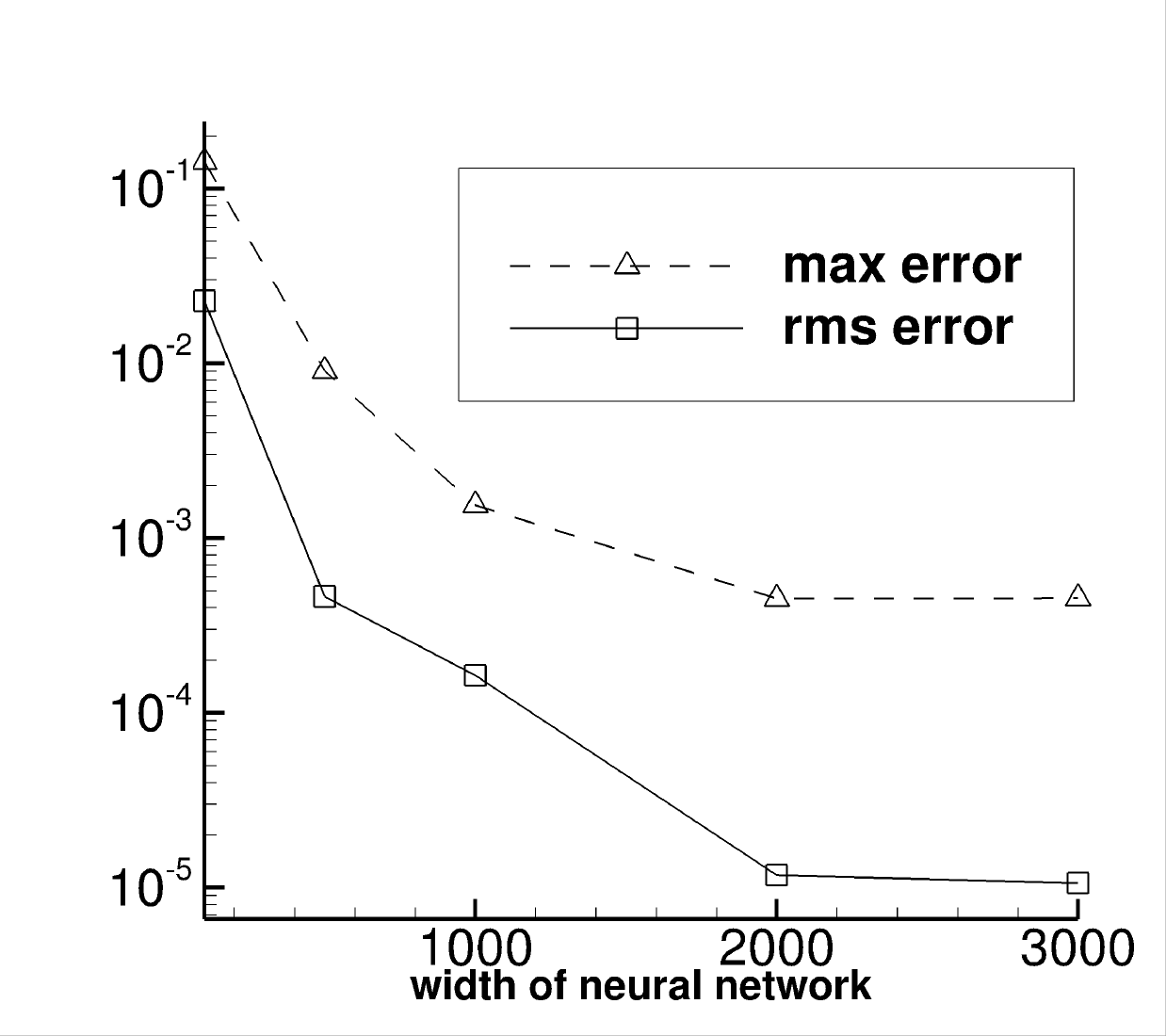}}
 }
	\caption{\small Nonlinear Poisson equation: $e_{l^{\infty}}$ and $e_{l^2}$ errors versus the number of training parameters ($M$) for dimensions (a) $d=3$, (b) $d=5$, and (c) $d=9$. 
 NN architecture: $[d,M,1]$ in (a,b,c);
 $(N_{\bcd},N_{\ins})=(100,100)$ in (a,b) and (100,220) in (c).
 }\label{fg_6}
\end{figure}

Figure~\ref{fg_6} illustrates the convergence behavior of the ELM errors with respect to the number of training parameters. It shows the $e_{l^{\infty}}$ and $e_{l^2}$ errors as a function of the number of training parameters ($M$) for three problem dimensions ($d=3$, $5$, $9$). The number of training collocation points is fixed, and their values are provided in the figure caption. The ELM errors decrease dramatically (quasi-exponentially) with increasing $M$ initially, but gradually plateau when $M$ becomes large. The $e_{l^2}$ error reaches a level around $10^{-11}$ for $d=3$, around $10^{-9}$ for $d=5$, and around $10^{-5}$ for $d=9$ in the range of parameters tested here.

\begin{figure}[tb]
	\centerline{
	\subfigure[d=3]
	{ \includegraphics[width=0.3\textwidth]{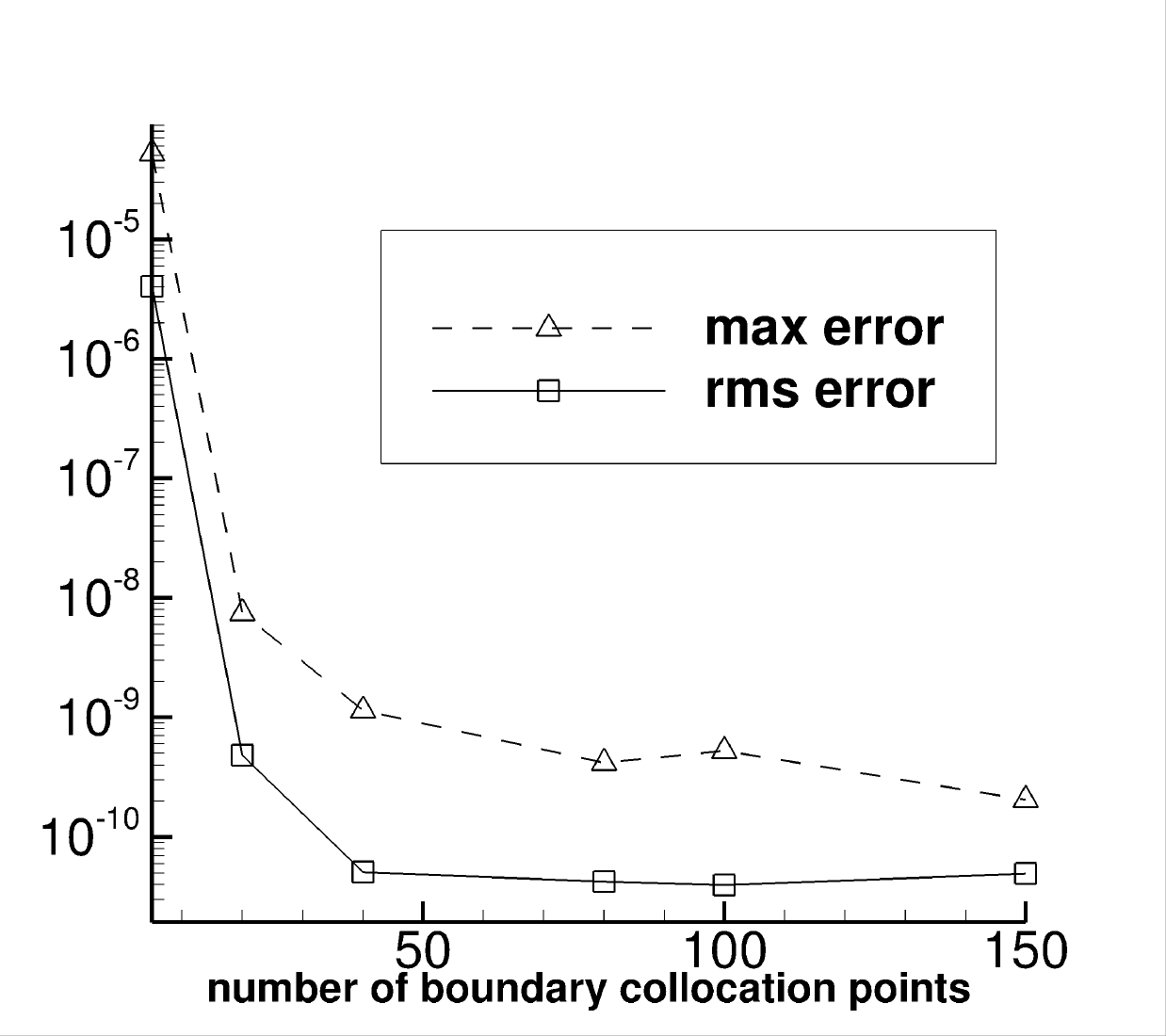}}
	\subfigure[d=5]
	{ \includegraphics[width=0.3\textwidth]{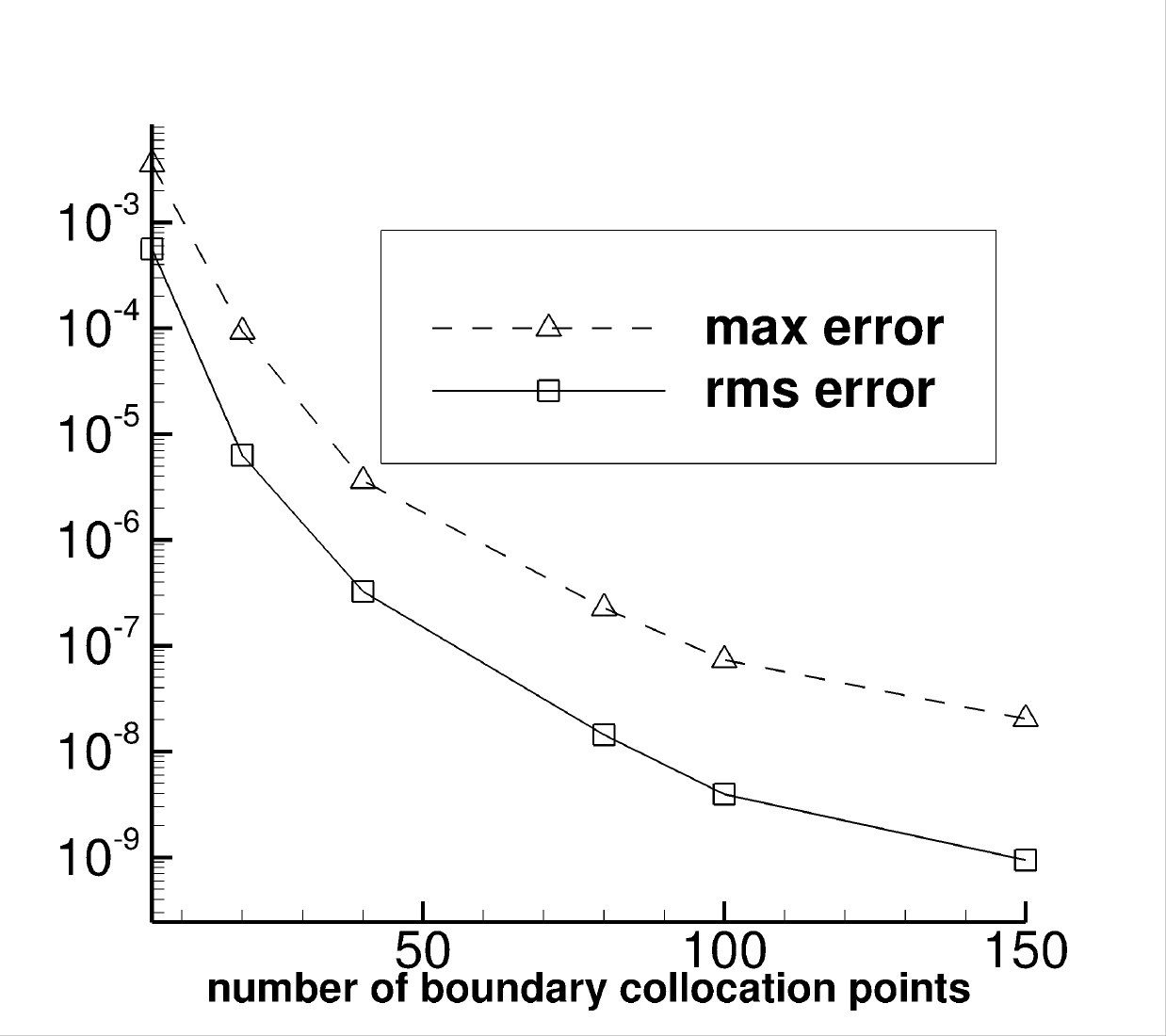}}
	\subfigure[d=9]
	{ \includegraphics[width=0.3\textwidth]{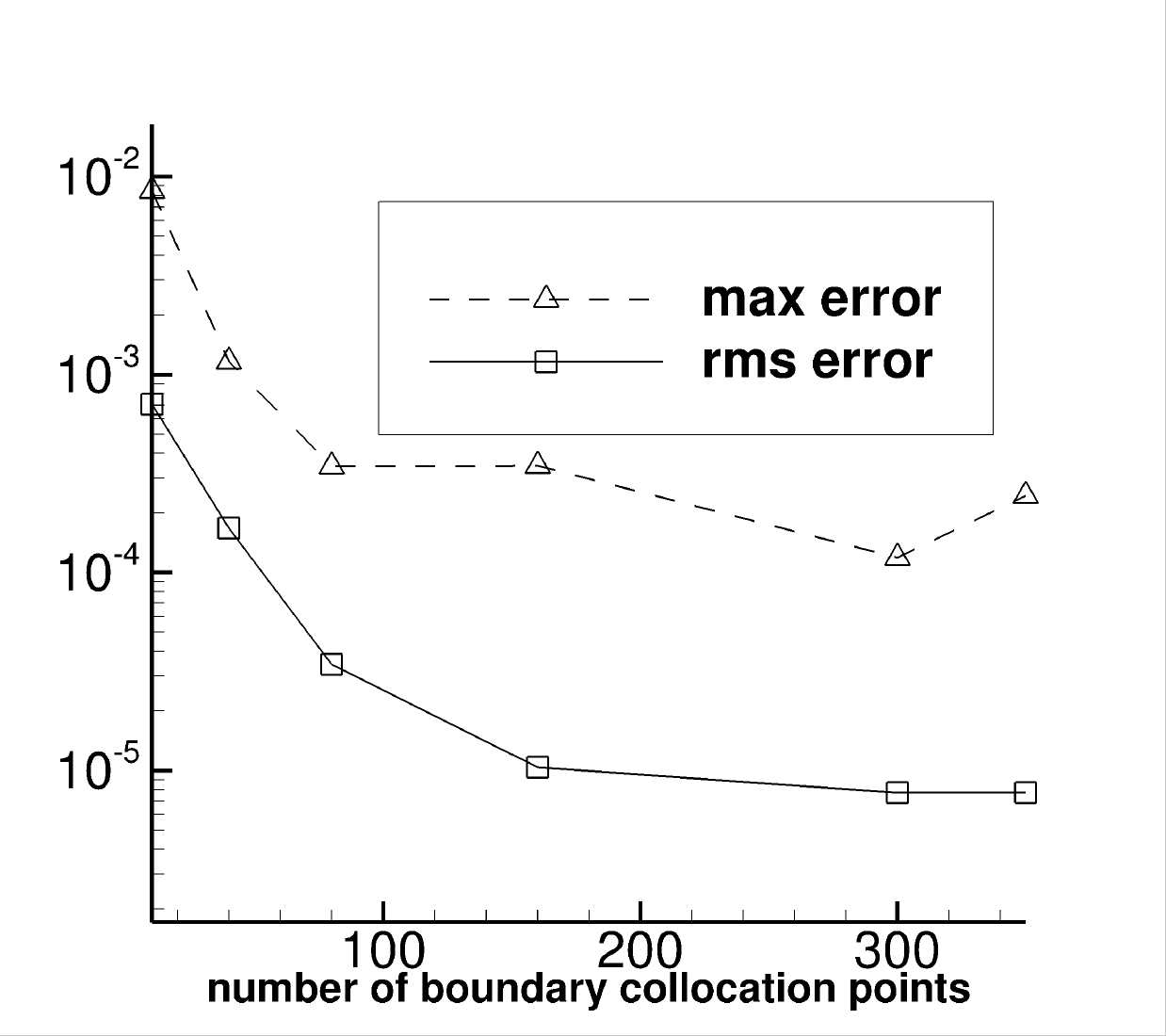}}
 }
 \caption{\small Nonlinear Poisson equation:
 $e_{l^{\infty}}$ and $e_{l^2}$ errors versus the number of  collocation points on each boundary ($N_{\bcd}$). 
 NN architecture: $[d,M,1]$ with $M=500$ in (a), $M=2500$ in (b), and $M=3000$ in (c). $N_{\ins}=100$ in (a,b,c).
 }
 \label{fg_7}
\end{figure}

\begin{figure}[tb]
\centerline{
	\subfigure[d=3]
	{ \includegraphics[width=0.3\textwidth]{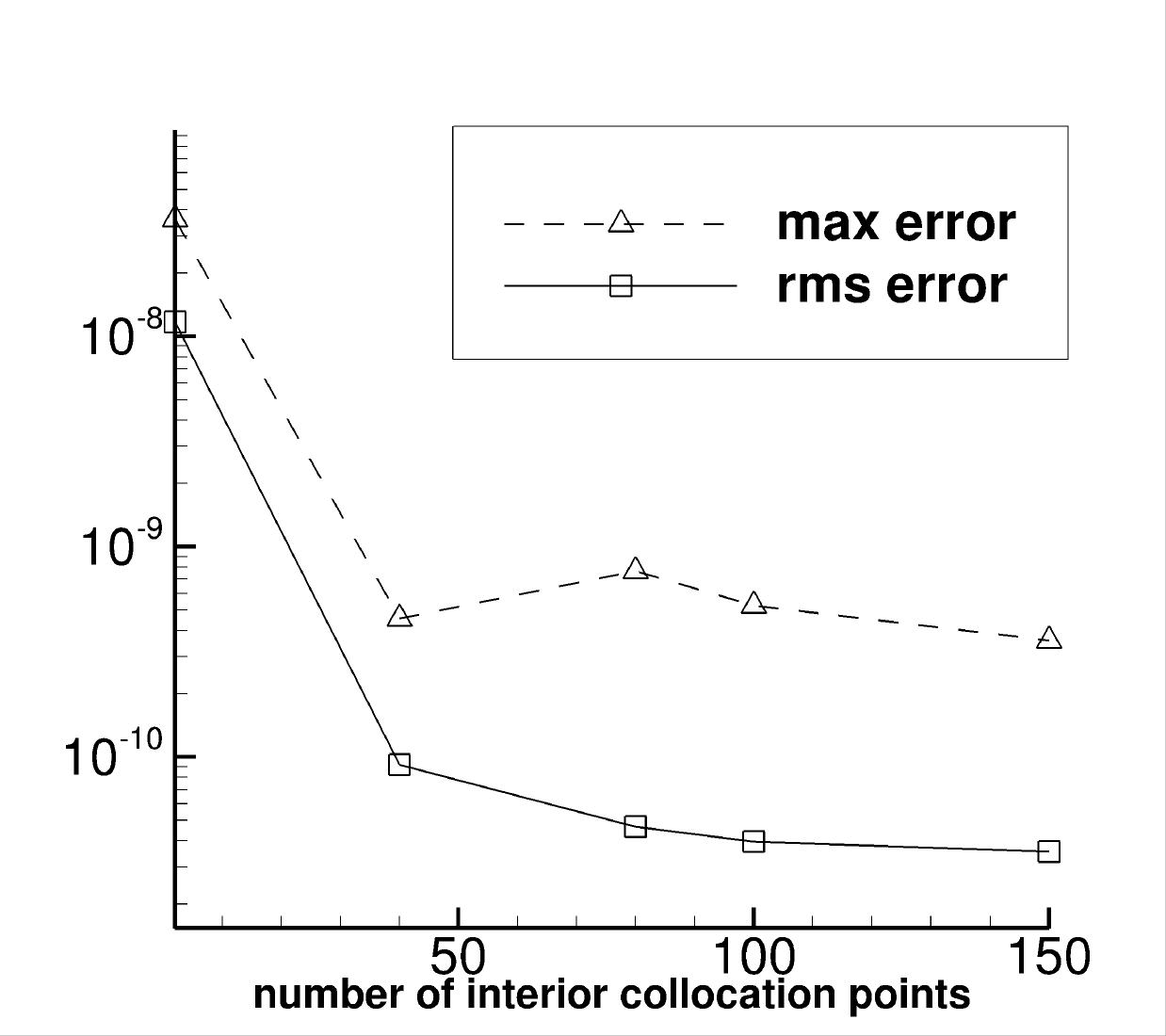}}
	\subfigure[d=5]
	{ \includegraphics[width=0.3\textwidth]{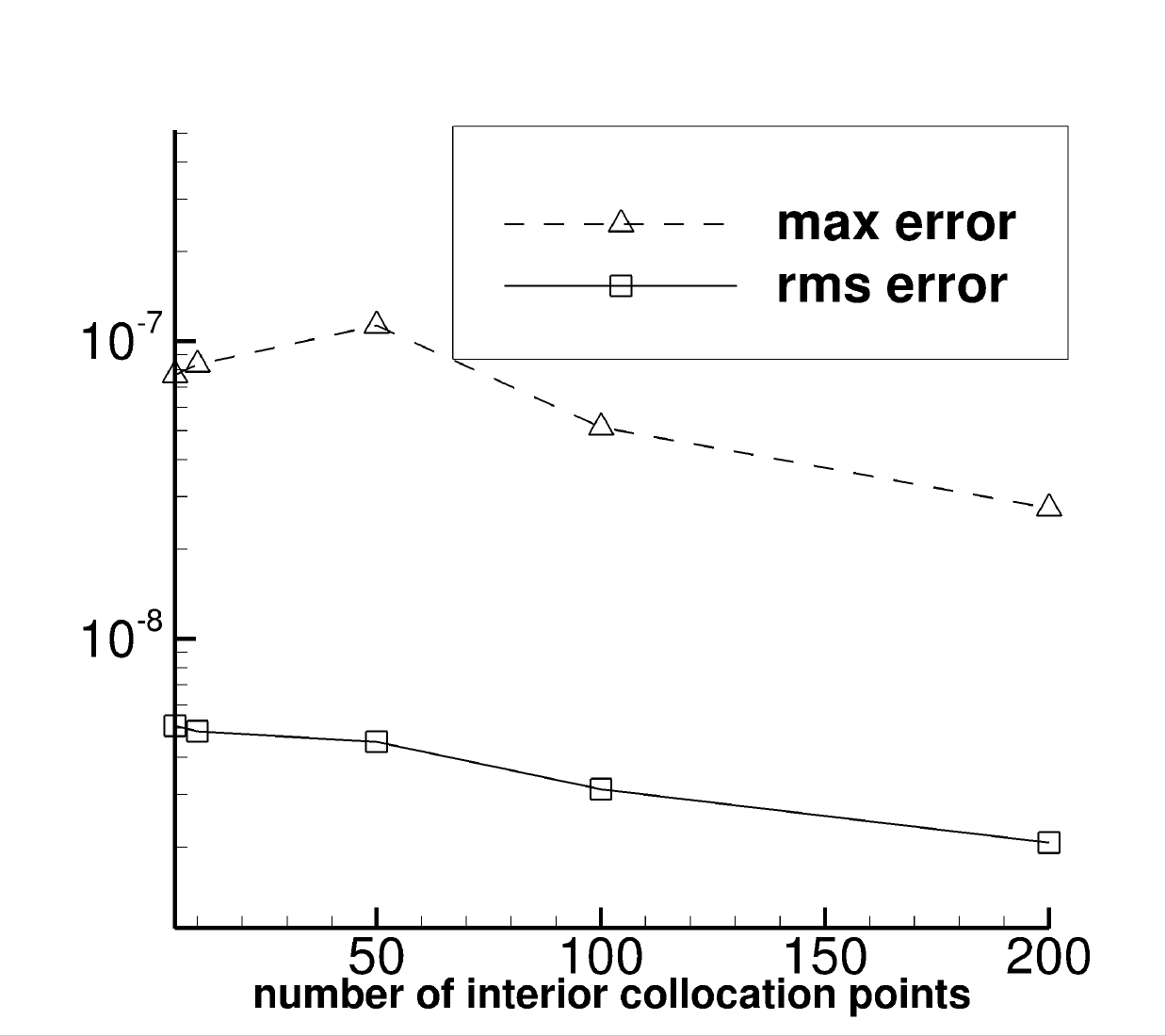}}
	\subfigure[d=9]
	{ \includegraphics[width=0.3\textwidth]{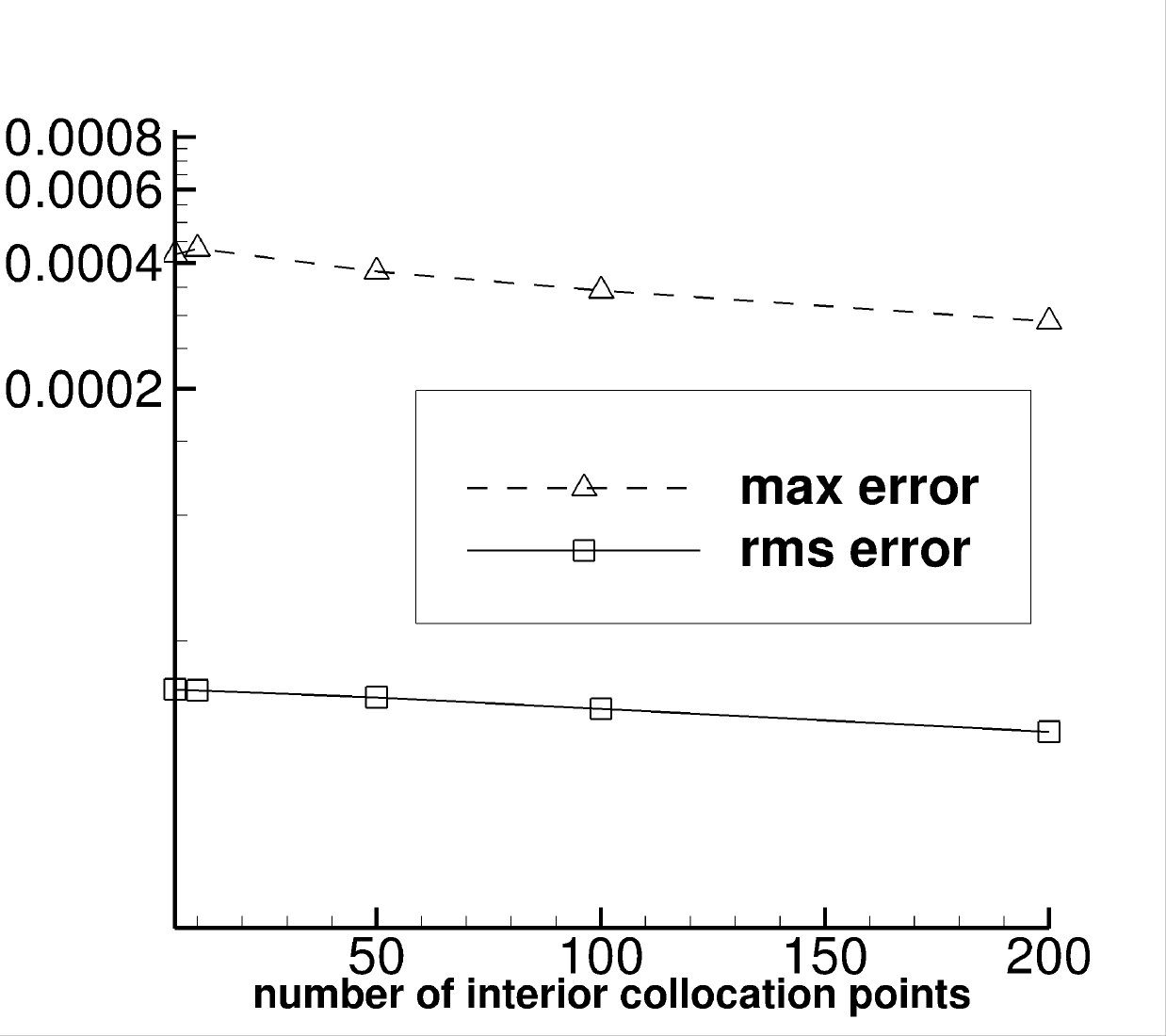}}
 }
	\caption{\small Nonlinear Poisson equation: $e_{l^{\infty}}$ and $e_{l^2}$ errors versus the number of interior collocation points ($N_{\ins}$).
 NN architecture: $[d,M,1]$ with $M=500$ in (a), $M=2500$ in (b), and $M=3000$ in (c). $N_{\bcd}=100$ in (a,b) and $80$ in (c).
 }
 \label{fg_8}
\end{figure}

The convergence of the ELM errors with respect to the number of collocation points is illustrated by Figures~\ref{fg_7} and~\ref{fg_8}, which show the $e_{l^{\infty}}$ and $e_{l^2}$ errors as a function of the boundary collocation points ($N_{\bcd}$) and the interior collocation points ($N_{\ins}$), respectively. The other crucial simulation parameters are provided in the captions of these figures. The ELM errors decrease significantly (approximately exponentially initially) with increasing number of boundary collocation points. On the other hand, increasing the number of interior collocation points in general only slightly improves  the error  in dimensions $d=5$ and $d=9$, and the error reduction is more significant in the lower dimension $d=3$. These behaviors are similar to what has been observed with the Poisson equation in the previous subsection.


\begin{figure}[!htb]
	\centerline{
\includegraphics[width=0.26\textwidth]{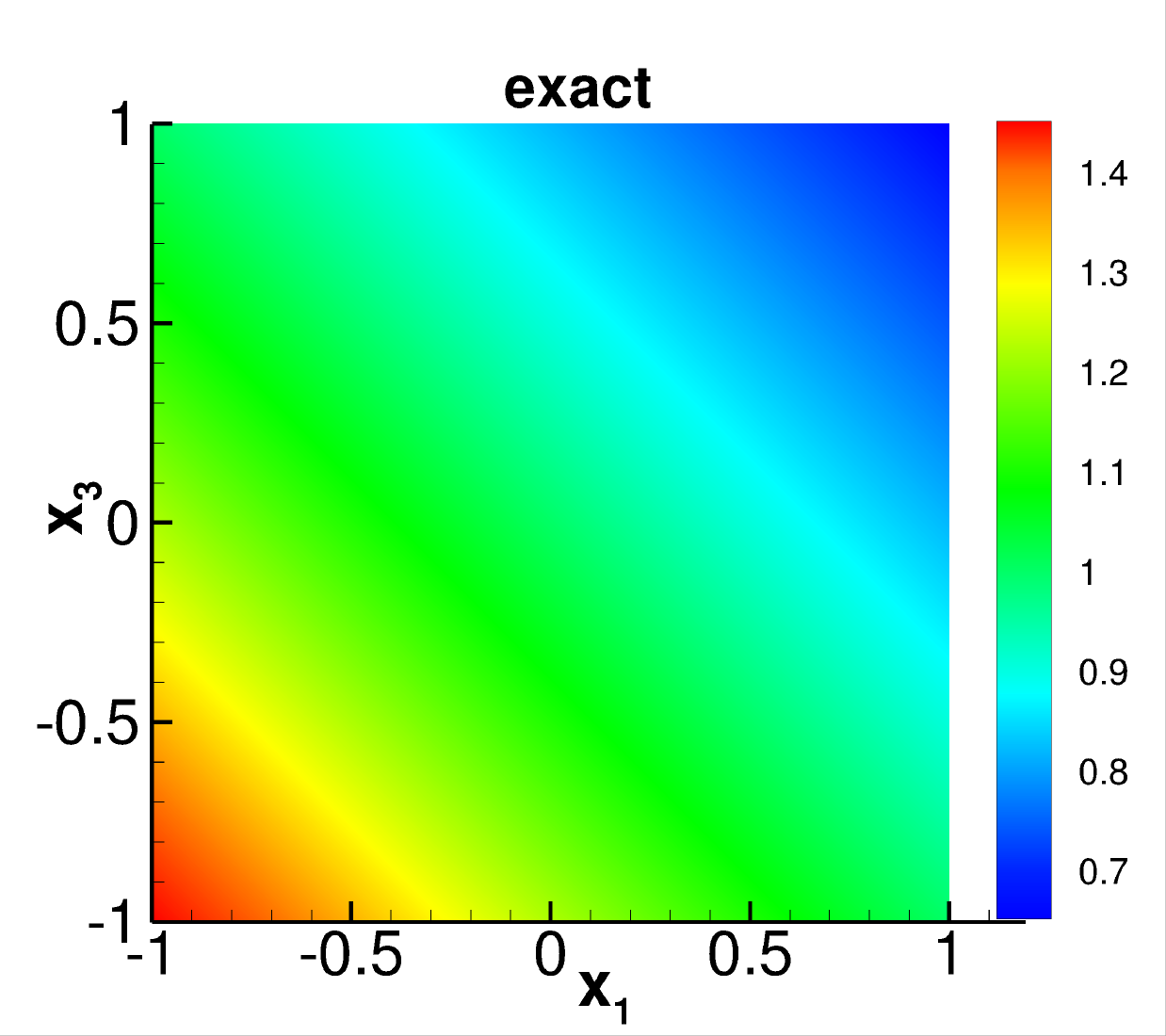}(a)
\includegraphics[width=0.26\textwidth]{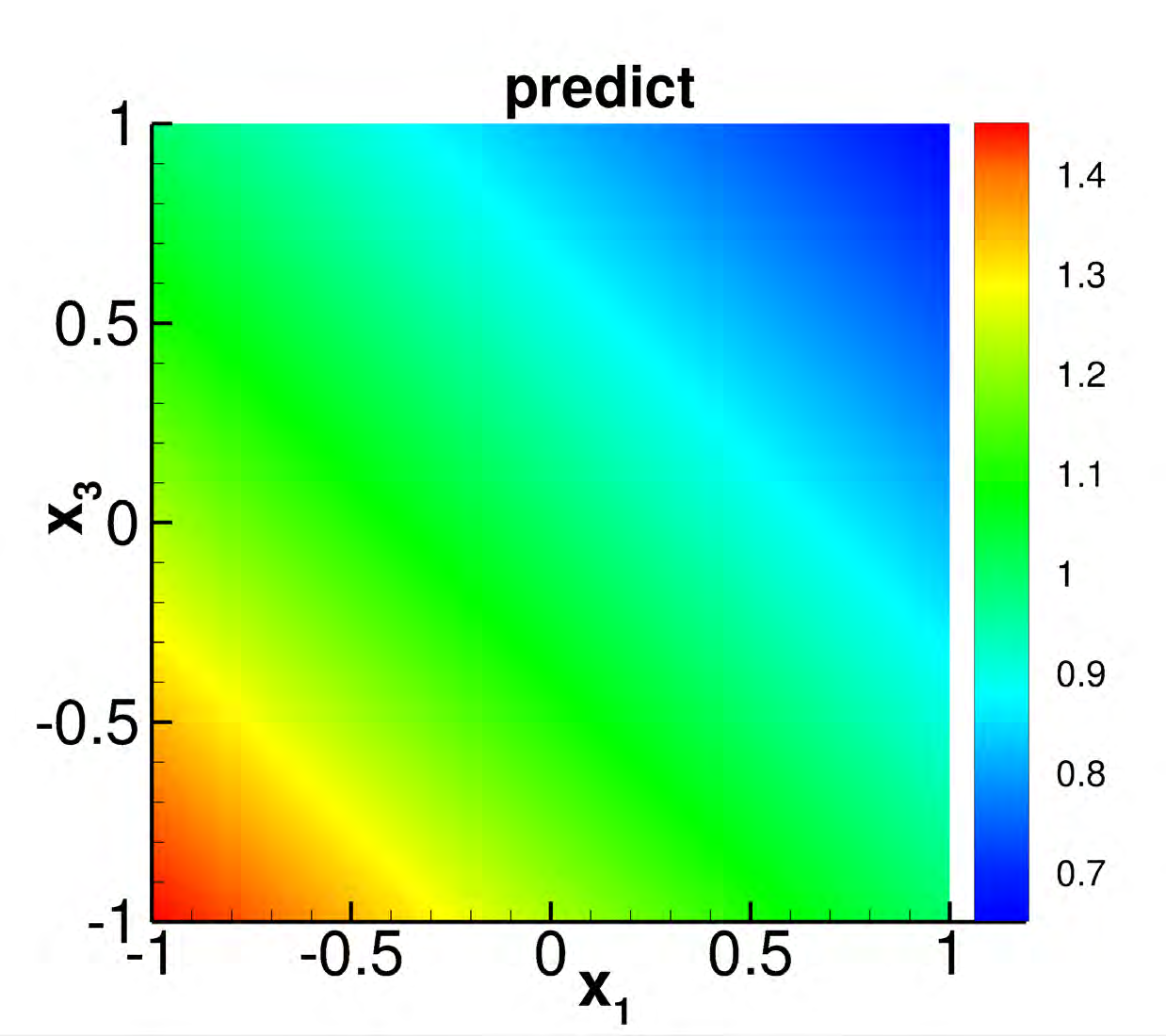}(b)
\includegraphics[width=0.26\textwidth]{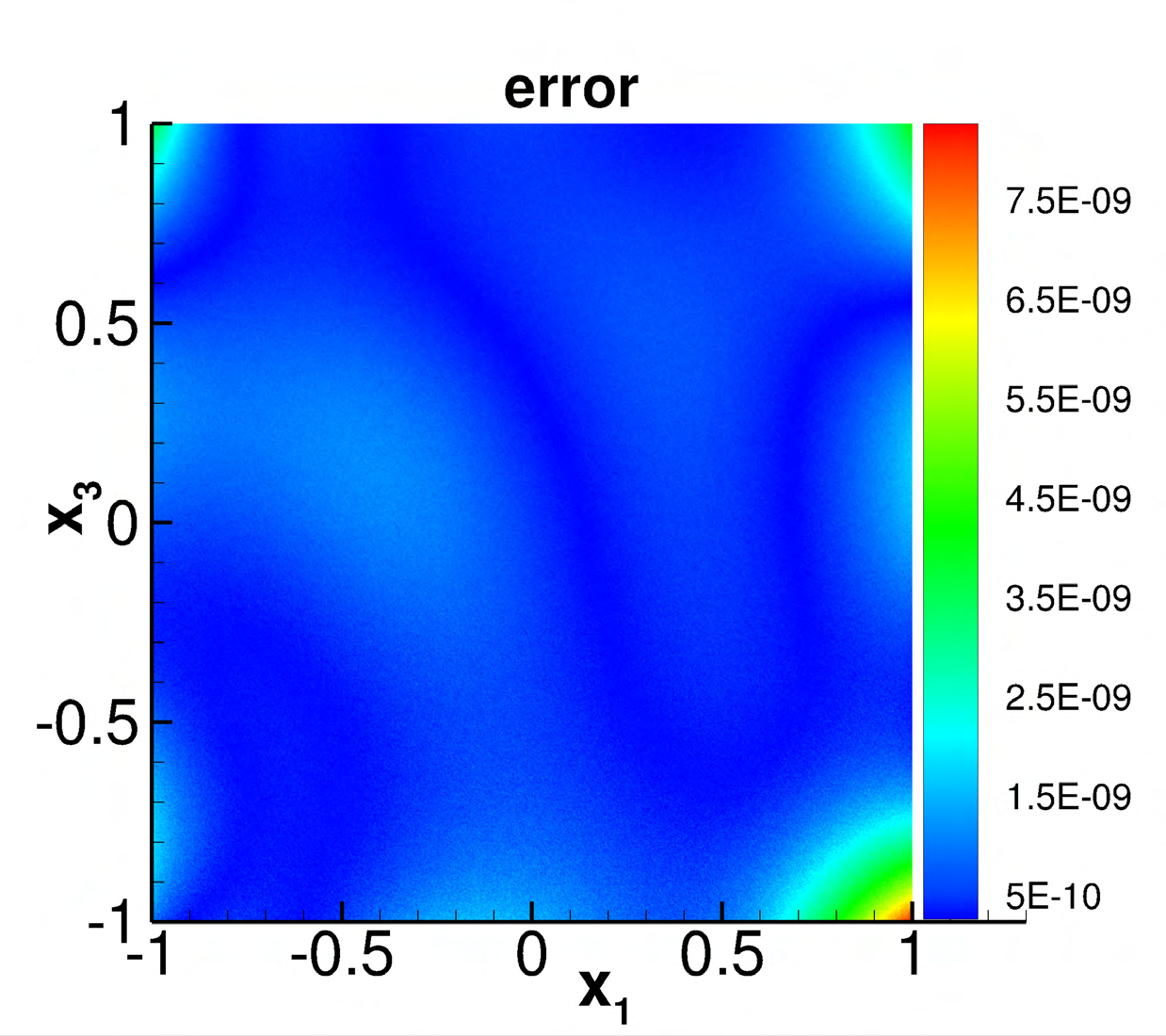}(c)
 }
 \centerline{
\includegraphics[width=0.26\textwidth]{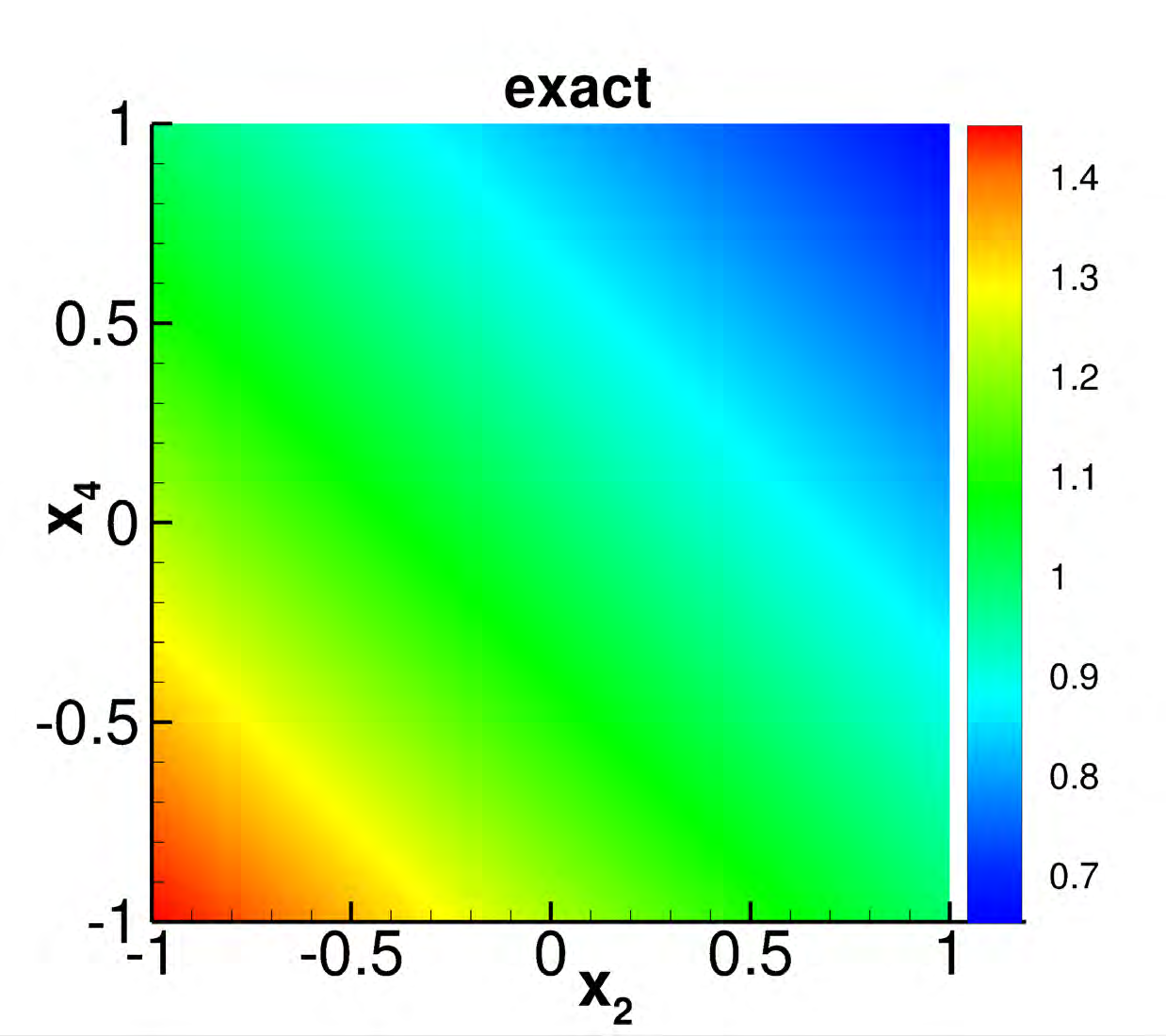}(d)
\includegraphics[width=0.26\textwidth]{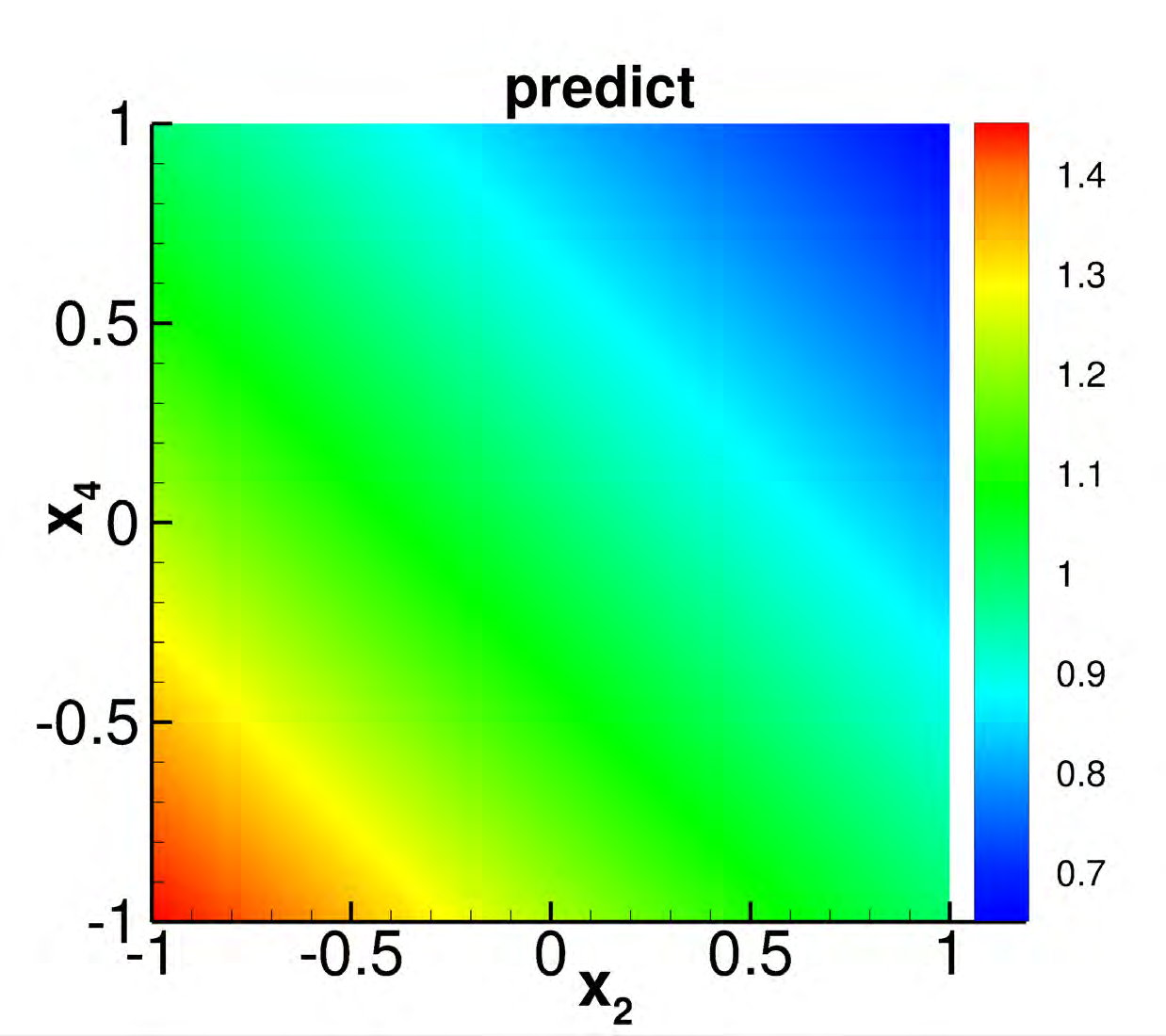}(e)
\includegraphics[width=0.26\textwidth]{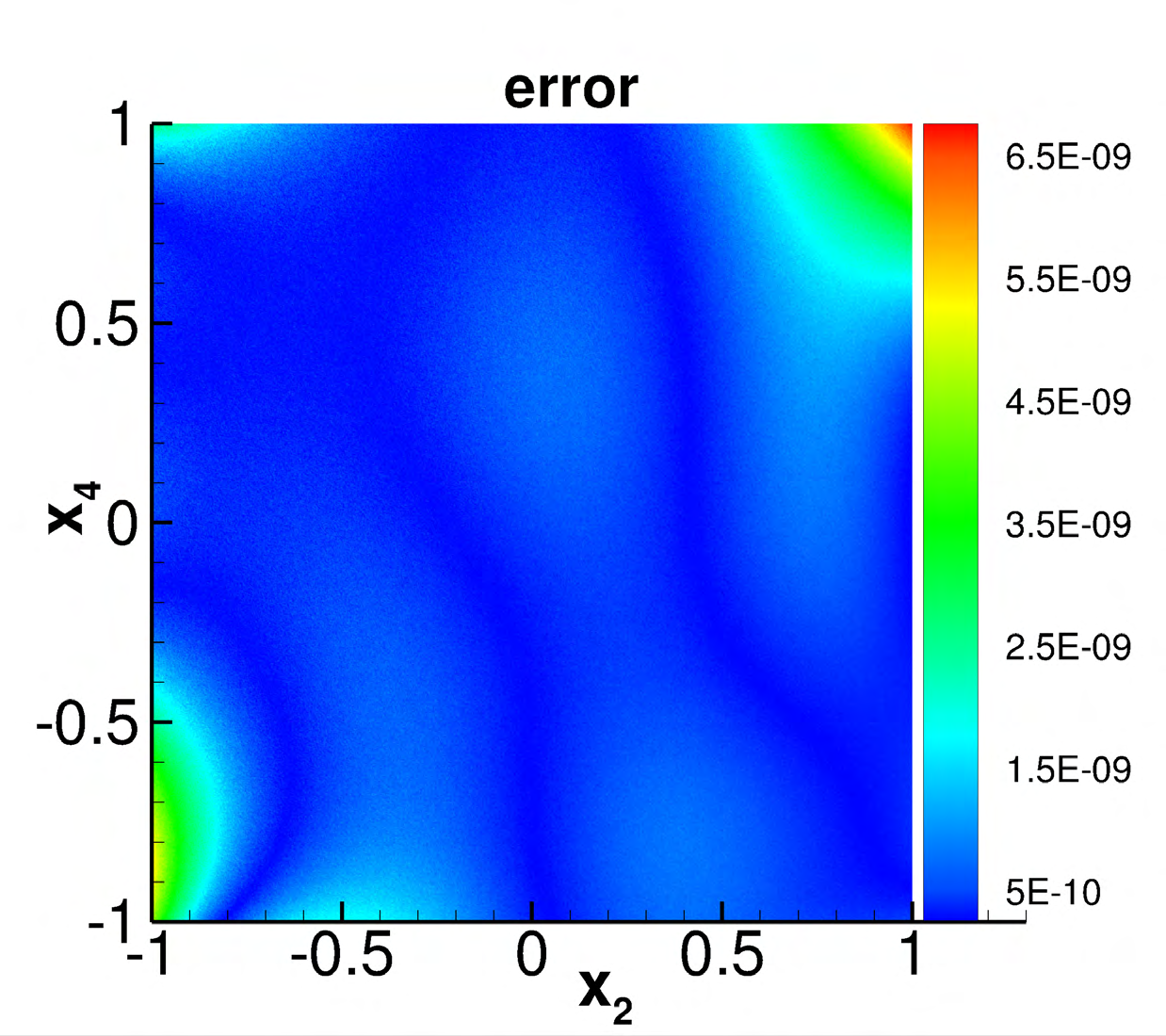}(f)
 }
 \centerline{
\includegraphics[width=0.26\textwidth]{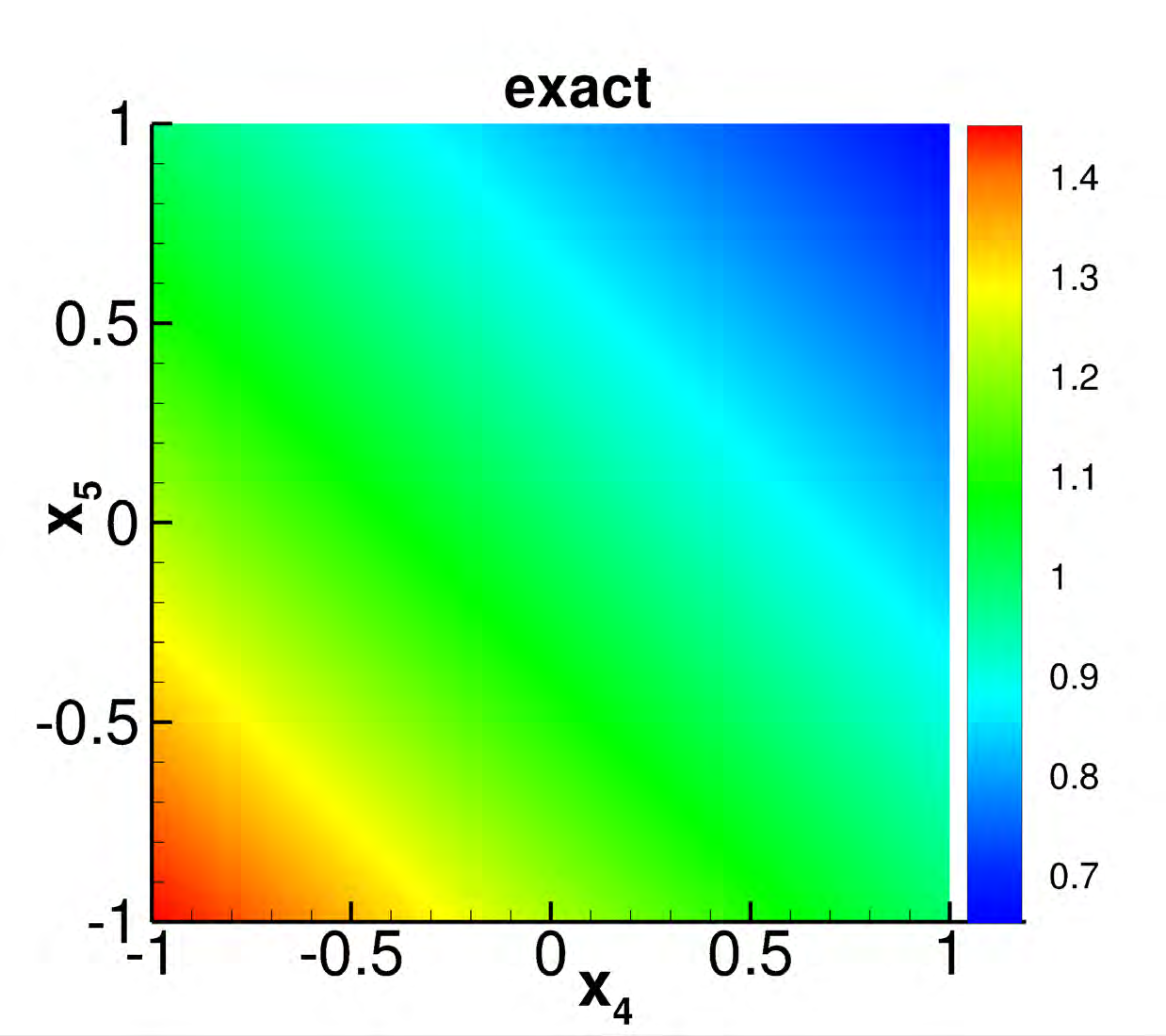}(g)
\includegraphics[width=0.26\textwidth]{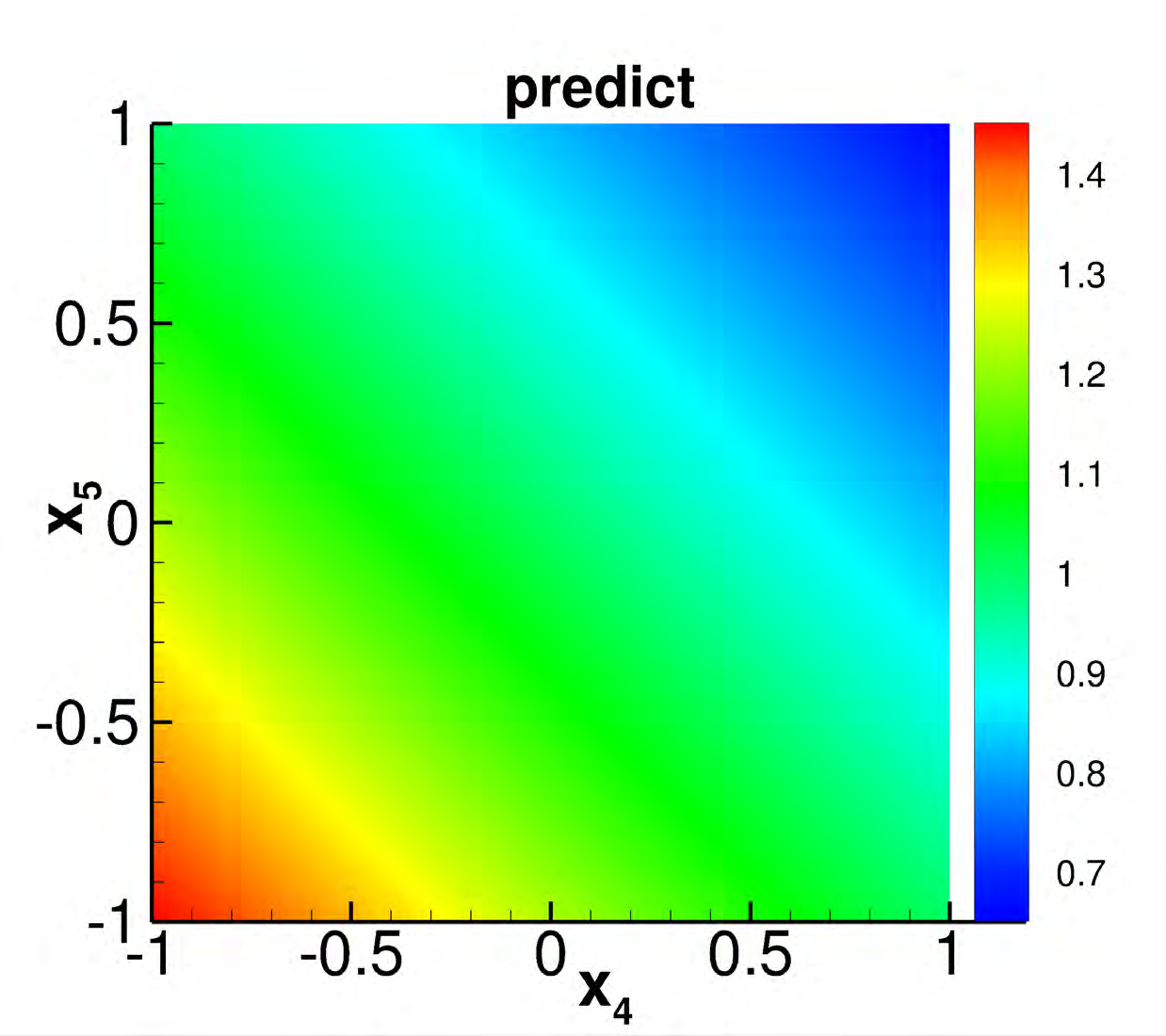}(h)
\includegraphics[width=0.26\textwidth]{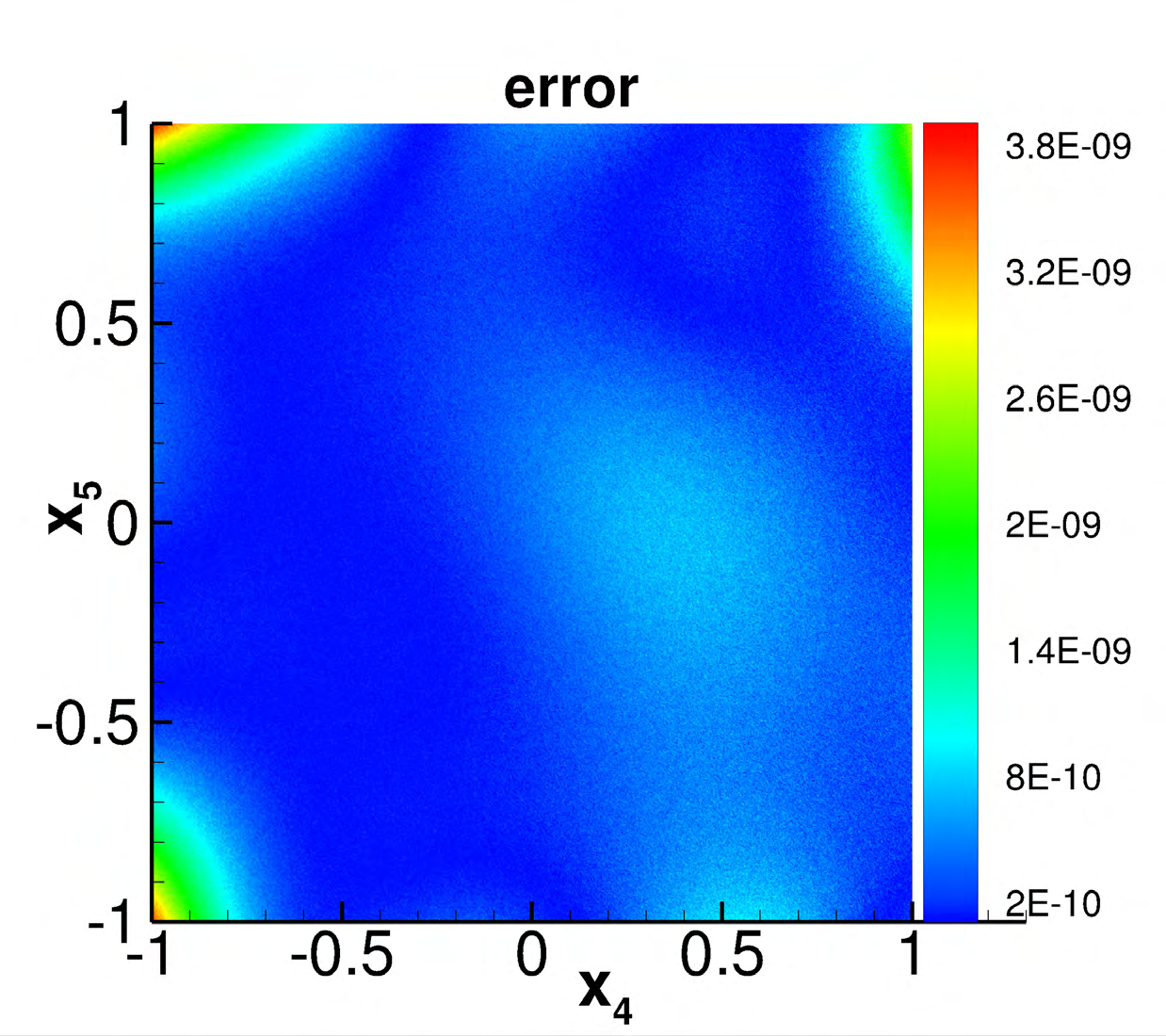}(i)
 }
	\caption{\small Nonlinear Poisson equation ($d=5$): Distributions of the exact solution (left column), the ELM solution (middle column), and the ELM point-wise absolute error (right column) in the $x_1$-$x_3$ plane (top row), $x_2$-$x_4$ plane (middle row), and the $x_4$-$x_5$ plane (bottom row).
 NN architecture: $[d,2500,1]$; $(N_{\bcd},N_{\ins})=(100,100)$.
 }
 \label{nl_poisson}
\end{figure}

Finally, Figure~\ref{nl_poisson} shows distributions of the ELM solution, the exact solution, and the ELM point-wise absolute error on three cross sections of the domain for $d=5$: the $x_1$-$x_3$, $x_2$-$x_4$, and $x_4$-$x_5$ planes. For each plane, the other coordinates of the plane are set to zero. The ELM simulation parameters are listed in the figure caption, and the distributions are plotted on a set of $800\times 800$ uniform grid points in these planes. It is evident that the ELM method has captured the solution very accurately.

\subsubsection{Advection Diffusion Equation}
\label{sec_314}

We next test the ELM method using the high-dimensional  advection-diffusion equation. Let $\Omega=[-1,1]^d$ and $T=1$.  Consider the  initial boundary value problem on the spatial-temporal domain $(x,t)\in\Omega\times [0,T]$,
\begin{subequations}\label{eq_45}
\begin{align}
\partial_t u-\nabla^2 u+R\cdot\nabla u &=f(x,t), \quad (x,t)\in \Omega \times [0,T],\\
u(x,t)&=g(x,t),\quad (x,t)\in \partial \Omega \times [0,T], \label{eq_45b}\\
u(x,0)&=h(x),\qquad x\in\Omega,
\end{align}
\end{subequations}
where 
$R=\frac{1}{d}(1,1,\ldots, 1)^T$, $f(x,t)=\frac{1}{d}\cos\left(\frac{1}{d}\sum_{i=1}^{d}x_i\right)\exp(-\frac{t}{d})$ on $\partial \Omega \times [0,T]$, and $h(x)=\sin\left(\frac{1}{d}\sum_{i=1}^{d}x_i\right)$. We employ the following analytic solution for this problem, $u(x,t)=\sin(\frac{1}{d}\sum_{i=1}^{d}x_i)\exp(-\frac{t}{d})$.
$g(x,t)$ in~\eqref{eq_45b} is set according to this expression.


To simulate this problem with ELM, we treat the time variable $t$ in the same way as the spatial coordinate $x$. We employ a neural network with an architecture, $\mbs M_{\text{arch}}=[d+1, M, 1]$, in which the $(d+1)$ input nodes represent $(x,t)$ and the single output node represents the solution $u(x,t)$. So the problem has been effectively treated as a $(d+1)$-dimensional problem in the simulations. We enforce the initial condition on $N_{t_0}=1000$ random collocation points on $\Omega$ at $t=0$, and enforce the boundary condition on $N_{\bcd}$ random collocation points on each of the boundaries $\partial\Omega\times[0,T]$. $N_{\ins}$ denotes the number of  random collocation points on the interior of $\Omega\times[0,T]$. After the neural network is trained, the $e_{l^{\infty}}$ and $e_{l^2}$ errors are computed on a set of random test points from $\Omega\times\{T\}$  characterized by $(N_{\bcd}^{(v)},N_{\ins}^{(v)})=(100,7000)$. Here $N_{\ins}^{(v)}$ and $N_{\bcd}^{(v)}$  denote the number of random test points from the interior of $\Omega\times\{T\}$ and from each of the boundaries of $\Omega\times\{T\}$, respectively.  


\begin{table}[tb]
\centering\small
\begin{tabular}{ccccccc}
	\hline
$R_m$	&5E-3&1E-2&5E-2&1E-1&5E-1&1 	\\ 
	$e_{l^\infty}$&9.79E-4&9.33E-5&6.21E-8&2.15E-8&1.92E-5&4.82E-4\\
	$e_{l^{2}}$ &2.28E-4&8.17E-6&7.01E-9&1.38E-9&1.98E-6&8.54E-5\\
	\hline
\end{tabular}
\caption{\small Advection diffusion equation: Determining $R_{m0}$ based on the procedure from Remark~\ref{rem_a6} for $d=3$.
NN architecture: $[d+1,1000,1]$; $(N_{\bcd},N_{\ins},N_{t_0})=(100,100,1000)$.
}
\label{tab_5}
\end{table}

\begin{table}[tb]
\centering\small
\begin{tabular}{l| ccc}
\hline
$d$ & 3 & 6 & 10 \\ \hline
$R_{m0}$ & 0.1 & 0.05 & 0.05 \\
\hline
\end{tabular}
\caption{\small Advection diffusion equation: $R_{m0}$ determined by the procedure from Remark~\ref{rem_a6} for several problem dimensions.
The simulation parameters  are the same as those in Table~\ref{tab_5}.
}
\label{tab_6}
\end{table}

Table~\ref{tab_5} shows the determination of  $R_{m0}$ using the procedure from Remark~\ref{rem_a6} for dimension $d=3$, leading to $R_{m0}\approx 0.1$. Table~\ref{tab_6} lists the $R_{m0}$ values corresponding to several problem dimensions for the advection diffusion equation. In these tests for determining $R_{m0}$, we have employed an NN architecture $[d+1,1000,1]$, and the number of  collocation points is characterized by $(N_{\bcd},N_{\ins},N_{t_0})=(100,100,1000)$.
In  subsequent tests of this section, the hidden-layer coefficients are set to uniform random values from  $[-R_{m0},R_{m0}]$.

\begin{figure}[!htb]
	\centerline{
\includegraphics[width=0.26\textwidth]{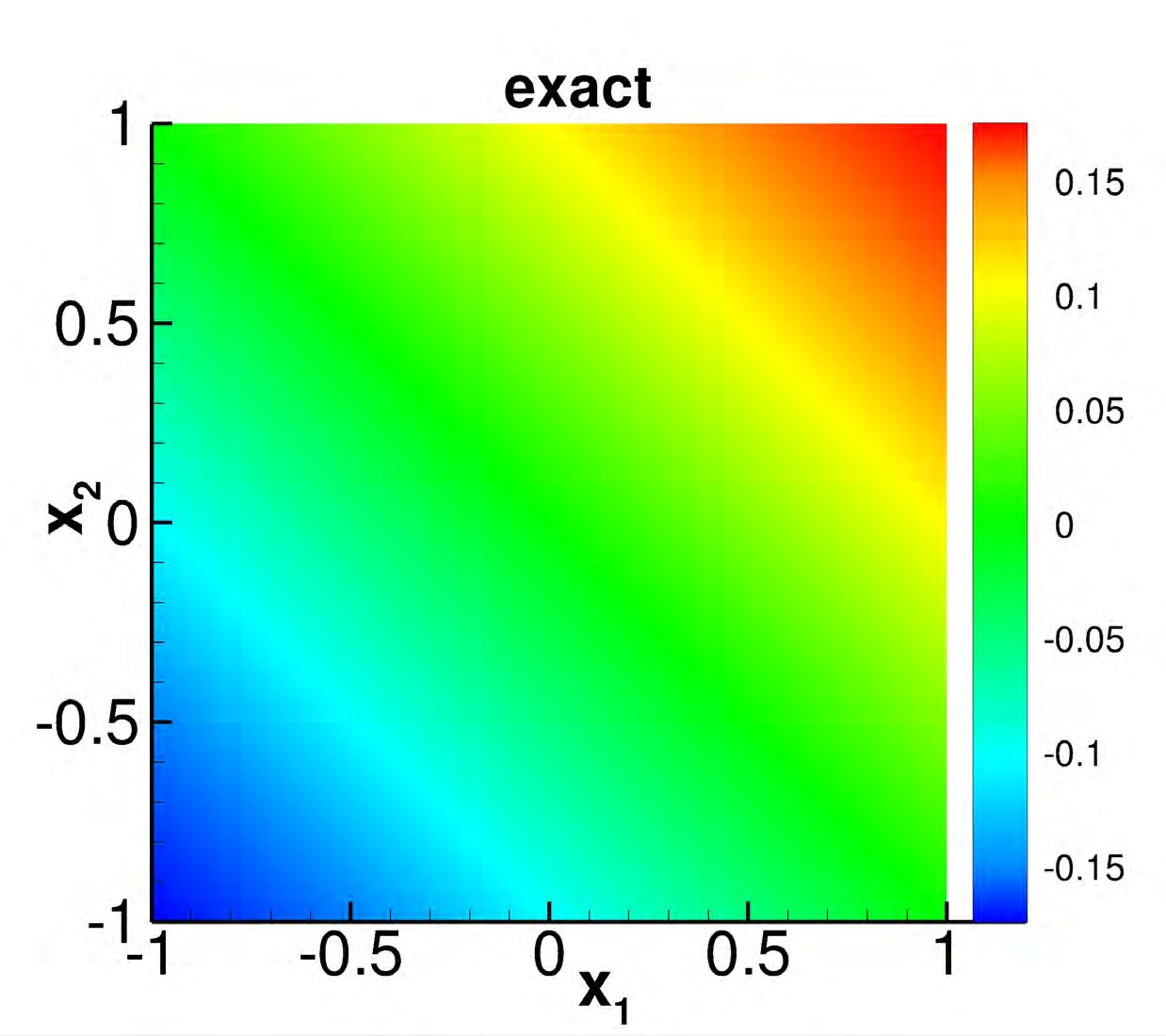}(a)
\includegraphics[width=0.26\textwidth]{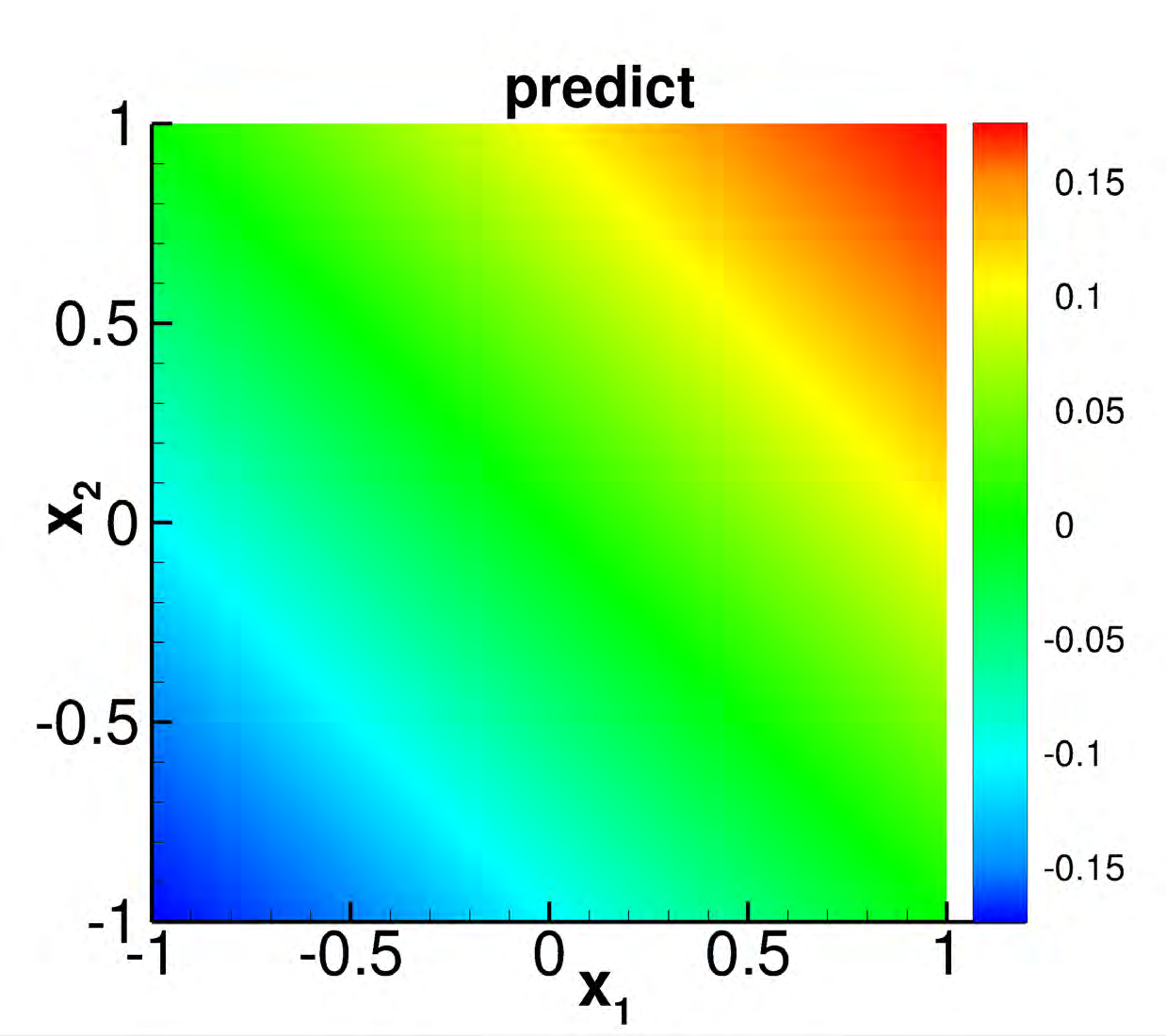}(b)
\includegraphics[width=0.26\textwidth]{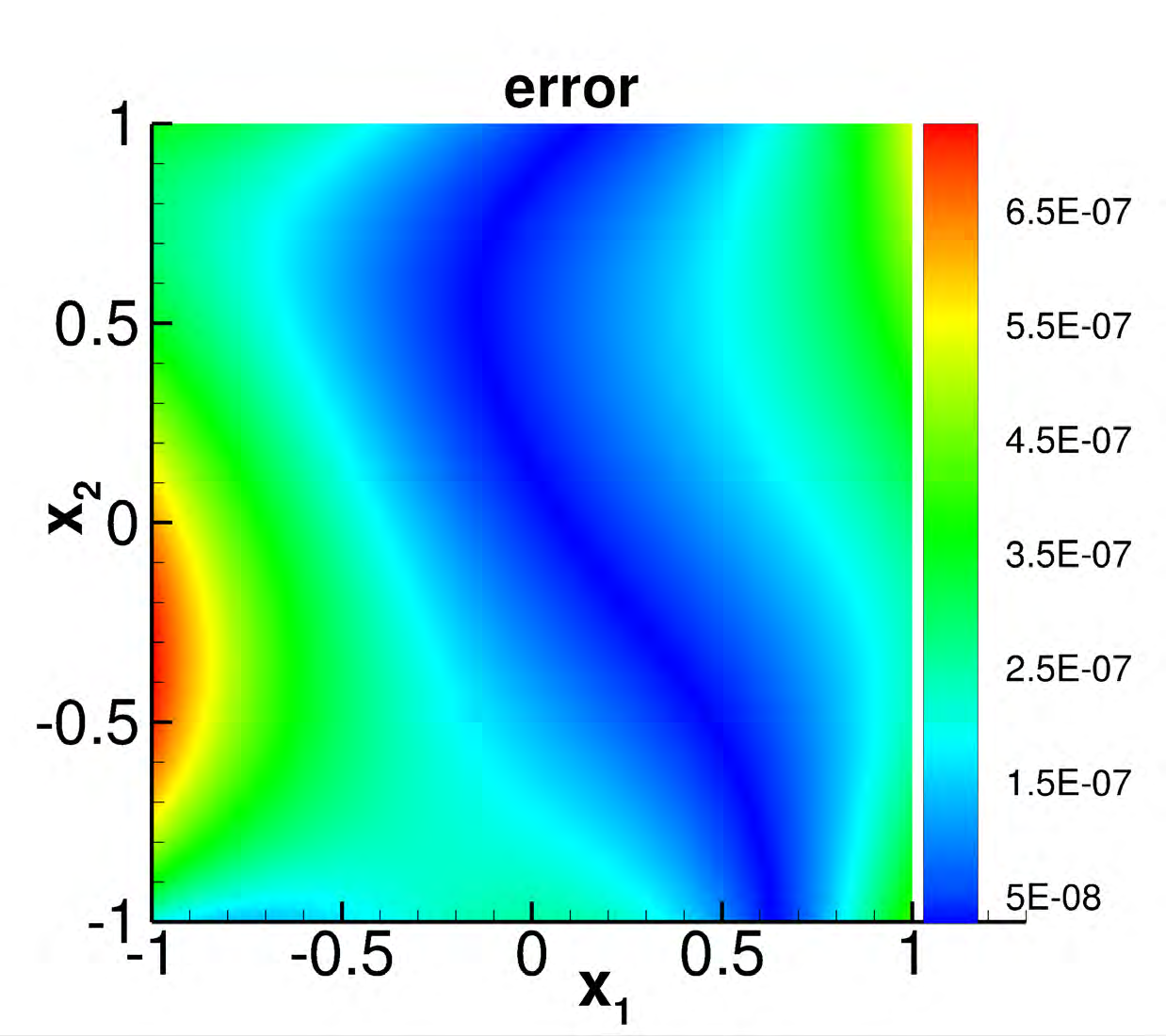}(c)
}
\centerline{
\includegraphics[width=0.26\textwidth]{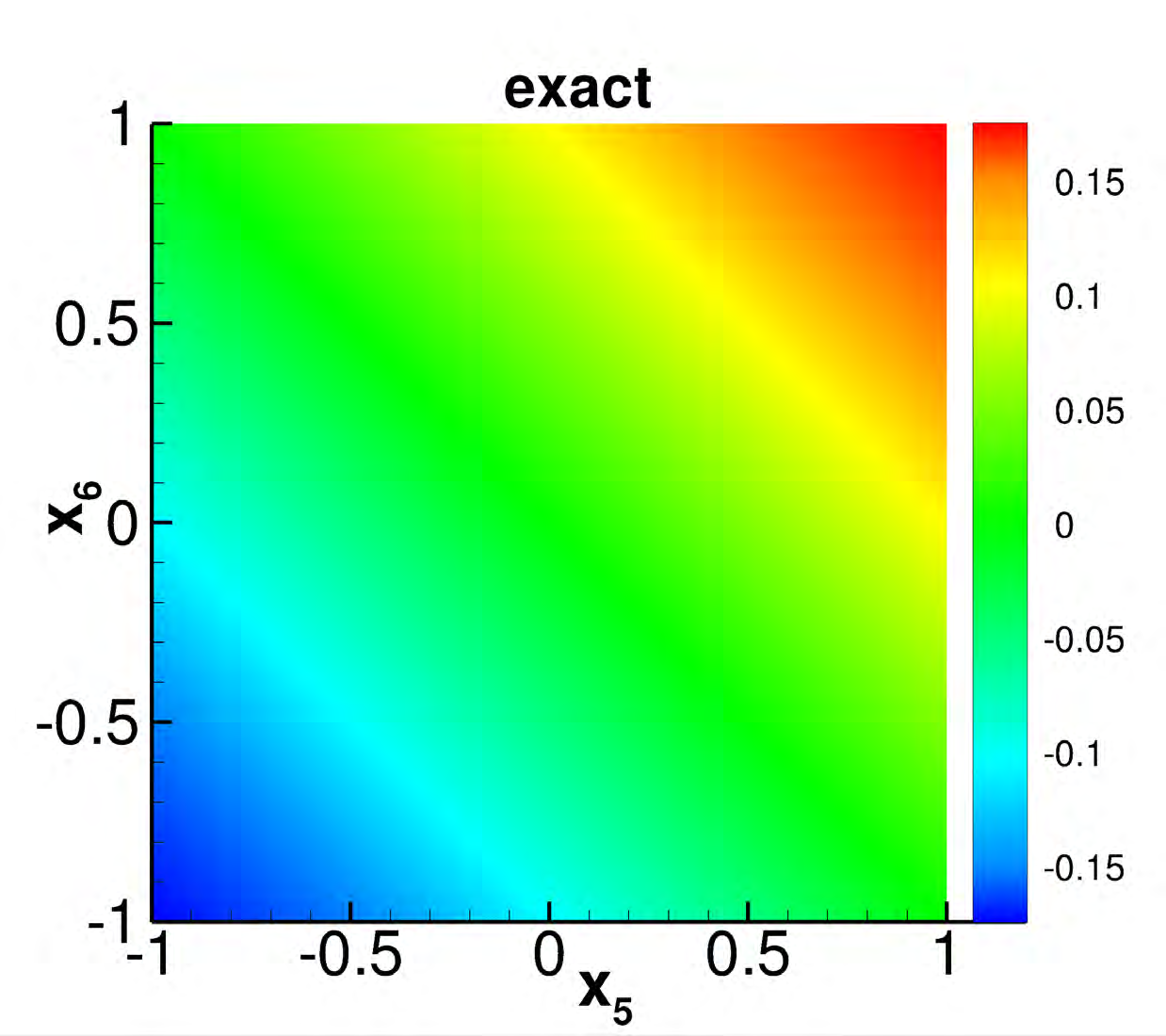}(d)
\includegraphics[width=0.26\textwidth]{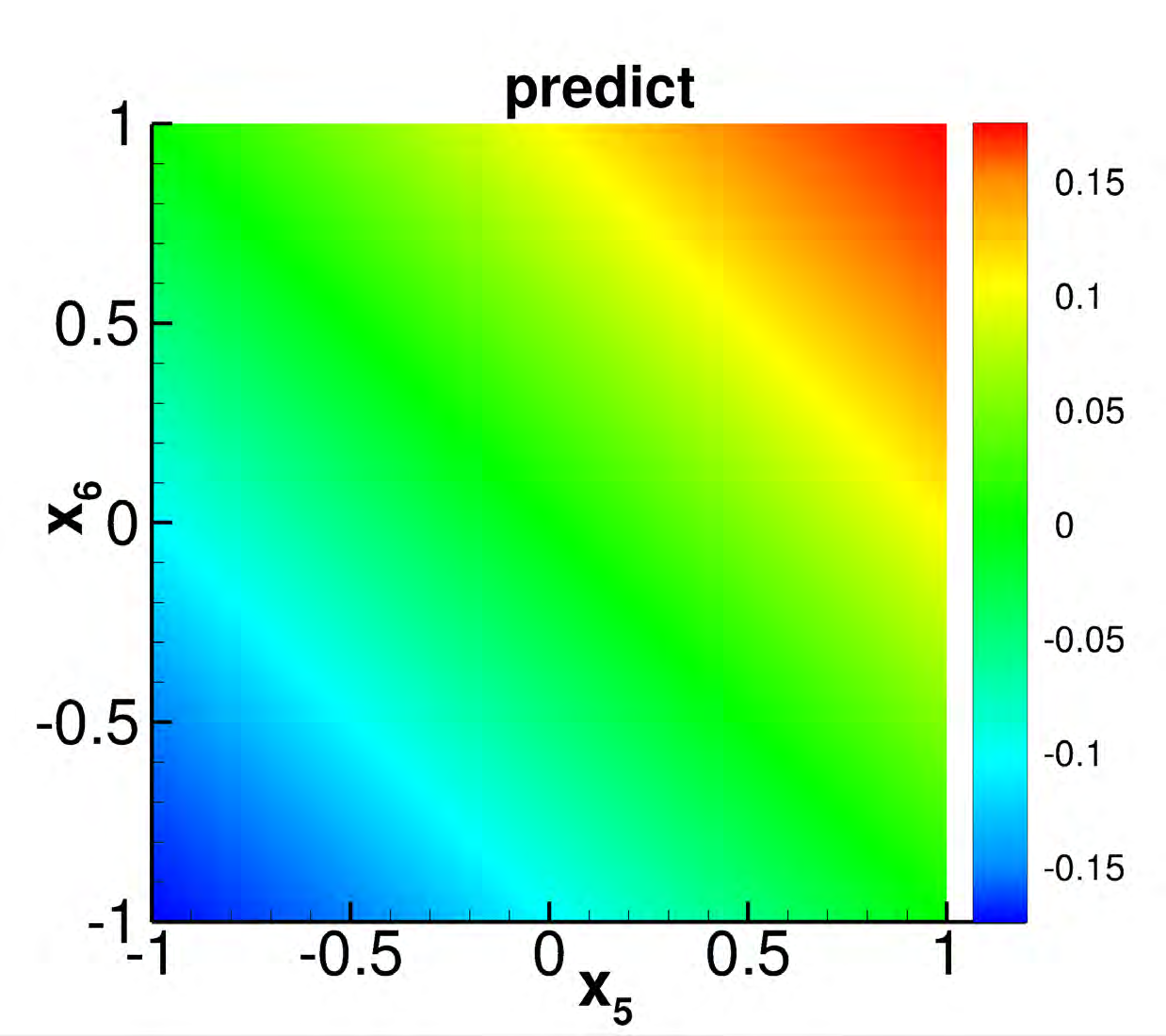}(e)
\includegraphics[width=0.26\textwidth]{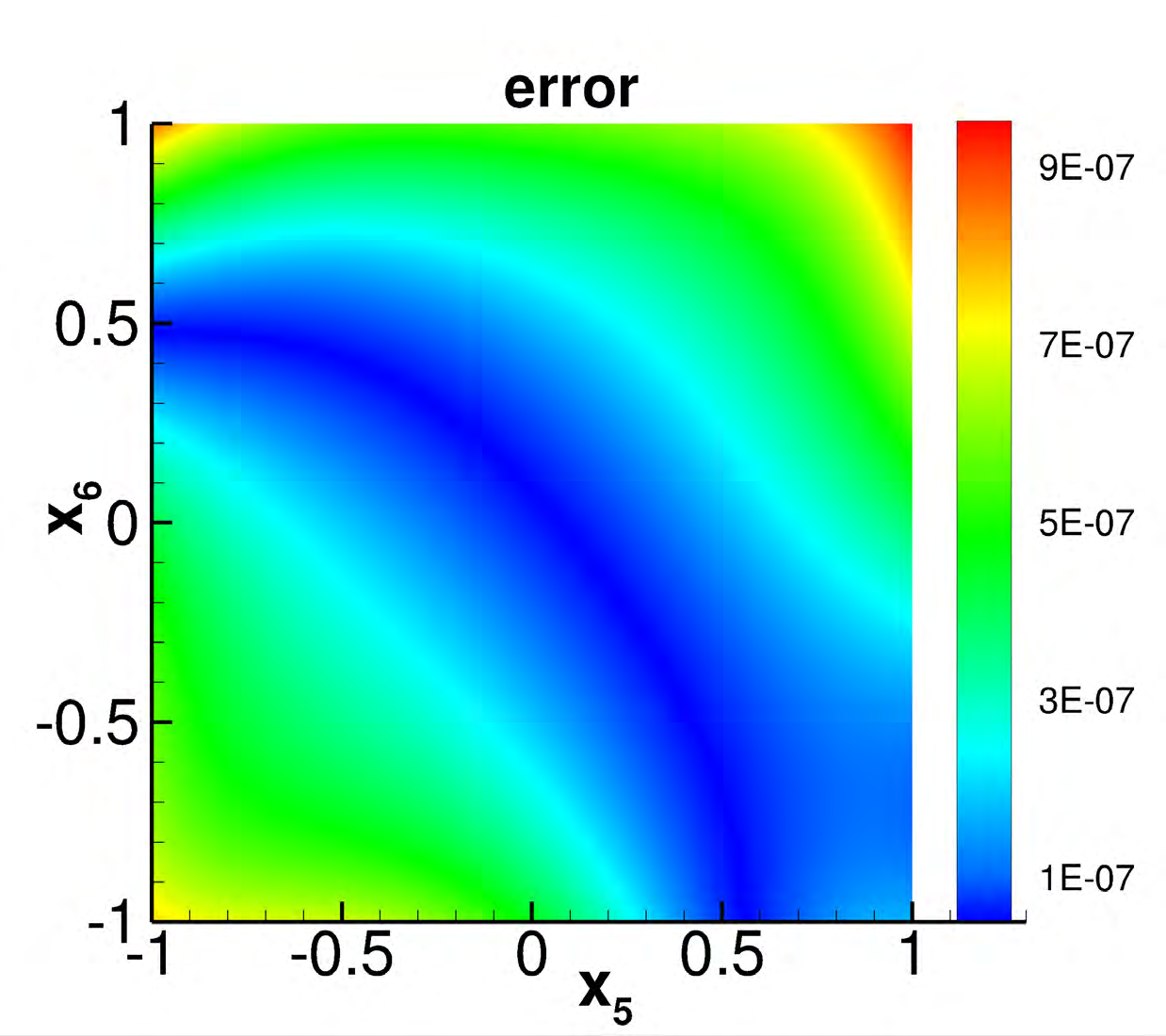}(f)
}
\centerline{
\includegraphics[width=0.26\textwidth]{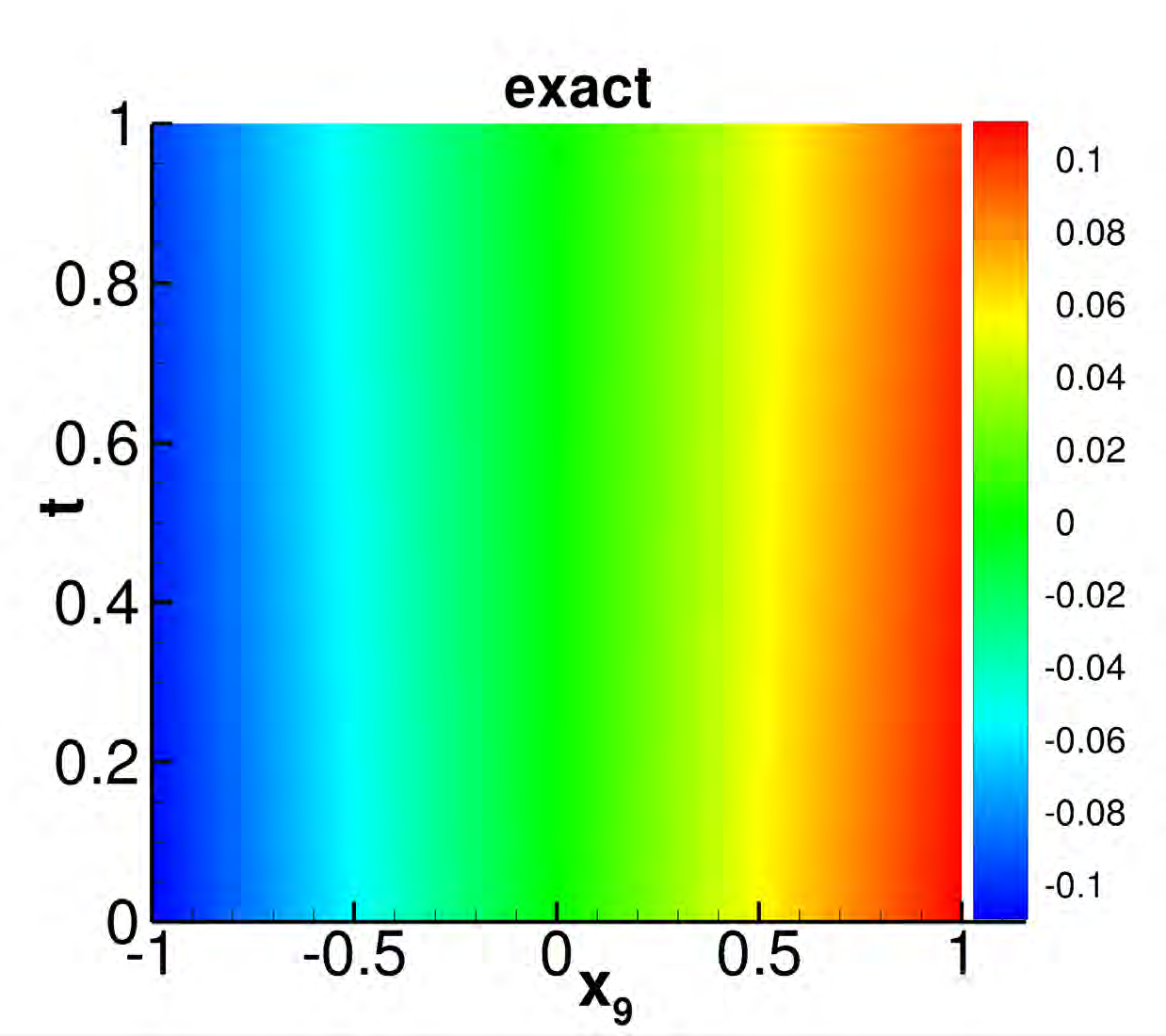}(g)
\includegraphics[width=0.26\textwidth]{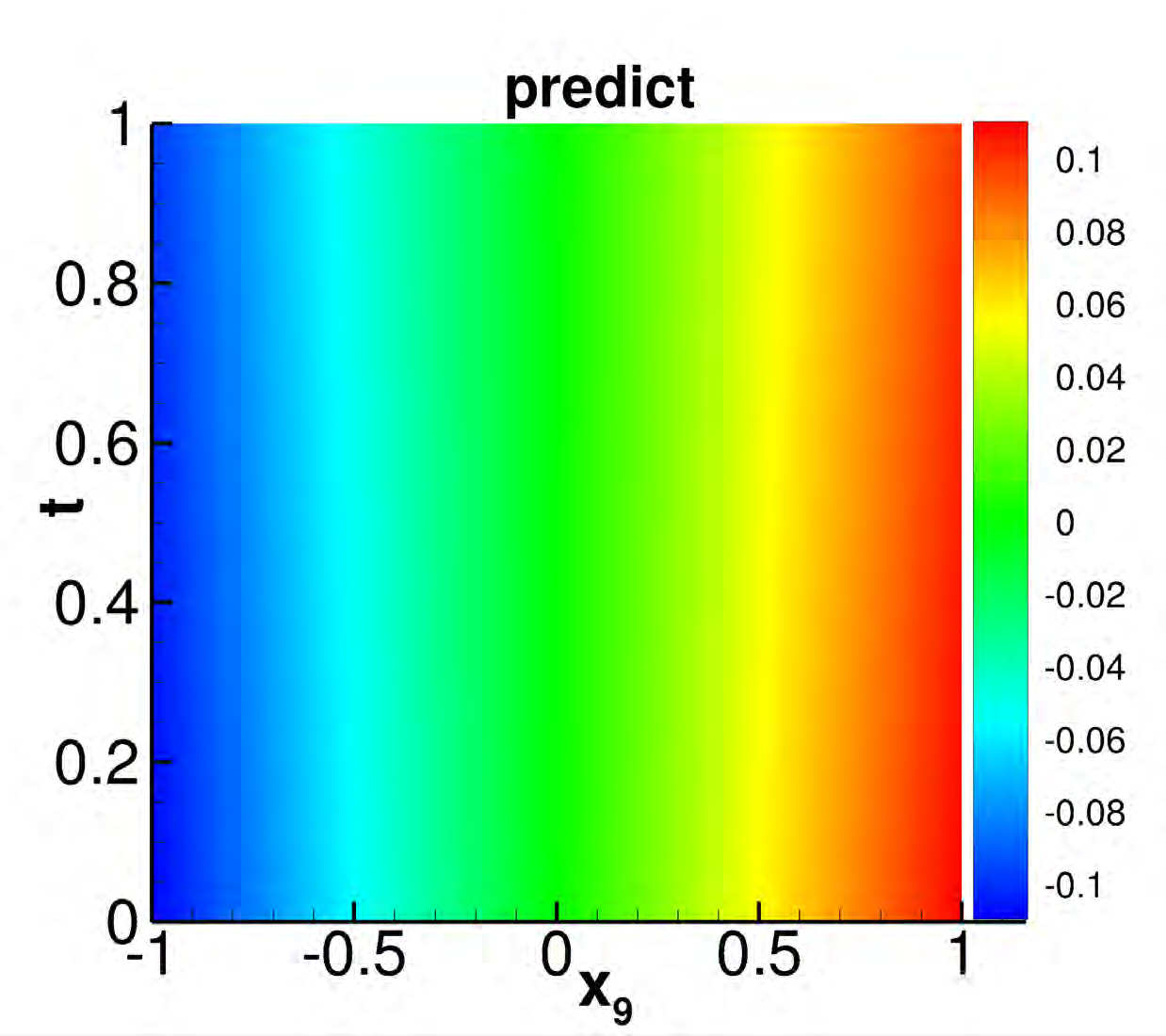}(h)
\includegraphics[width=0.26\textwidth]{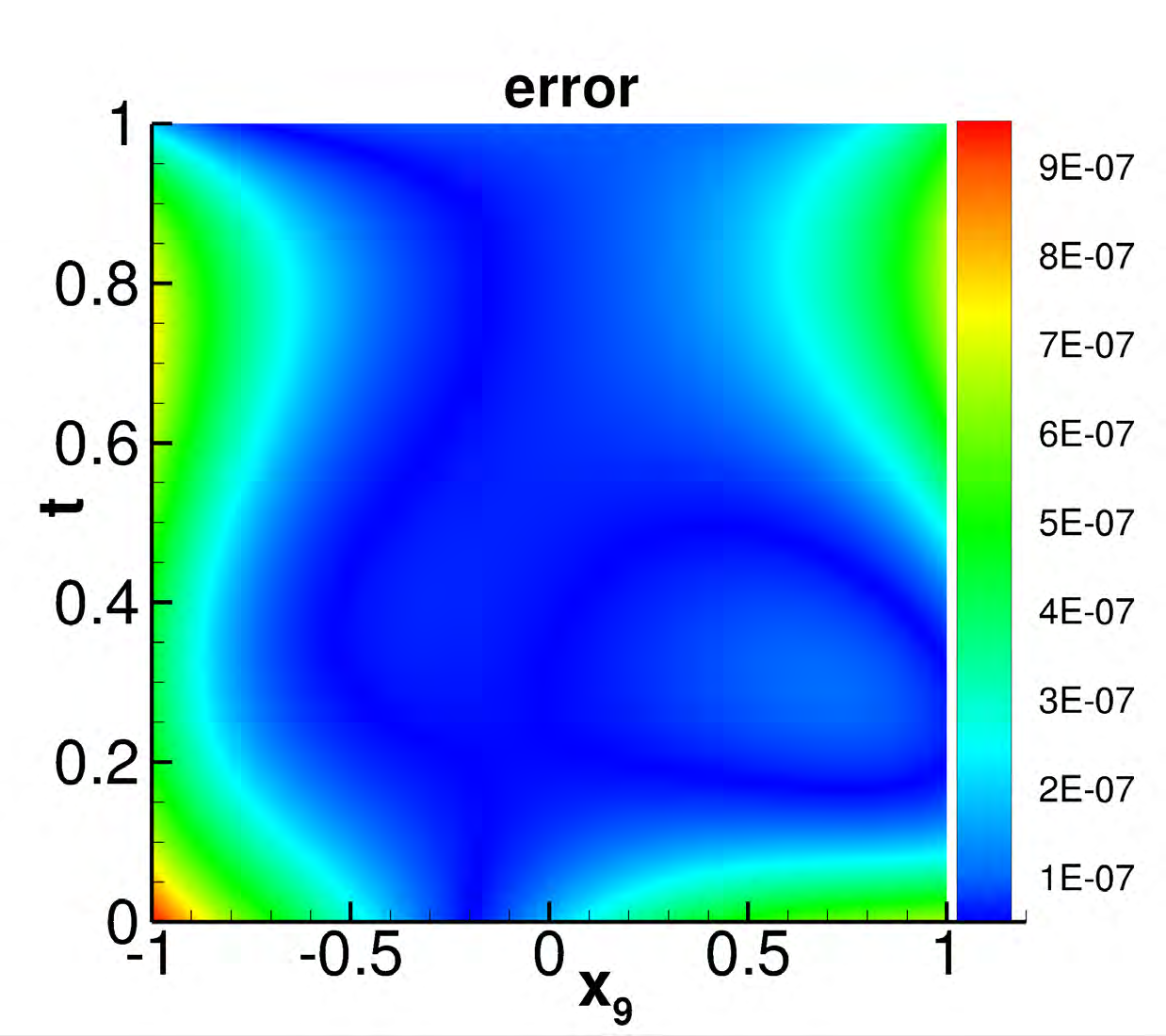}(i)
}
	\caption{\small Advection diffusion equation ($d=10$): Distributions of the exact solution (left column), the ELM solution (middle column), and the point-wise absolute error of ELM (right column) in selected cross sections of the domain, the $x_1$-$x_2$ plane (first row), the $x_5$-$x_6$ plane (middle row), and the $x_9$-$t$ plane (bottom row).
 NN architecture: $[d+1,5000,1]$; $(N_{\bcd},N_{\ins},N_{t_0})=(160,10,1000)$.
 For each cross section, the other coordinates of the plane are in the middle of the domain in each direction.
 }
 \label{fg_10}
\end{figure}

Figure~\ref{fg_10} shows distributions of the exact solution and the ELM solution, as well as the point-wise absolute error of the ELM solution, to the advection diffusion equation in dimension $d=10$ on several  cross sections of the spatial-temporal domain. These cross sections are the $x_1$-$x_2$ plane, the $x_5$-$x_6$ plane, and the $x_9$-$t$ plane. Each plane  is located in the middle of the spatial-temporal domain with respect to the rest of the coordinates. For example, the $x_i$-$x_j$ plane is characterized by $t=T/2$ and $x_k=0$ ($k\neq i,j$) for $1\leqslant k\leqslant d$, and the $x_i$-$t$ plane is characterized by $x_k=0$ ($k\neq i$) for $1\leqslant k\leqslant d$.
The network architecture and the other simulation parameters are provided in figure caption.
It is evident that ELM has captured the solution quite accurately, with the maximum error on the order of $10^{-7}$ on these cross sections.

\begin{figure}[!htb]
	\centering
	\subfigure[d=3]
	{ \includegraphics[width=0.3\textwidth]{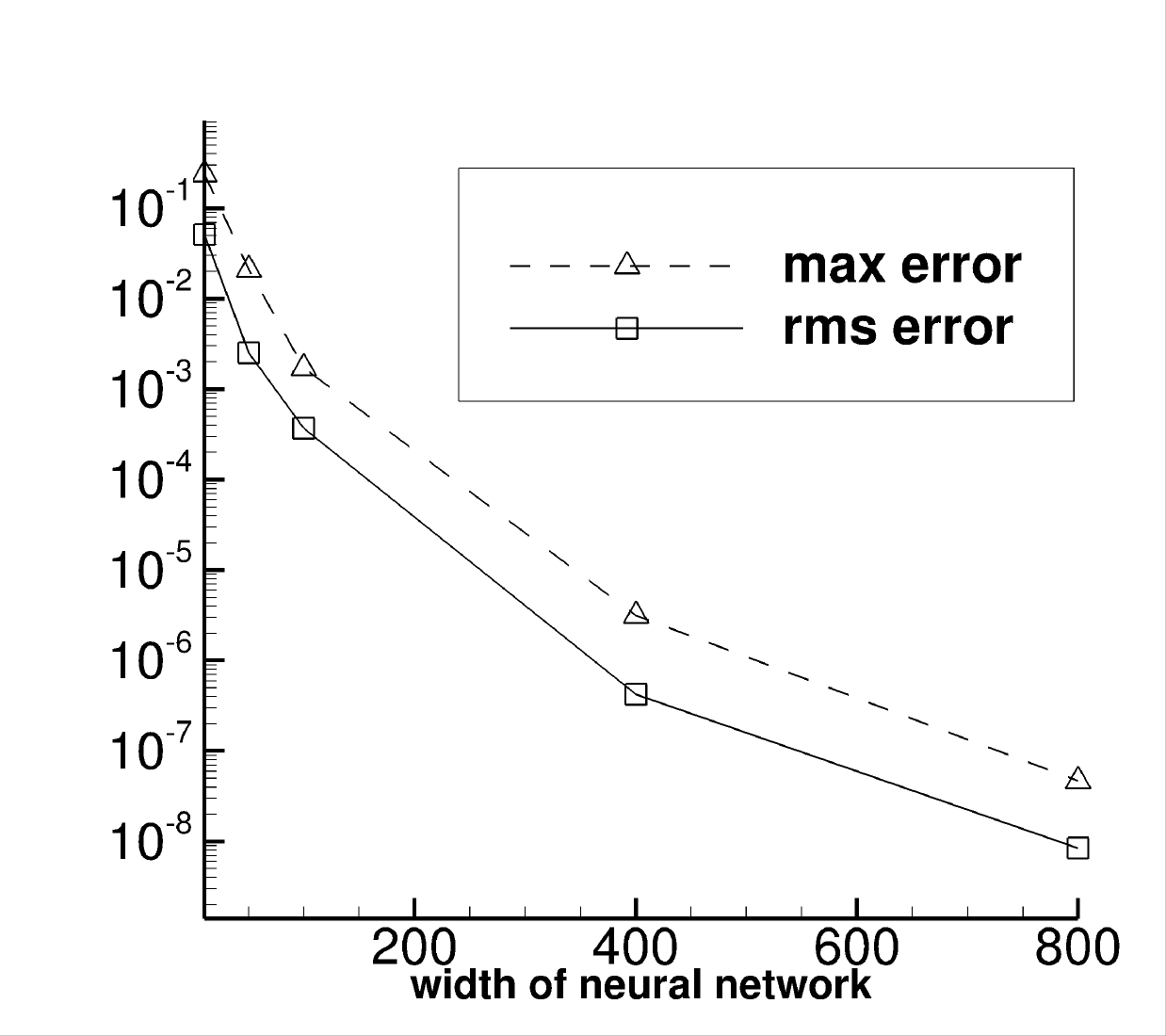}}
	\subfigure[d=6]
	{ \includegraphics[width=0.3\textwidth]{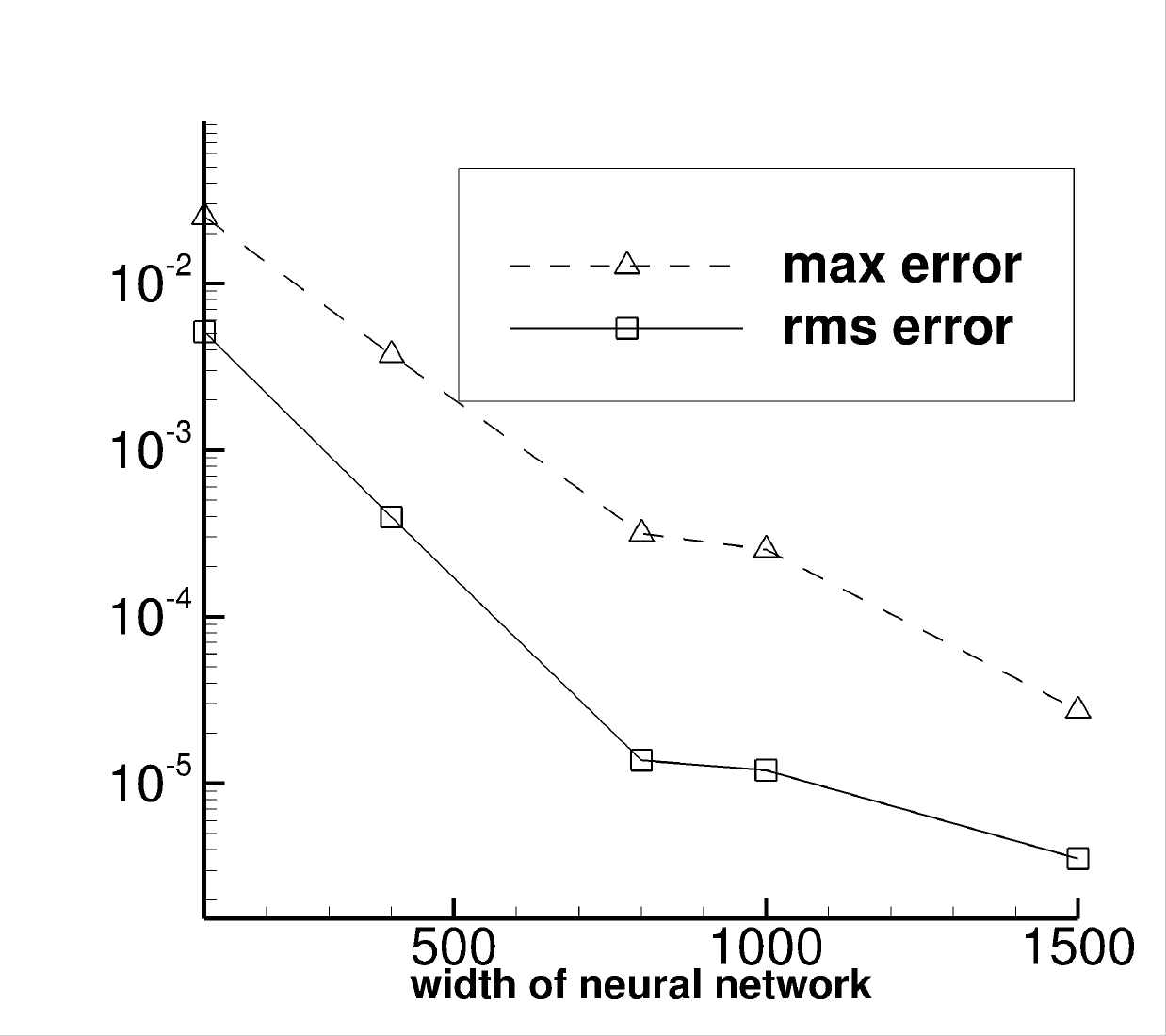}}
	\subfigure[d=10]
	{ \includegraphics[width=0.3\textwidth]{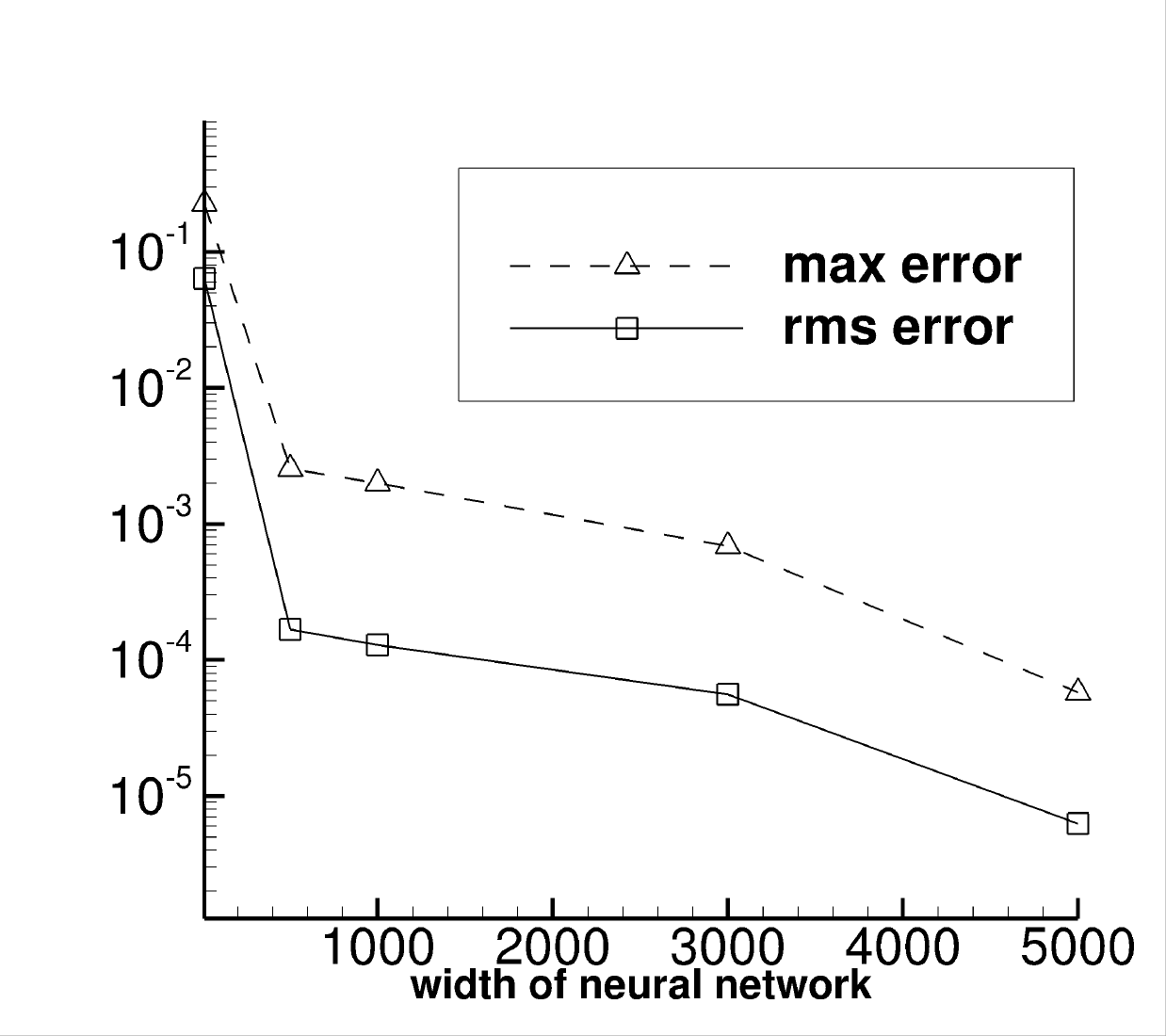}}
	\caption{\small Advection diffusion equation: $e_{l^{\infty}}$ and $e_{l^2}$ errors versus the number of training parameters ($M$) for problem dimensions (a) $d=3$, (b) $d=6$, and (c) $d=10$.
 NN architecture: $[d+1,M,1]$; $(N_{\bcd},N_{\ins},N_{t_0})=(200,10,1000)$ in (a), $(80,10,1000)$ in (b), and $(160,10,1000)$ in (c).
 $M$ is varied in (a,b,c).
 }
 \label{adv_err_width}
\end{figure}

\begin{figure}[!htb]
	\centering
	\subfigure[d=3]
	{ \includegraphics[width=0.3\textwidth]{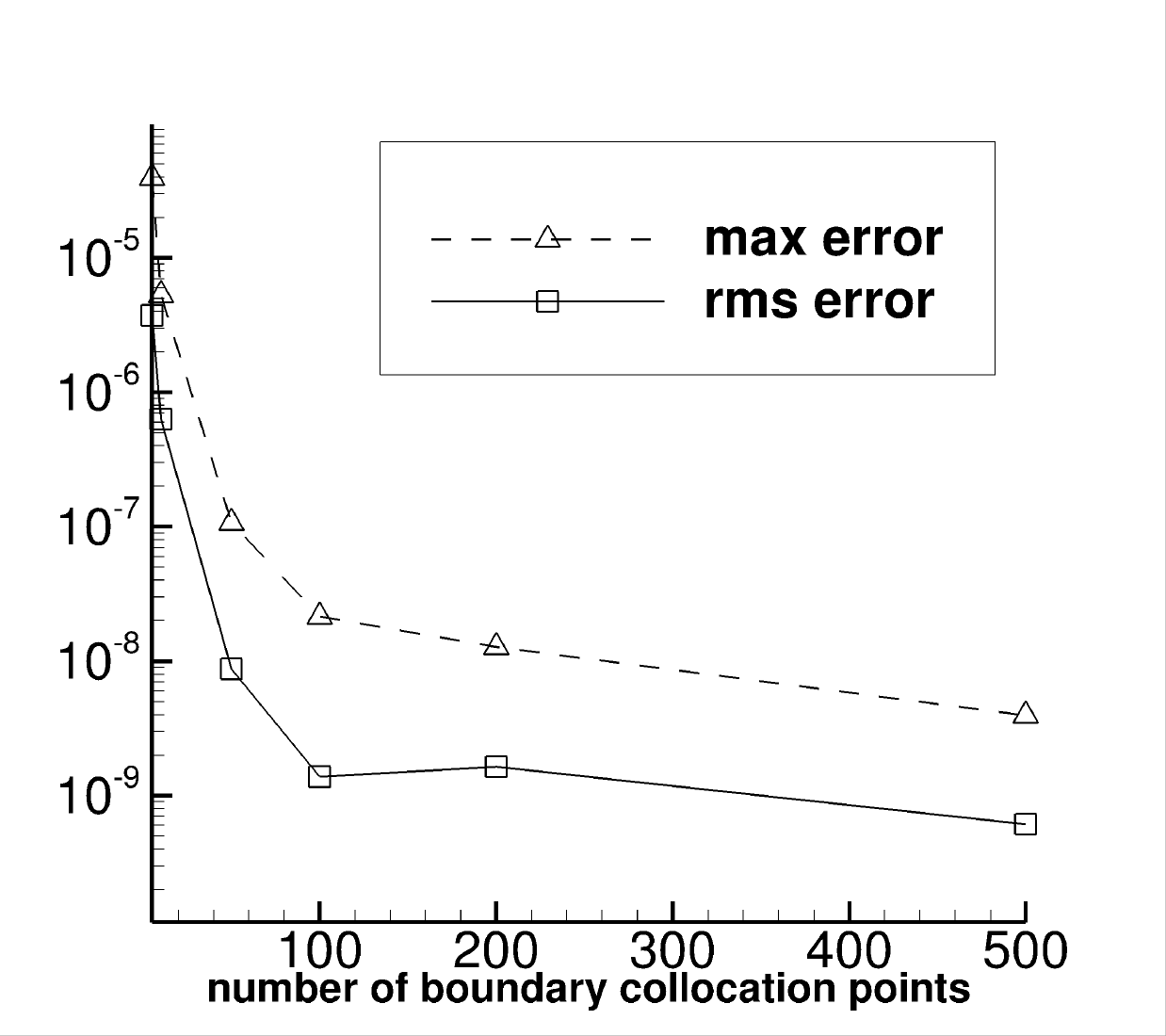}}
	\subfigure[d=6]
	{ \includegraphics[width=0.3\textwidth]{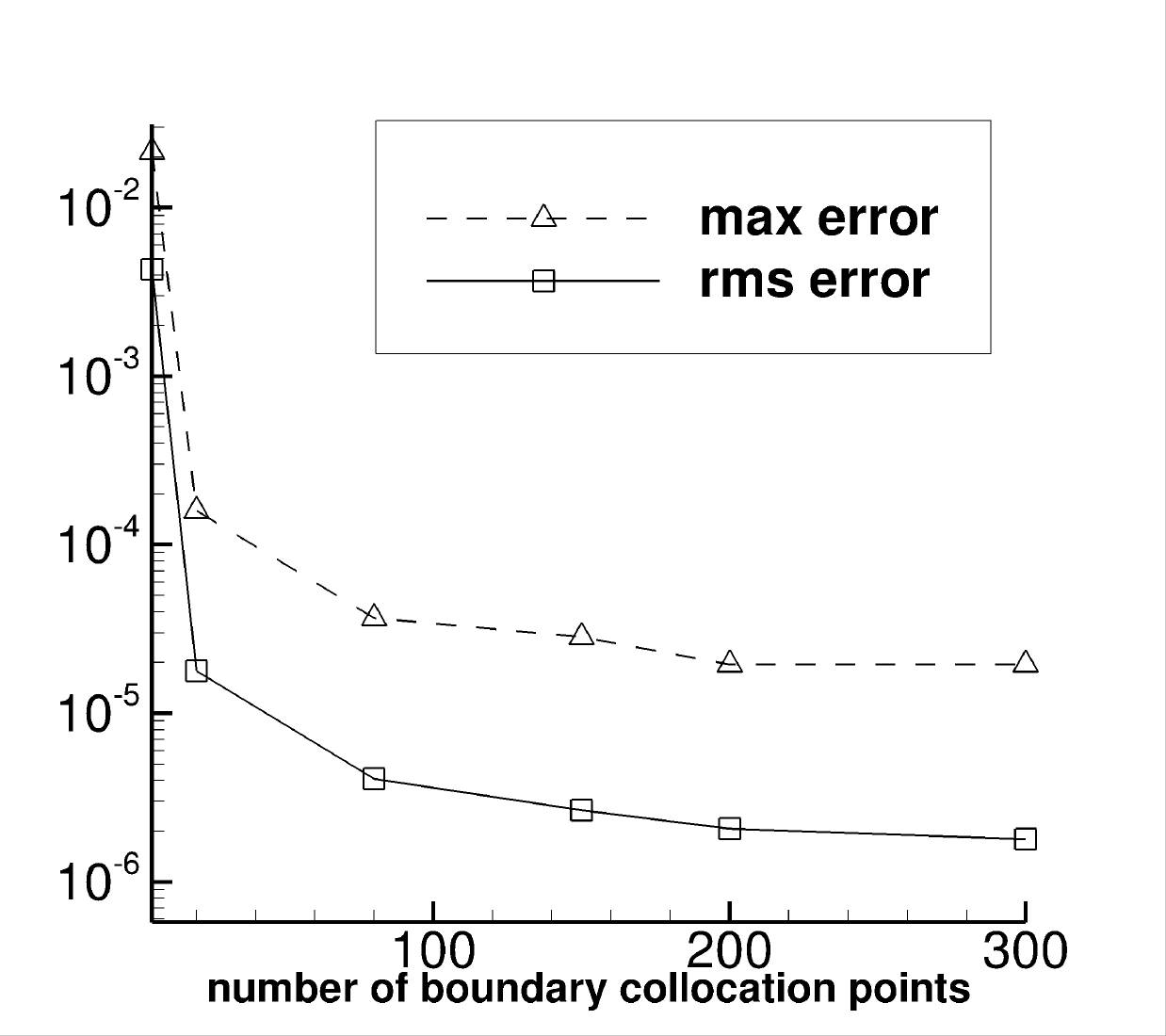}}
	\subfigure[d=10]
	{ \includegraphics[width=0.3\textwidth]{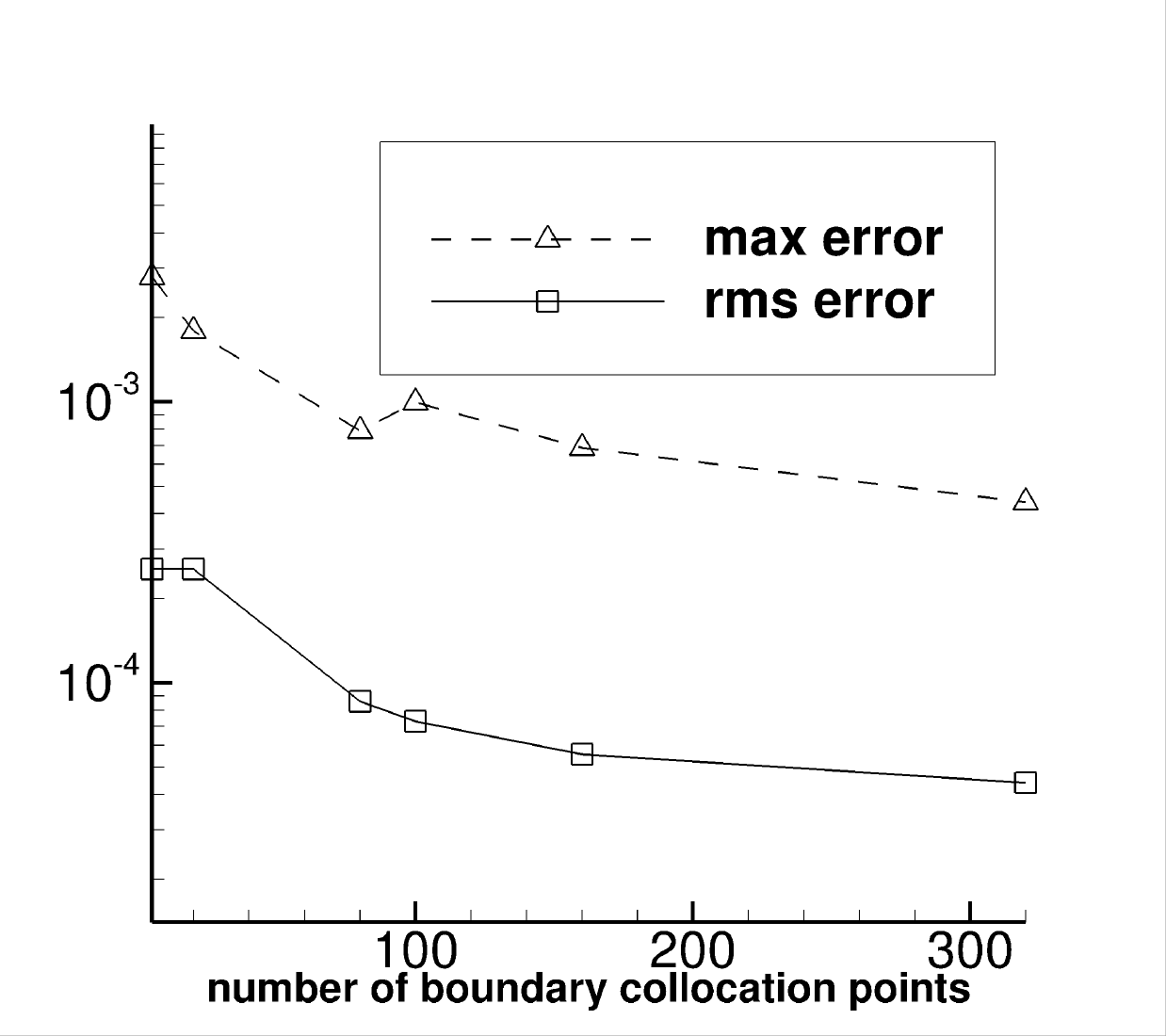}}
	\caption{\small Advection diffusion equation: $e_{l^{\infty}}$ and $e_{l^2}$ errors versus the number of boundary collocation points  ($N_{\bcd}$) for  dimensions (a) $d=3$, (b) $d=6$, and (c) $d=10$.
 NN architecture: $[d+1,M,1]$ with $M=1000$ in (a), $M=2000$ in (b), and $M=3000$ in (c). 
 $(N_{\ins},N_{t_0})=(10,1000)$ in (a,b,c). $N_{\bcd}$ is varied in (a,b,c). 
 }\label{adv_err_nc}
\end{figure}

Figure~\ref{adv_err_width} illustrates the effect of the trainable parameters on the ELM accuracy. Here we show the $e_{l^{\infty}}$ and $e_{l^2}$ errors of ELM versus the number of training parameters $M$ in the neural network for solving the advection-diffusion equation in dimensions $d=3$, $6$ and $10$. The network architecture is given by $[d+1,M,1]$, where $M$ is varied in the tests. The other simulation parameters are listed in the figure caption. 
The ELM errors can be observed to decrease dramatically (close to exponential rate) with increasing number of training parameters. The $e_{l^2}$ (rms) error levels are on the order of $10^{-9}$ (for $d=3$), and $10^{-6}$ (for $d=6$ and $10$) for the range of parameters tested here.

The effect of the boundary collocation points on the ELM accuracy is illustrated in Figure~\ref{adv_err_nc}. The number of interior collocation points, on the other hand, has little (or much less) influence on the ELM results compared with the boundary points. 
Figure~\ref{adv_err_nc} shows the $e_{l^{\infty}}$ and $e_{l^2}$ errors versus $N_{\bcd}$ (collocation points on each boundary) for dimensions $d=3$, $6$ and $10$. As $N_{\bcd}$ increases, the errors appear to decrease approximately exponentially, and then level off when $N_{\bcd}$ reaches a certain level.



\subsubsection{Korteweg-De Vries Equation}

In this subsection we consider the Korteweg-De Vries (KdV) equation,
\begin{subequations}\label{eq_46}
\begin{align}
\partial_t u +\sum_{i=1}^d \partial_{x_ix_ix_i}^3 u &= f(x,t), \quad (x,t)\in \Omega \times [0,T],\\
u(x,t)&=g(x,t),\quad (x,t)\in \partial \Omega \times [0,T],\\
u(x,0)&=h(x),\qquad x\in \Omega,
\end{align}
\end{subequations}
where $\Omega=[-1,1]^d$ and $T=1$. In these equations $f(x,t)=-\frac{1}{d^2}\left[\sin\left(\frac{1}{d}\sum_{i=1}^{d}x_i\right)+\cos\left(\frac{1}{d}\sum_{i=1}^dx_i\right)\right]\exp\left(-\frac{t}{d^2}\right)$, $g(x,t)=\sin\left(\frac{1}{d}\sum_{i=1}^dx_i\right)\exp\left(-\frac{t}{d^2}\right)$ on $\partial \Omega\times [0,T]$, and $h(x)=\sin\left(\frac{1}{d}\sum_{i=1}^dx_i\right)$ in $\Omega$. This problem has an exact solution $u(x,t)=\sin\left(\frac{1}{d}\sum_{i=1}^dx_i\right)\exp\left(-\frac{t}{d^2}\right)$.
The notations below follow those of the previous subsections.

\begin{table}[tb]
\centering\small
\begin{tabular}{cccccc}
	\hline
	$R_m$&1E-3&5E-3&1E-2&5E-2&1E-1	\\ 
	$e_{l^\infty}$&3.15E-4&1.48E-5&5.58E-6&2.51E-7&1.05E-6\\
	$e_{l^{2}}$ &2.31E-5&1.90E-6&2.97E-7&3.55E-8&1.71E-7\\
	\hline
\end{tabular}
\caption{\small KdV equation: Determining $R_{m0}$ based on the procedure of Remark~\ref{rem_a6} for $d=5$.
NN architecture: $[d+1,2000,1]$; $(N_{\ins},N_{\bcd},N_{t_0})=(100,100,1000)$.\\
}
\label{tab_7}
\end{table}

\begin{table}[tbp]
\centering\small
\begin{tabular}{l| ccc}
\hline
$d$ & 3 & 5 & 10 \\ \hline
$R_{m0}$ & 0.05 & 0.05 & 0.05 \\
\hline
\end{tabular}
\caption{\small KdV equation: $R_{m0}$ determined by the procedure from Remark~\ref{rem_a6} for several problem dimensions. In these tests
the NN architecture and the number of collocation points are the same as those in Table~\ref{tab_7} for $d=3$ and $5$. For $d=10$, the NN architecture is $[d+1,3000,1]$, and $(N_{\ins},N_{\bcd},N_{t_0})=(10,100,1000)$.
}
\label{tab_8}
\end{table}

Tables~\ref{tab_7} documents the tests for determining the $R_{m0}$ for $d=5$ using the procedure from Remark~\ref{rem_a6}, and Table~\ref{tab_8} lists the resultant $R_{m0}$ values corresponding to the dimensions $d=3$, $5$ and $10$ from this procedure. The simulation parameters employed in these tests are provided in the captions of these figures. We use $R_m=R_{m0}$ for generating the random hidden-layer coefficients in the following tests with ELM.

\begin{figure}[!tb]
	\centerline{
\includegraphics[width=0.26\textwidth]{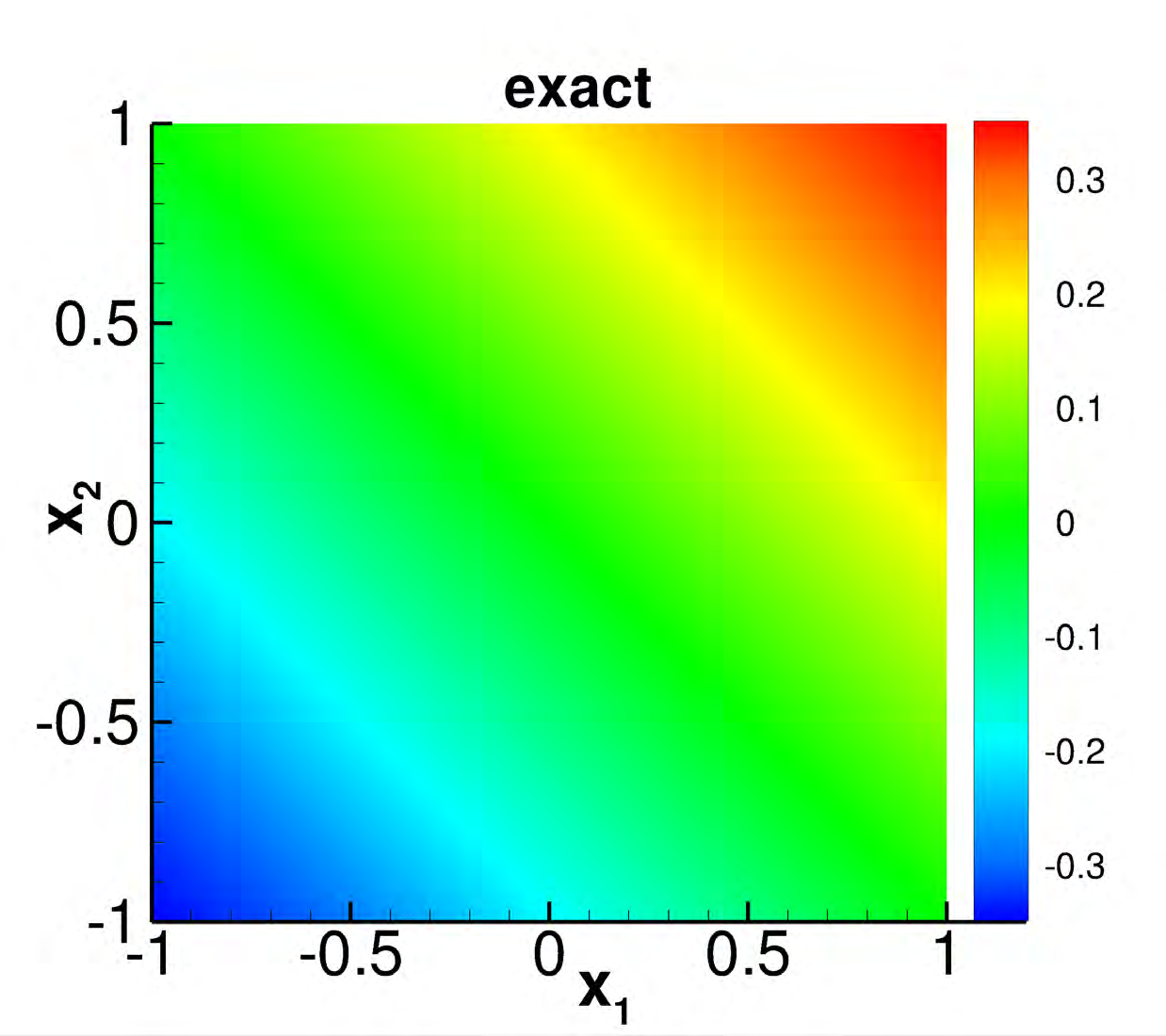}(a)
\includegraphics[width=0.26\textwidth]{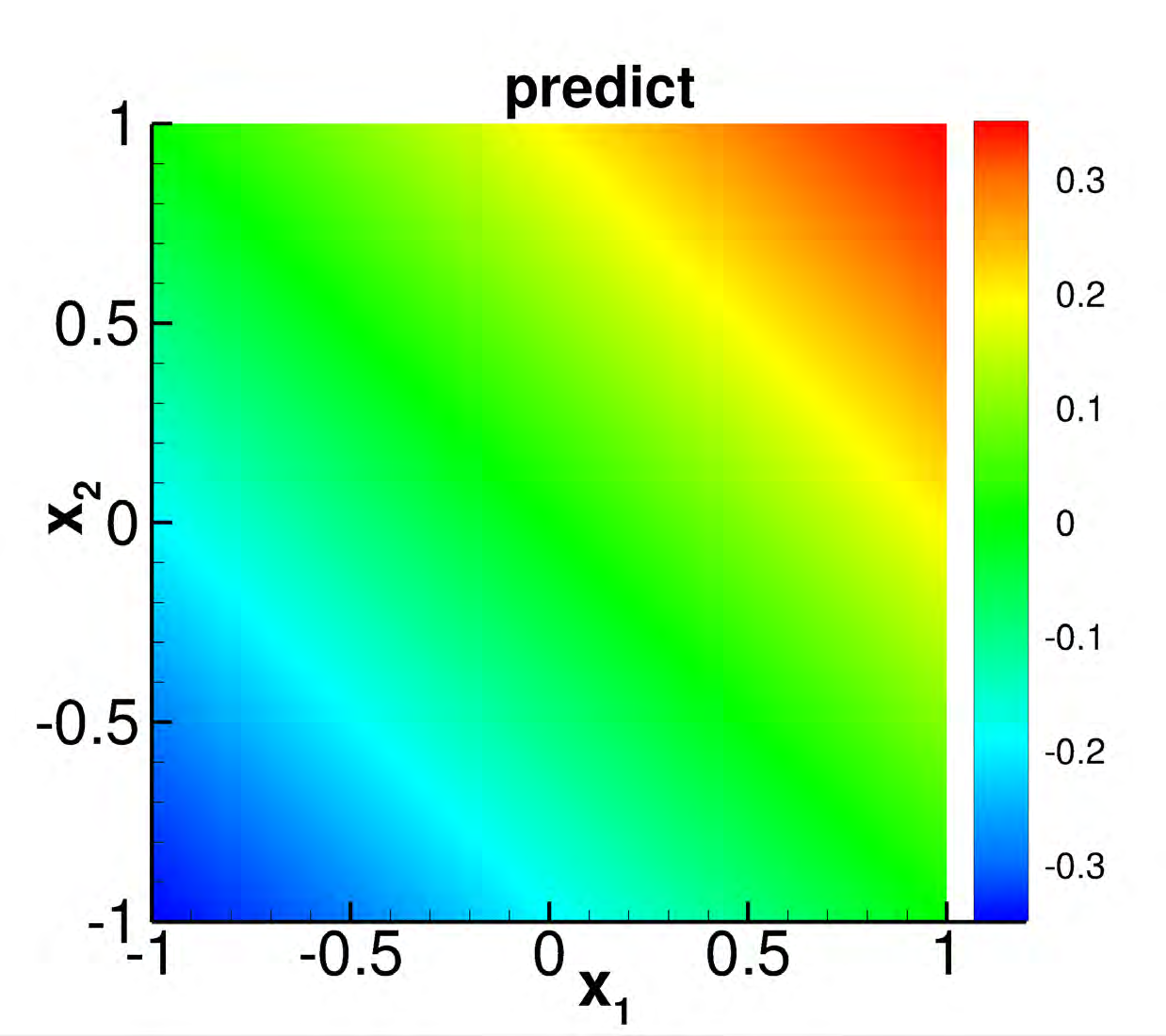}(b)
\includegraphics[width=0.26\textwidth]{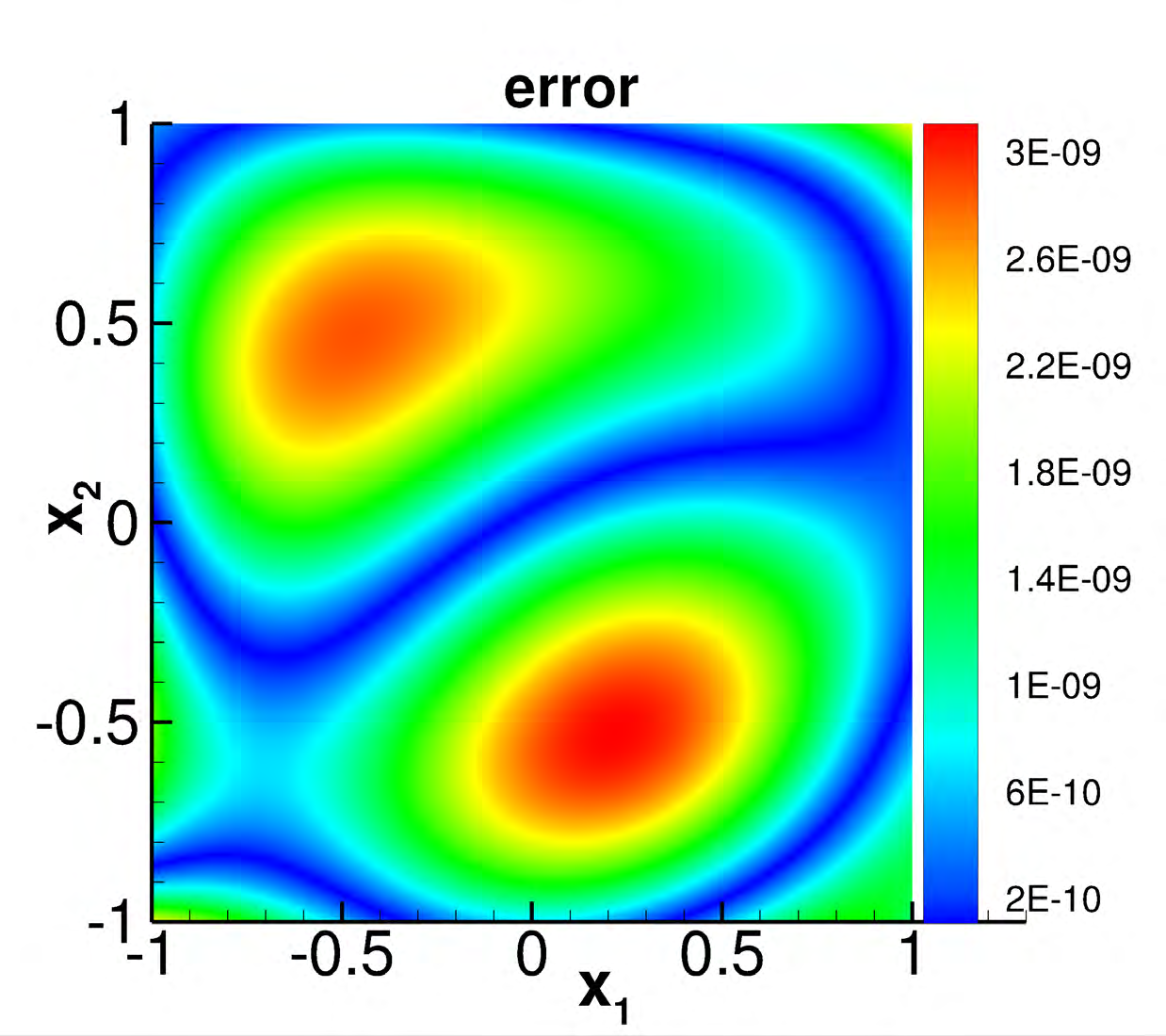}(c)
}
\centerline{
\includegraphics[width=0.26\textwidth]{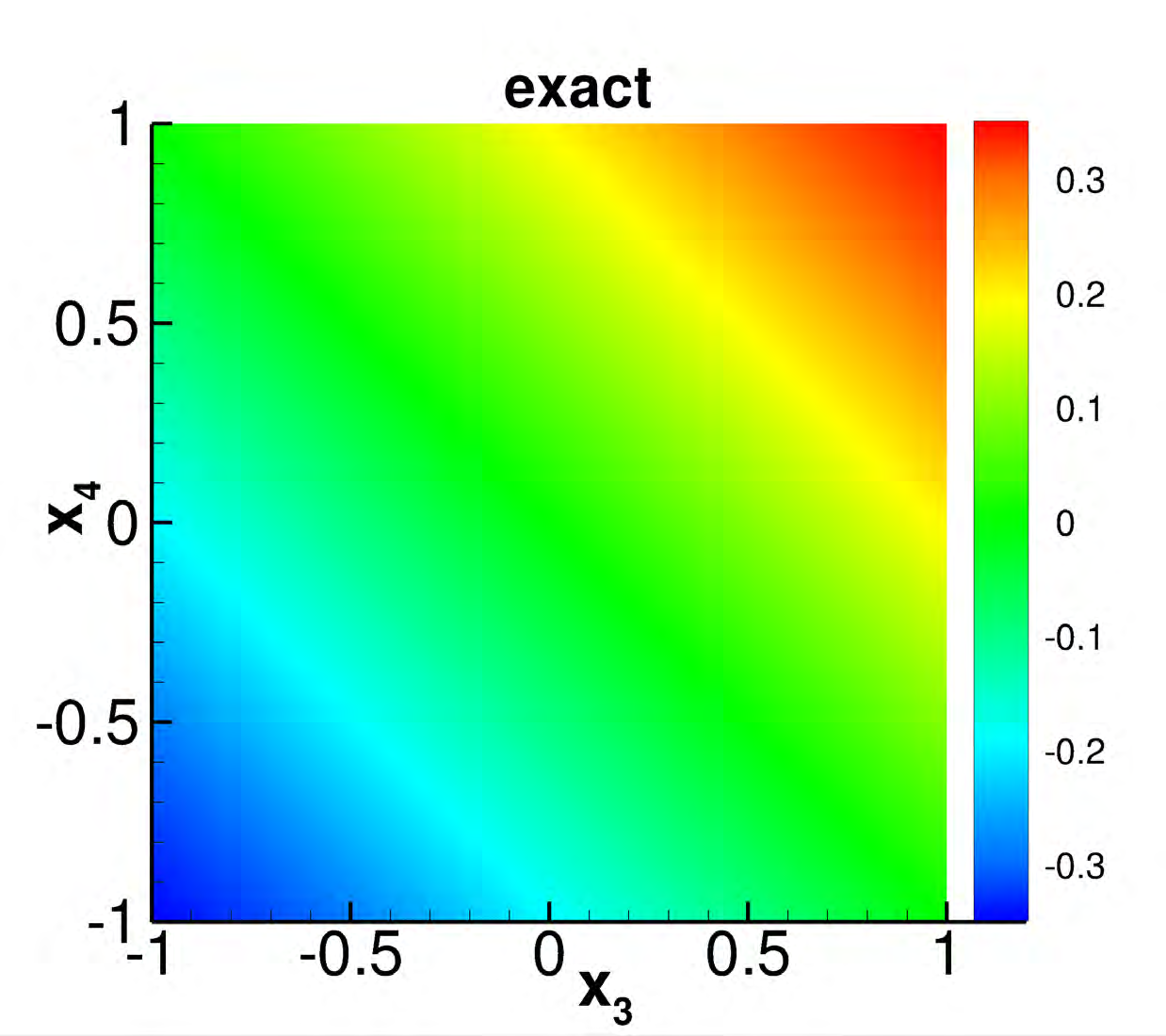}(d)
\includegraphics[width=0.26\textwidth]{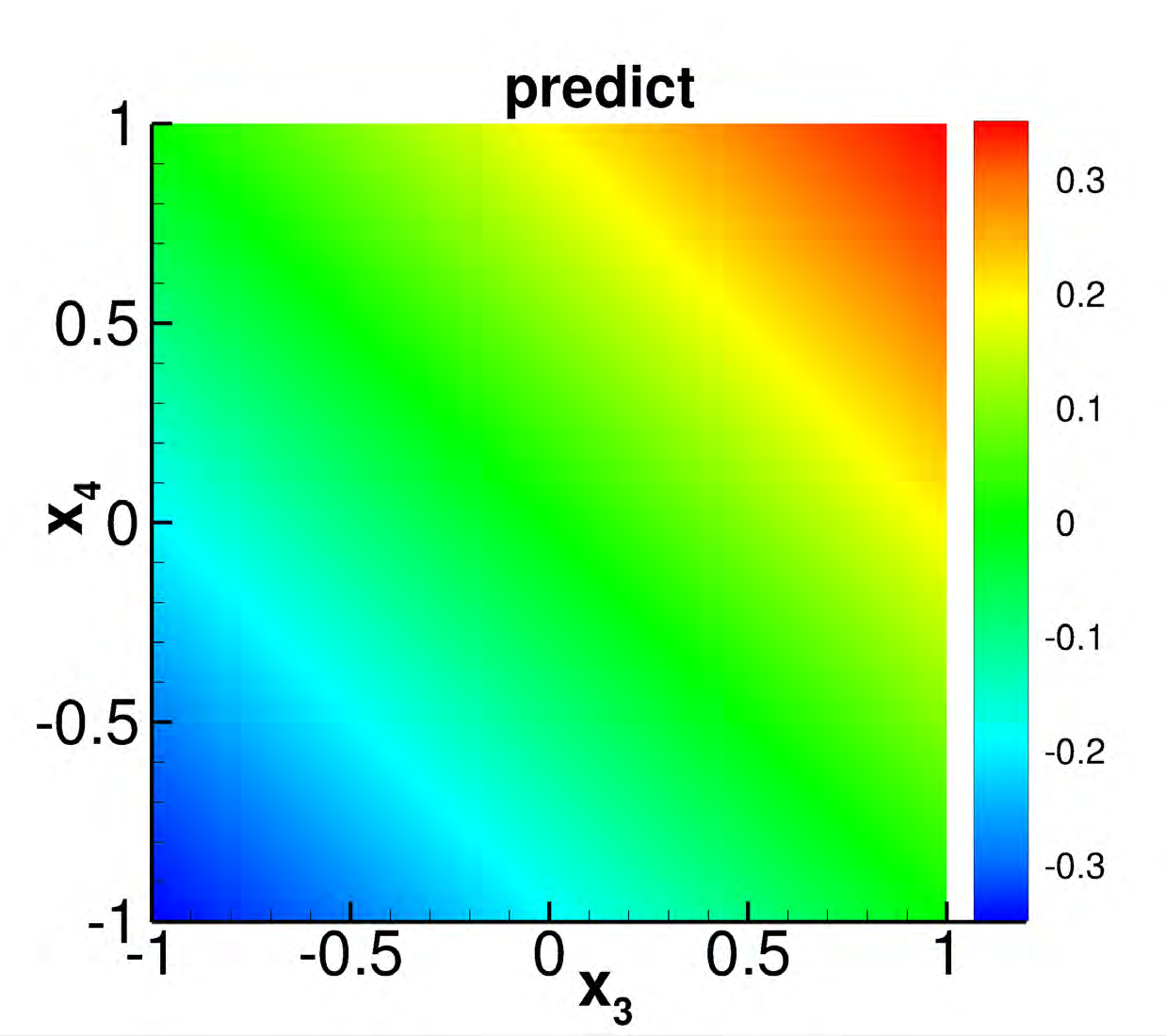}(e)
\includegraphics[width=0.26\textwidth]{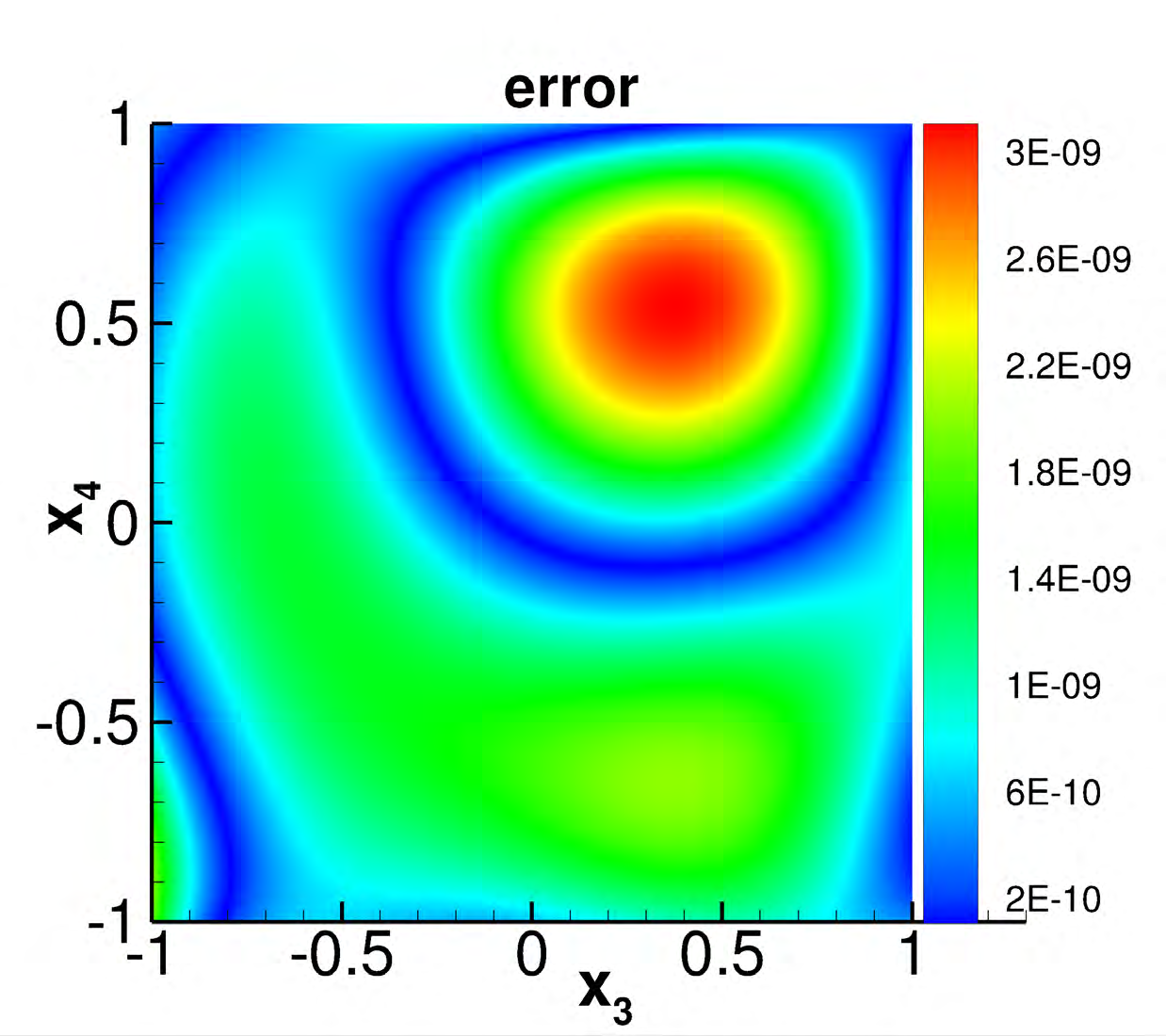}(f)
}
\centerline{
\includegraphics[width=0.26\textwidth]{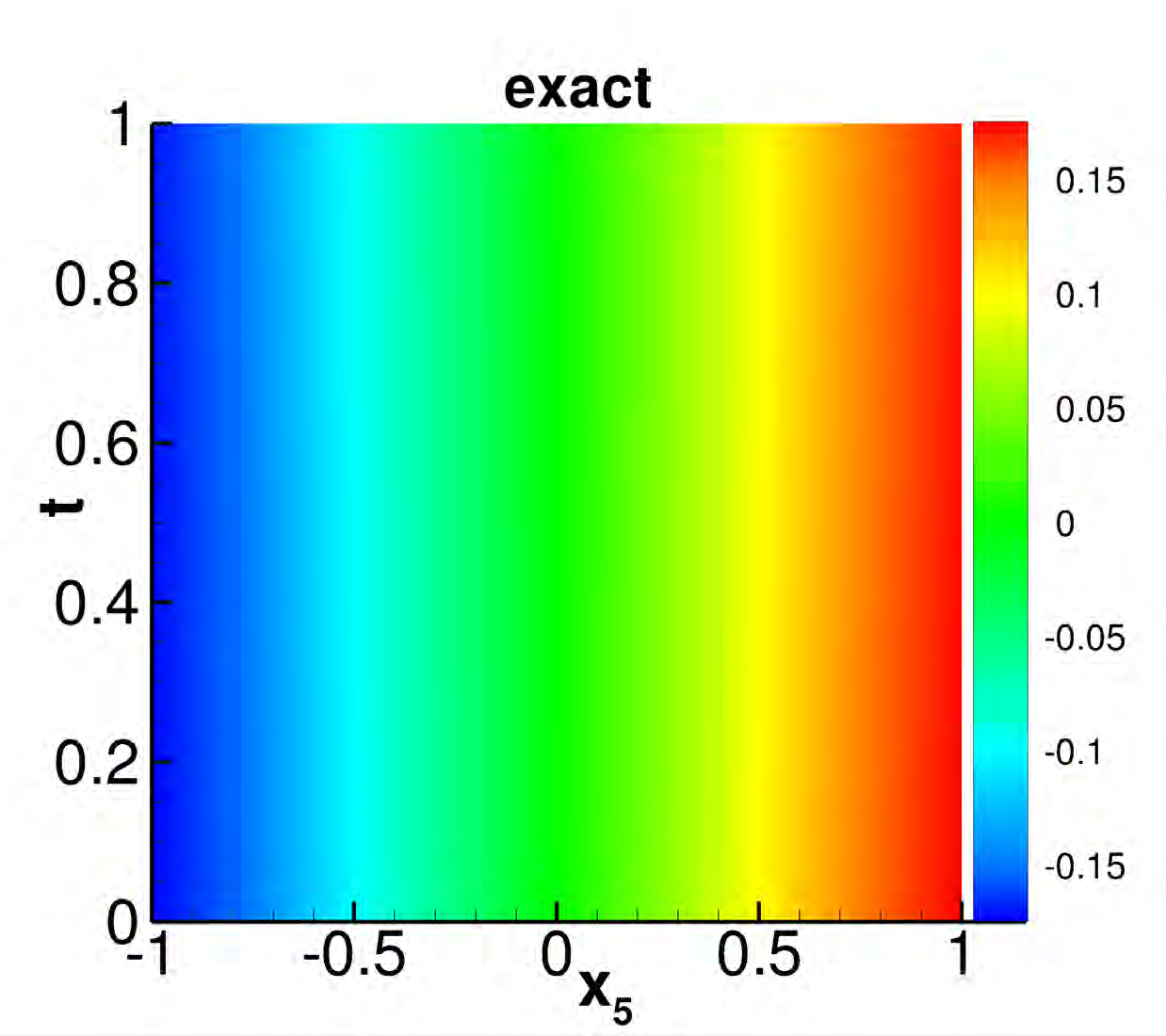}(g)
\includegraphics[width=0.26\textwidth]{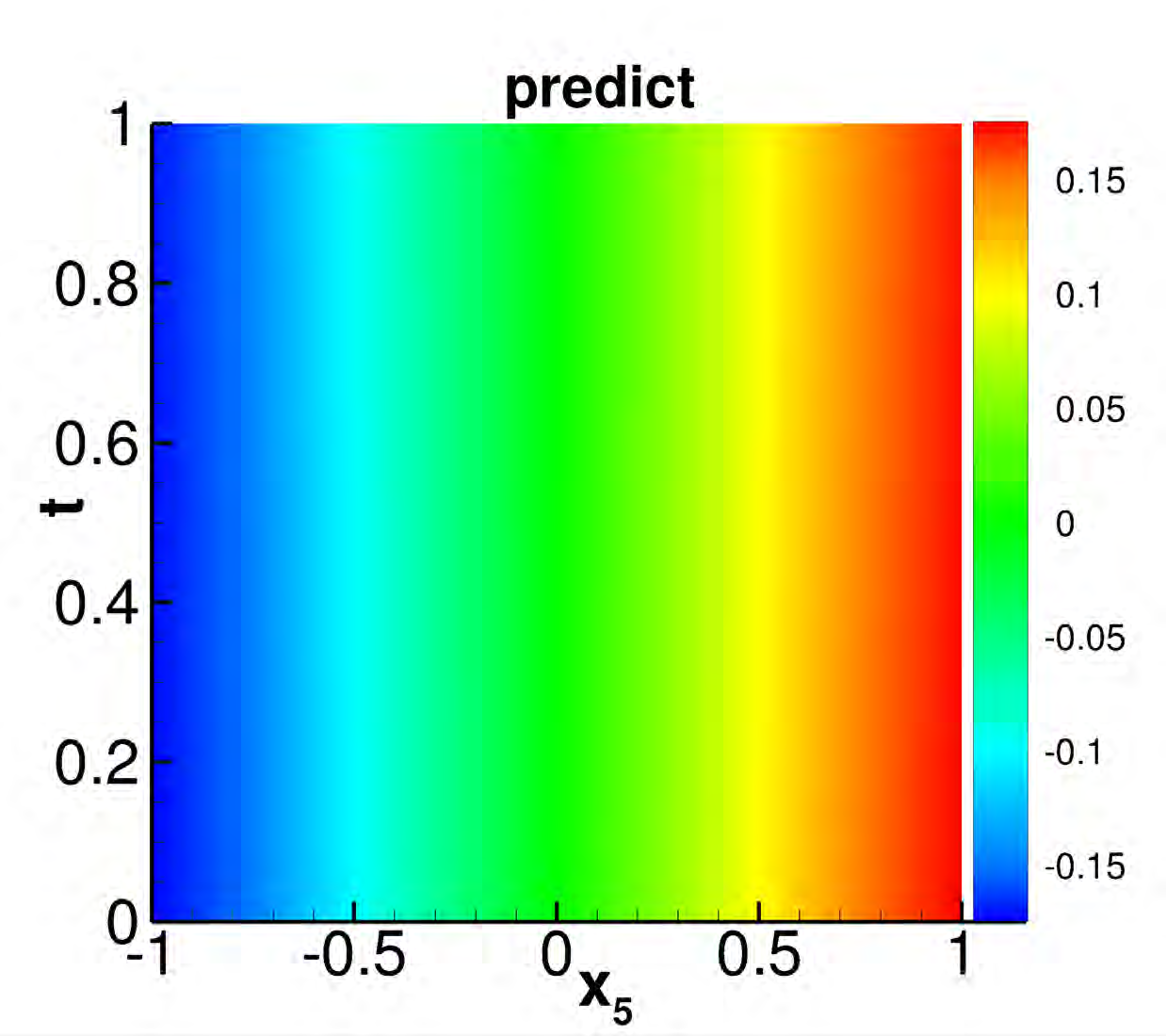}(h)
\includegraphics[width=0.26\textwidth]{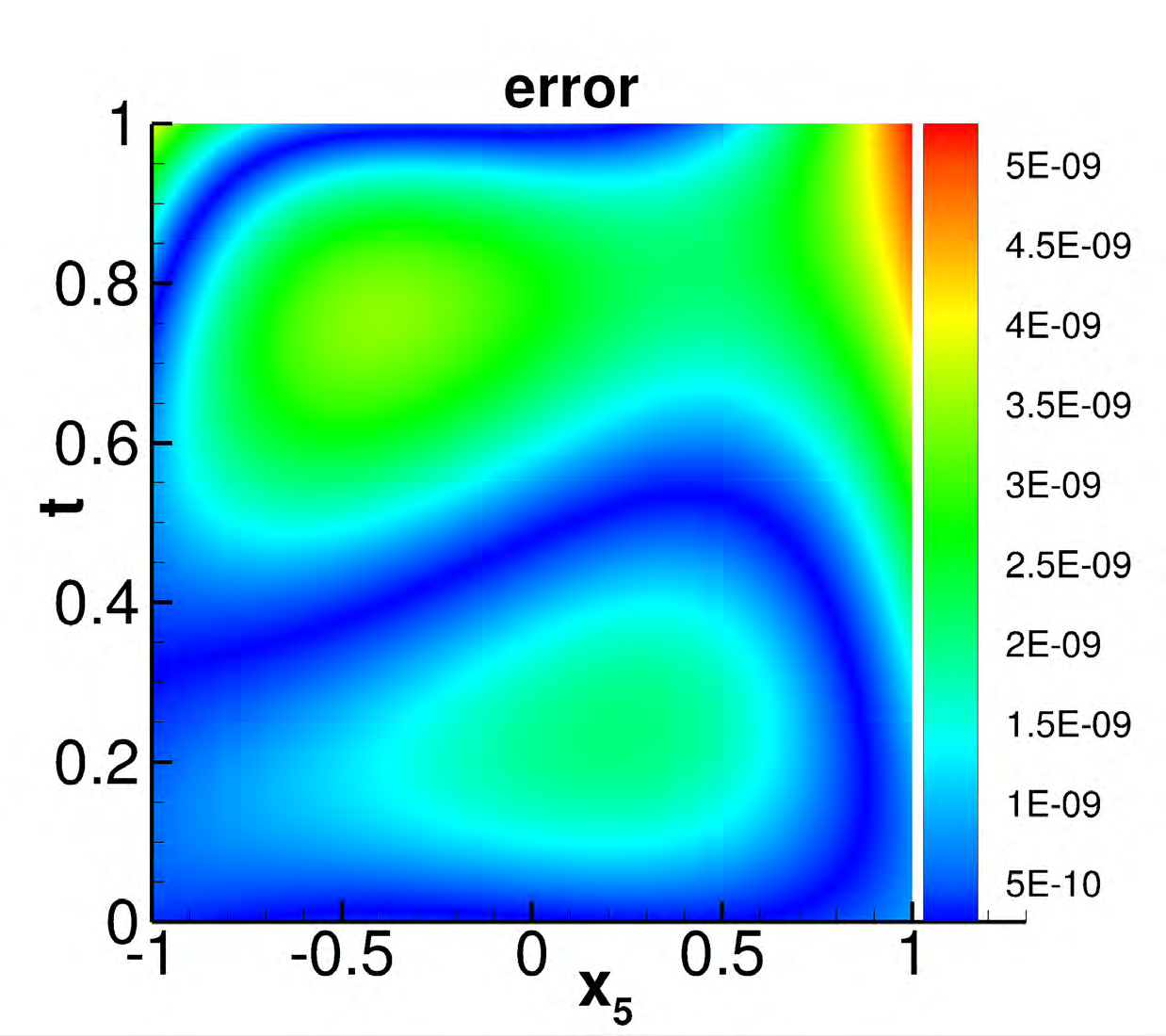}(i)
}
	\caption{\small KdV equation ($d=5$): Distributions of exact solution (left column), the ELM solution (middle column), and the point-wise absolute error of ELM (right column) on cross sections of the spatial-temporal domain: the $x_1$-$x_2$ plane (top row), $x_3$-$x_4$ plane (middle row), and $x_5$-$t$ plane (bottom row).
 NN architecture: $[d+1,1000,1]$; $(N_{\ins},N_{\bcd},N_{t_0})=(10,150,1000)$.
 }\label{kdv}
	\end{figure}

Figure~\ref{kdv} provides an overview of distributions of the exact solution, the ELM solution, and the ELM point-wise absolute error in several 2D cross sections ($x_1$-$x_2$ plane, $x_3$-$x_4$ plane, $x_5$-$t$ plane) of the spatial-temporal domain for the KdV equation in dimension $d=5$.
These cross sections are located in the middle of the domain with regard to the rest of the coordinates. The main simulation parameters for these results are listed in the figure caption. The ELM method has evidently captured the solution  accurately, with an absolute error on the level of $10^{-9}$ in these cross sections. 

\begin{figure}[tb]
	\centering
	\subfigure[d=3]
	{ \includegraphics[width=0.3\textwidth]{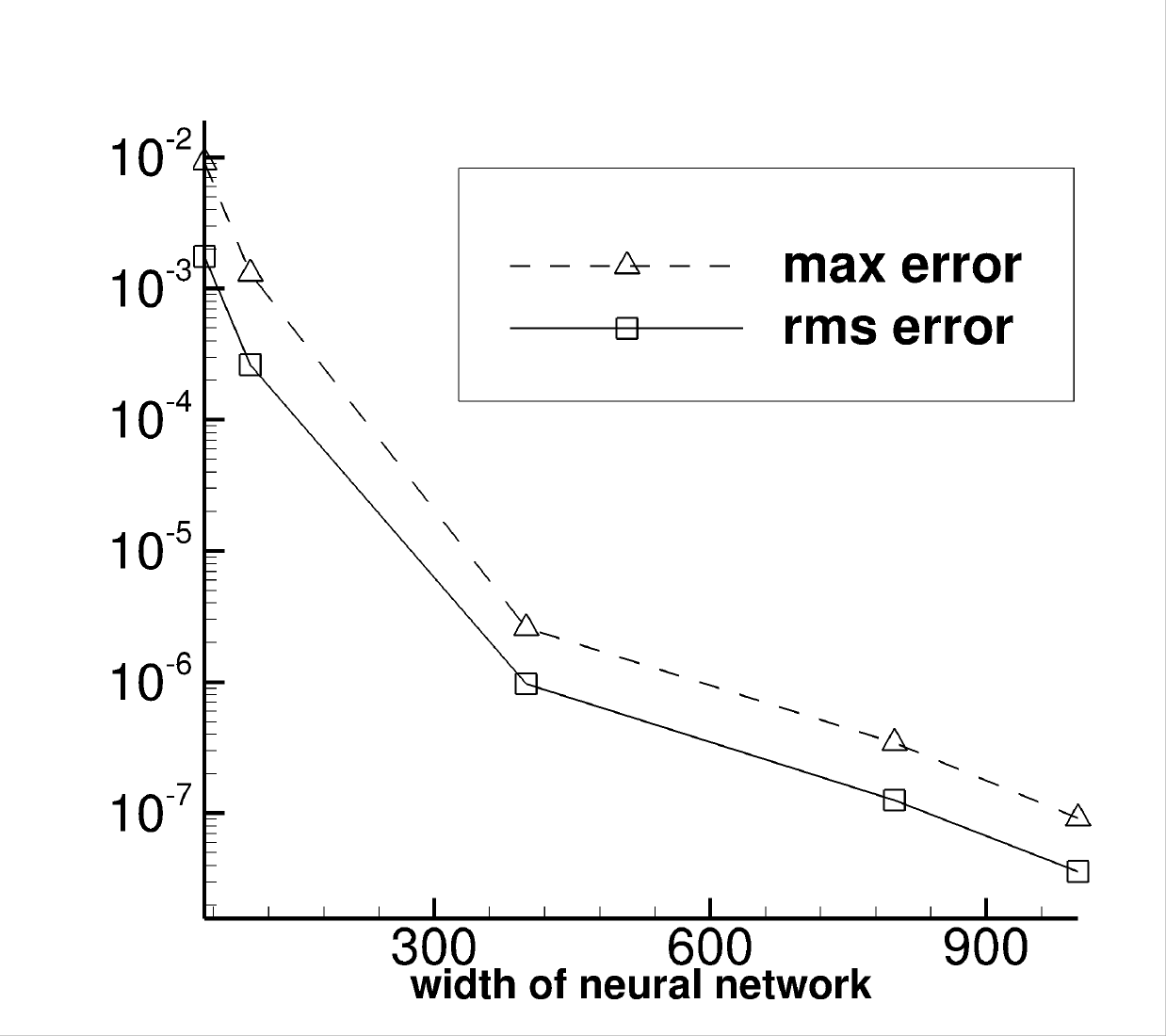}}
	\subfigure[d=5]
	{ \includegraphics[width=0.3\textwidth]{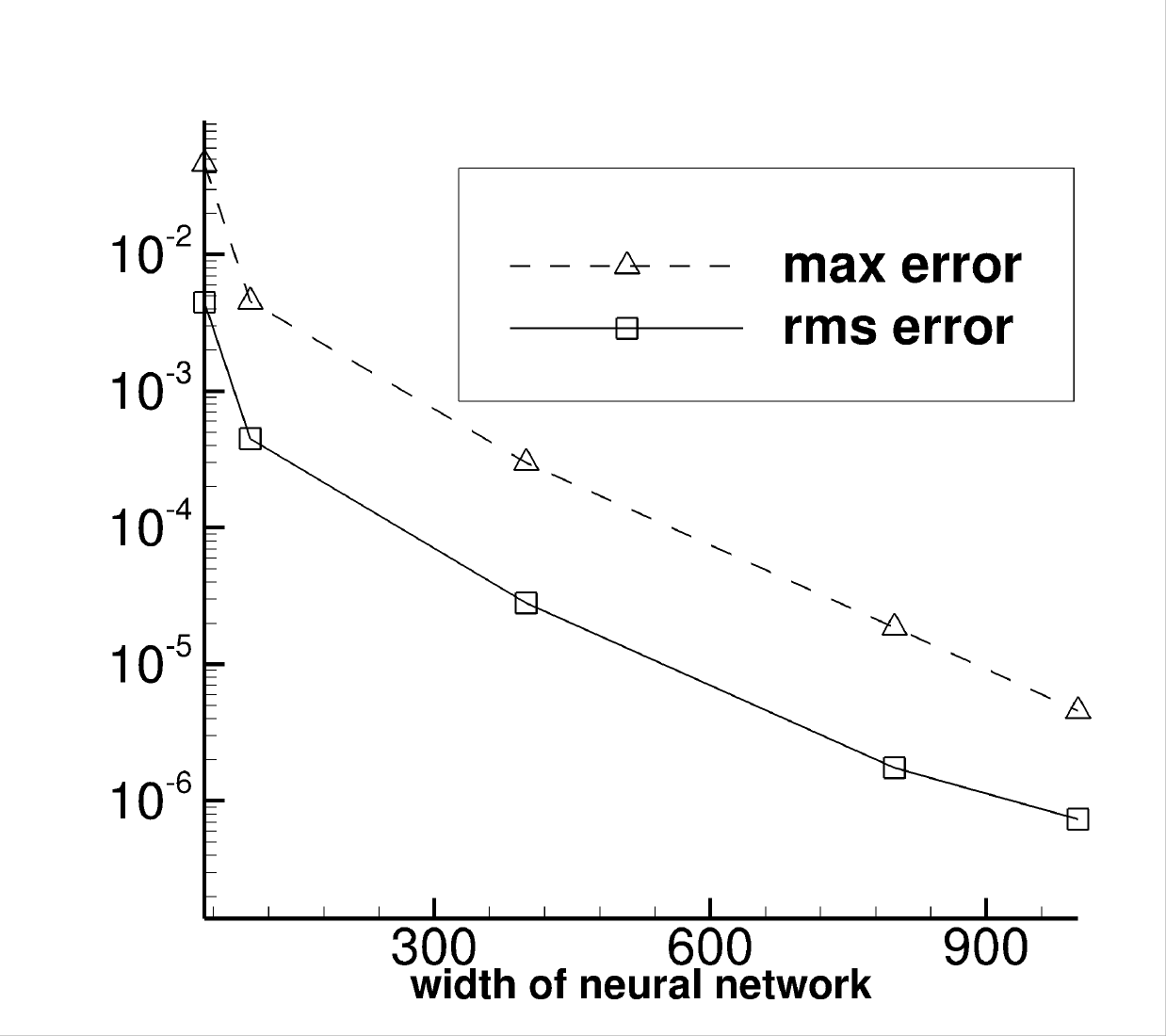}}
	\subfigure[d=10]
	{ \includegraphics[width=0.3\textwidth]{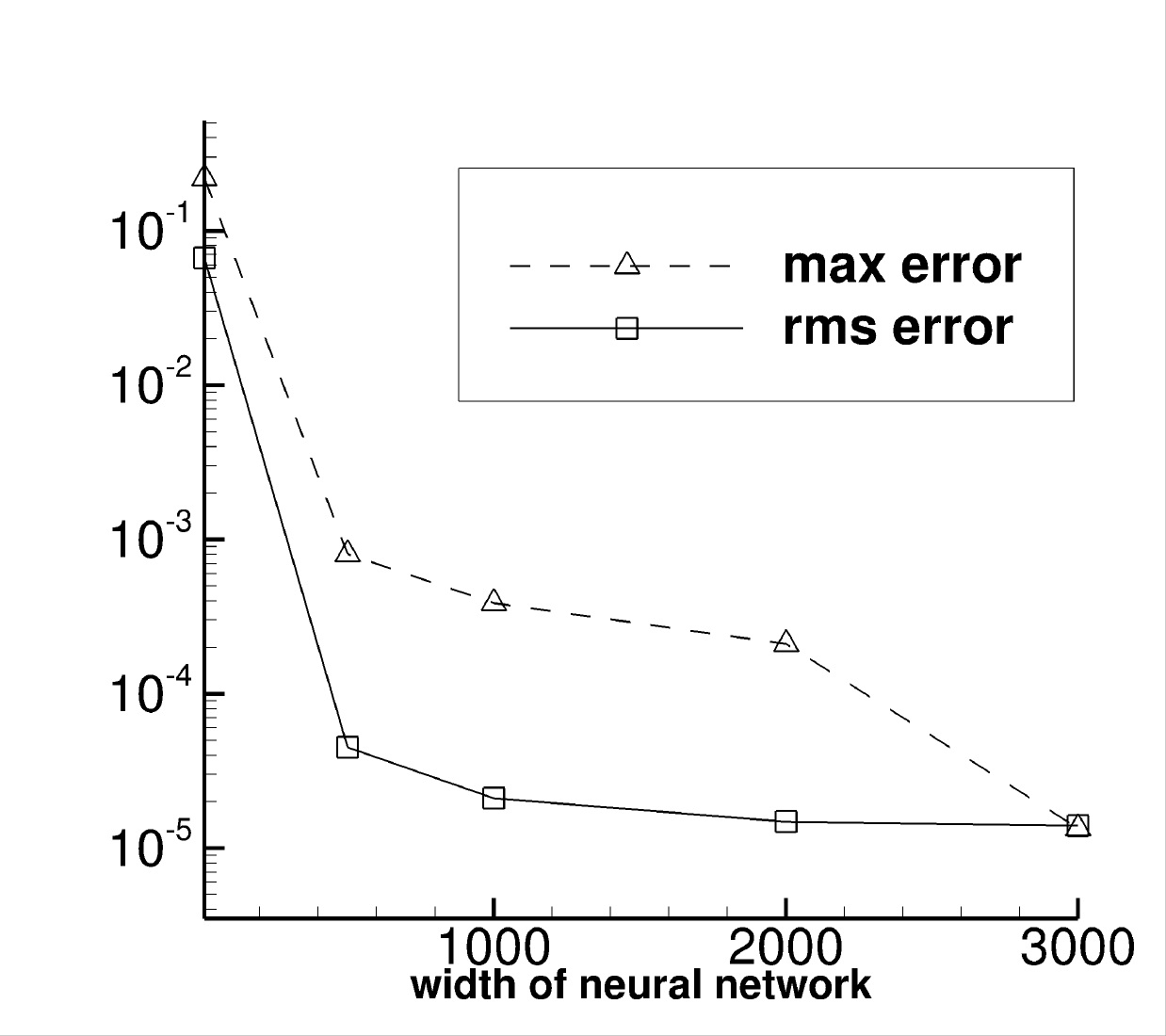}}
	\caption{\small KdV equation: $e_{l^{\infty}}$ and $e_{l^2}$ errors versus the number of training parameters ($M$) for dimensions (a) $d=3$, (b) $d=5$, and (c) $d=10$.
 NN architecture: $[d+1,M,1]$; $(N_{\ins},N_{t_0})=(10,1000)$ in (a,b,c); $N_{\bcd}=200$ in (a), $150$ in (b), and $100$ in (c); $M$ is varied in (a,b,c).
 }
	\label{kdv_err_width}
\end{figure}

\begin{figure}[tb]
	\centering
	\subfigure[d=3]
{\includegraphics[width=0.3\textwidth]{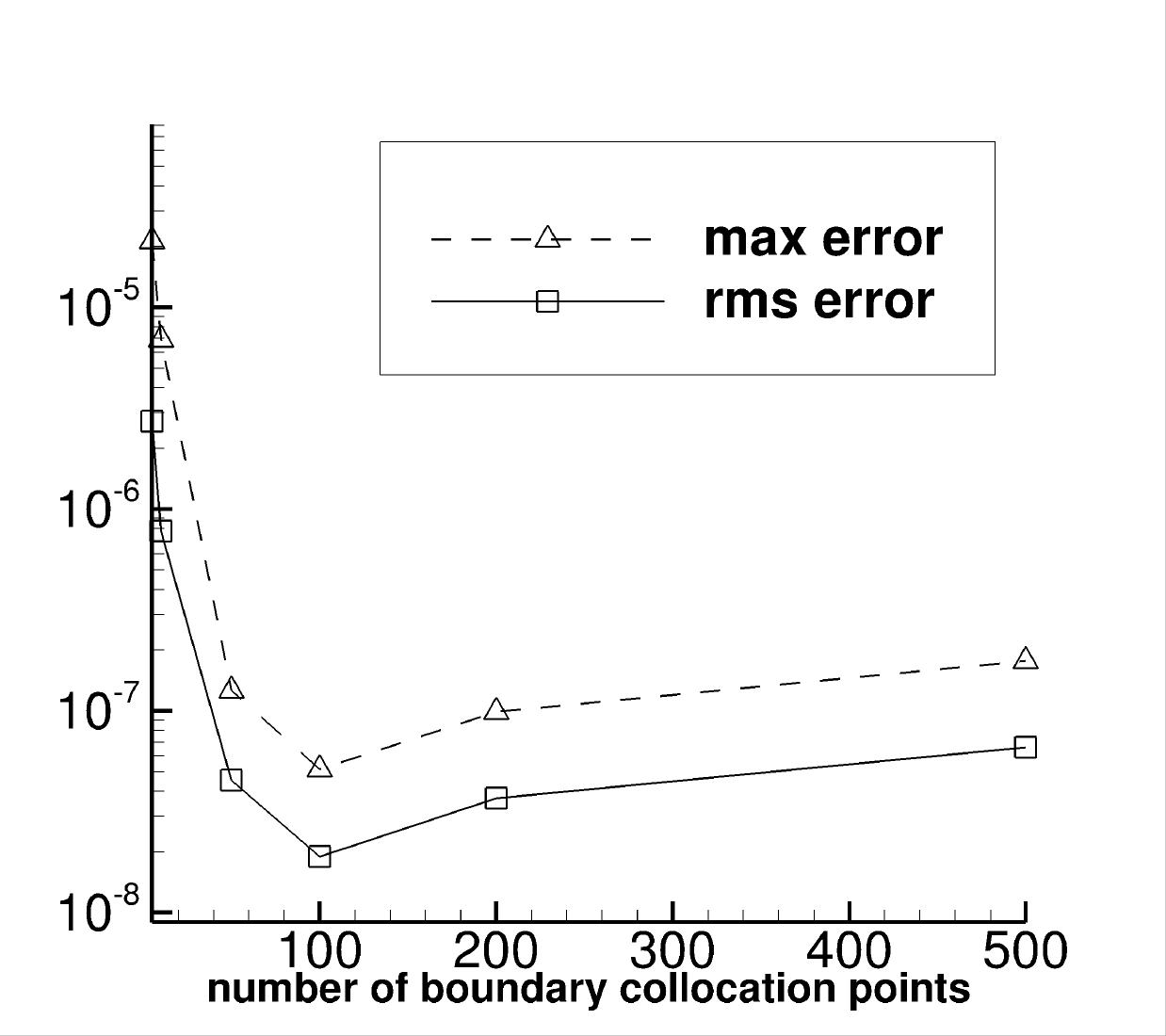}}
	\subfigure[d=5]
	{ \includegraphics[width=0.3\textwidth]{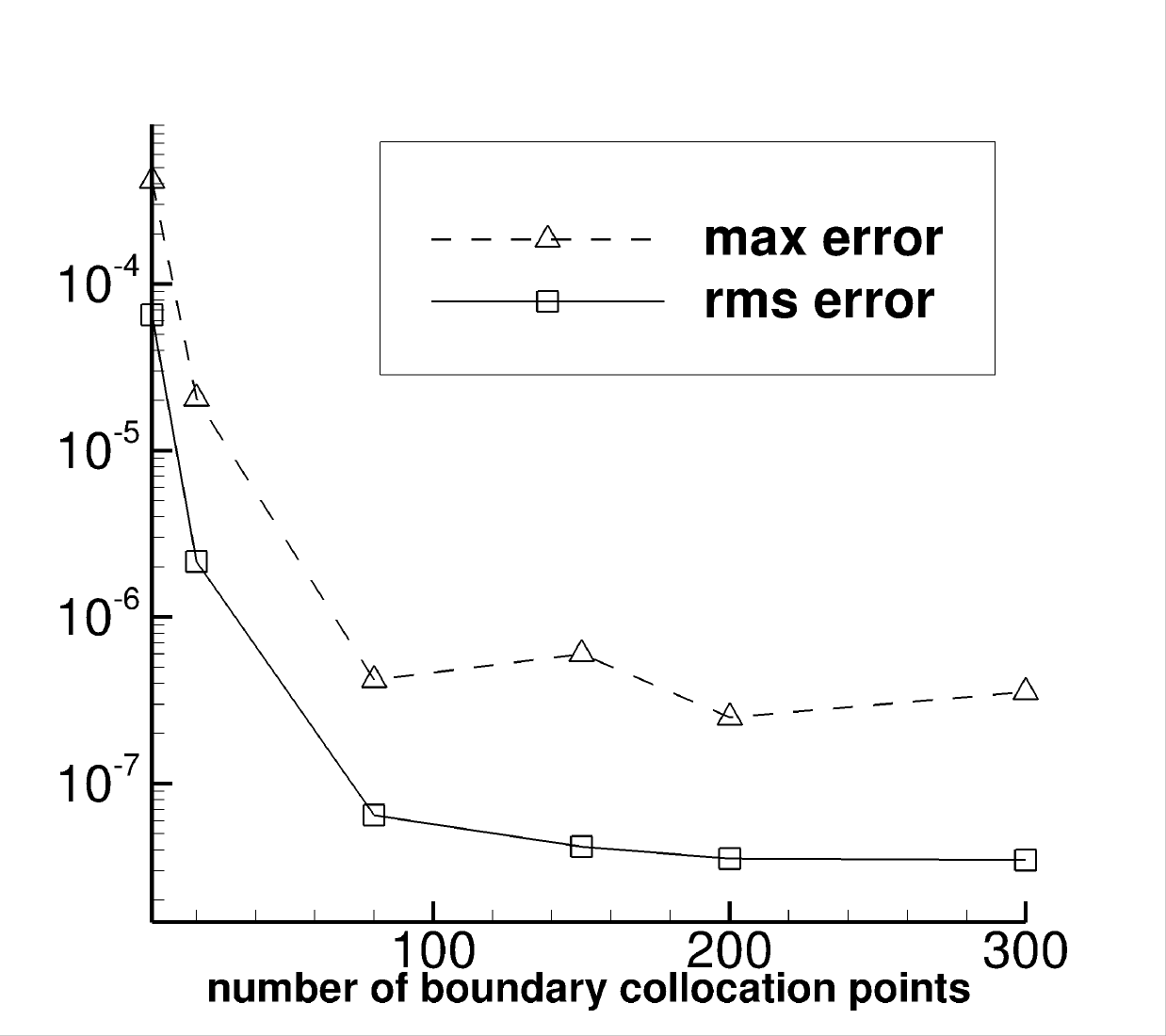}}
	\subfigure[d=10]
	{ \includegraphics[width=0.3\textwidth]{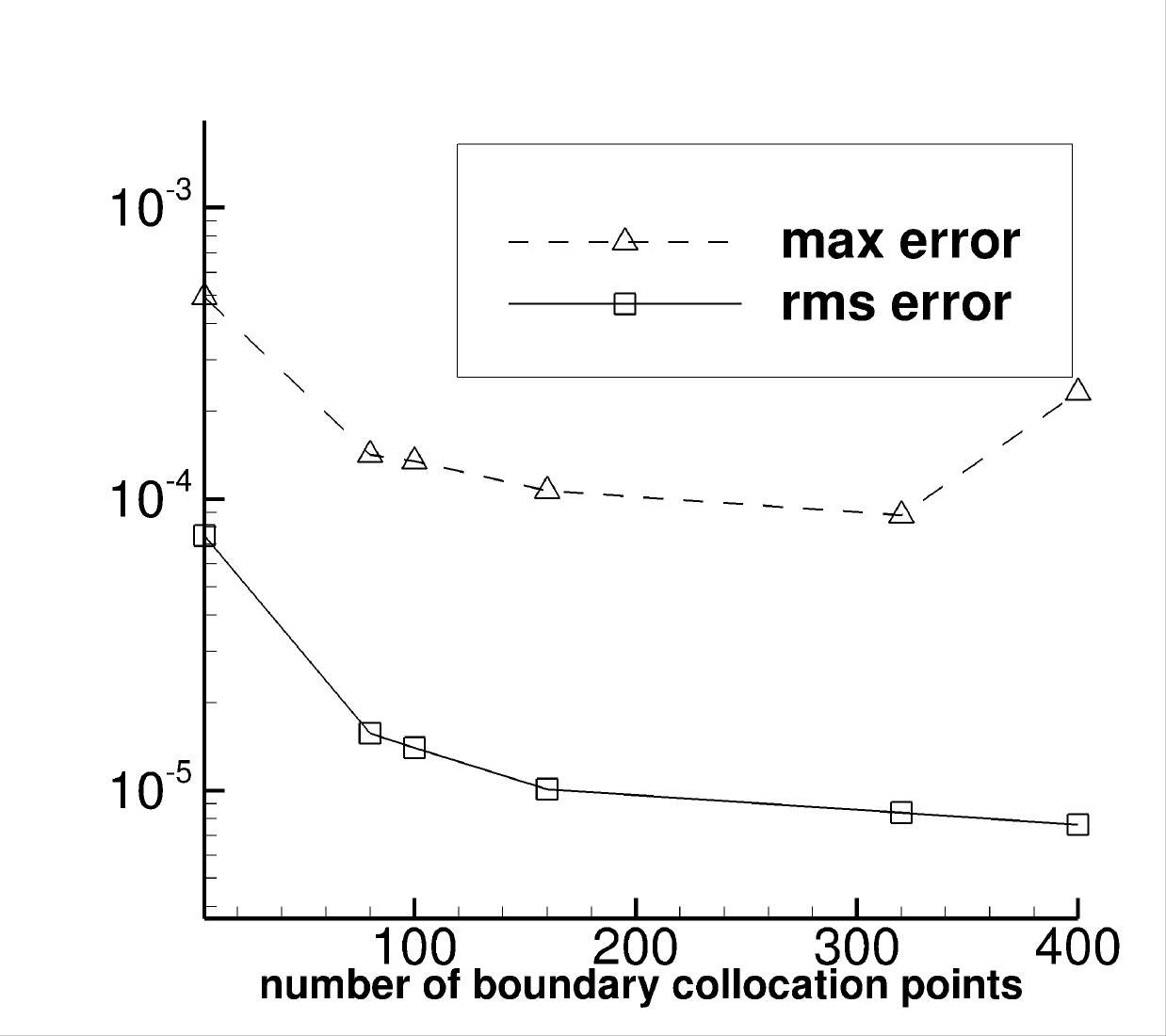}}
	\caption{\small KdV equation: $e_{l^{\infty}}$ and $e_{l^2}$ errors versus the number of collocation points on each  boundary ($N_{\bcd}$) for several problem dimensions.
 NN architecture: $[d+1,M,1]$; $M=1000$ in (a), $2000$ in (b), and $3000$ in (c). $(N_{\ins},N_{t_0})=(100,1000)$ in (a), $(10,1000)$ in (b,c). $N_{\bcd}$ is varied in (a,b,c).
 }\label{kdv_err_nc}
\end{figure}

The convergence behavior of the ELM method has been investigated and the test results are documented in Figures~\ref{kdv_err_width} and~\ref{kdv_err_nc}. These figures depict the $e_{l^{\infty}}$ and $e_{l^2}$ errors for three problem dimensions ($d=3$, $5$ and $10$) with respect to number of training parameters and the number of boundary collocation points (per boundary), respectively. The crucial simulation parameters in the tests are listed in the captions of these figures.
With increasing number of training parameters in the network, the ELM errors decrease approximately exponentially.
With respect to the number of boundary collocation points ($N_{\bcd}$), the ELM errors initially decreases approximately exponentially and gradually stagnates as $N_{\bcd}$ reaches a certain level for $d=3$ and $5$. But for $d=10$, the reduction in the ELM errors is not as significant as for the lower dimensions with increasing $N_{\bcd}$.

\subsection{Numerical Tests with the ELM/A-TFC Method }
\label{atfc_test}

In this subsection we test the performance of the combined ELM/A-TFC method from Section~\ref{sec:atfc} using several high-dimensional linear/nonlinear PDEs.


\subsubsection{Poisson Equation}

\begin{table}[tb]
\centering\small
	\begin{tabular}{cccccc}
		\hline
		$R_m$&1E-3&5E-3&0.01&0.05&0.1\\  $e_{l^\infty}$&2.65E-4&6.06E-5&2.94E-5&4.40E-5&6.64E-4\\
		$e_{l^{2}}$ &3.19E-5&5.19E-6&2.51E-6&4.11E-6&6.06E-5\\
		\hline
	\end{tabular}
	\caption{\small Poisson equation ($d=7$): determining $R_{m0}$ for the ELM/A-TFC method using the procedure from Remark~\ref{rem_a6}.
 NN architecture: $[d,3000,1]$; $(N_{\bcd},N_{\ins})=(100,1000)$.
 }
	\label{tab_9}
\end{table}

\begin{table}[tb]
\centering\small
\begin{tabular}{l| cc}
\hline
$d$ & 3 & 7  \\ \hline
$R_{m0}$ & 0.1 & 0.01 \\
\hline
\end{tabular}
\caption{\small Poisson equation: $R_{m0}$ for ELM/A-TFC determined by the procedure from Remark~\ref{rem_a6}. The simulation parameters (NN architecture, collocation points) here follow those of Table~\ref{tab_9}, except that for $d=3$ the NN architecture is $[d,1000,1]$.
}
\label{tab_10}
\end{table}

We employ the same Poisson problem as in Section~\ref{sec_poisson} to test the ELM/A-TFC method. the governing equations  are given by equations~\eqref{eq_42}. The problem settings here follow those of Section~\ref{sec_poisson}.
The notations below follow those of the test problems in Section~\ref{elm_test}.

Let us first determine the $R_{m0}$ using the procedure from Remark~\ref{rem_a6} for generating the random hidden-layer coefficients  with the ELM/A-TFC method. 
Table~\ref{tab_9} shows the test using this procedure for   dimension $d=7$, leading to $R_{m0}\approx 0.01$. The $R_{m0}$ values for different dimensions are listed in Table~\ref{tab_10}, which we will use for generating the random hidden-layer coefficients with ELM/A-TFC  in subsequent tests.

\begin{figure}[!tb]
	\centerline{
	\subfigure
	{ \includegraphics[width=0.26\textwidth]{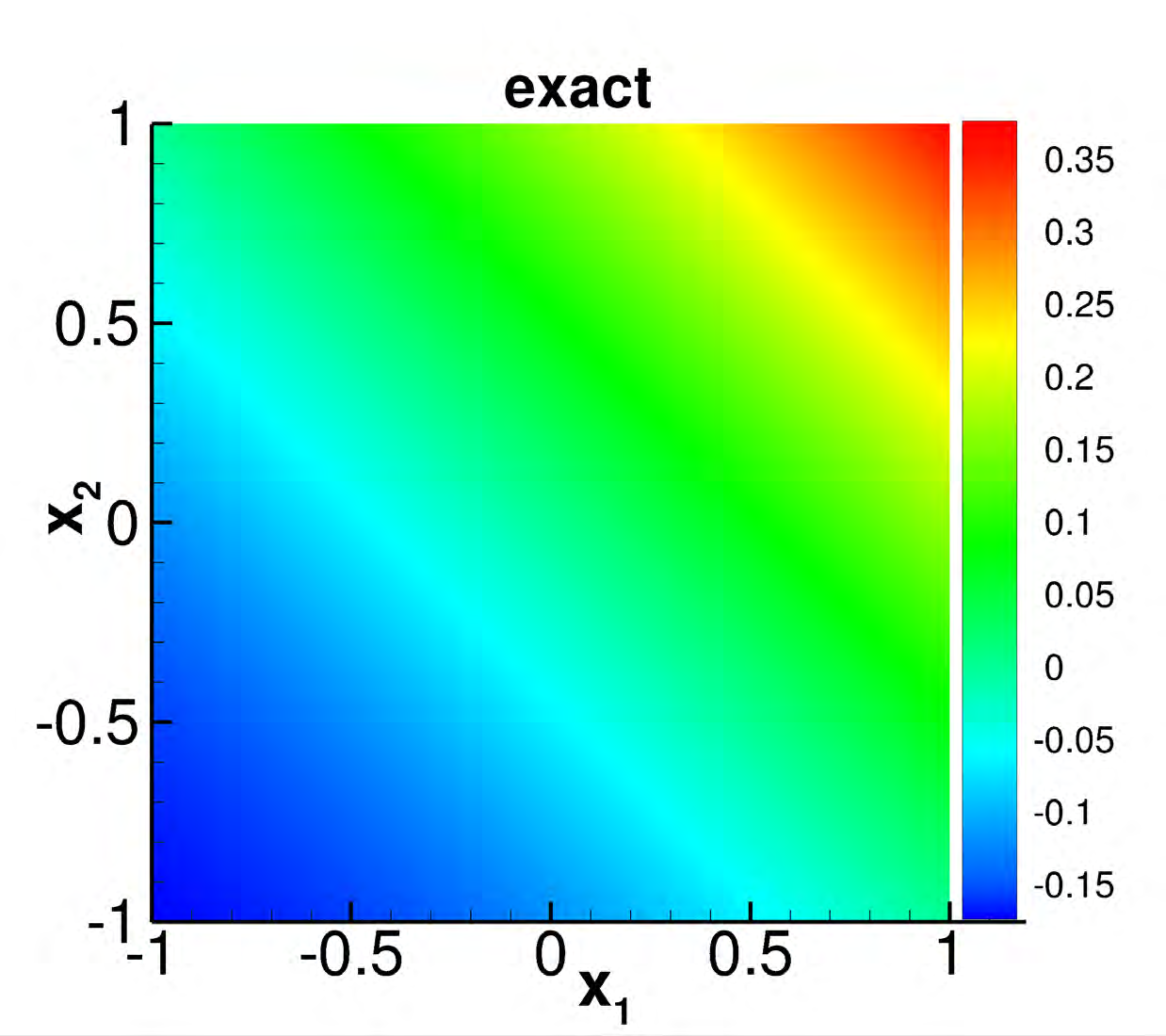}(a)}
 \subfigure
	{ \includegraphics[width=0.26\textwidth]{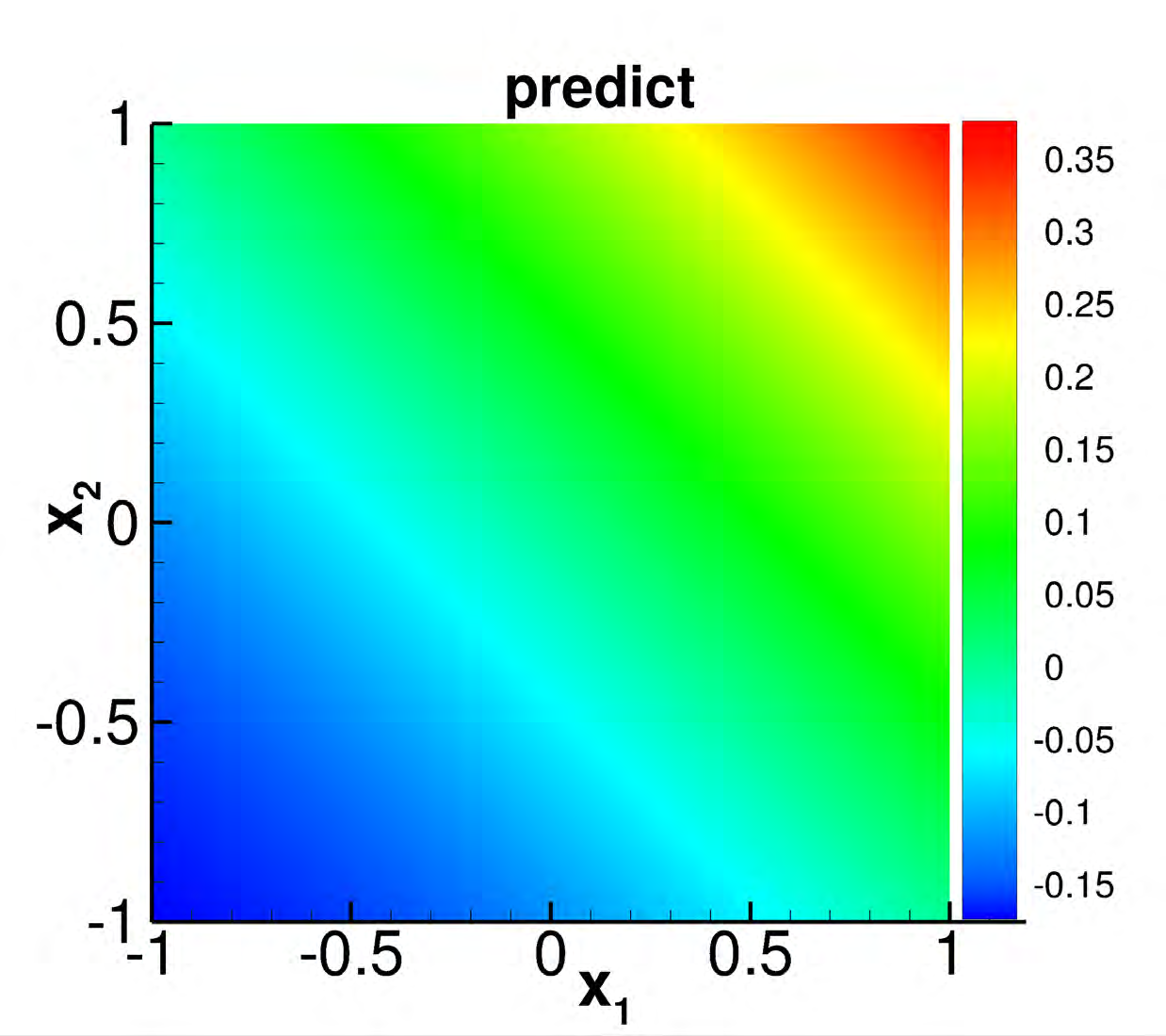}(b)}	
	\subfigure
	{ \includegraphics[width=0.26\textwidth]{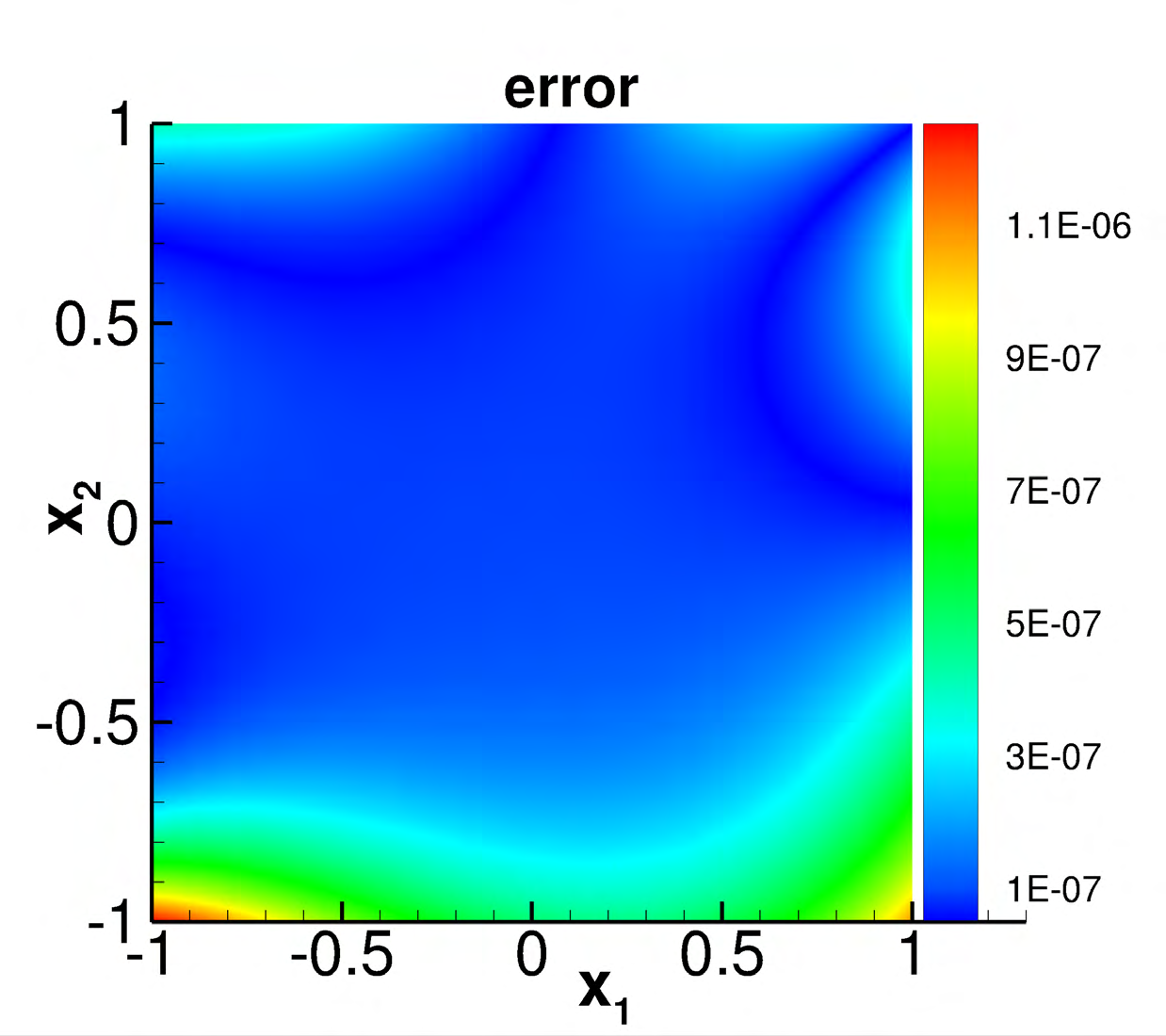}(c)}
 }
 \centerline{
	\subfigure
	{ \includegraphics[width=0.26\textwidth]{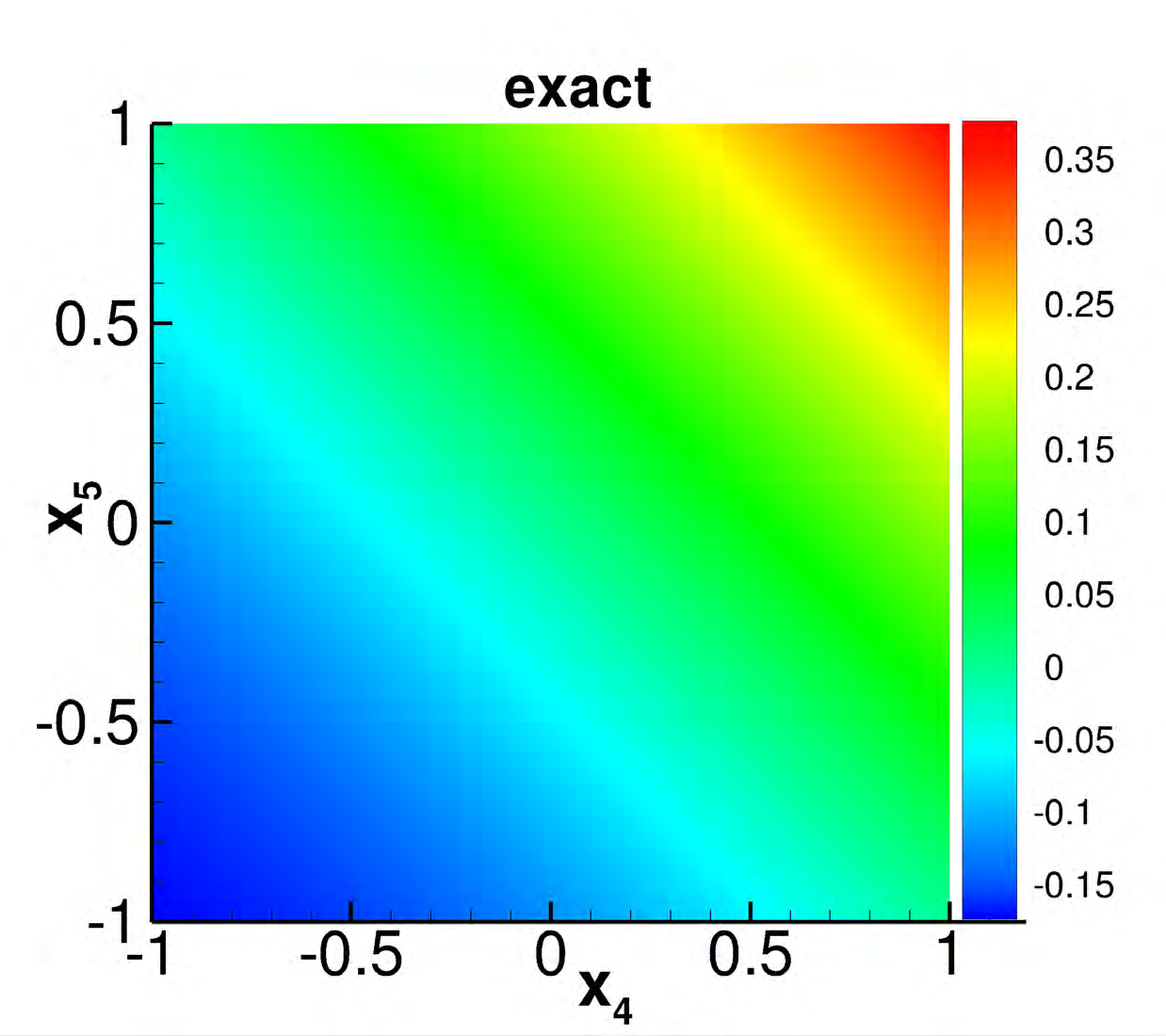}(d)}
	\subfigure
	{ \includegraphics[width=0.26\textwidth]{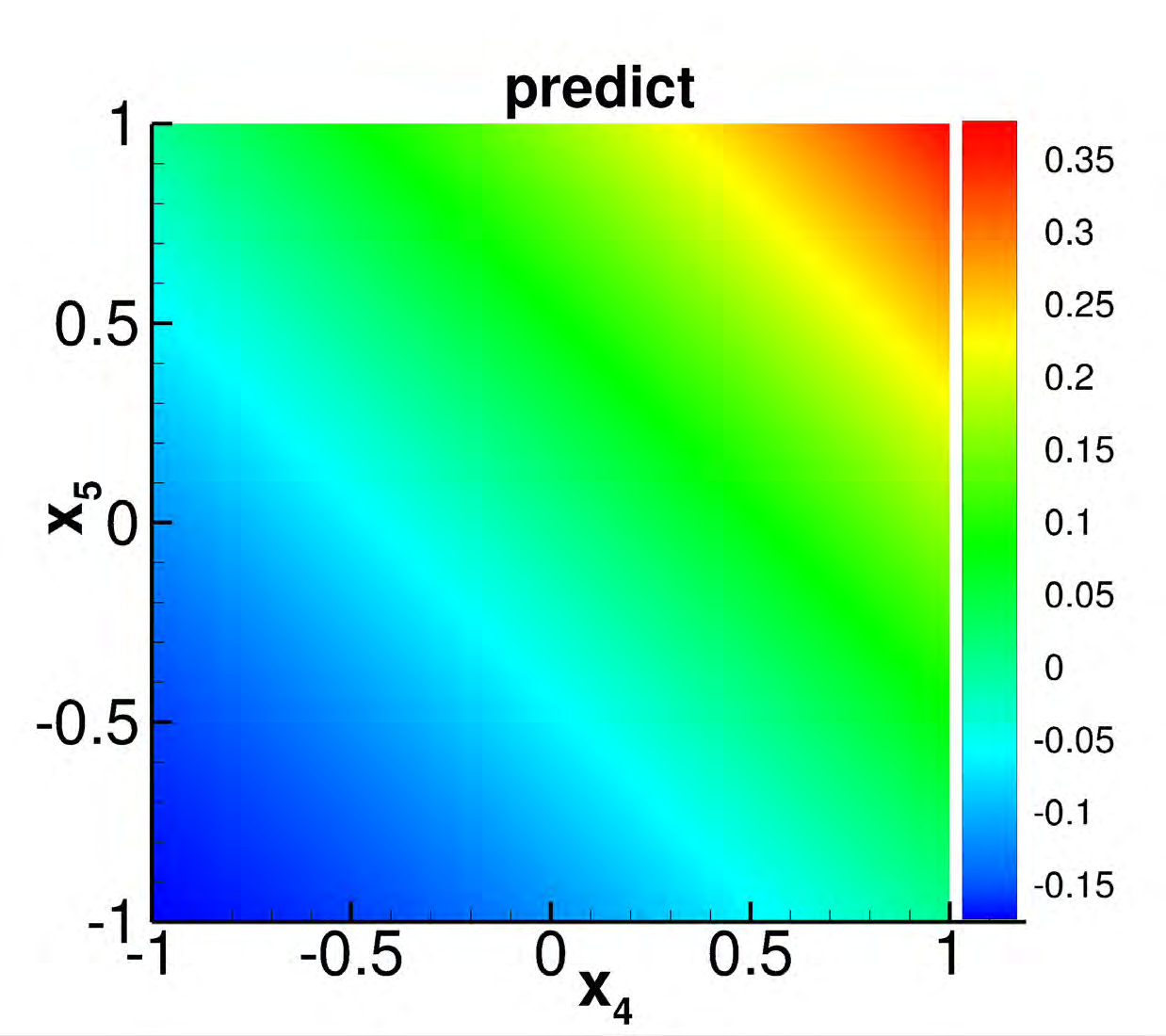}(e)}
	\subfigure
	{ \includegraphics[width=0.26\textwidth]{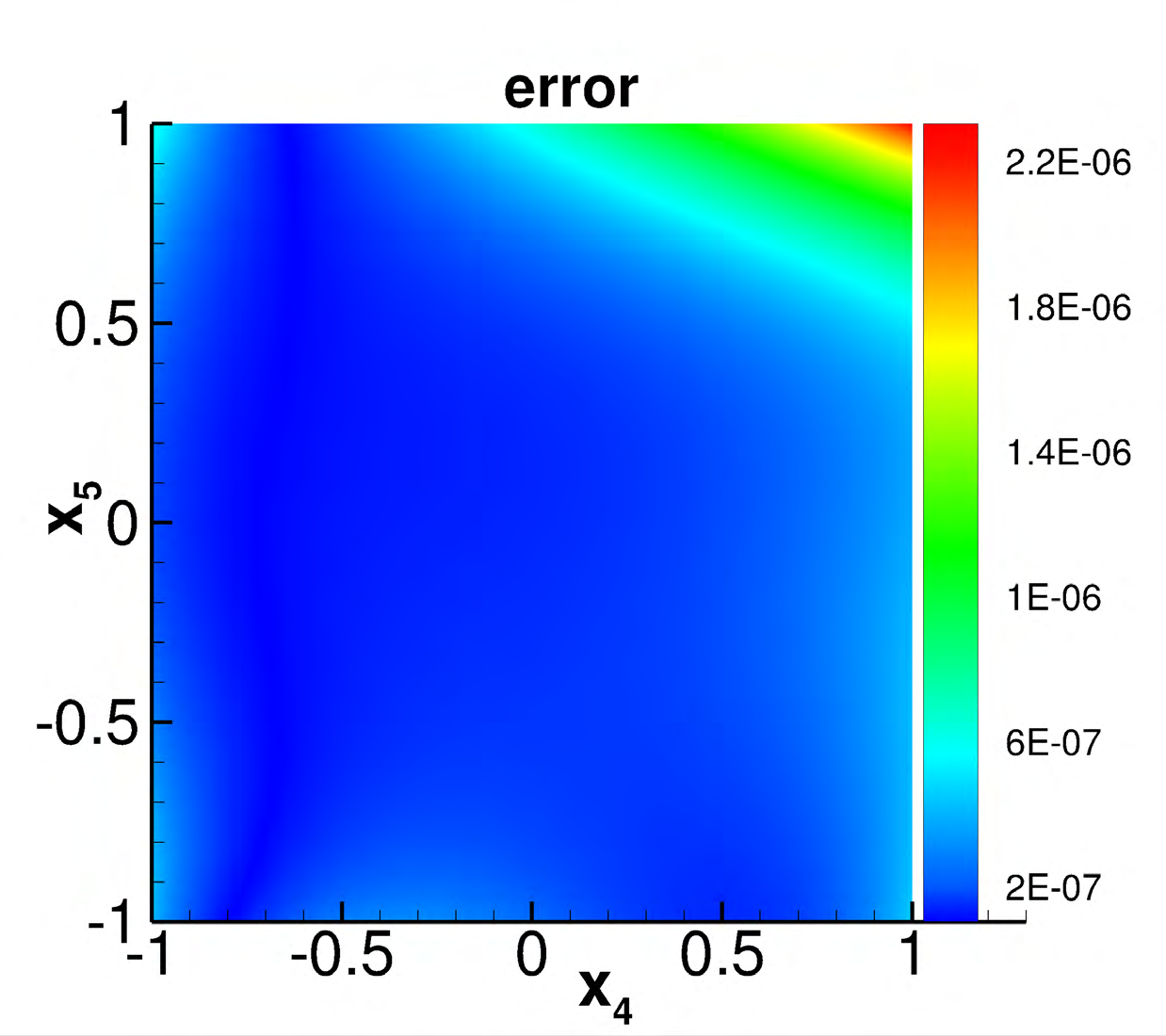}(f)}
 }
\centerline{
	\subfigure
	{ \includegraphics[width=0.26\textwidth]{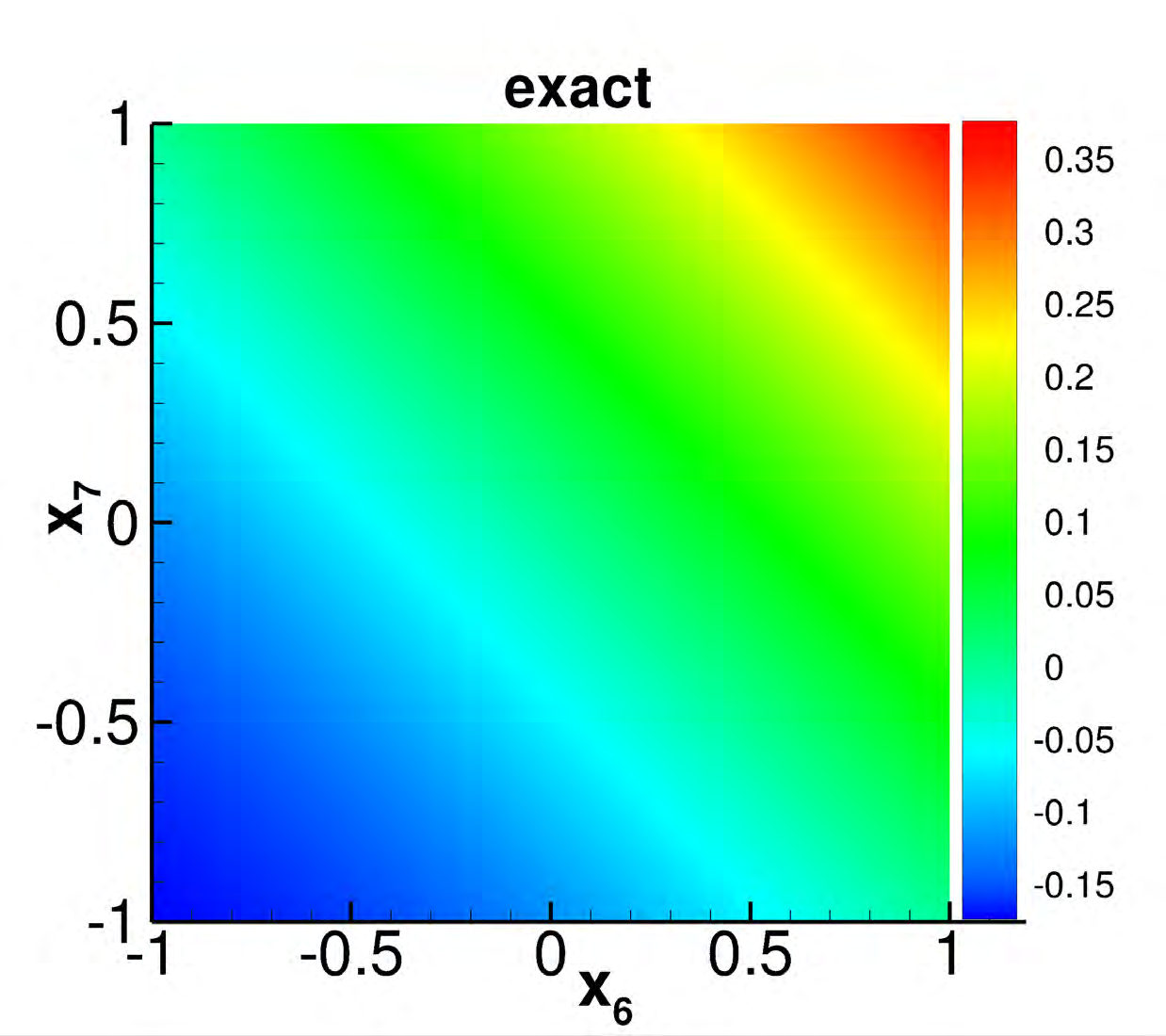}(g)}
 \subfigure
	{ \includegraphics[width=0.26\textwidth]{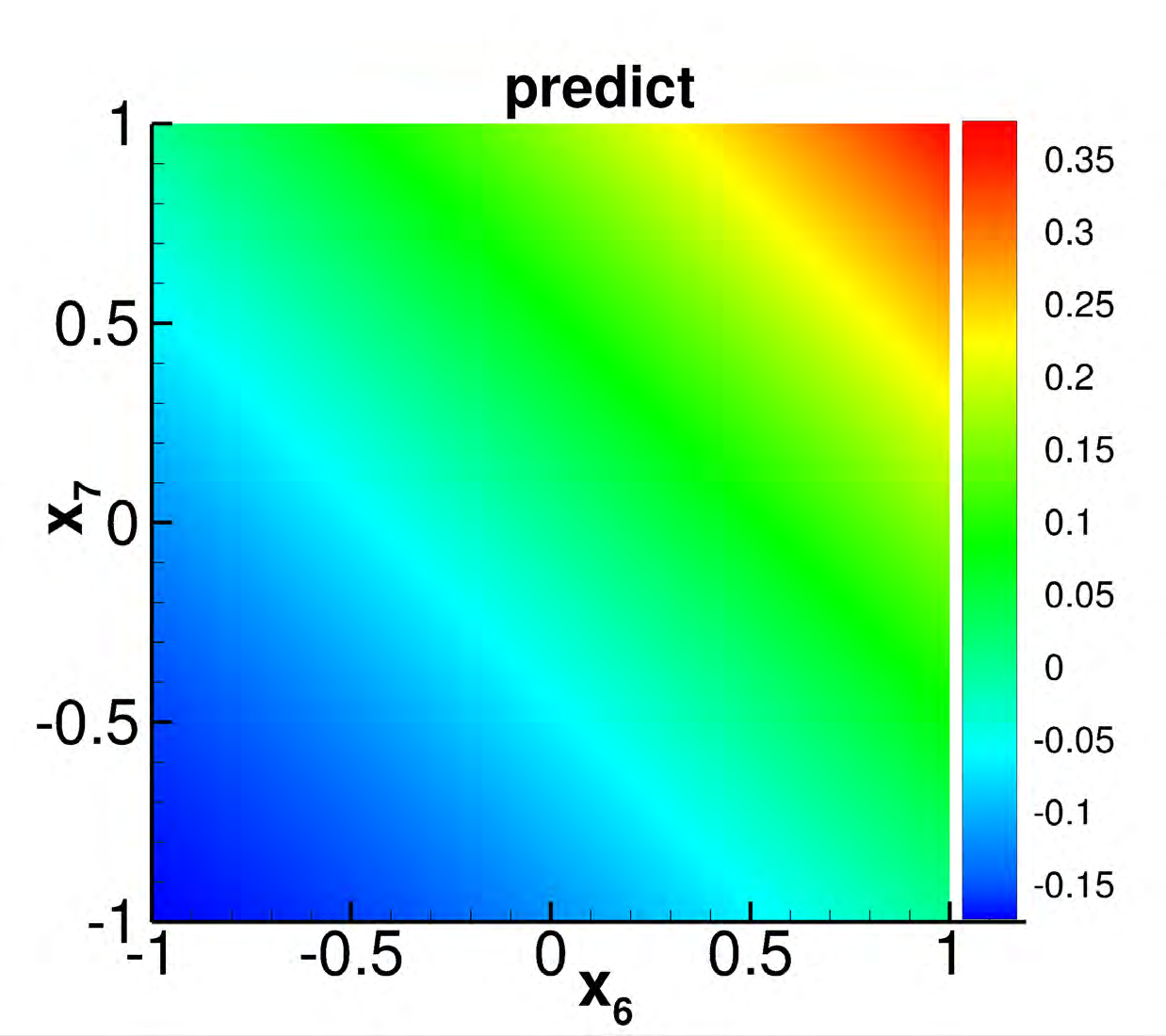}(h)}
	\subfigure
	{ \includegraphics[width=0.26\textwidth]{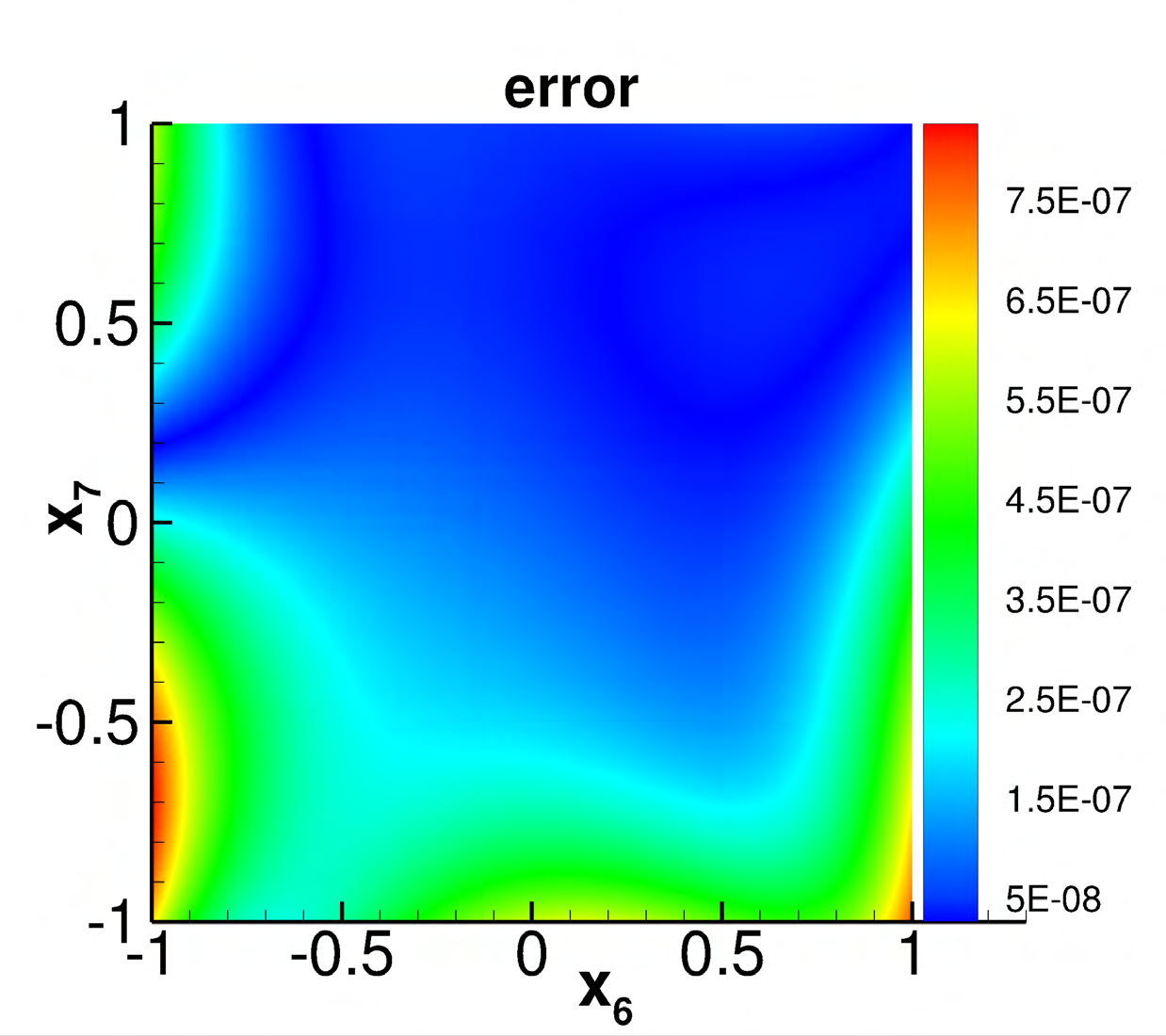}(i)}
 }
	\caption{\small Poisson equation ($d=7$): Distributions of the exact solution (left column), the ELM/A-TFC solution (middle column), and the point-wise absolute error of ELM/A-TFC (right column) in several 2D cross-sections of the domain (top row: $x_1$-$x_2$ plane; middle row: $x_4$-$x_5$ plane; bottom row: $x_6$-$x_7$ plane). 
 NN architecture: $[d,3000,1]$; $(N_{\bcd},N_{\ins})=(200,10)$.
 }
 \label{tfc_poi}
\end{figure}

An illustration of the distributions of the exact solution, the ELM/A-TFC solution, and the point-wise absolute error of the ELM/A-TFC solution for $d=7$ is provided in Figure~\ref{tfc_poi} for several 2D cross-sections of the domain (the $x_1$-$x_2$, $x_4$-$x_5$, $x_6$-$x_7$ planes). These planes are located in the middle of the domain with respect to the rest of the coordinates.
These results are obtained using a network architecture $\mbs M_{\text{arch}}=[d,3000,1]$, with the collocation points characterized by $(N_{\bcd},N_{\ins})=(200,10)$. The ELM/A-TFC results are observed to be quite accurate, with the maximum errors on the order of $10^{-6}$ or $10^{-7}$ in these cross sections.

\begin{figure}[tb]
	\centering
	\subfigure[d=3]
	{ \includegraphics[width=0.3\textwidth]{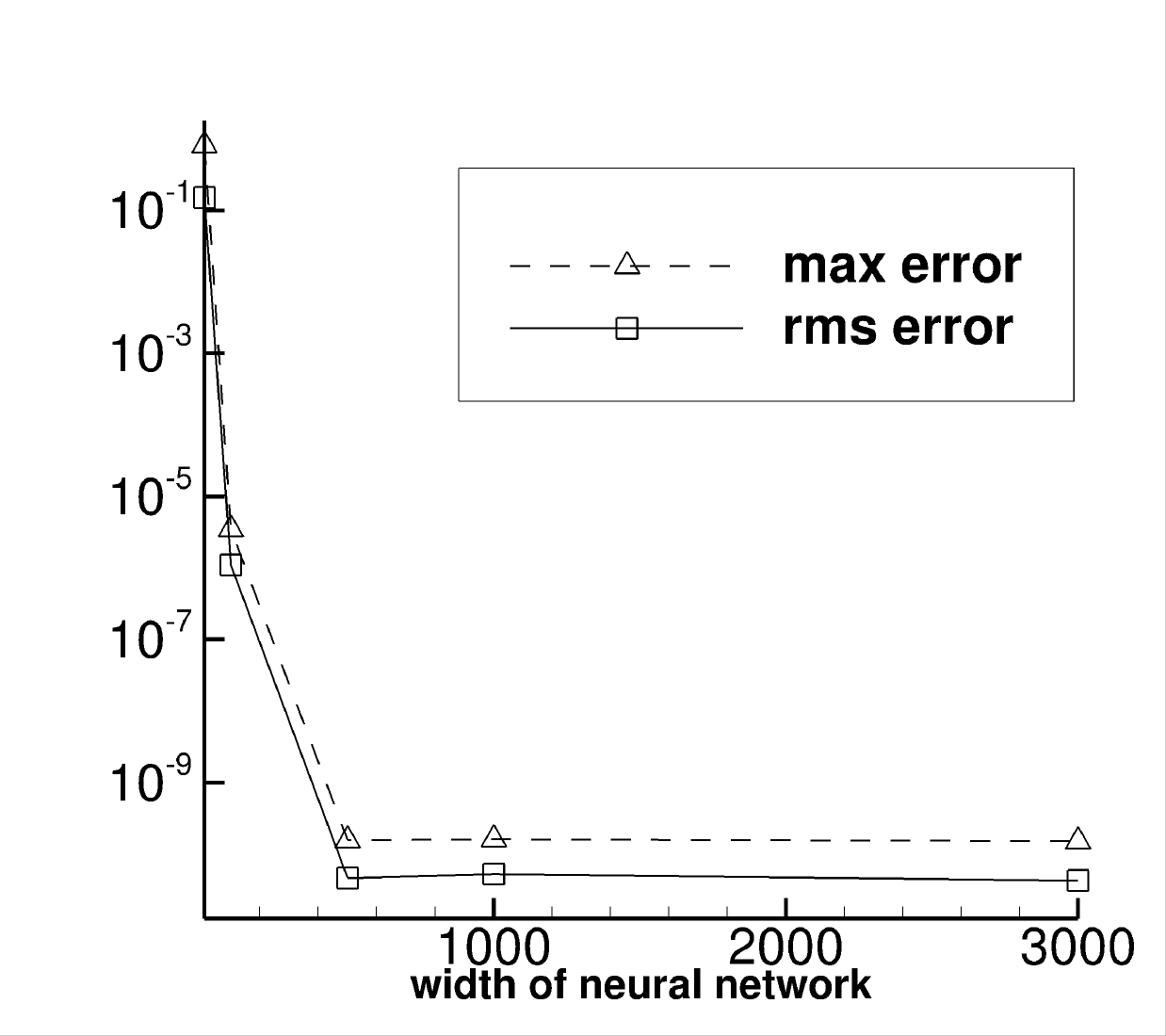}}
	\subfigure[d=7]
	{ \includegraphics[width=0.3\textwidth]{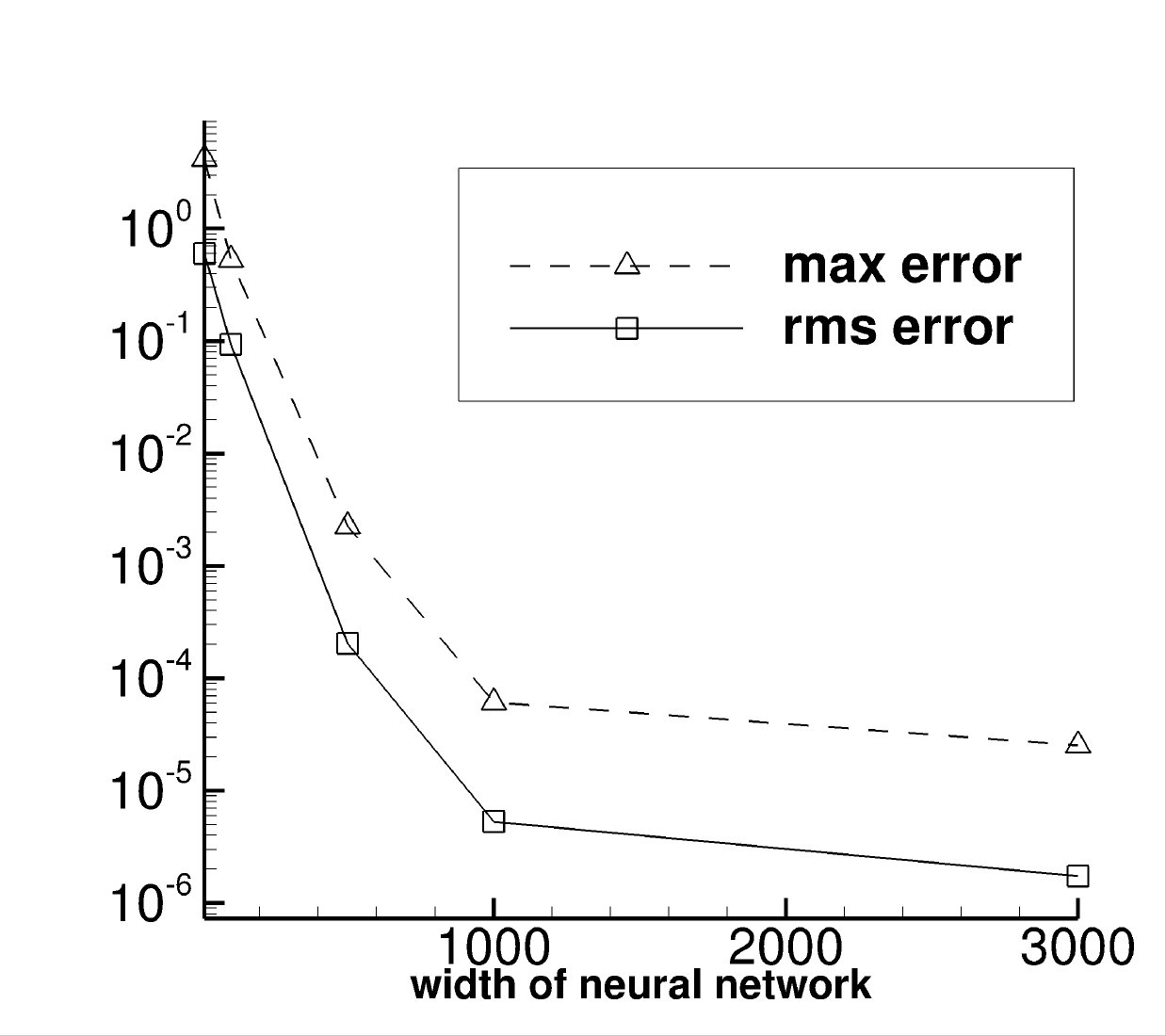}}
	\caption{\small Poisson equation: $e_{l^{\infty}}$ and $e_{l^2}$ errors versus the number of training parameters ($M$) for dimensions (a) $d=3$ and (b) $d=7$.
 NN architecture: $[d,M,1]$ ($M$ varied); $(N_{\bcd},N_{\ins})=(100,1000)$ in (a,b).
 }\label{fg_17}
\end{figure}

\begin{figure}[tb]
	\centering
	\subfigure[dim=3]
	{ \includegraphics[width=0.3\textwidth]{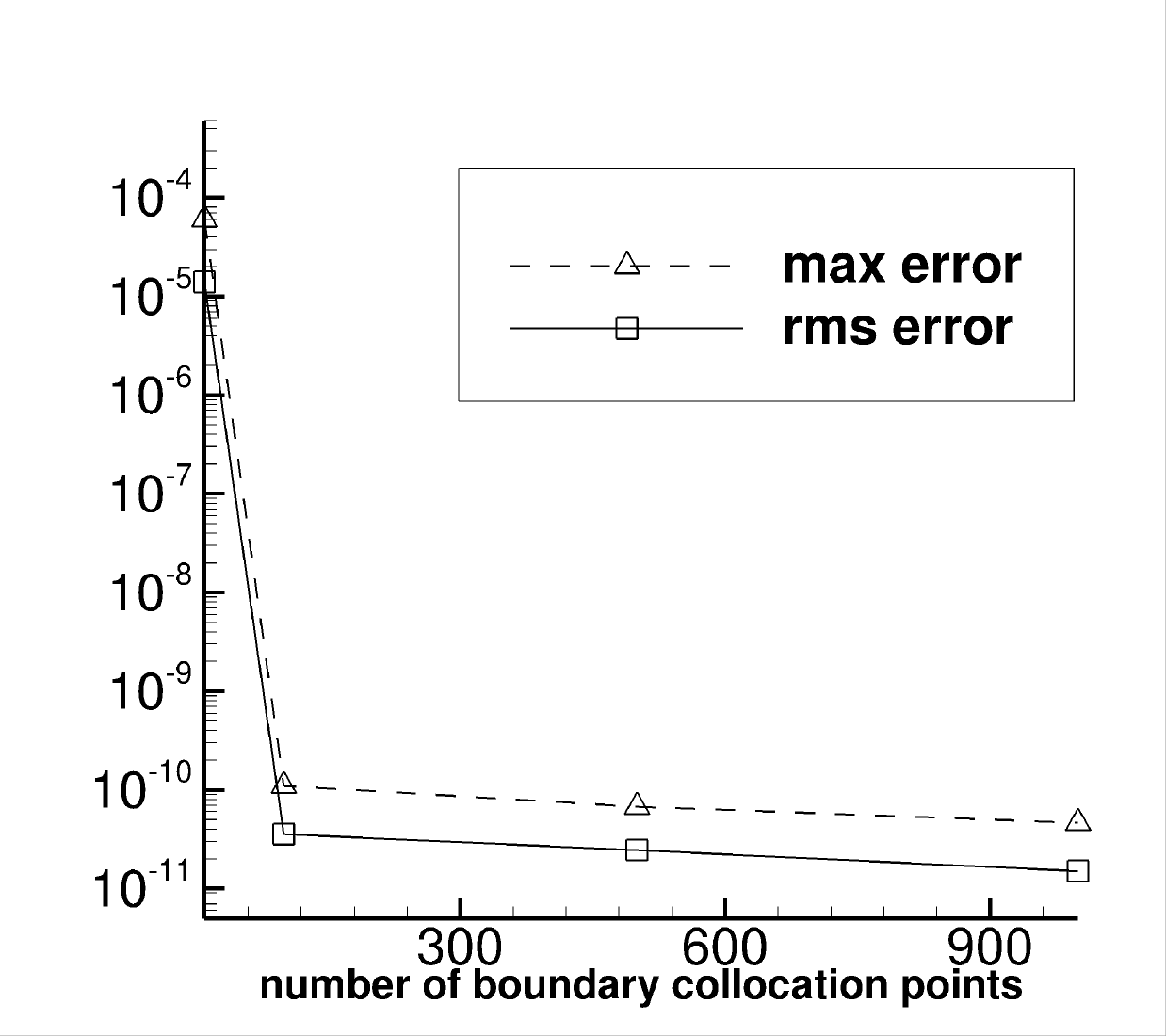}}
	\subfigure[dim=7]
	{ \includegraphics[width=0.3\textwidth]{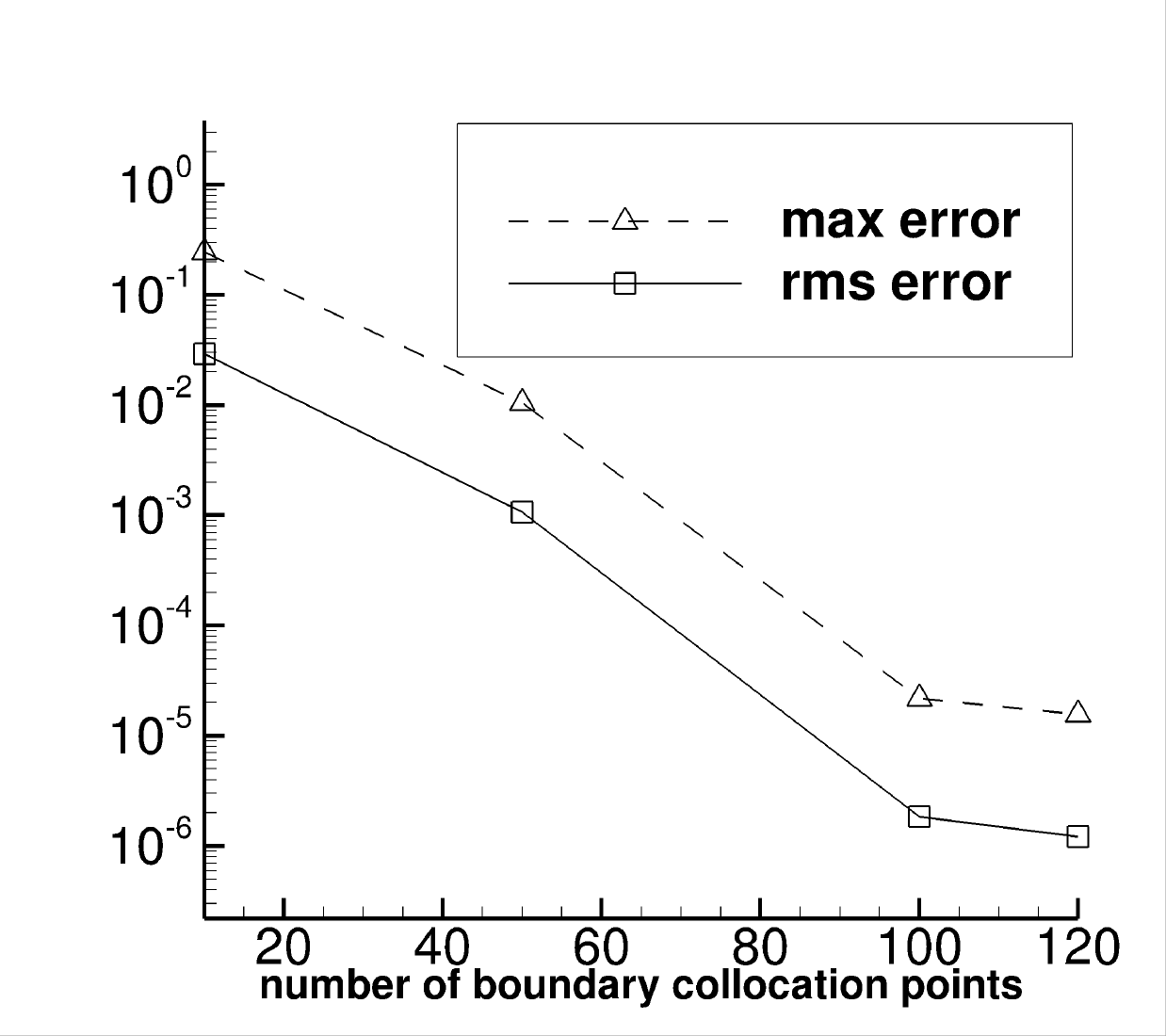}}
	\caption{\small Poisson equation: $e_{l^{\infty}}$ and $e_{l^2}$ errors versus the number of boundary collocation points ($N_{\bcd}$). 
 NN architecture: $[d,M,1]$, with $M=1000$ in (a) and $3000$ in (b). $N_{\ins}=1000$ in (a,b), while $N_{\bcd}$ is varied.
 }
 \label{fg_18}
\end{figure}

\begin{figure}[tb]
	\centering
	\subfigure[dim=3]
	{ \includegraphics[width=0.3\textwidth]{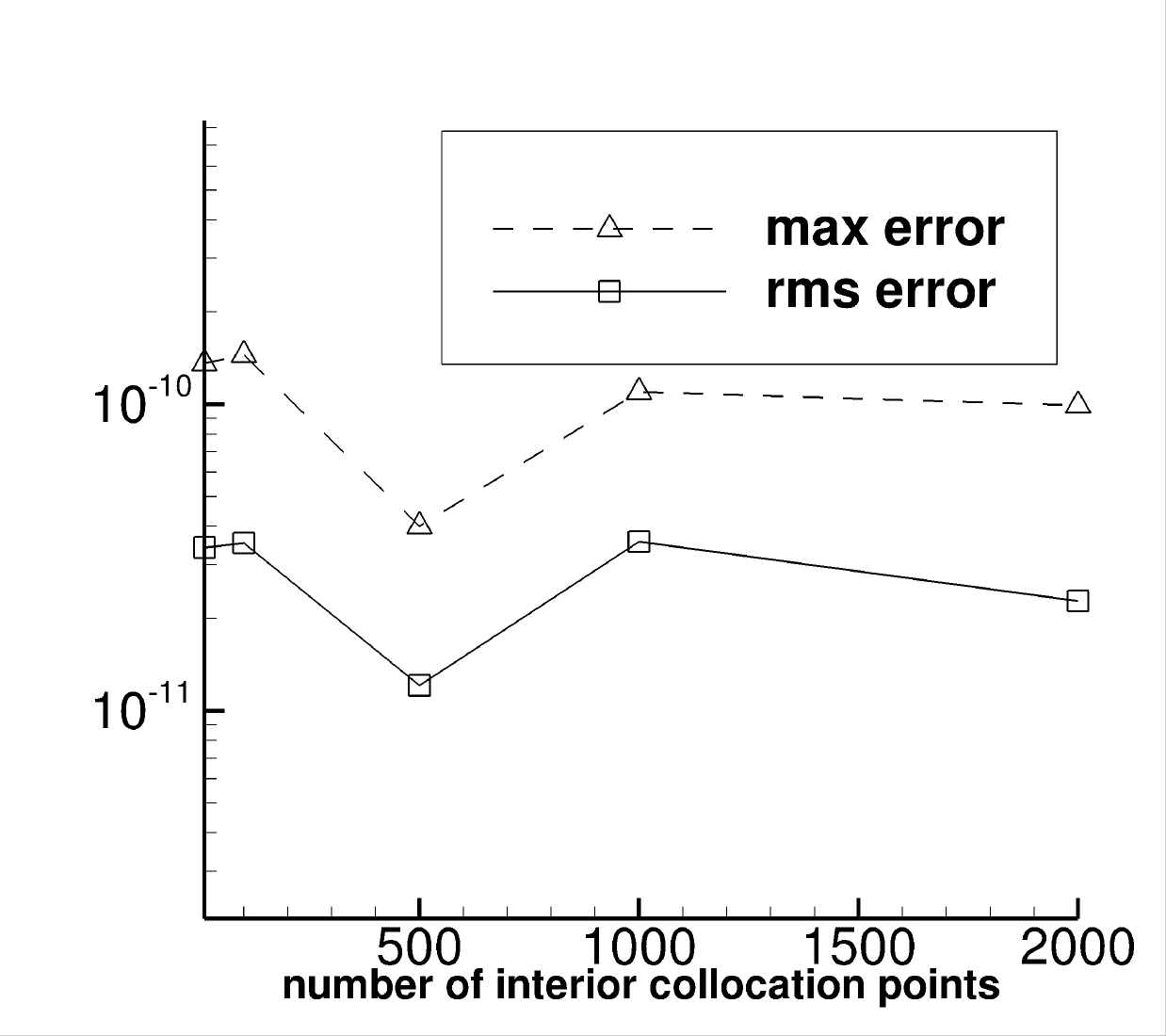}}
	\subfigure[dim=7]
	{ \includegraphics[width=0.3\textwidth]{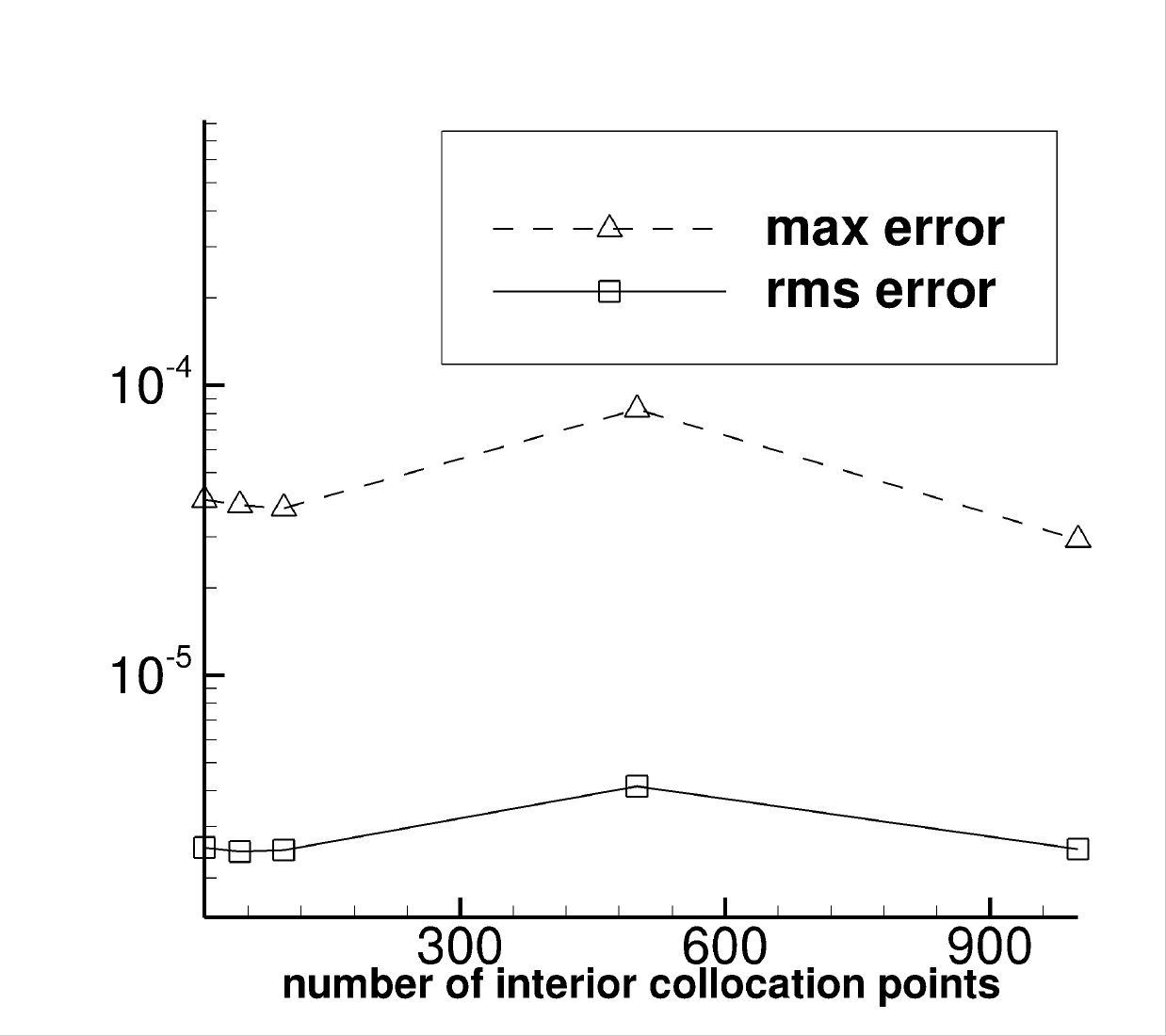}}
	\caption{\small Poisson equation: $e_{l^{\infty}}$ and $e_{l^2}$ errors versus the number of interior collocation points ($N_{\ins}$).
 NN architecture: $[d,M,1]$, with $M=1000$ in (a) and $3000$ in (b). $N_{\bcd}=100$ in (a,b), while $N_{\ins}$ is varied.
 }
 \label{fg_19}
\end{figure}

The convergence behavior of ELM/A-TFC with respect to the number of training parameters is illustrated in Figure~\ref{fg_17} for  dimensions $d=3$ and $d=7$. Here the width of the single hidden layer in the network  is varied, while the numbers of boundary/interior collocation points are fixed and listed in the figure caption. We observe an initial exponential decrease in the $e_{l^{\infty}}$ and $e_{l^2}$ errors with increasing number of training parameters. Then the errors stagnate as the number of training parameters reaches a certain level.

Figures~\ref{fg_18} and~\ref{fg_19} demonstrate the convergence behavior of ELM/A-TFC with respect to the number of boundary and interior collocation points, respectively. Here the neural network architecture is $\mbs M_{\text{arch}}=[d,M,1]$, with a fixed  $M=1000$ for  $d=3$ and $M=3000$ for $d=7$. We observe an exponential decrease in the $e_{l^{\infty}}$ and $e_{l^2}$ errors (before saturation), as the number of boundary collocation points ($N_{\bcd}$) increases. On the other hand,  the number of interior collocation points ($N_{\ins}$) appears to have little effect on the  ELM/A-TFC accuracy.  These behaviors  are similar to what have been observed with the ELM method in Section~\ref{elm_test}.

\begin{table}[tb]
\centering\small
\begin{tabular}{c|l|ll|ll}
\hline
 & & ELM & & ELM/A-TFC & \\ \cline{3-6}
Dimension & $N_{\bcd}$ & $e_{l^{\infty}}$ & $e_{l^2}$ & $e_{l^{\infty}}$ & $e_{l^2}$ \\ \hline
$d=3$ & 10 & 1.16E-4 & 1.11E-5 & 6.08E-5& 1.40E-5\\
& 100 & 1.03E-9 & 7.62E-11 & 1.10E-10 &3.57E-11 \\
& 500 & 2.96E-10 & 3.83E-11 & 6.76E-11&2.47E-11\\
& 1000 & 2.62E-10 & 3.58E-11 & 4.69E-11&1.50E-11\\ \hline
$d=7$ & 80 & 1.12E-5 & 1.09E-6 & 9.81E-5&7.24E-6\\
& 90 & 9.26E-6 & 7.30E-7 &6.85E-5 &4.64E-6 \\
& 100 & 8.16E-6 & 6.64E-7 & 2.96E-5&2.51E-6\\
& 110 & 7.30E-6 & 6.04E-7 & 2.48E-5&2.21E-6\\
\hline
\end{tabular}
\caption{\small Poisson equation: Comparison of $e_{l^{\infty}}$ and $e_{l^2}$ errors obtained with ELM and ELM/A-TFC corresponding to a set of boundary collocation points ($N_{\bcd}$). 
NN architecture: $[d,M,1]$, with $M=1000$ for $d=3$ and with $M=3000$ for $d=7$.
$N_{\ins}=1000$ in all cases.
}
\label{tab_11}
\end{table}

Finally we shows a comparison between the ELM method and the ELM/A-TFC method for solving the Poisson equation. Table~\ref{tab_11} lists the $e_{l^{\infty}}$ and $e_{l^2}$ errors obtained by ELM and ELM/A-TFC corresponding to a set of $N_{\bcd}$ values for two problem dimensions ($d=3$ and $7$). In these tests, the NN architecture is $\mbs M_{\text{arch}}=[d,M,1]$ with $M=1000$ for $d=3$ and $M=3000$ for $d=7$, and the interior collocation points is fixed at $N_{\ins}=1000$. It is observed that the accuracy with ELM and ELM/A-TFC is generally comparable, and the ELM/A-TFC method appears to be slightly more accurate for lower dimensions. This can be attributed to the fact that  the  A-TFC resembles the full TFC more closely in lower dimensions. Therefore A-TFC enforces the boundary conditions more accurately (closer to  TFC) in lower dimensions.
On the oter hand, we  note that the computational effort and cost involved in ELM/A-TFC is generally higher than that of ELM, because of the computations associated with the A-TFC terms.

\subsubsection{Nonlinear Poisson Equation}

In this subsection we test the ELM/A-TFC method using the  nonlinear Poisson problem from Section~\ref{sec_312}, under the same problem settings and parameters.

\begin{table}[tb]
\centering\small
	\begin{tabular}{cccccc}
		\hline
	$R_m$	&0.01&0.05&0.1&0.5&1.0 \\ 
 $e_{l^\infty}$&5.82E-6&5.21E-9&2.92E-9&3.47E-6&4.85E-1 \\ 
		$e_{l^{2}}$ &1.64E-6&9.78E-10&5.72E-10&7.34E-7&6.71E-2 \\ 
		\hline
	\end{tabular}
	\caption{\small Nonlinear Poisson equation ($d=3$): Determining $R_{m0}$ for the ELM/A-TFC method based on the procedure from Remark~\ref{rem_a6}.
 NN architecture: $[d,500,1]$; $(N_{\bcd},N_{\ins})=(100,10)$.
 }
	\label{tab_12}
\end{table}

\begin{table}[tb]
\centering\small
\begin{tabular}{l| cc}
\hline
$d$ & 3 & 7  \\ \hline
$R_{m0}$ & 0.1 & 0.01 \\
\hline
\end{tabular}
\caption{\small Nonlinear Poisson equation: $R_{m0}$ for ELM/A-TFC determined by the procedure from Remark~\ref{rem_a6}. The simulation parameters here for $d=3$ follow those of Table~\ref{tab_10}. For $d=7$, NN architecture: $[d,3000,1]$, $(N_{\bc},N_{\ins})=(300,100)$. 
}
\label{tab_13}
\end{table}

Tables~\ref{tab_12} and~\ref{tab_13} show the tests for determining the $R_{m0}$ using the procedure from Remark~\ref{rem_a6}. The results lead to $R_{m0}\approx 0.1$ for dimension $d=3$ and $R_{m0}\approx 0.01$ for $d=7$. These values are employed for generating the hidden-layer coefficients in the neural network in the subsequent simulations.

 \begin{figure}[!htb]
	\centerline{
	\subfigure
	{ \includegraphics[width=0.26\textwidth]{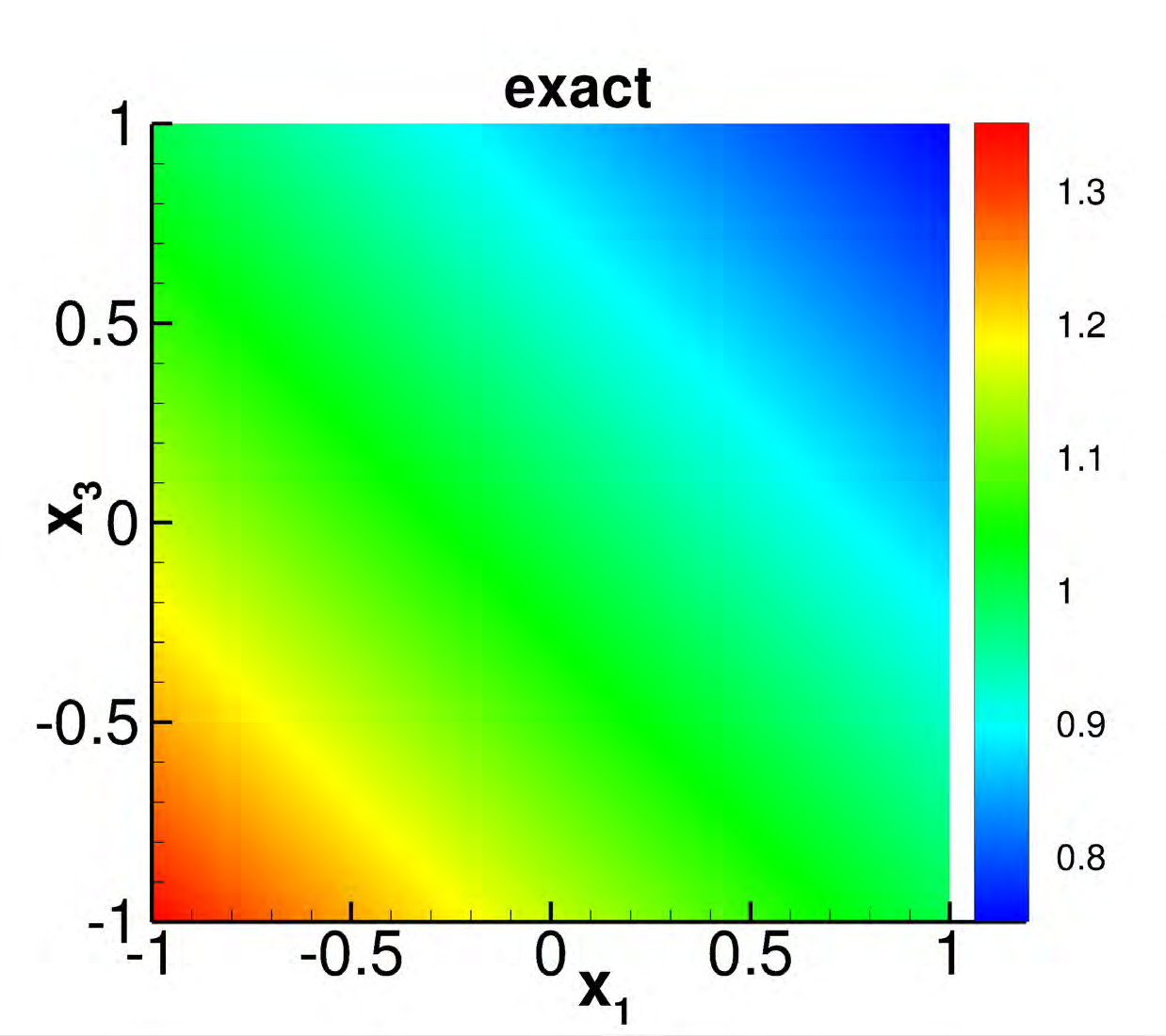}(a)}
	\subfigure
	{ \includegraphics[width=0.26\textwidth]{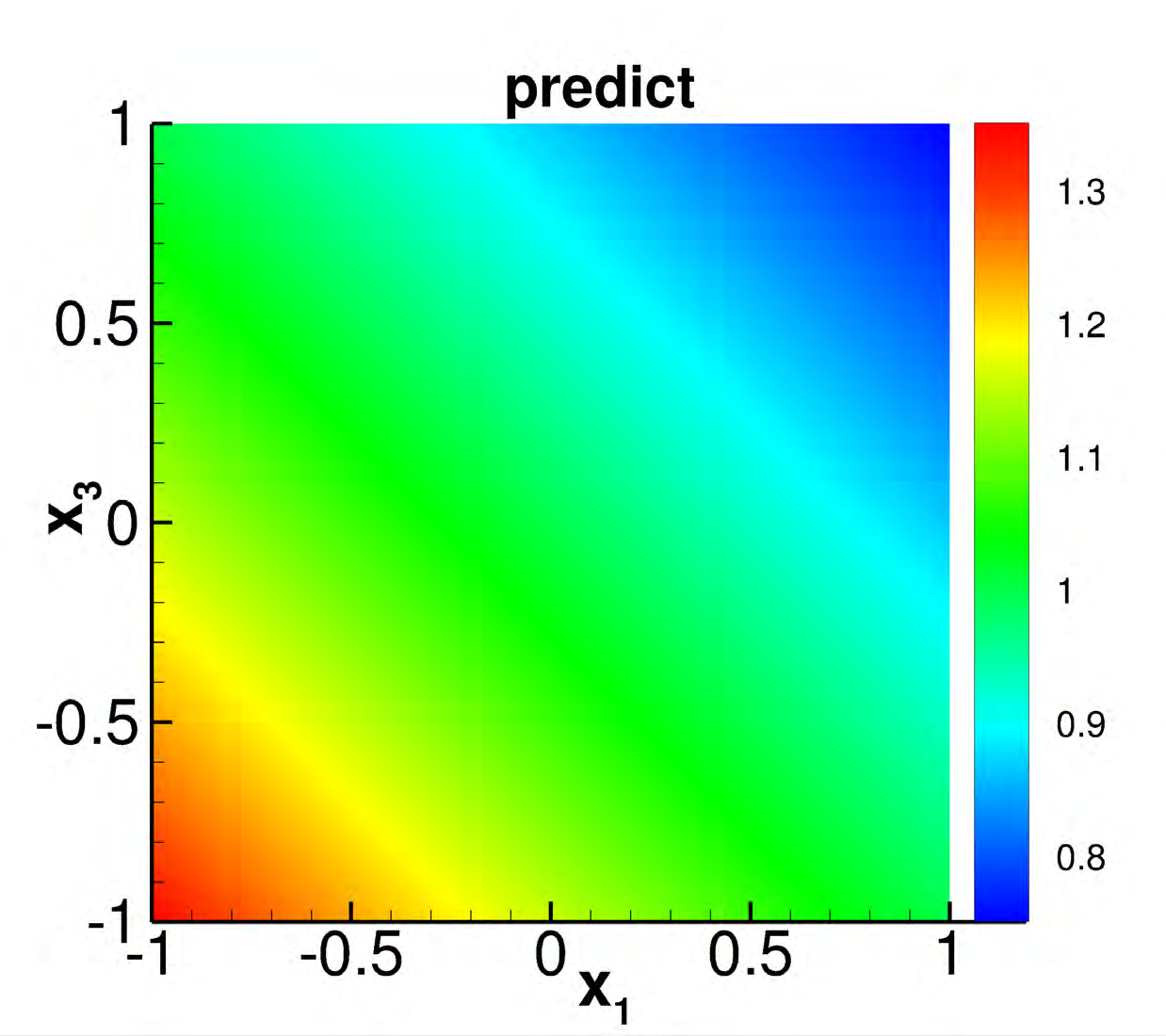}(b)}
	\subfigure
	{ \includegraphics[width=0.26\textwidth]{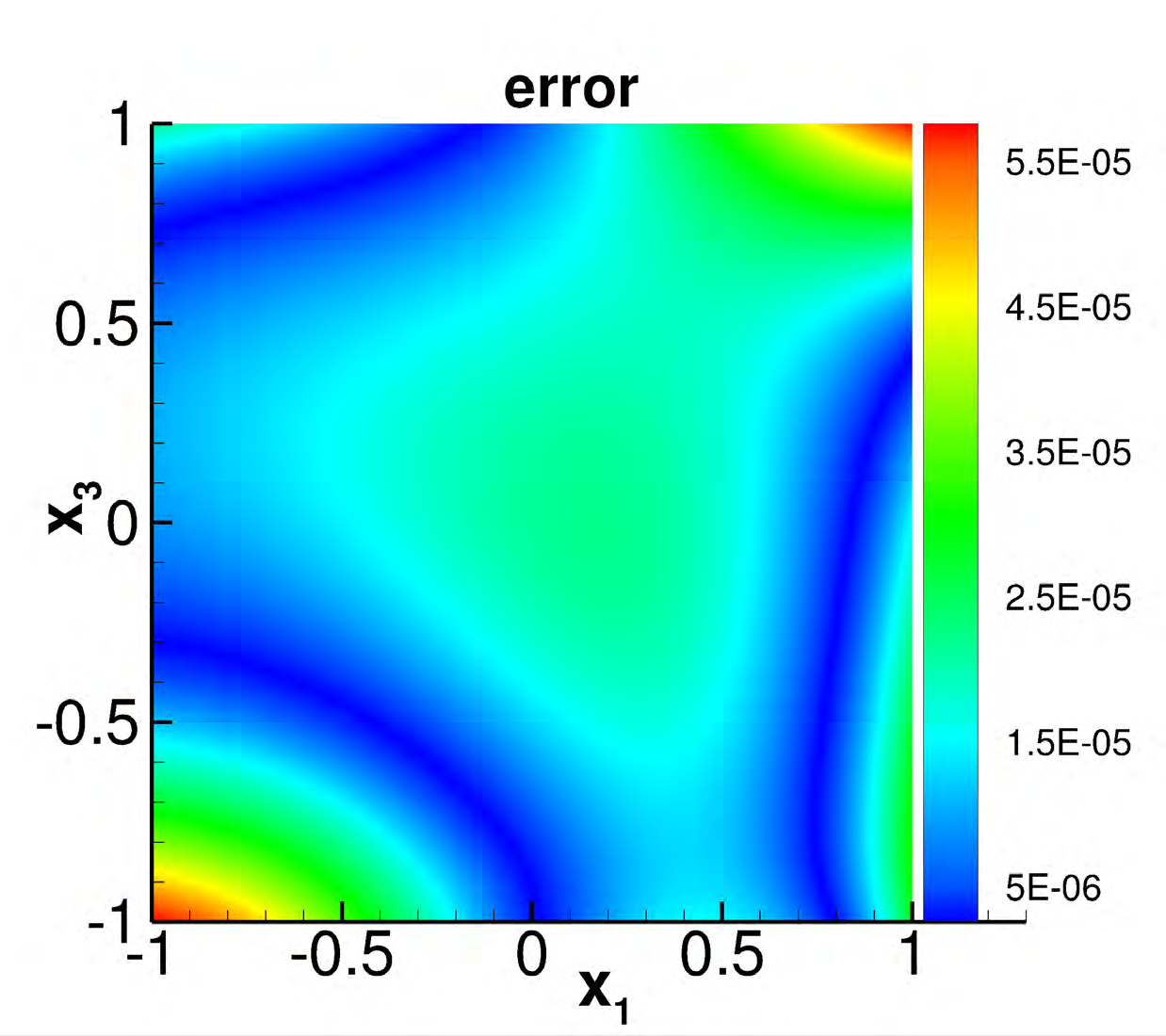}(c)}
 }
 \centerline{
	\subfigure
	{ \includegraphics[width=0.26\textwidth]{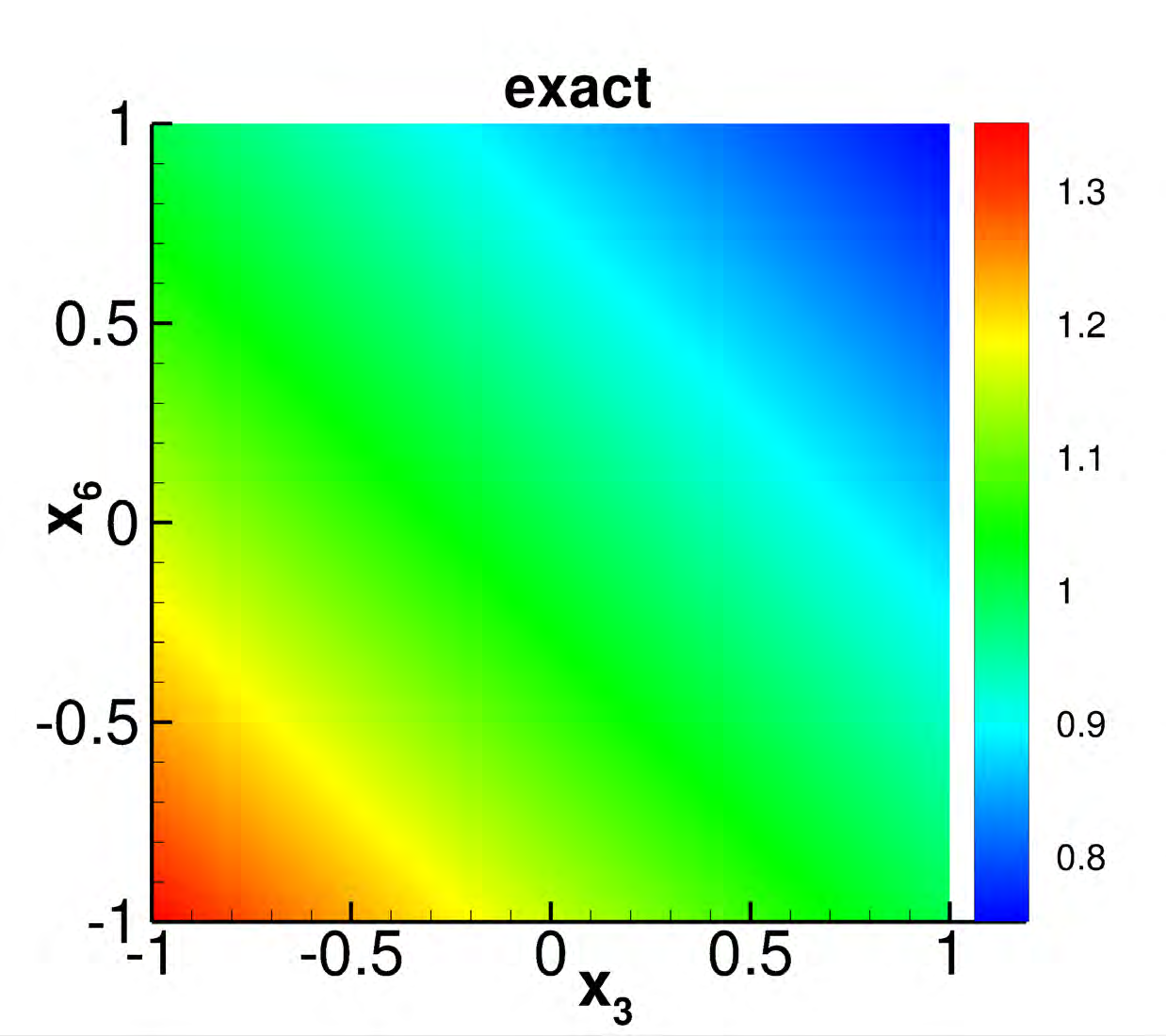}(d)}
	\subfigure
	{ \includegraphics[width=0.26\textwidth]{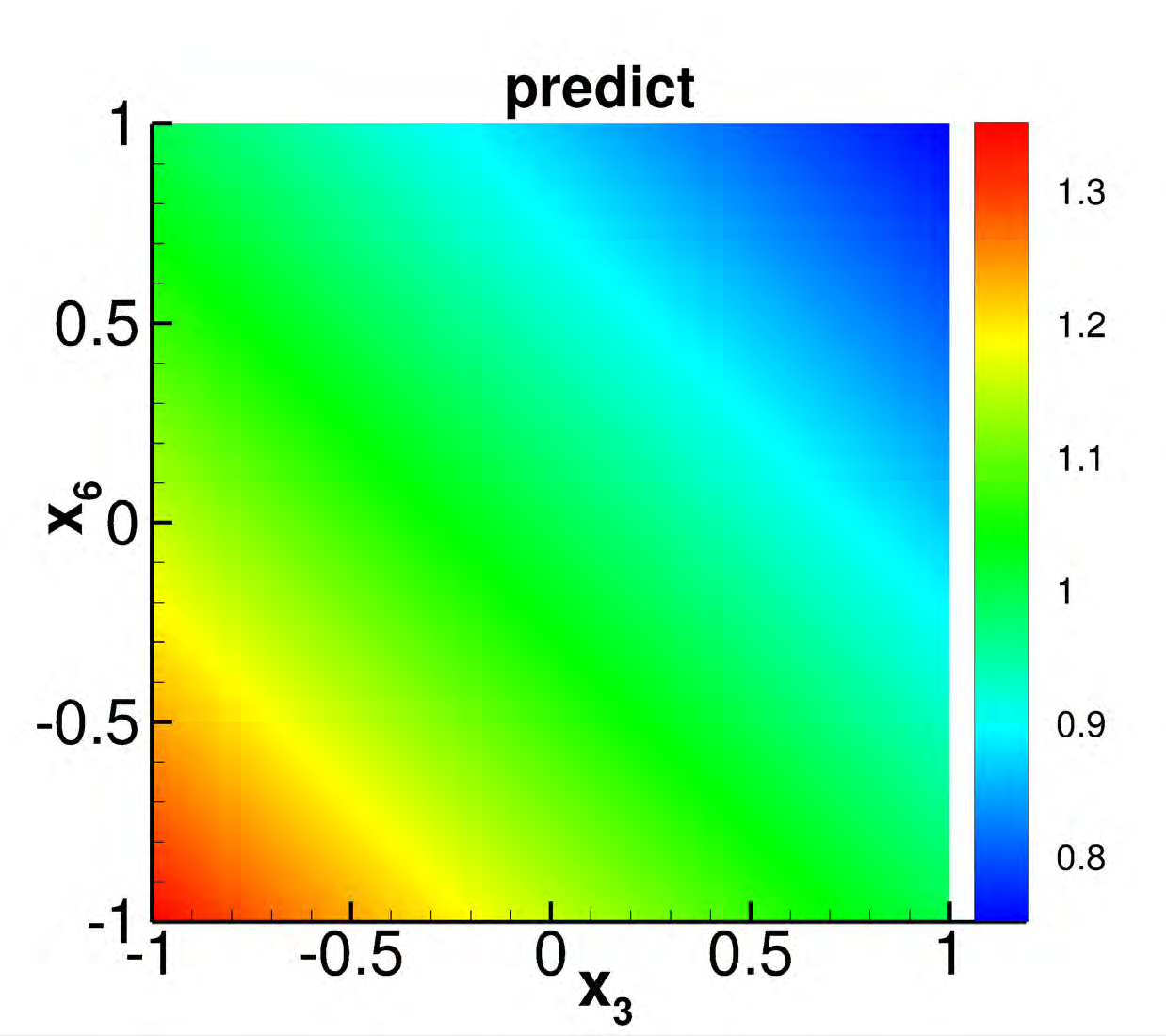}(e)}
	\subfigure
	{ \includegraphics[width=0.26\textwidth]{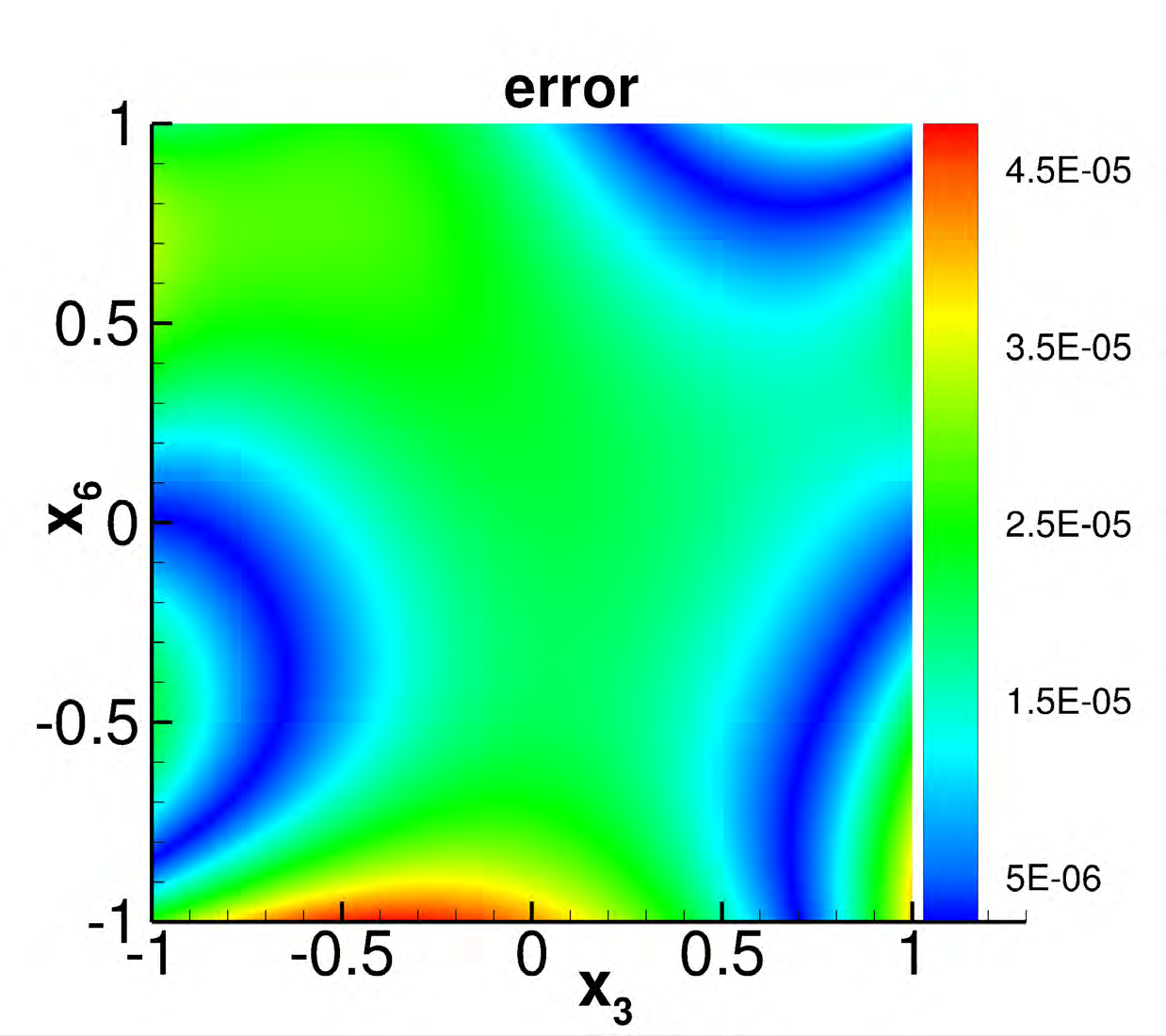}(f)}
 }
 \centerline{
	\subfigure
	{ \includegraphics[width=0.26\textwidth]{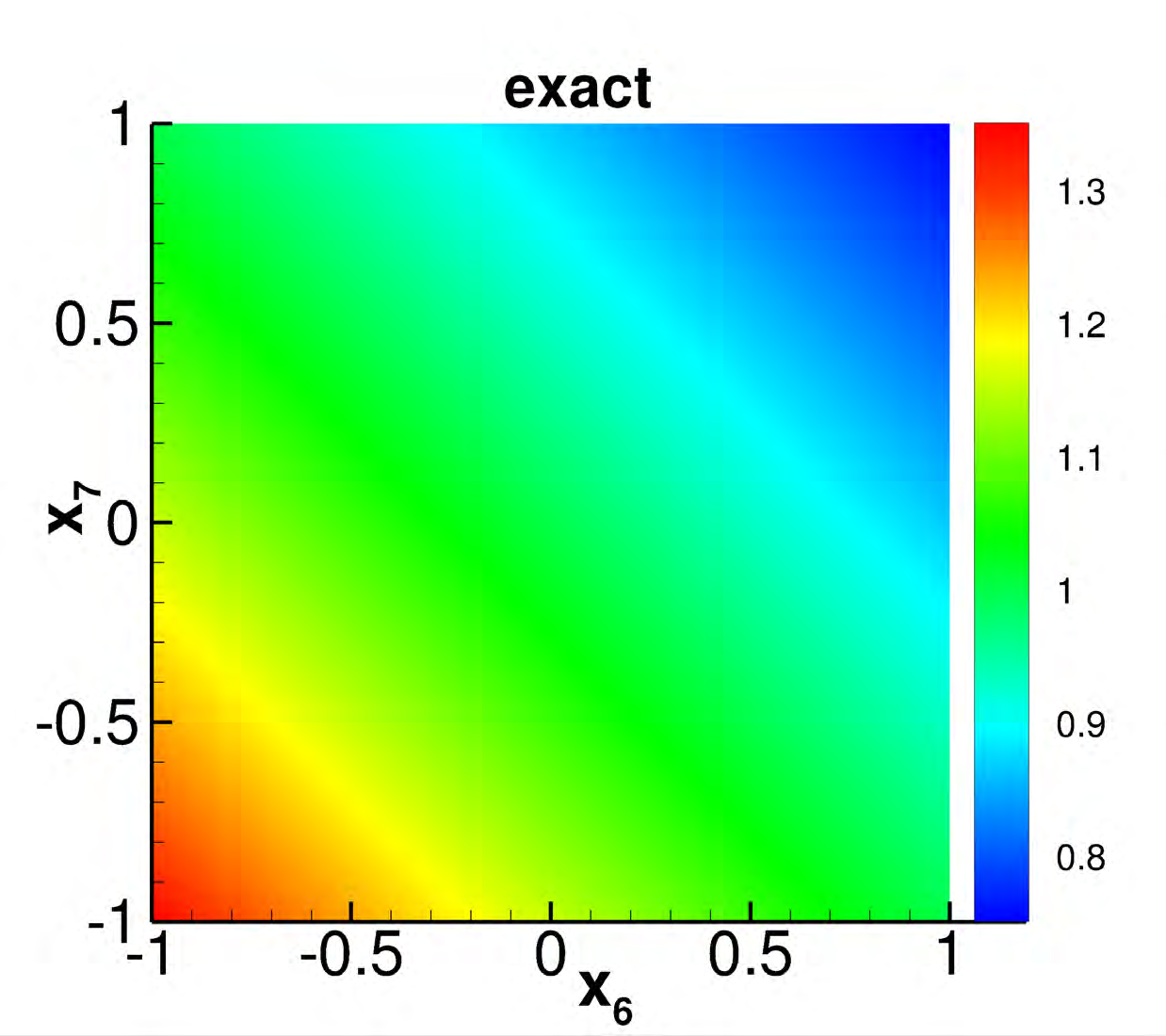}(g)}
	\subfigure
	{ \includegraphics[width=0.26\textwidth]{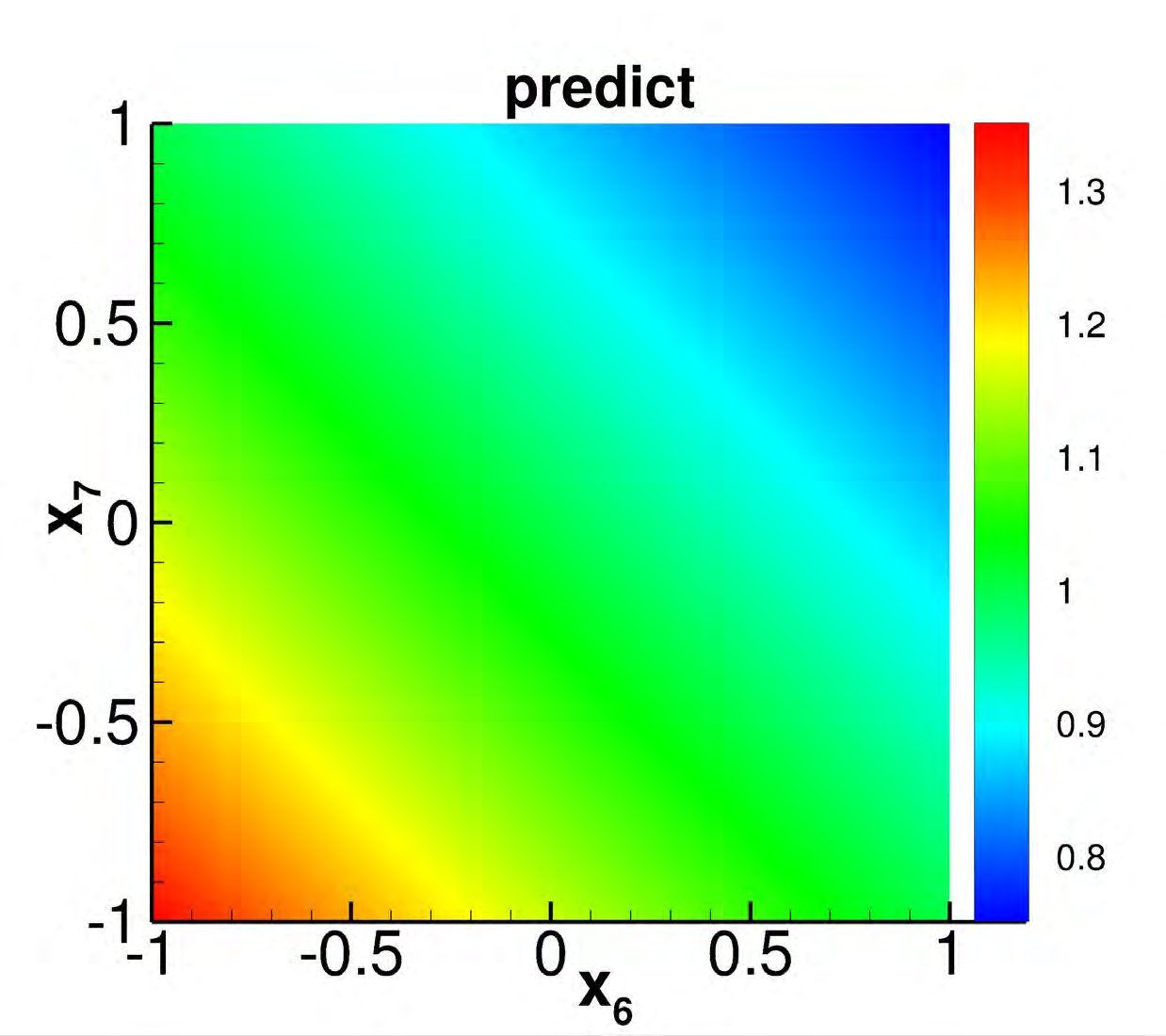}(h)}
	\subfigure
	{ \includegraphics[width=0.26\textwidth]{tfc_poi/67error-eps-converted-to.pdf}(i)}
 }
	\caption{\small Nonlinear Poisson equation ($d=7$): Distributions of the exact solution (left column), the ELM/A-TFC solution (middle column), and the point-wise absolute error of ELM/A-TFC in several cross sections (top row: $x_1$-$x_3$ plane; middle row: $x_3$-$x_6$ plane; bottom row: $x_6$-$x_7$ plane)  of the domain. These cross sections are located in the middle of the domain with respect to the rest of the coordinates.
 NN architecture: $[d,1000,1]$, $(N_{\bcd},N_{\ins})=(50,10)$.
 }
 \label{fg_20}
\end{figure}

Figure~\ref{fg_20} shows distributions of the exact solution, the ELM/A-TFC solution, and the point-wise absolute error of ELM/A-TFC in several cross sections of the domain for the nonlinear Poisson equation in dimension $d=7$. The ELM/A-TFC results are obtained with an NN  architecture $\mbs M_{\text{arch}}=[d,1000,1]$ and the random collocation points characterized by $(N_{\bcd},N_{\ins})=(50,10)$. The ELM/A-TFC method has captured the solution accurately, with the maximum error on the order $10^{-5}$ in the $x_1$-$x_3$ and $x_3$-$x_6$ planes and on the order of $10^{-7}$ in the $x_6$-$x_7$ plane.

\begin{figure}[tb]
	\centering
	\subfigure[d=3]
	{ \includegraphics[width=0.3\textwidth]{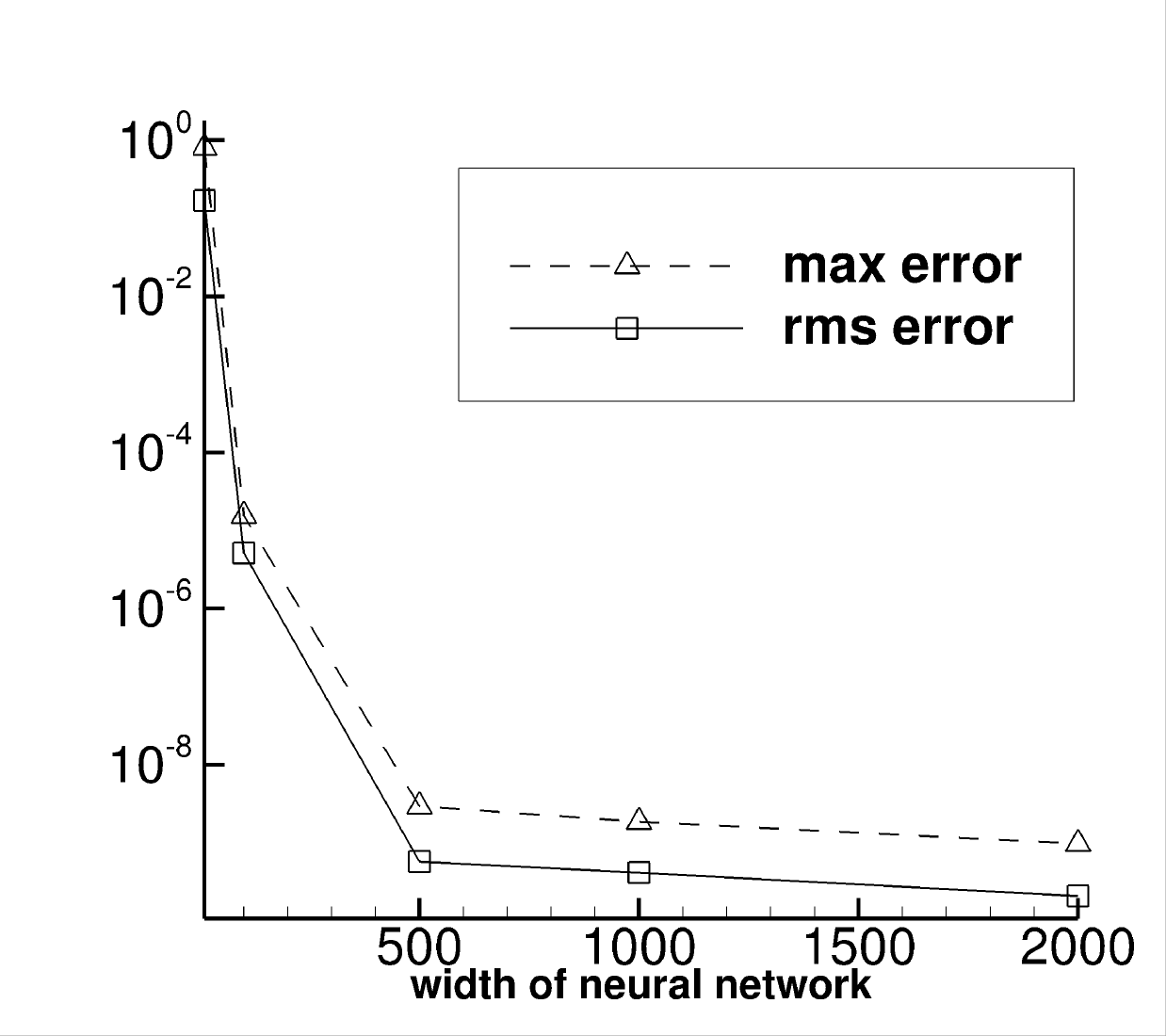}}
	\subfigure[dim=7]
	{ \includegraphics[width=0.3\textwidth]{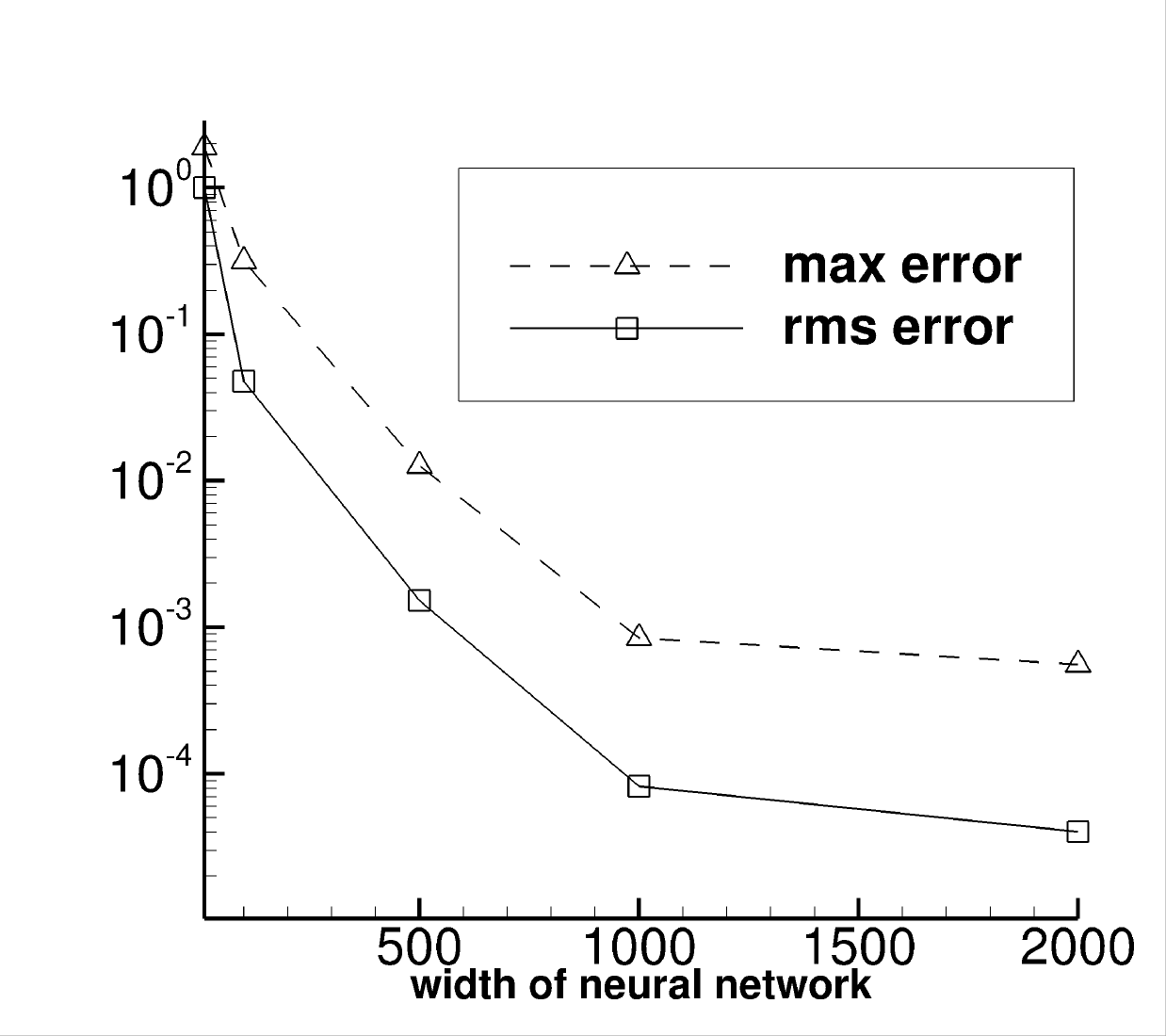}}
	\caption{\small Nonlinear Poisson equation: $e_{l^{\infty}}$ and $e_{l^2}$ errors versus the number of training parameters ($M$) in the ELM/A-TFC  network.
 NN architecture: $[d,M,1]$ (with $M$ varied); $(N_{\bcd},N_{\ins})=(100,10)$ in (a) and $(50,100)$ in (b).
 }
 \label{tfc_nlpoi_err_width}
\end{figure}

\begin{figure}[tb]
	\centering
	\subfigure[d=3]
	{ \includegraphics[width=0.3\textwidth]{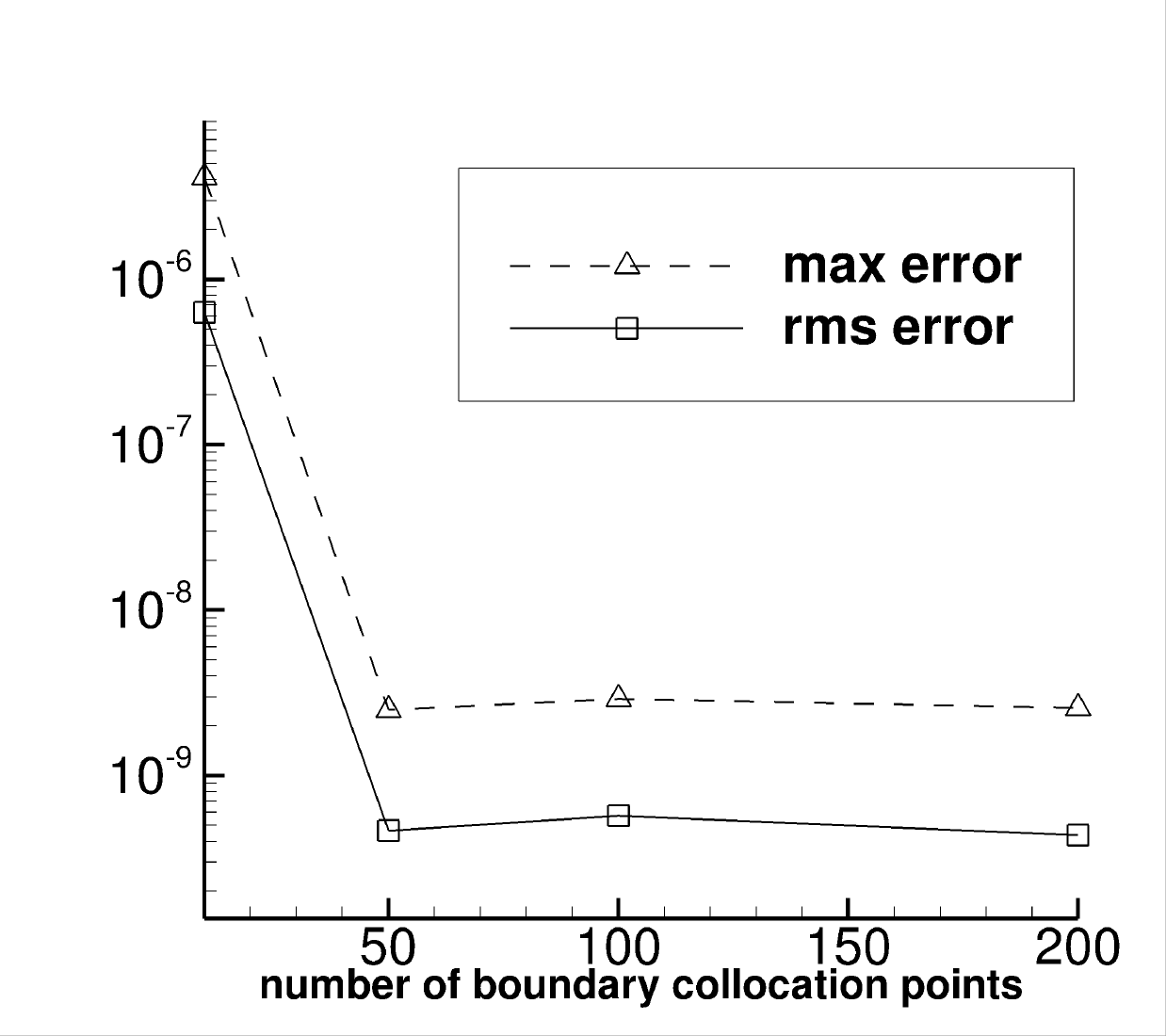}}
	\subfigure[d=7]
	{ \includegraphics[width=0.3\textwidth]{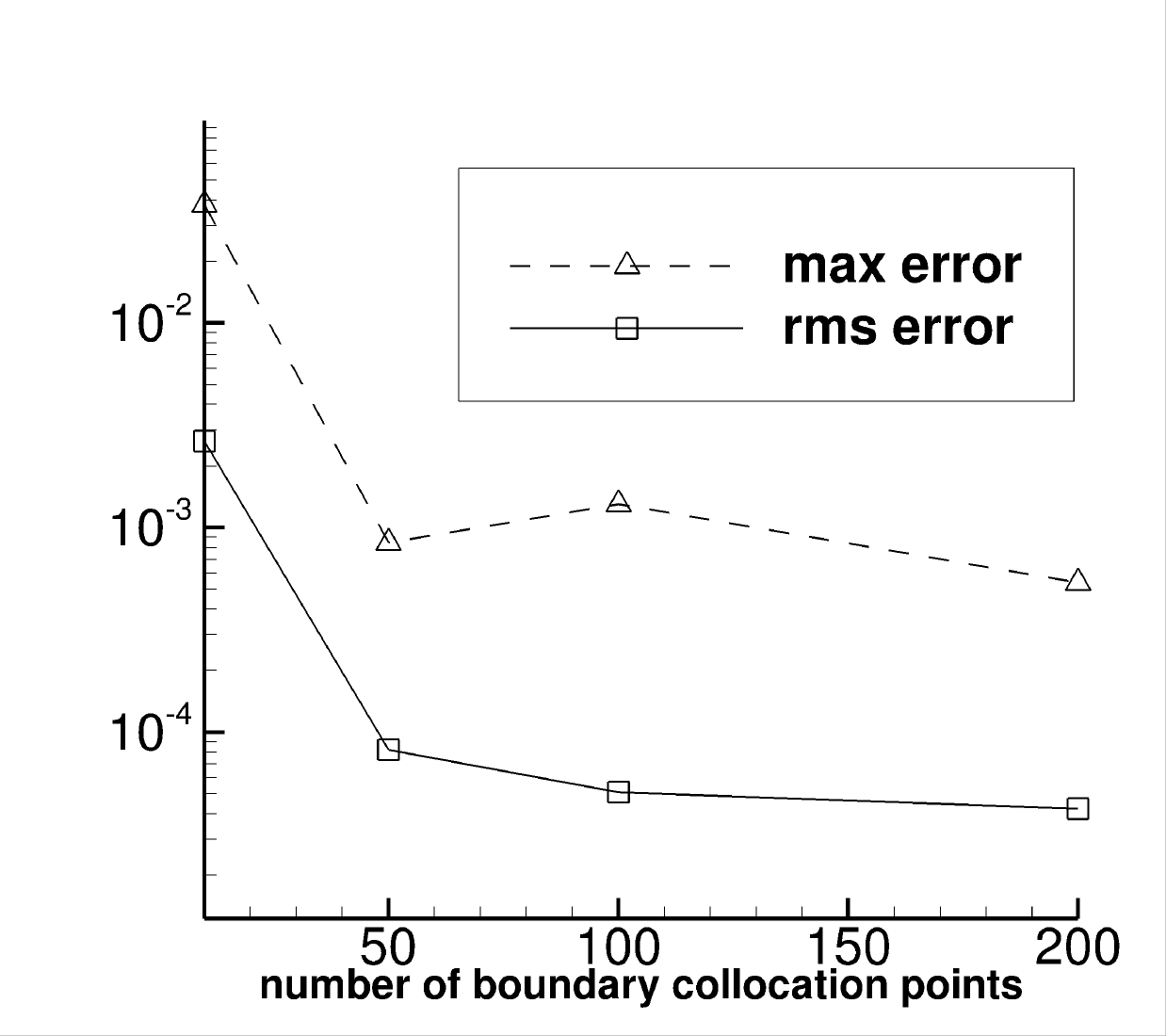}}
	\caption{\small Nonlinear Poisson equation: $e_{l^{\infty}}$ and $e_{l^2}$ errors of ELM/A-TFC versus the number of boundary collocation points ($N_{\bcd}$). 
 NN architecture: $[d,M,1]$, with $M=500$ in (a) and $M=1000$ in (b); $N_{\ins}=10$ in (a) and $N_{\ins}=100$ in (b); $N_{\bcd}$ varied in (a,b).
 }\label{tfc_nlpoi_err_nc}
\end{figure}

The convergence behavior of the ELM/A-TFC method is illustrated by Figures~\ref{tfc_nlpoi_err_width} and~\ref{tfc_nlpoi_err_nc} for problem dimensions $d=3$ and $d=7$. 
Figure~\ref{tfc_nlpoi_err_width} shows the $e_{l^{\infty}}$ and $e_{l^2}$ errors as a function of the number of training parameters in the neural network. In this set of tests, the number of boundary/interior collocation points is fixed while the number of training parameters is varied.
We can observe a rapid decrease (approximately exponential) in the ELM/A-TFC errors as the number of training parameters increases (before saturation).
Figure~\ref{tfc_nlpoi_err_nc} shows the $e_{l^{\infty}}$ and $e_{l^2}$ errors as a function of the number of boundary collocation points ($N_{\bcd}$). In this set of tests the number of training parameters and the number of interior collocation points are fixed, while the number of boundary collocation points ($N_{\bcd}$) is varied systematically.
The ELM/A-TFC errors initially decrease rapidly as $N_{\bcd}$ increases, and then level off when $N_{\bcd}$ increases beyond a certain level.
We would like to further point out that the number of interior collocation points has little effect on the ELM/A-TFC accuracy (result not shown here), similar to what has been observed with the ELM method.

\subsubsection{Heat Equation}

We next consider the domain $\Omega=[-1,1]^d$ and the heat equation on $\Omega\times[0,T]$ (with $T=1$), 
\begin{align}
\partial_t u-\Delta u&=f(x,t), \quad(x,t)\in \Omega \times [0,T],\\
u(x,t) &= g(x,t), \quad (x,t)\in \partial \Omega \times [0,T],\\
u(x,0)&=h(x), \qquad x \in \Omega,
\end{align}
where  $f(x,t)=(\frac{1}{d}-1)\cos\left(\frac{1}{d}\sum_{i=1}^{d}x_i\right)\exp(-t)$, $g(x,t)=\cos\left(\frac{1}{d}\sum_{i=1}^{d}x_i\right)\exp(-t)$, and $h(x)=\cos\left(\frac{1}{d}\sum_{i=1}^{d}x_i\right)$. This problem has the  exact solution $u(x,t)=\cos\left(\frac{1}{d}\sum_{i=1}^{d}x_i\right)\exp(-t)$.

The simulation settings and the notations here follow those of Section~\ref{sec_314}. We employ a neural network architecture $\mbs M_{\text{arch}}=[d+1,M,1]$, where the $(d+1)$ input nodes denote $(x,t)$ and $M$ is the number of training parameters in the network. $N_{\bcd}$, $N_{\ins}$ and $N_{t_0}$ denote the number of random collocation points on each of the boundary of $\partial\Omega\times[0,T]$, on the interior of $\Omega\times[0,T]$, and on $\Omega$ at $t=0$, respectively.

\begin{table}[tb]
\centering\small
	\begin{tabular}{ccccccc}
		\hline
		$R_m$&1E-3&1E-2&5E-2&1E-1&5E-1&1\\  
		$e_{l^\infty}$&6.05E-3&4.15E-4&3.56E-6&5.32E-7&2.64E-5&8.21E-4\\
		$e_{l^{2}}$ &1.16E-3&6.18E-5&6.37E-7&9.35E-8&6.03E-6&1.54E-4\\
		\hline
	\end{tabular}
	\caption{\small Heat equation: Determination of $R_{m0}$ for ELM/A-TFC using the procedure from Remark~\ref{rem_a6} for $d=3$.
 NN architecture: $[d+1,1000,1]$; $(N_{\bcd},N_{\ins},N_{t_0})=(100,100,1000)$.
 }
	\label{tab_16_RM}
\end{table}

\begin{table}[tb]
\centering\small
\begin{tabular}{l| cc}
\hline
$d$ & 3 & 7  \\ \hline
$R_{m0}$ & 0.1 & 0.005 \\
\hline
\end{tabular}
\caption{\small Heat equation: $R_{m0}$ for ELM/A-TFC determined by the procedure from Remark~\ref{rem_a6}. The simulation parameters here follow those of Table~\ref{tab_16_RM}. 
}
\label{tab_a17}
\end{table}

Table~\ref{tab_16_RM} summarizes the tests for determining the $R_{m0}$ with the ELM/A-TFC method using the procedure from Remark~\ref{rem_a6} for $d=3$, which lead to $R_{m0}\approx 0.1$. Table~\ref{tab_a17} lists the $R_{m0}$ values for the problem dimensions we have considered for this problem. The hidden-layer coefficients are set to uniform random values generated on $[-R_{m0},R_{m0}]$ in the subsequent simulations.

\begin{figure}[!tb]
	\centerline{
	\subfigure
	{ \includegraphics[width=0.26\textwidth]{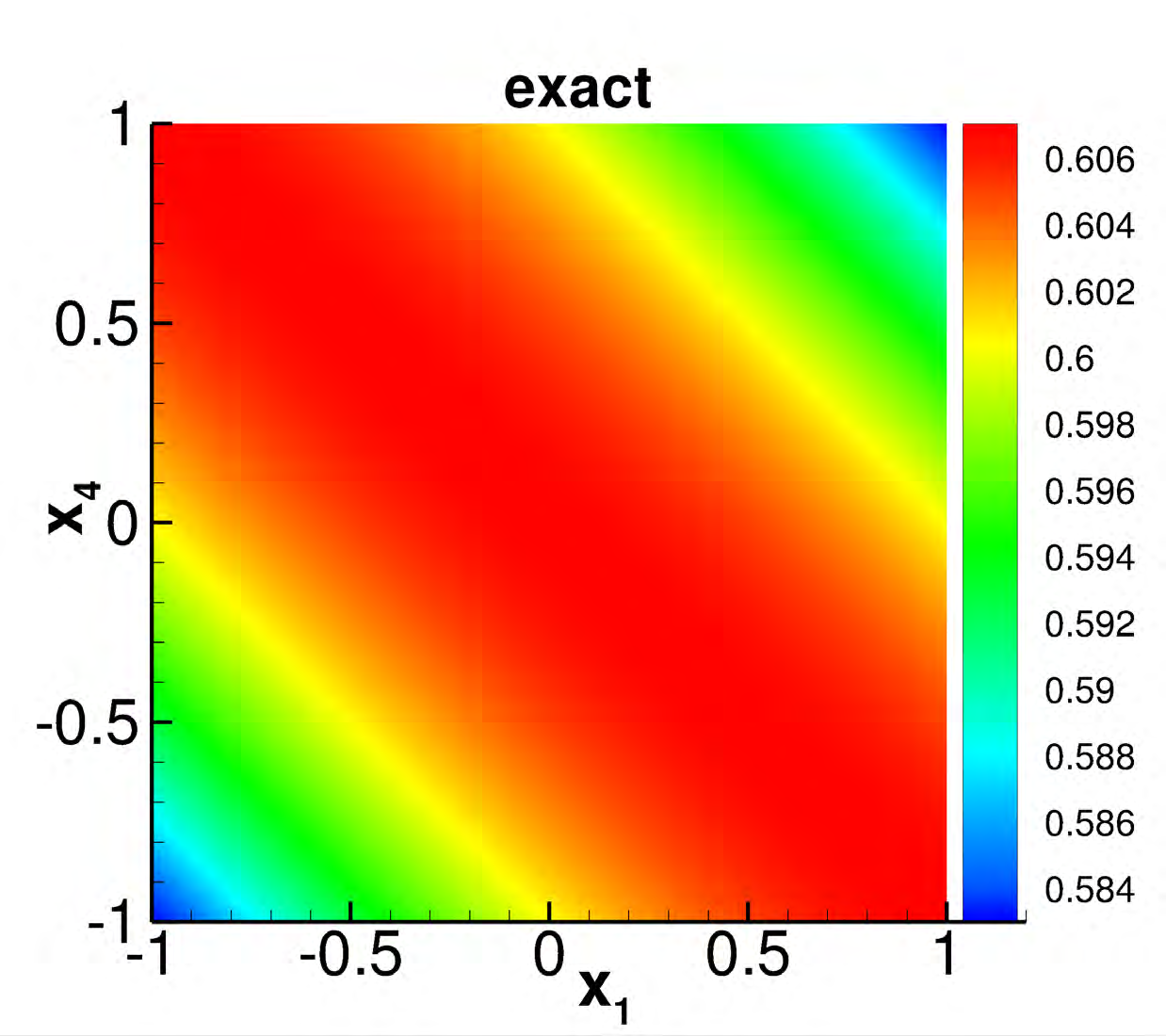}(a)}
	\subfigure
	{ \includegraphics[width=0.26\textwidth]{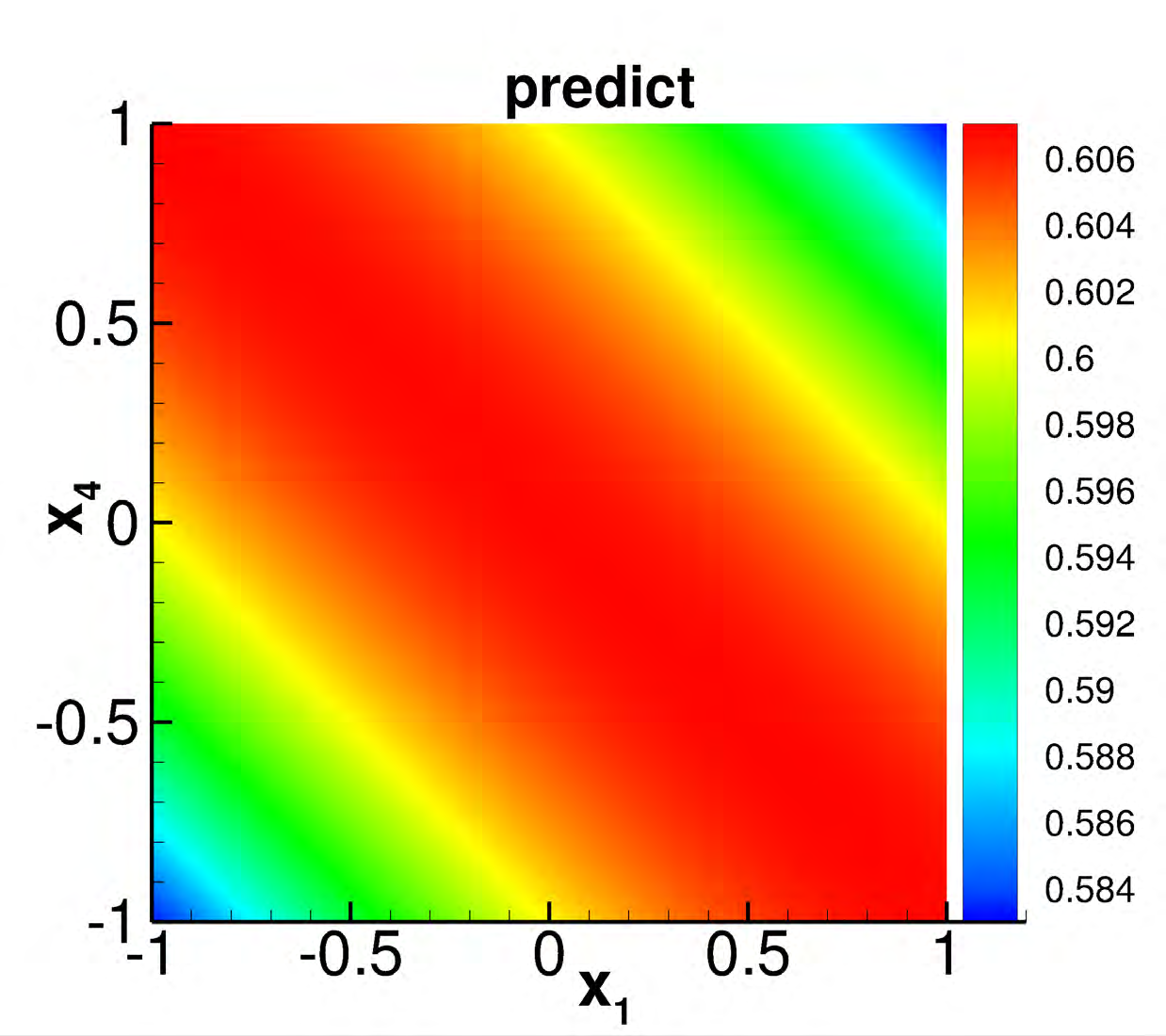}(b)}
	\subfigure
	{ \includegraphics[width=0.26\textwidth]{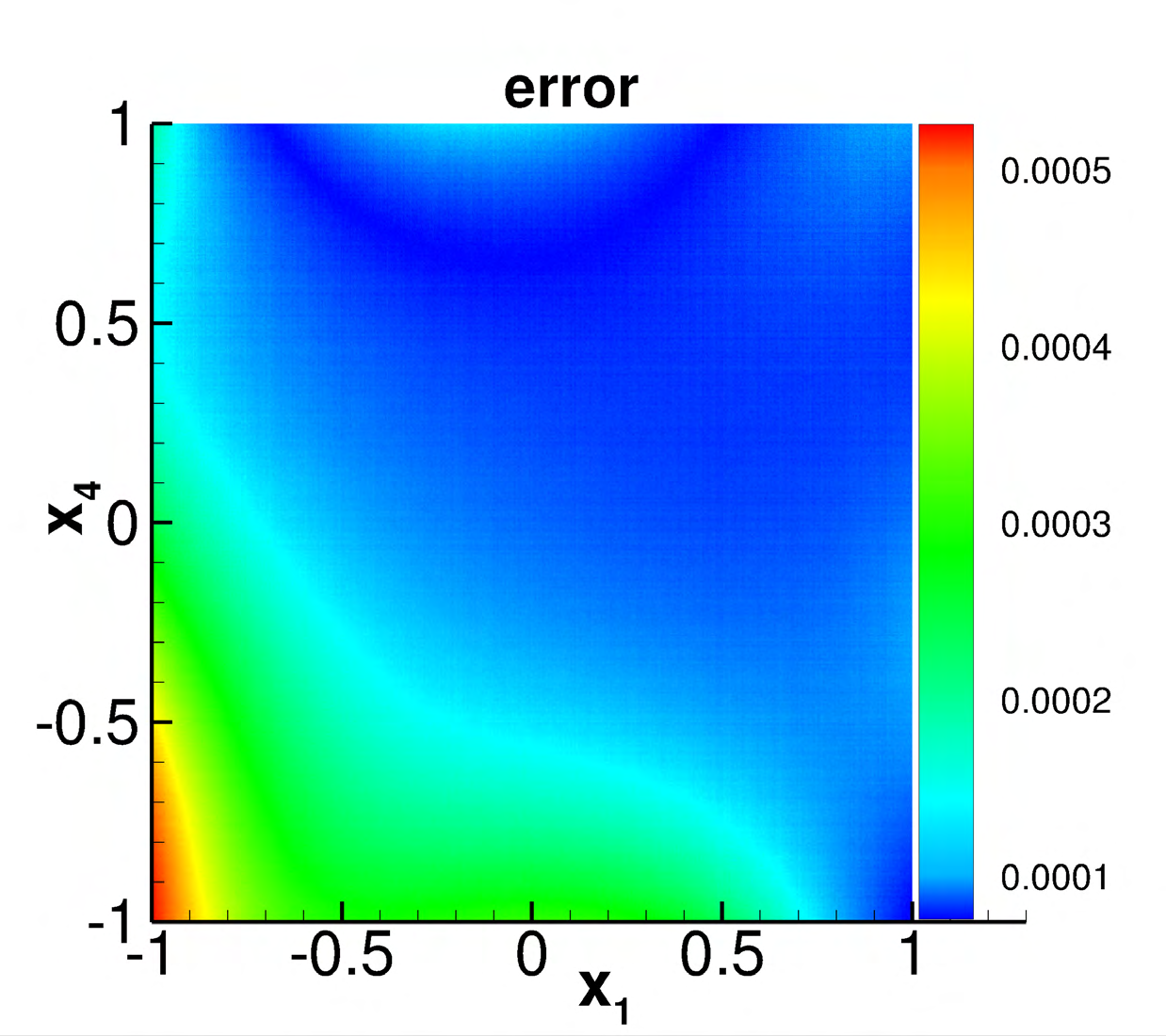}(c)}
 }
 \centerline{
	\subfigure
	{ \includegraphics[width=0.26\textwidth]{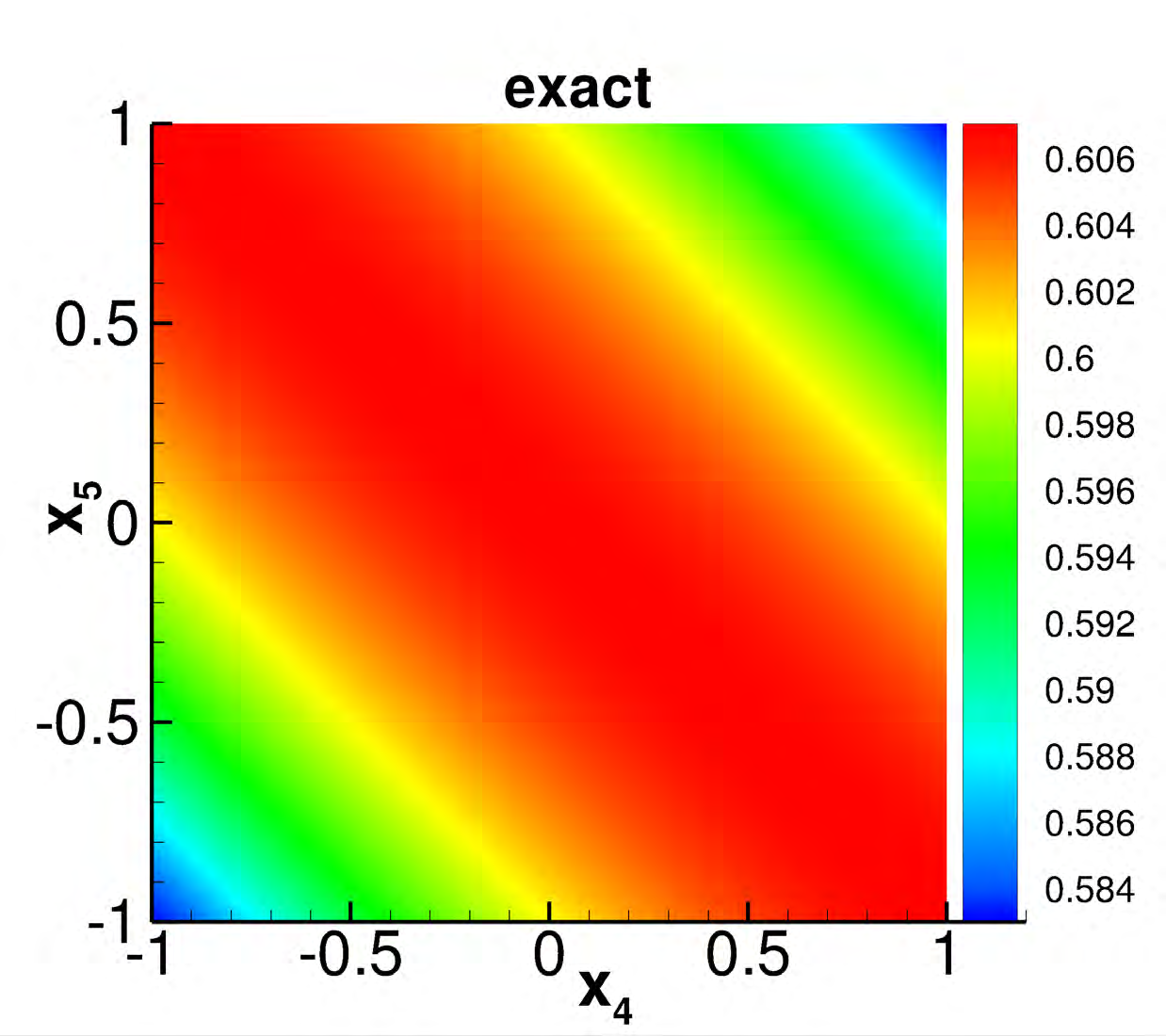}(d)}
	\subfigure
	{ \includegraphics[width=0.26\textwidth]{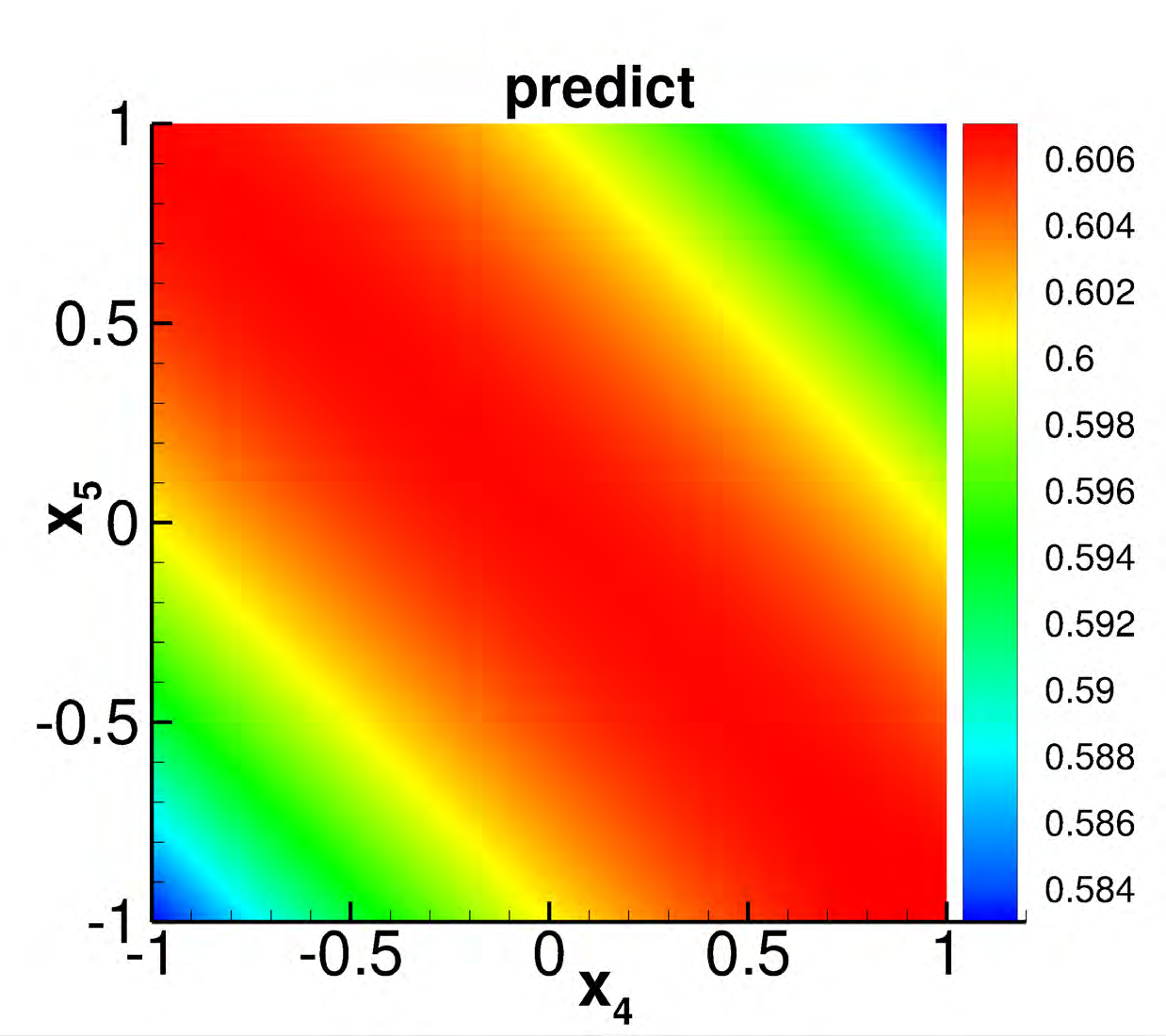}(e)}
	\subfigure
	{ \includegraphics[width=0.26\textwidth]{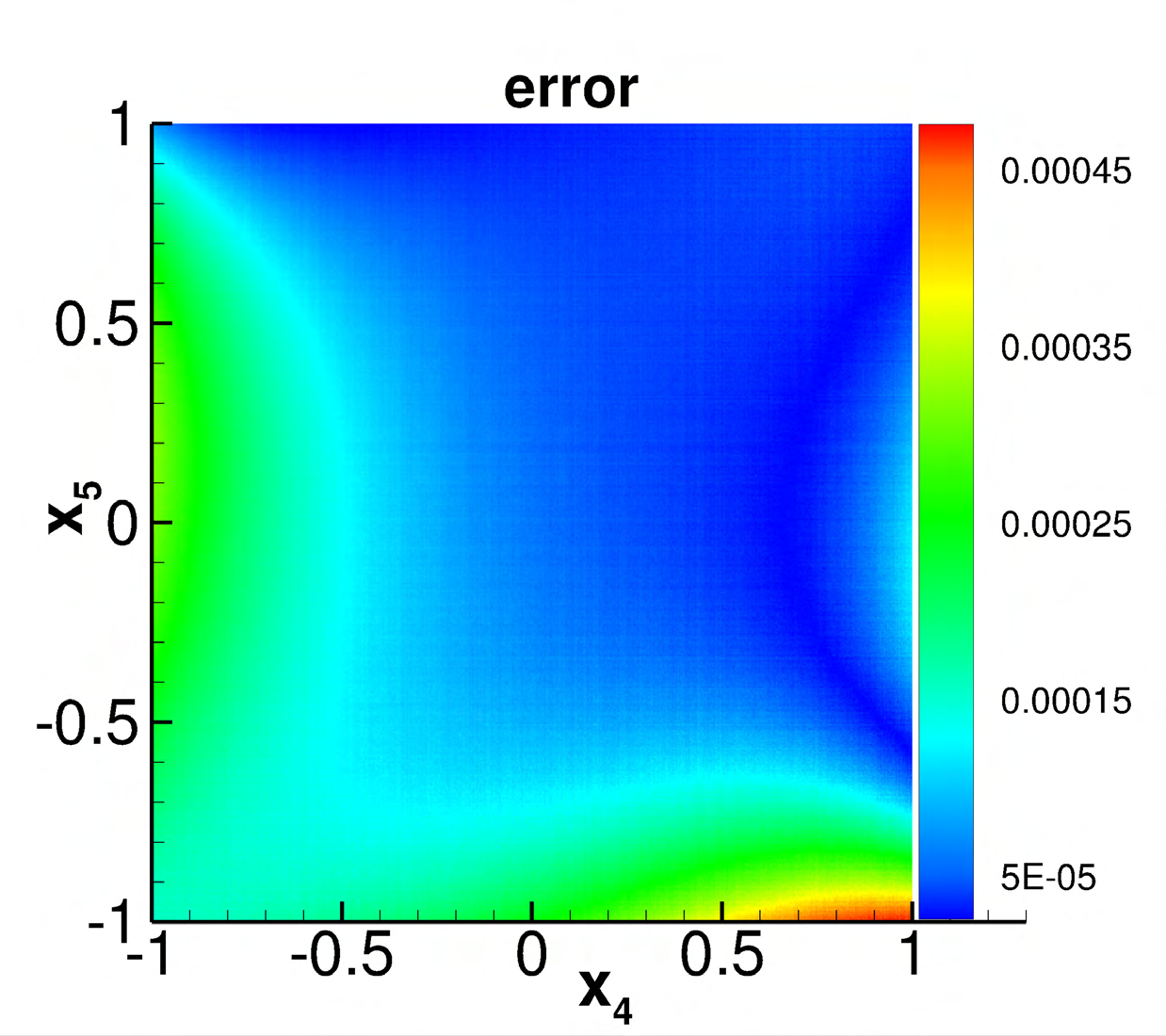}(f)}
 }
 \centerline{
	\subfigure
	{ \includegraphics[width=0.26\textwidth]{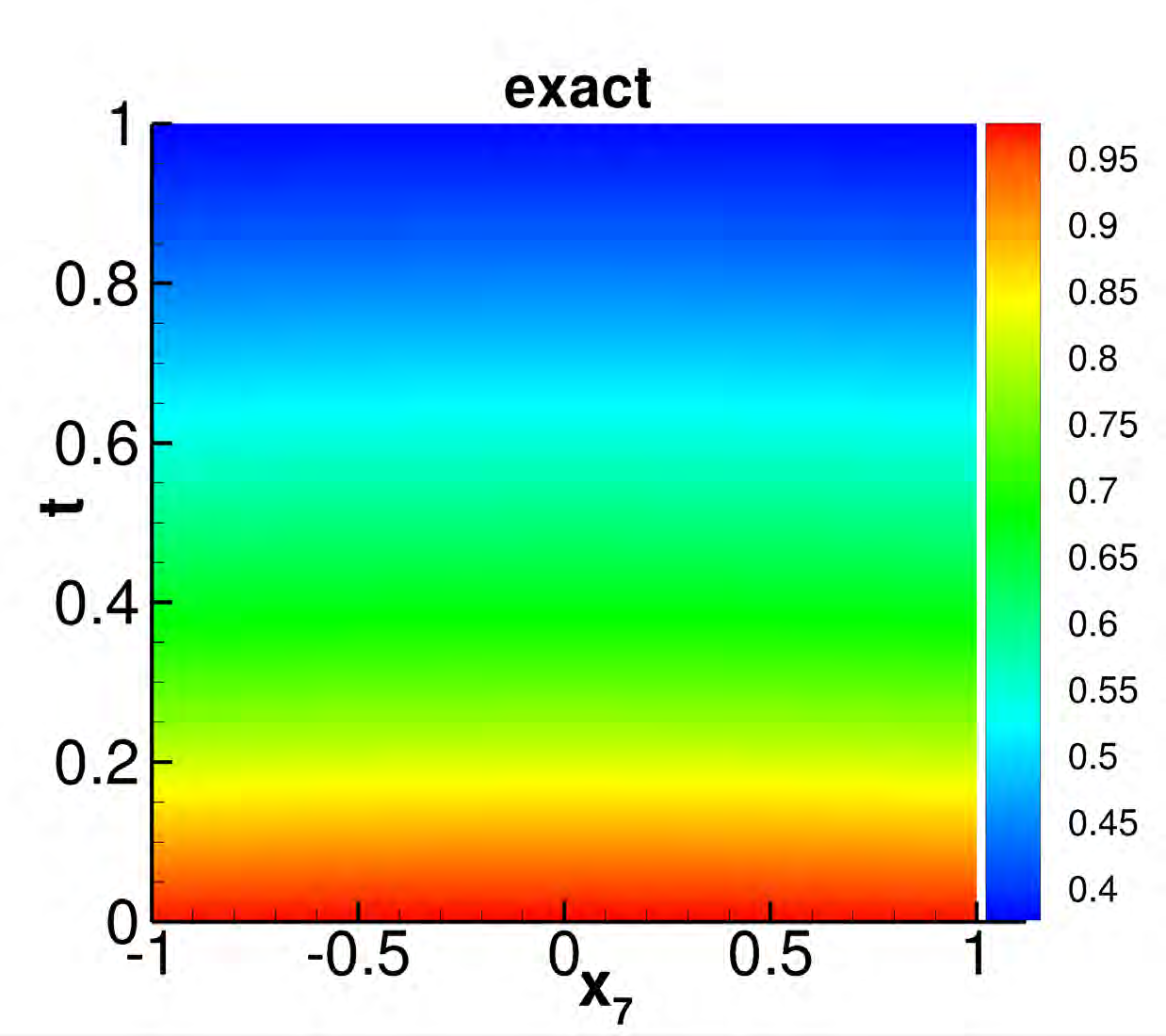}(g)}
	\subfigure
	{ \includegraphics[width=0.26\textwidth]{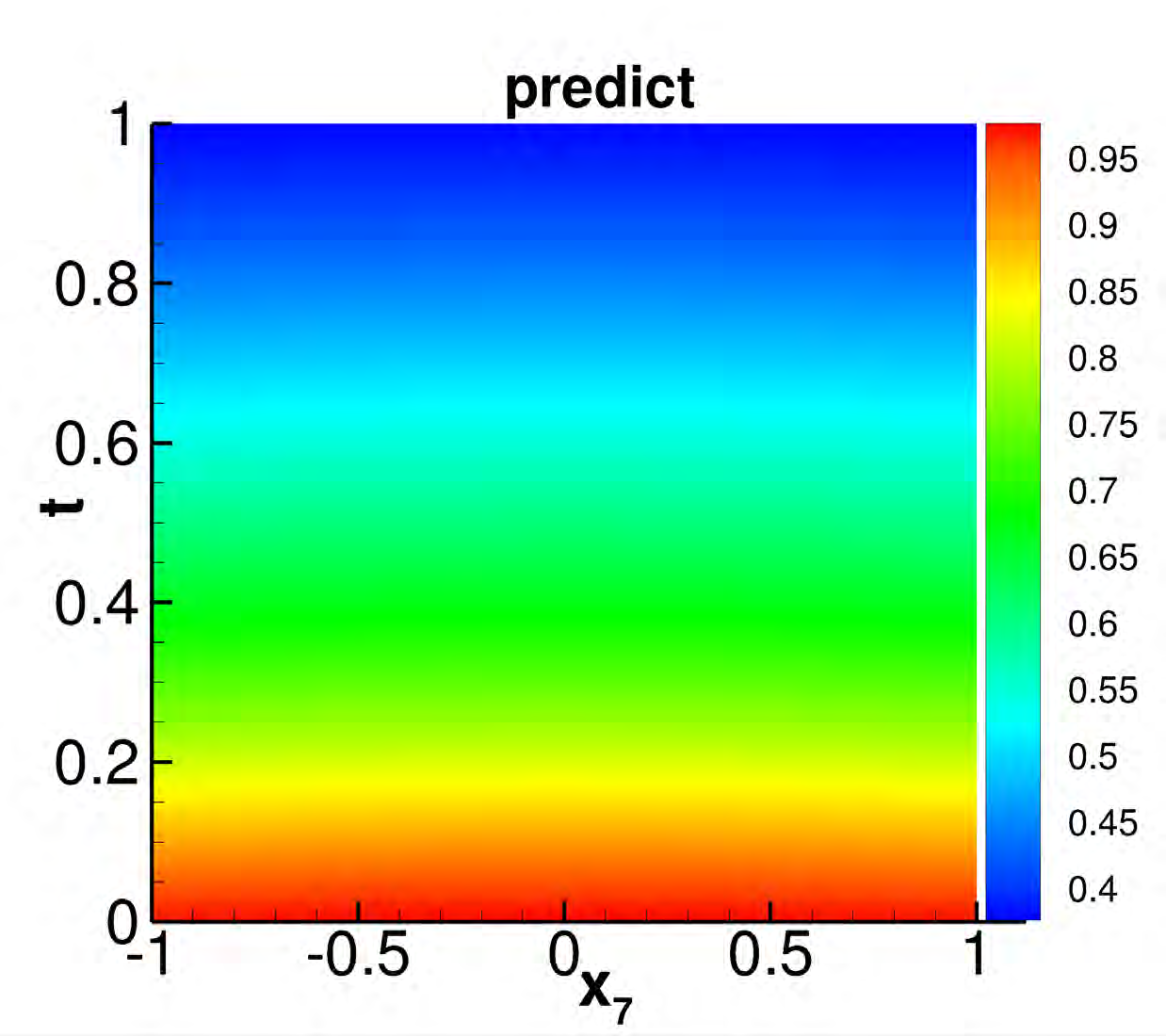}(h)}
	\subfigure
	{ \includegraphics[width=0.26\textwidth]{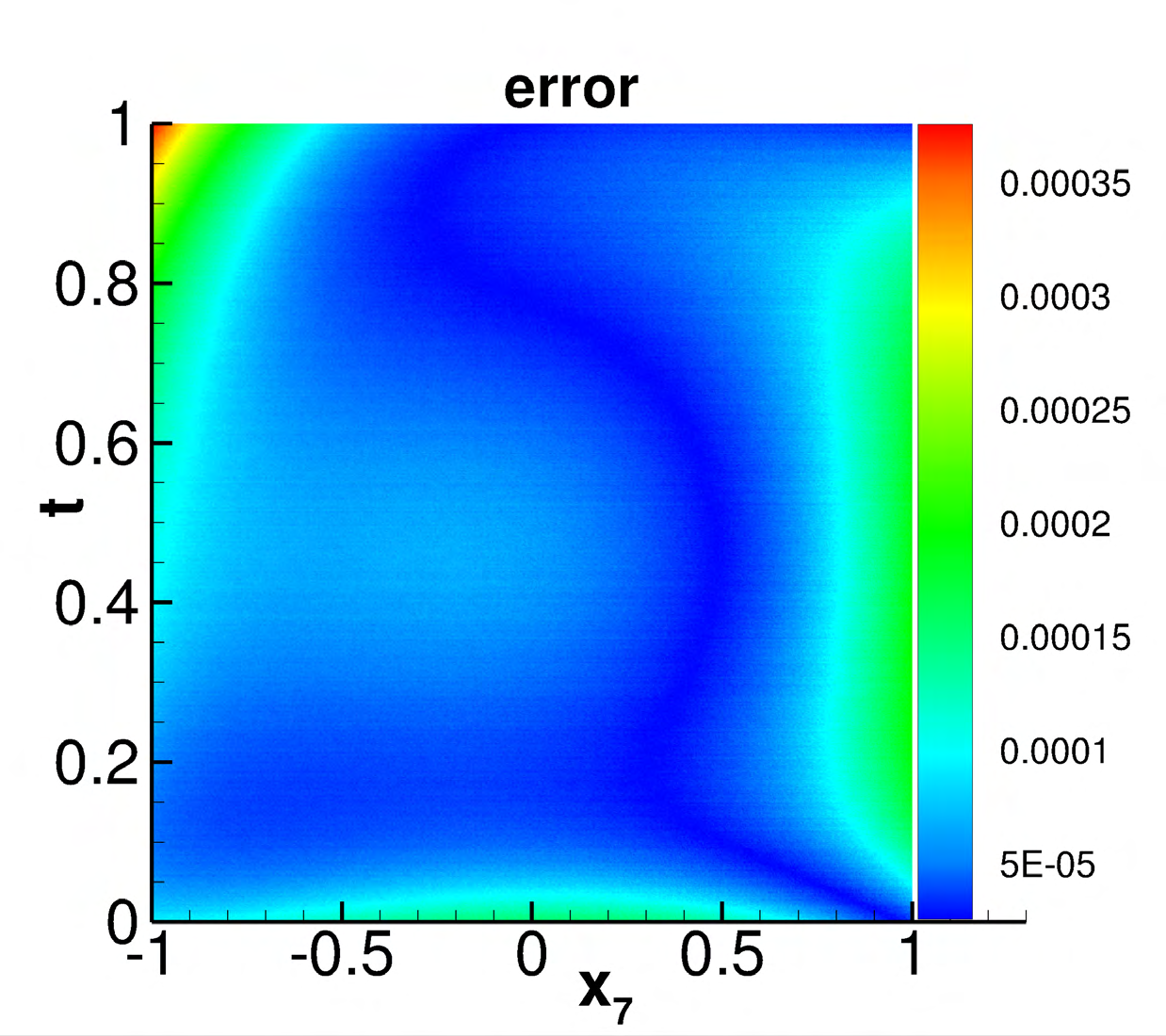}(i)}
 }
	\caption{\small Heat equation: ($d=7$): Distributions of the exact solution (left column), the ELM/A-TFC solution (middle column), and the point-wise absolute error of ELM/A-TFC in several cross sections (top row: $x_1$-$x_4$ plane; middle row: $x_4$-$x_5$ plane; bottom row: $x_7$-$t$ plane) of the spatial-temporal domain $\Omega\times[0,T]$. These cross sections are located in the middle of the spatial-temporal domain with respect to the rest of coordinates.
 NN architecture: $[d+1,3000,1]$; $(N_{\bcd},N_{\ins},N_{t_0})=(120,10,1000)$.
 }\label{tfc_heat}
\end{figure}

Figure~\ref{tfc_heat} illustrates the distributions of the exact solution, the ELM/A-TFC solution, and the point-wise absolute error of ELM/A-TFC in several cross sections of the spatial-temporal domain for the problem dimension $d=7$. It is observed that the ELM/A-TFC method has captured the solution fairly accurately, with the maximum error on the order of $10^{-4}$ in these cross sections.

\begin{figure}[tb]
	\centering
	\subfigure[d=3]
	{ \includegraphics[width=0.3\textwidth]{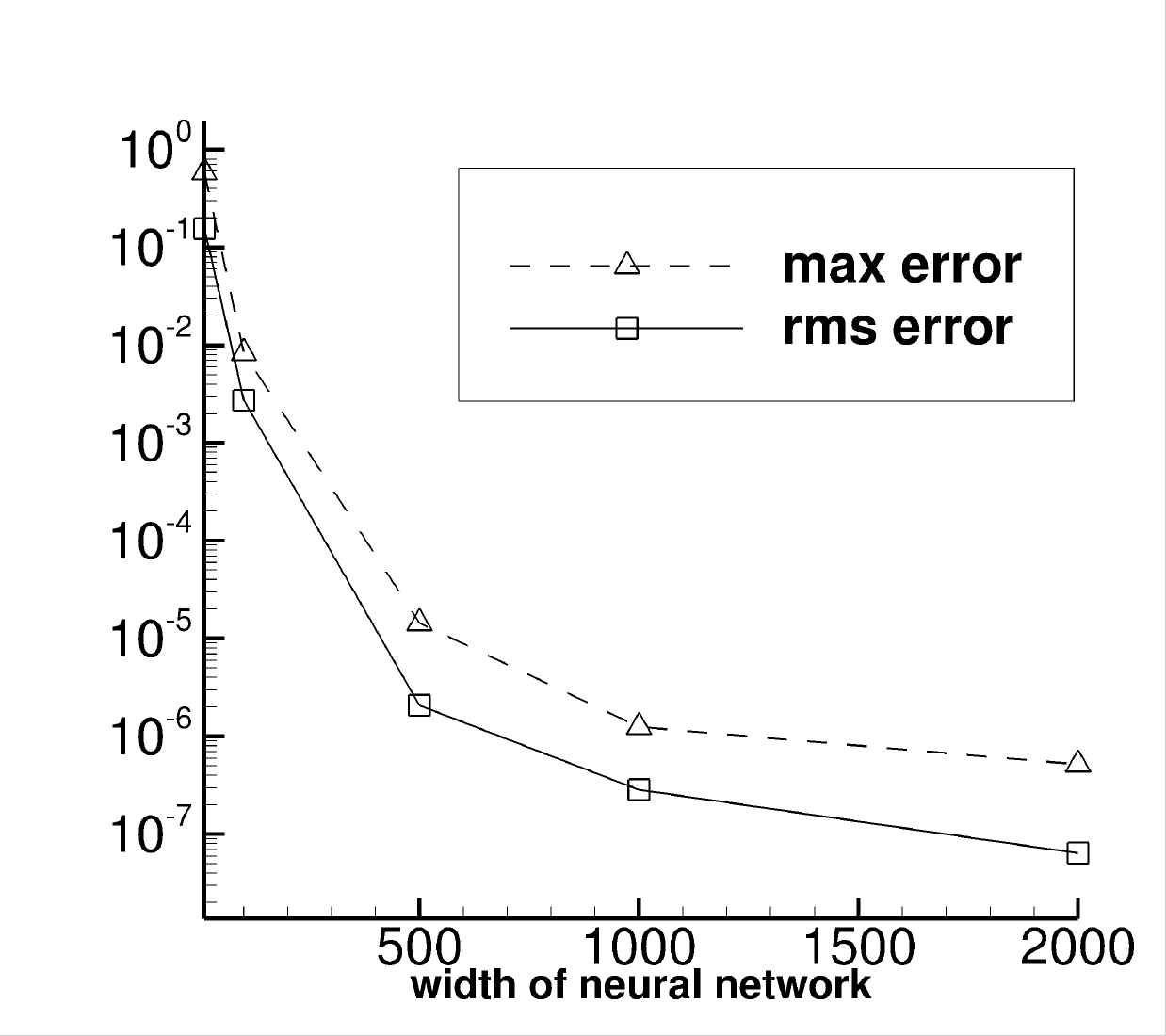}}
	\subfigure[d=7]
	{ \includegraphics[width=0.3\textwidth]{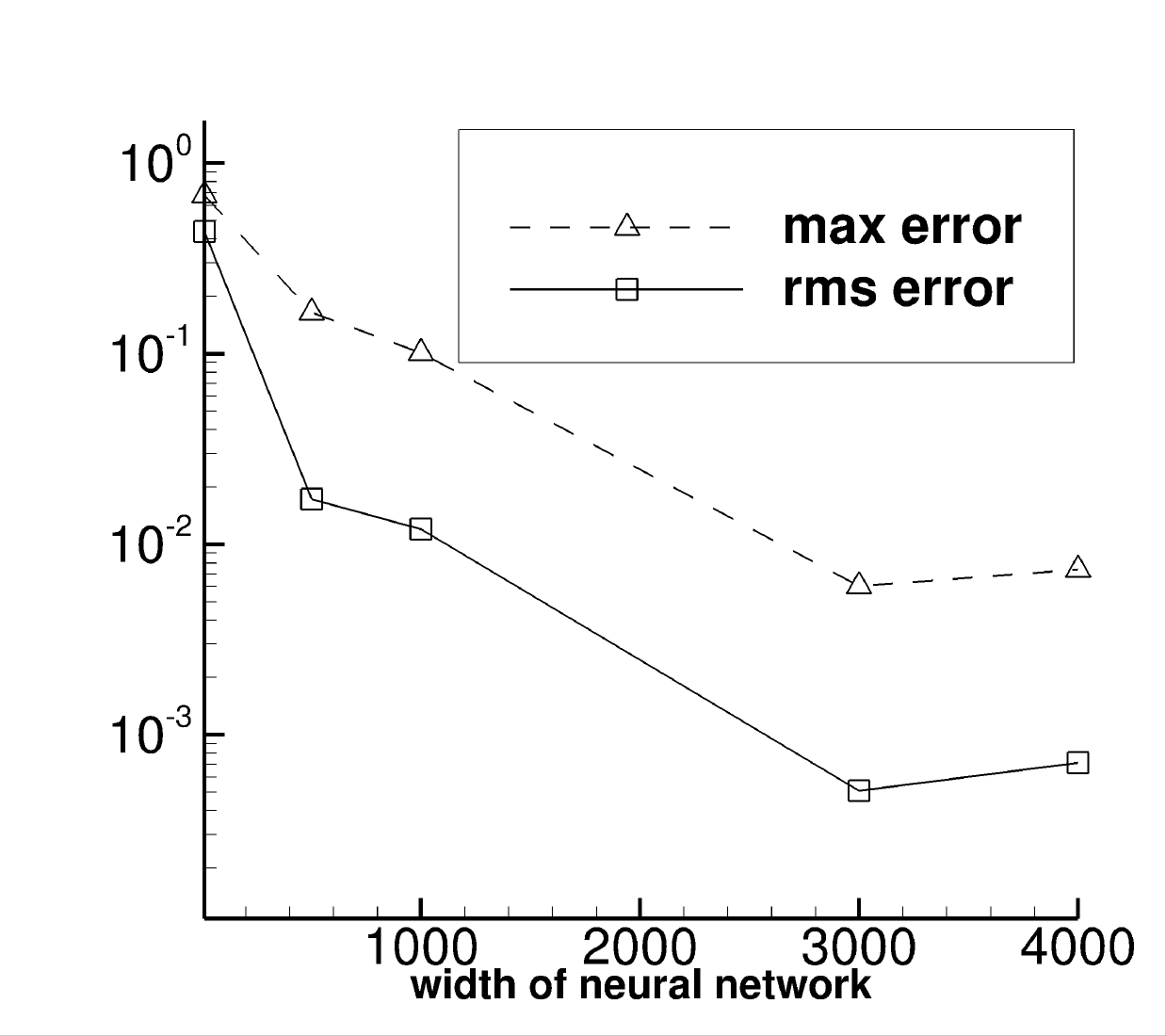}}
	\caption{\small Heat equation: $e_{l^{\infty}}$ and $e_{l^2}$ errors of ELM/A-TFC versus the number of training parameters ($M$).
 NN architecture: $[d+1,M,1]$ (with $M$ varied); $(N_{\bcd},N_{\ins},N_{t_0})=(100,100,1000)$ in (a) and $(100,10,1000)$ in (b).
 }\label{fg_27}
\end{figure}

\begin{figure}[tb]
	\centering
	\subfigure[dim=3]
	{ \includegraphics[width=0.3\textwidth]{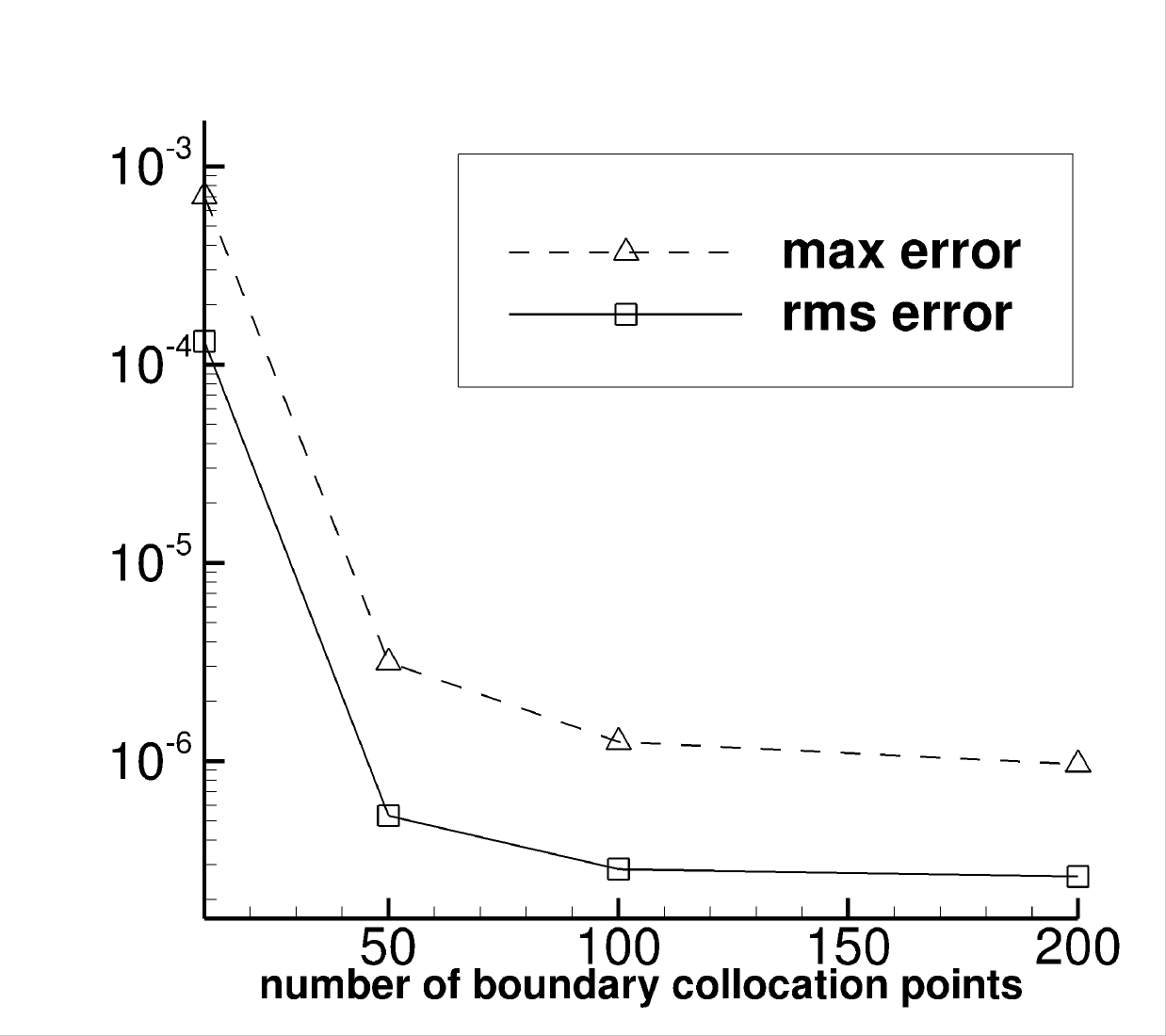}}
	\subfigure[dim=7]
	{ \includegraphics[width=0.3\textwidth]{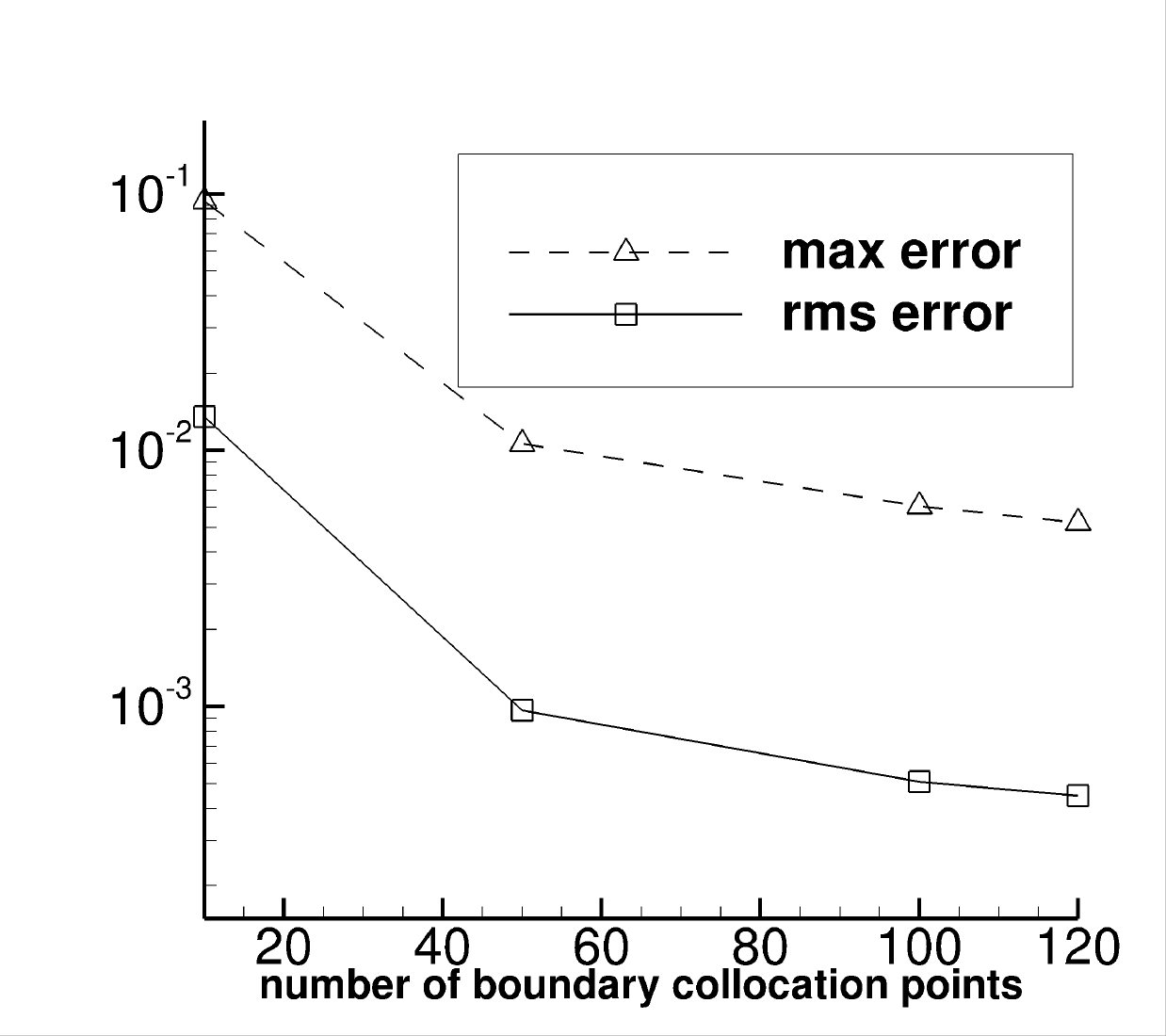}}
	\caption{\small Heat equation: $e_{l^{\infty}}$ and $e_{l^2}$ errors of ELM/A-TFC versus the number of boundary collocation points ($N_{\bcd}$).
 NN architecture: $[d+1,M,1]$, with $M=1000$ in (a) and $M=3000$ in (b). $(N_{\ins},N_{t_0})=(100,1000)$ in (a) and $(10,1000)$ in (b). $N_{\bcd}$ is varied in (a,b).
 }\label{fg_28}
\end{figure}

Figures~\ref{fg_27} and~\ref{fg_28} demonstrate the convergence behavior of the ELM/A-TFC method with respect to the number of training parameters  and the number of boundary collocation points. The simulation parameters in these tests have been provided in the figure captions. The characteristics are similar to what have been observed for other test problems in previous subsections.
With regard to the number of interior collocation points, we again observe that it has little influence on the accuracy of ELM/A-TFC (result  not shown here).

\subsection{Comparison with PINN}
\label{sec_pinn}

\begin{table}[tb]
\centering\small
\begin{tabular}{c|lll|lll}
\hline \multirow{2}{*}{ Dimension } & \multicolumn{3}{c|}{PINN} & \multicolumn{3}{c}{Current (ELM)}\\
\cline { 2 - 7 } & $e_{l^{\infty}}$ & $e_{l^{2}}$  & train-time(sec)& $e_{l^{\infty}}$ & $e_{l^{2}}$ & train-time(sec)  \\
\hline $d=3$  & 3.84E-3 &6.90E-4 &326.8 & 2.34E-10&4.42E-12& 2.6 \\
$d=5$ & 1.25E-2  & 1.45E-3 & 552.5 & 1.30E-8  &5.20E-10 &3.5 \\ 
$d=7$ & 8.12E-2&2.85E-3&688.1 & 8.16E-6 & 6.64E-7 & 11.2 \\
$d=9$ & 8.56E-1 & 1.16E-1 & 964.0 & 3.95E-5&5.26E-6&14.6 \\ \hline
\end{tabular}
\caption{\small Poisson equation: Comparison of $e_{l^{\infty}}$ and $e_{l^2}$ errors and the network training time between PINN and ELM. See the text for simulation parameters.
}
\label{tab_16}
\end{table}

\begin{table}[tb]
\centering\small
\begin{tabular}{c|lll|lll}
\hline \multirow{2}{*}{ Dimension } & \multicolumn{3}{c|}{PINN} & \multicolumn{3}{c}{Current (ELM)} \\
\cline { 2 - 7 } & $e_{l^{\infty}}$ & $e_{l^{2}}$ &train-time(sec) & $e_{l^{\infty}}$ & $e_{l^{2}}$ &train-time(sec) \\
\hline $d=3$  & 2.23E-3&4.77E-4& 1977.4 & 8.10E-10&1.20E-11&29.8 \\
$d=5$ & 2.52E-3&3.75E-4 &2526.3 &3.20E-8&1.13E-9&112.5 \\
$d=7$ & 1.89E-2&1.24E-3&3517.9 & 3.96E-5&3.23E-6& 220.9 \\
$ d=9$ & 4.781E-2&1.33E-2&5036.0 & 1.61E-4&1.27E-5&678.8 \\
\hline 
\end{tabular}
\caption{\small Nonlinear Poisson equation: Comparison of $e_{l^{\infty}}$ and $e_{l^2}$ errors and the network training time between PINN and ELM. See the text for simulation parameters.
}
\label{tab_17}
\end{table}

We next compare the current ELM method  with the physics-informed neural network (PINN) method~\cite{RaissiPK2019} for the  Poisson equation of Section~\ref{sec_poisson} and the nonlinear Poisson equation of Section~\ref{sec_312} for a range of problem dimensions. In PINN the loss function consists of two terms, the term for the PDE residual and the one for the residual of the boundary conditions. We employ the penalty coefficients  $(1-p_{\text{bc}})$ and  $p_{\text{bc}}$ in front of the loss terms for the PDE and the boundary conditions, respectively, where $p_{\text{bc}}\in(0,1)$ is a constant. The Adam optimizer is used to train the neural network in PINN. Our PINN implementation is  also based on the Tensorflow and Keras libraries.  


With PINN, we have varied the random initialization of the weight/bias coefficients, the neural network architecture, the learning rate, the learning rate schedule, and the penalty coefficient $p_{\text{bc}}$ systematically for training the neural network. The PINN results reported below are the best we have obtained in these tests. It should be noted that much poorer PINN results (not used for comparison or shown here) have been obtained in these tests.

Table~\ref{tab_16} compares the $e_{l^{\infty}}$ and $e_{l^2}$ errors, as well as the network training time (in seconds), obtained with PINN and with the current ELM method for solving the Poisson problem from Section~\ref{sec_poisson} in dimensions ranging from $d=3$ to $d=9$. The PINN results are obtained with a network architecture $[d,30,30,30,30,30,30,1]$ ($\tanh$ activation function), a penalty coefficient $p_{\text{bc}}=0.99$, and the random collocation points characterized by $(N_{\bcd},N_{\ins})=(100,3000)$ for the boundary and interior of the domain for dimensions $d=3$, $5$ and $7$ and $(N_{\bcd},N_{\ins})=(200,5000)$ for dimension $d=9$. A staircase learning rate schedule has been employed with PINN, starting with a learning rate $0.01$ and decaying by a rate $0.5$ every $500$ epochs. The PINN has been trained for a total of $5000$ epochs for the Poisson problem. 
The ELM results are obtained using a network architecture $[d,2000,1]$ for dimensions $d=3$ and $5$ and $[d,3000,1]$ for dimensions $d=7$ and $9$. The ELM random collocation points are characterized by $(N_{\bcd},N_{\ins})=(100,200)$ for $d=3$ and $5$ and $(N_{\bcd},N_{\ins})=(100,1000)$ for $d=7$ and $9$ for the boundary and interior of the domain. For generating the ELM random hidden-layer coefficients we employ $R_m=0.5$ for $d=3$,  $R_m=0.05$ for $d=5$ and $d=7$, and $R_m=0.001$ for $d=9$.  Compared with PINN, the current ELM method produces significantly more accurate results with a much smaller training time. For example, for dimension $d=5$ the PINN method produces an error on the order of $10^{-3}\sim 10^{-2}$ with a network training time close to $600$ seconds. In contrast, for this case the ELM method produces an error on the order of $10^{-10}\sim 10^{-8}$ with a network training time around $3$ or $4$ seconds.

Table~\ref{tab_17} compares the $e_{l^{\infty}}$ and $e_{l^2}$ errors, as well as the network training time, obtained by PINN and ELM for  the nonlinear Poisson problem from Section~\ref{sec_312} for dimensions ranging from $d=3$ to $d=9$. The PINN results correspond to a network architecture $[d,30,30,30,30,30,30,1]$ ($\tanh$ activation function), a penalty coefficient $p_{\text{bc}}=0.99$, a staircase learning rate schedule starting with a learning rate $0.01$ and decaying at a rate $0.8$ every $1000$ epochs, a set of random collocation points characterized by $(N_{\bc},N_{\ins})=(100,3000)$ on the boundary and interior of the domain, and a total of $20000$ training epochs. The ELM results correspond to a network architecture $[d,2000,1]$, and a set of random collocation points characterized by $(N_{\bcd},N_{\ins})=(100,200)$ for dimensions $d=3$, $5$ and $9$ and $(N_{\bcd},N_{\ins})=(100,100)$ for $d=7$ on the boundary and interior of the domain. We employ $R_m=0.5$ for $d=3$,  $R_m=0.05$ for $d=5$ and $d=7$, and $R_m=0.001$ for $d=9$ for generating the ELM random hidden-layer coefficients. The results here signify the considerably higher accuracy and less network training cost of ELM, when compared with PINN, for the nonlinear problem. For example, for dimension $d=5$ the PINN method achieves an error level on the order of $10^{-4}\sim 10^{-3}$ with a  training time around $2500$ seconds, while the ELM method achieves an error on the order of $10^{-9}\sim 10^{-8}$ with a network training time around $110$ seconds.

\section{Concluding Remarks}
\label{sec:summary}


In this paper we have presented two methods for computing high-dimensional PDEs  based on  randomized  neural networks.  These methods are motivated by the  theoretical result established in the literature that the ELM-type randomized NNs  can effectively approximate high-dimensional functions, with a rate of convergence independent of the function  dimension in the sense of expectations.

The first method extends the ELM approach, and its local variant locELM, developed in a previous work for low-dimensional problems to linear/nonlinear PDEs in high dimensions. We represent the solution field to the high-dimensional PDE problem by a randomized NN, with its hidden-layer coefficients assigned to  random values  and fixed and its output-layer coefficients trained. Enforcing the PDE problem on a set of  collocation points randomly distributed on the  interior/boundary of the domain  leads to an algebraic system of equations, which is linear for linear PDE problems and nonlinear for nonlinear PDE problems, about the ELM trainable parameters. By seeking a least squares solution to this algebraic system, attained by either a linear or a nonlinear least squares method, we can determine the values for the training parameters and complete the network training. 
ELM can be combined with domain decomposition and local randomized NNs for solving high-dimensional PDEs, leading to a local variant of this method. 
In this case,  domain decomposition is performed  along a maximum of two  designated directions for a $d$-dimensional problem, and the PDE problem, together with  appropriate continuity conditions, is enforced on the random collocation points on each sub-domain and the shared sub-domain boundaries.

Compared with the ELM for low-dimensional problems, the  difference of the method here for high-dimensional PDEs lies in at least two aspects. First, the collocation points employed for training the ELM  network for high-dimensional PDEs are randomly generated on the interior and the boundaries of the domain (or the sub-domains), and the number of interior collocation points has little (essentially no) effect on the ELM accuracy in high dimensions. In contrast, for low-dimensional PDE problems the ELM neural network is trained largely on grid-based collocation points (e.g.~uniform grid points, or quadrature points), and the number of interior collocation points  critically influences the ELM accuracy. Second, with the local variant of ELM (plus domain decomposition) for solving high-dimensional PDEs, the domain is only decomposed along a maximum of $\mathcal M$ directions, where $\mathcal M$ is a prescribed  small integer ($\mathcal M=2$ in this paper), so as for the method to be feasible in high dimensions. This is an issue not present for low-dimensional PDEs. 

The second method (ELM/A-TFC) combines the ELM approach and an approximate variant of TFC (A-TFC) for solving high-dimensional PDEs. While TFC provides a systematic approach to enforce  the boundary/initial conditions, 
the number of terms involved in  TFC constrained expression grows exponentially as the problem dimension increases, rendering it  infeasible for high-dimensional problems. By noting that the TFC constrained expression can be decomposed into a hierarchical form, we introduce the A-TFC by retaining only the   dominant terms in the constrained expression. A-TFC avoids the exponential growth in the number of terms of TFC and is feasible for high-dimensional PDEs. On the other hand, the A-TFC constrained expression 
does not unconditionally satisfy the boundary/initial conditions  for an arbitrary free function in the expression. However, the conditions that the free function in the A-TFC constrained expression needs to fulfill, in order to satisfy the boundary/initial conditions, involve functions of simpler forms, which in some sense can be considered as an effective linearization of those involved in the original boundary/initial conditions. A-TFC carries a  level of benefit of TFC for enforcing the boundary/initial conditions and is simultaneously suitable for high-dimensional problems.
With the ELM/A-TFC method, we reformulate the high-dimensional PDE problem using the A-TFC constrained expression, and attain a transformed problem about the free function involved in the A-TFC expression. We  represent this free function by ELM, and determine the ELM trainable parameters by the linear or nonlinear least squares method in a fashion analogous to  the first method. 
After the free function is determined by the ELM network, the solution field to the original high-dimensional PDE problem is then computed by  the A-TFC constrained expression. 


The two methods have been tested numerically using a number of linear/nonlinear stationary/dynamic PDEs for a range of problem dimensions. The  method has also been compared with the PINN method. We have the following observations from these numerical results: 
\begin{itemize}
\item
Both the ELM method and the ELM/A-TFC method produce accurate solutions to high-dimensional PDEs, in particular with their errors  reaching levels not far from the machine accuracy for relatively lower dimensions.  

\item
Both methods exhibit a clear sense of convergence with respect to the number of trainable parameters and the number of boundary collocation points. Their errors decrease rapidly (exponentially or nearly exponentially) for an initial range of parameter values (before saturation).

\item
The number of interior collocation points appears to have a minimal (essentially no) effect on the accuracy of ELM and ELM/A-TFC for high-dimensional PDEs. 

\item 
For a given PDE, the problem becomes more challenging to compute with increasing  dimension, in the sense that the errors of both methods in higher dimensions generally appear somewhat worse than in lower dimensions, at least with the range of  parameter values tested in this work.  

\item 
The error levels obtained by the ELM method and the ELM/A-TFC method are generally comparable, with  ELM/A-TFC  appearing slightly better in lower dimensions. On the other hand, the ELM/A-TFC method generally involves a larger computational effort and cost than ELM, due to the A-TFC constrained expression.

\item 
The current method exhibits a clear advantage compared with  PINN  for solving high-dimensional PDEs, and achieves a significantly better accuracy under markedly smaller training time than the latter.


\end{itemize}
The simulation results signify that the ELM-based methods developed herein are effective for computational PDEs in high dimensions.

\section*{Acknowledgment}

This work was partially supported by the US National Science Foundation (DMS-2012415).

 \bibliographystyle{plain}
	\bibliography{references,dnn1,dnn,elm,elm1,mypub,varpro}
\end{document}